\documentstyle[epsf,12pt]{article}
\def\Bbb R{{\rm \bf R}}
\def\proclaim#1{\vskip2mm{\bf #1}\em}
\def\endproclaim{\em \vskip2mm}
\def\tag#1{\eqno(#1)}
\def\gathered{\begin{array}{c}}
\def\endgathered{\end{array}}
\def\text{\mbox}

\begin{document}

\title {{\bf Extracting discontinuity using the probe and enclosure methods}}
\author{Masaru IKEHATA\footnote{
Laboratory of Mathematics,
Graduate School of Advanced Science and Engineering,
Hiroshima University,
Higashihiroshima 739-8527, JAPAN}}
\maketitle

\begin{abstract}
This is a review article on the development of the probe and enclosure methods from past to present, focused on their central ideas together with various applications.

\noindent
AMS: 35R30, 35L05

\end{abstract}


\tableofcontents

\section{Introduction}

More than twenty years have been passed since the {\it probe} and {\it enclosure} methods were introduced by the author.
Those are analytical methods based on the governing equations in so-called {\it inverse obstacle problems}.
The purpose of solving inverse obstacle problems is to extract various information about unknown {\it discontinuity}
embedded in a given medium or material, such as cavities, inclusions, obstacles, etc., from the data observed
at an accessible part.  The observation data depend on the used signal which propagates inside the medium and are described by using the solutions of the governing equation of the signal.  The common idea of both methods, constructs the so-called {\it indicator function} by using the observed data.
And from the asymptotic behaviour one extracts information about unknown discontinuity.  In early times the object of the methods 
is only the stationary, or frequency domain data, however, now we have succeeded to extend the range of applications
of the methods to time dependent observation data.  In this paper, we make an order of the development of the both methods, in particular,
focused on the main idea together with several results, however, not to overlap from the previous survey papers 
\cite{PSurvey0, ESurvey, ITwomethods, ISugakuEnglish, ISA, IRomania} with a different point of view.
Besides some of still unsolved questions and remarks on the methods themselves or in their applications
are also mentioned.







\section{The Probe Method}

In 1997 April in a train to Tokyo the author suddenly got a general idea to extract 
information about unknown discontinuity embedded in a known
background medium from the Dirichlet-to-Neumann map on the boundary of the medium.
The author named the idea the {\it probe method} because of its geometric style.

In the same year Eighth international colloquium on differential equations was held at Plovdiv, 
Bulgaria, 1997 August 18th-23rd.  In the colloquium originally the author planned to make a presentation 
of some uniqueness results in elasticity, however, asked the organizer Prof. Drumi Bainov to change the contents of the talk
and title.  He accepted this proposal and the probe method was presented.  This means that the year when the probe method
was firstly presented in front of the third parsons abroad is 1997.

Later, in 1998-1999 two applications of the probe method were published.
The first one \cite{IProbe} is an application to an important version of the inverse conductivity problem raised by A. P. Calder\'on \cite{Cal} and 
another one \cite{IScattering0} together with \cite{IScattering1} an application
to the inverse obstacle scattering problem at a fixed frequency.
These are briefly explained in the author's review paper \cite{PSurvey0} about early development of the probe method.
See also \cite{ISA} for compact explanations on the proofs of the results \cite{IScattering0, IScattering1} in inverse obstacle scattering problems.

\subsection{Inverse crack problem}

The probe method is a mathematical method of {\it probing} inside a given domain by using a {\it virtual needle}
inserted from the boundary of the domain.
Here we introduce the probe method by making a review of a result \cite{IProbeCrack} for an inverse crack problem.

Assume that we have a body $\Omega$ which is a bounded and connected open subset of $\Bbb R^m$ ($m=2,3$) with a smooth boundary.
We denote by $\Sigma$ the unknown crack occurred in the body.  
We assume that $\Sigma$ is given by a $(m-1)$-dimensional closed submanifold of $\Bbb R^m$ of class $C^0$ with boundary.  
The $\Sigma$ is divided into two parts: the interior and the boundary denoted by $\text{Int}\,\Sigma$ and $\partial\Sigma$, respectively
(cf. \cite{L}).  Note that $\text{Int}\,\Sigma=\Sigma\setminus\partial\Sigma$.

We say that $\Sigma$ is {\it extendable} of class $C^{0,1}$, if $\Sigma$ admits the existence of an open set $D$ of $\Bbb R^m$ 
with Lipschitz
boundary, having finitely many connected components and
satisfying the following:
$$
(\star)\left\{
\begin{array}{l}
\displaystyle
\text{$\overline D\subset\Omega$;}\\
\\
\displaystyle
\text{$\Omega\setminus\overline D$ is connected;}\\
\\
\displaystyle \text{$\Sigma\subset\partial D$.}
\end{array}
\right.
$$
\noindent
Of course, there should be infinitely many $D$ satisfying ($\star$)
for given extendable $\Sigma$.  In this paper, we always
assume that $\Sigma$ is extendable of class $C^{0,1}$ unless otherwise
specified.

We denote by $\nu$ the unit outward
normal relative to $D$ and $\Omega$ unless otherwise specified.
Set $\partial D=\Gamma$.  Let
$\Omega_{+}=\Omega\setminus\overline D$ and write $D=\Omega_{-}$.
For a function $v\in L^2(\Omega)$ set
$$\displaystyle
v_{+}=v\vert_{\Omega_{+}},\,v_{-}=v\vert_{\Omega_{-}}.
$$

\noindent We introduce some Hilbert spaces. Define
$$\begin{array}{c}
\displaystyle
X(\Omega\setminus\Sigma;D)
=\{v\in L^2(\Omega)\,\vert\,
v_{+}\in H^1(\Omega_{+})\,, v_{-}\in H^1(\Omega_{-}),
v_{+}\vert_{\Gamma}(x)
=v_{-}\vert_{\Gamma}(x)\,\text{a.e.}\,x\in\Gamma\setminus\Sigma\};\\
\\
\displaystyle
\Vert v\Vert_{X(\Omega\setminus\Sigma;D)}
=\Vert v_{+}\Vert_{H^1(\Omega_{+})}
+\Vert v_{-}\Vert_{H^1(\Omega_{-})}.
\end{array}
$$
$X(\Omega\setminus\Sigma;D)$ is complete
with respect to the norm $\Vert\,\cdot\,\Vert_{X(\Omega\setminus\Sigma;D)}$.

\noindent
Define
$$\displaystyle
X_0(\Omega\setminus\Sigma;D)=\{v\in\,X(\Omega\setminus\Sigma;D)\,\vert\,v\vert_{\partial\Omega}(x)=0\,
\text{a.e.}\,x\in\partial\Omega\}.
$$
This space is a closed subspace of $X(\Omega\setminus\Sigma;D)$.

We clarify what we mean by the {\it weak solution}.

{\bf\noindent Definition 2.1.}
Given $f\in H^{1/2}(\partial\Omega)$ we say that $u\in X(\Omega\setminus\Sigma;D)$ is
a weak solution of the elliptic problem
$$\left\{
\begin{array}{ll}
\displaystyle
\Delta u=0, & x\in\Omega\setminus\Sigma,\\
\\
\displaystyle
\frac{\partial u}{\partial\nu}=0, & x\in\Sigma,\\
\\
\displaystyle
u=f, & x\in\partial\Omega
\end{array}
\right.
\tag {2.1}
$$
if $u$ satisfies $u=f$ on $\partial\Omega$ in the sense of trace
and, for all $\varphi\in X_0(\Omega\setminus\Sigma;D)$
$$\displaystyle
\int_{\Omega\setminus\Sigma}\nabla u\cdot\nabla\varphi dy=0.
$$

The existence and uniqueness of the weak solution of (2.1)
and the invariance of the solution with respect to the choice of $D$ as long as $D$ satisfies condition ($\star$)
are summarized as

\proclaim{\noindent Proposition 2.1(\cite{IProbeCrack}).}
For each fixed $D$ satisfying ($\star$) there exists a unique weak solution of (2.1).
Moreover the solution does not depend on the choice of $D$.
\endproclaim

The proof was omitted in \cite{IProbeCrack}, however, later given in Appendix A.1 in \cite{IO2}.

{\bf\noindent Definition 2.2.}
For each fixed $D$ satisfying ($\star$), define the bounded linear functional
$\Lambda_{\Sigma}f$ on $H^{1/2}(\partial\Omega)$
by the formula
$$\displaystyle
\int_{\partial\Omega}\,(\Lambda_{\Sigma}f)h\,dS
=\int_{\Omega\setminus\Sigma}\nabla u\cdot\nabla vdy,\,\,h\in H^{1/2}(\partial\Omega),
\tag {2.2}
$$
where $u$ is the weak solution of (2.1) and $v\in X(\Omega\setminus\Sigma;D)$ is an arbitrary function with
$v=h$ on $\partial\Omega$ in the sense of the trace. Note that the integral on the left hand side on (2.2)
is a formal expression for the dual pairing $<\,,\,>$ between $H^{-1/2}(\partial\Omega)$ and $H^{1/2}(\partial\Omega)$.

The map $\Lambda_{\Sigma}:f\longmapsto \Lambda_{\Sigma}f$ is called the {\it Dirichlet-to-Neumann map}.
We set $\Lambda_{\Sigma}=\Lambda_0$ in the case when $\Sigma=\emptyset$.
Hereafter we formally write 
$$\begin{array}{ll}
\displaystyle
\Lambda_{\Sigma}f=\frac{\partial u}{\partial\nu}, & x\in\partial\Omega
\end{array}
$$

Since the solution does not depend on the choice of $D$ as long as $D$ satisfies condition ($\star$),
the Dirichlet-to-Neumann map $\Lambda_{\Sigma}$ also does not depend on $D$.

The pair $(f, \Lambda_{\Sigma}f)$ which is a point on the graph of $\Lambda_{\Sigma}$
is viewed as a possible pair of the voltage and current density in the continuum model of electric impedance tomography (EIT). 

In \cite{IProbeCrack} we considered the following inverse problem.

$\quad$

{\bf\noindent Problem 2.1.}  Extract information about the shape and location of $\Sigma$ from $\Lambda_{\Sigma}$.

$\quad$

On uniqueness, Eller \cite{E} had proven that the full knowledge of $\Lambda_{\Sigma}$ uniquely determines $\Sigma$ itself.
His proof employs a contradiction argument and is based on Isakov's {\it method of singular solutions} developed in \cite{Is0}.
See \cite{IkNaSurvey} for a detailed review of his method.
Alessandrini-DiBenedetto \cite{AD} considered the case when unknown cracks are planar (flat) in a given three dimensional body
and established a uniqueness theorem in a different formulation.  
See also the references therein and a classical uniqueness result \cite{FV} in two dimensions.

Problem 2.1 asks us to find an {\it exact procedure} for extracting information about the geometry of $\Sigma$.
For this we have already two approaches \cite{IN} and \cite{IProbeCrack} using the probe method: one is an application of 
the original version of the probe method which is nothing but the version introduced in \cite{IProbe}
and another a new formulation of the probe method published in \cite{IProbeNew}.


\subsubsection{Needle and needle sequence}

The probe method starts with introducing the notion of the needle.

{\bf\noindent Definition 2.3.}  Given a point $x\in\Omega$ let $N_x$ denote the set of all piecewise linear curves 
$\sigma:[0,\,1]\rightarrow\overline\Omega$ such that

$\bullet$  $\sigma(0)\in\partial\Omega$, $\sigma(1)=x$ and $\sigma(t)\in\Omega$ for all $t\in\,]0,\,1[$.

$\bullet$  $\sigma$ is injective.

We call $\sigma\in N_x$ a {\it needle with tip at} $x$.  We write also $\sigma([0,\,1])=\sigma$.

The next important step is a design of the special sequence of functions $\{v_n\}$.

In this paper, a set $V$ of $\Bbb R^m$ is called a {\it finite cone} with vertex $x\in\Bbb R^m$ if $V$ has the expression
$$\displaystyle
V=B_{\rho}(x)\cap\text{int}\,C_x(\omega,\alpha),
$$
where $\omega$ is a unit vector in $\Bbb R^m$, $\rho>0$, $\alpha\in\,]0,\,1[$,
$B_{\rho}(x)=\{y\in\Bbb R^m\,\vert\,\vert y-x\vert<\rho\}$,
$$\displaystyle
C_x(\omega,\alpha)
=\left\{y\in\Bbb R^m\,\left\vert\right.\,(y-x)\cdot\omega\ge\vert y-x\vert\cos\frac{\pi\alpha}{2}\,\right\}
$$
and $\text{int}\,C_x(\omega,\alpha)$ denotes its open kernel.  We call the closed set $C_x(\omega,\alpha)$
the $m$-space dimensional cone about $\omega$ of opening angle $\frac{\pi\alpha}{2}$ with vertex at $x$.

Let $G=G(y)$ be an arbitrary solution of the Laplace equation in $\Bbb R^m\setminus\{0\}$
such that for any finite cone $V$ with vertex at $0$
we have
$$\displaystyle
\int_V\,\vert\nabla G(y)\vert^2\,dy=\infty.
\tag {2.3}
$$
Hereafter we fix $G$.

{\bf\noindent Definition 2.4.}  Let $x\in\Omega$.
Let $\sigma\in N_x$.  We call the sequence $\xi=\{v_n\}$ of $H^1(\Omega)$ solutions of the Laplace equation in $\Omega$ 
a {\it needle sequence} for $(x,\sigma)$ if it satisfies, for each fixed compact set $K$ of $\Bbb R^m$
with $K\subset\Omega\setminus\sigma$
$$\displaystyle
\lim_{n\rightarrow\infty}\left(\Vert v_n(\,\cdot\,)-G(\,\cdot\,-x)\Vert_{L^2(K)}
+\Vert\nabla\{v_n(\,\cdot\,)-G(\,\cdot\,-x)\}\Vert_{L^2(K)}\right)=0,
$$
that is, $v_n\rightarrow G(\,\cdot\,-x)$ in $H^1_{\text{loc}}(\Omega\setminus\sigma)$.
Note that this yields the convergence of the higher order derivatives of $v_n$, that is,
$v_n\rightarrow G(\,\cdot\,-x)$ in $H^m_{\text{loc}}(\Omega\setminus\sigma)$ with $m\ge 2$
and thus all the derivatives converge to those of  $G(\,\cdot\,-x)$
compact uniformly in $\Omega\setminus\sigma$ since $\{v_n\}$ is a sequence of harmonic functions.

It is well known that the {\it existence} of the needle sequence for arbitrary given $(x,\sigma)$ has been ensured
by applying the Runge approximation property of the Laplace equation, which is a consequence of the unique continuation
property of the solution of the Laplace equation.  This property for elliptic equations 
has been used in the uniqueness issue of the Calder\'on problem for a class of
piecewise analytic conductivities in \cite{KVII} and for an inverse inclusion problem in \cite{Is0}. 
For the proof, for example, see \cite{KVII}.  Note also that one can add a constraint to the size of the support
of the trace of $v_n$ onto $\partial\Omega$.

\proclaim{\noindent Proposition 2.2.}  Let $B$ be an open ball with $B\cap\partial\Omega\not=\emptyset$.
Given $x\in\Omega$ and $\sigma\in N_x$ there exists a needle sequence $\xi$ for $(x,\sigma)$ with
$\text{supp}\,(v_n\vert_{\partial\Omega})\subset B\cap\partial\Omega$.

\endproclaim

See, \cite{Is0} and Appendix in \cite{IProbe2} for the proof of such type of
the Runge approximation property.  
In the present case the function $G(\,\cdot\,-x)$ is harmonic in $\Bbb R^m\setminus\{x\}$, one can apply the classical Runge approximation
property to an exhaustion of the set $\Omega\setminus\sigma$, provided, for example, $\Bbb R^m\setminus\overline\Omega$ is connected.
In this case one can choose each $v_n$ from harmonic polynomials, however, of course, no constraint on the size of
the trace of $v_n$ onto $\partial\Omega$.  See \cite{dP}.

The needle sequence has two important properties which are not pointed out in the original version of the probe method \cite{IProbe}.

\proclaim{\noindent Proposition 2.3.}  Let $x\in\Omega$ and $\sigma\in N_x$.  Let
$\xi=\{v_n\}$ be an arbitrary needle sequence for $(x,\sigma)$.  We have

$\bullet$  for any finite cone $V$ with vertex at $x$
$$\displaystyle
\lim_{n\rightarrow\infty}\int_{V\cap\Omega}\,\vert\nabla v_n(y)\vert^2\,dy=\infty.
$$

$\bullet$  for any point $z\in\,\sigma(]0,\,1[)$ and open ball $B$ centered at $z$
$$\displaystyle
\lim_{n\rightarrow\infty}\int_{B\cap\Omega}\,\vert\nabla v_n(y)\vert^2\,dy=\infty.
$$

\endproclaim

This is a special case of a result in \cite{IProbeNew} which covers also the case when the governing equation
of the needle sequence is replaced with the Helmholtz equation $(\Delta+k^2)v=0$ with an arbitrary fixed $k>0$.

From these and Definition 2.4 we can recover the needle $\sigma$ by the formula
$$\displaystyle
\sigma(]0,\,1])
=\left\{y\in\Omega\,\vert\,\text{for all open ball $B$ centered at $y$}\,
\lim_{n\rightarrow\infty}\int_{B\cap\Omega}\vert\nabla v_n(y)\vert^2\,dy=\infty\right\},
$$
where $\{v_n\}$ is an arbitrary needle sequence for $(x,\sigma)$.  This means that the needle sequence
is an {\it analytical realization} of the needle without loosing the geometry.

\subsubsection{Indicator sequence and function}

Using the needle sequence and Dirichlet-to-Neumann map we define another sequence which is called
the {\it indicator sequence}.

{\bf\noindent Definition 2.5.}  Given $x\in\Omega$, $\sigma\in N_x$ and needle sequence $\xi=\{v_n\}$ for $(x,\sigma)$ define
$$\begin{array}{ll}
\displaystyle
I(x,\sigma,\xi)_n=\int_{\partial\Omega}\{(\Lambda_0-\Lambda_{\Sigma})f_n\}f_n\,dS,
&
\displaystyle
n=1,2,\cdots,
\end{array}
$$
where $f_n(y)=v_n(y)$, $y\in\partial\Omega$.

To study the behaviour of the indicator sequence as $n\rightarrow\infty$ we prepare

{\bf\noindent Definition 2.6.}  Define the {\it indicator function} $I(x)$ of independent variable $x\in\Omega\setminus\Sigma$ by the formula
$$\displaystyle
I(x)=\int_{\Omega\setminus\Sigma}\vert\nabla w_x\vert^2\,dy,
\tag {2.4}
$$
where the function $w_x=w_x(y)$ is the {\it reflected solution} by $\Sigma$, that is, the unique weak solution in $X_0(\Omega\setminus\Sigma;D)$
of the variational problem:
find $w\in X_0(\Omega\setminus\Sigma;D)$ such that, 
for all $\Psi\in X_0(\Omega\setminus\Sigma;D)$
$$\displaystyle
\int_{\Sigma}\frac{\partial G(\,\cdot\,-x)}{\partial\nu}(\Psi_{+}-\Psi_{-})dS
=\int_{\Omega\setminus\Sigma}\nabla w\cdot\nabla\Psi dx.
$$
The solution $w_x$ is also independent of $D$ as long as $D$ satisfies ($\star$).

\subsubsection{Two sides of the probe method}

The following result proved in \cite{IProbeCrack} states
the way of calculating the value of the indicator function from the indicator sequence
and the properties of the calculated indicator function.

\proclaim{\noindent Theorem 2.1.}

$\bullet$  (A.1)  Given $x\in\Omega\setminus\Sigma$ let $\sigma\in N_x$ satisfy $\sigma\cap\Sigma=\emptyset$.
Then, for all needle sequence $\xi=\{v_n\}$ for $(x,\sigma)$ the indicator sequence $\{I(x,\sigma,\xi)_n\}$
converges to the value of indicator function $I$ at $x$.

$\bullet$  (A.2)  For each $\epsilon>0$ we have $\sup_{\text{dist}\,(x,\Sigma)>\epsilon}\,I(x)<\infty$.

$\bullet$  (A.3)  Given $a\in\text{Int}\,\Sigma$  we have $\lim_{x\rightarrow a}I(x)=\infty$.

\endproclaim

{\it\noindent Sketch of Proof.} 
We only mention assertion (A.1) since assertion (A.2) is a consequence of the well-posedness of the direct problem
and assertion (A.3) can be proved similarly to Theorem 2.2 mentioned later.

Given $x\in\Omega$, $\sigma\in N_x$ and needle sequence $\xi=\{v_n\}$ for $(x,\sigma)$ let 
$u_n\in X(\Omega\setminus\Sigma;D)$ be the weak solution of (2.1) for $f=v_n\vert_{\partial\Omega}$.
Here we have the representation formula of the indicator sequence
$$\displaystyle
I(x,\sigma,\xi)_n=
\int_{\Omega\setminus\Sigma}\vert\nabla w_n\vert^2dx,
\tag {2.5}
$$
where $w_n=u_n-v_n\in X_0(\Omega\setminus\Sigma;D)$ and satisfies
for all $\Psi\in X_0(\Omega\setminus\Sigma;D)$
$$\displaystyle
\int_{\Sigma}\frac{\partial v_n}{\partial\nu}(\Psi_{+}-\Psi_{-})dS
=\int_{\Omega\setminus\Sigma}\nabla w_n\cdot\nabla\Psi dx.
\tag {2.6}
$$
See Appendix A.2 in \cite{IO2} for the proof.  Besides, if $\sigma\cap\Sigma=\emptyset$, 
it is easy to see that $w_n\rightarrow w_x$ in $X_0(\Omega\setminus\Sigma;D)$ as $n\rightarrow\infty$.
Now from (2.4) one gets (A.1).

\noindent
$\Box$

Note that, in the case when $\text{Int}\,\Sigma$ is smooth, an equivalent form of Theorem 2.1 which is written in terms
of the original probe method \cite{IProbe} has been established in \cite{IN}.

After having Theorem 2.1, it is natural to ask the question that {\it what happens on the
indicator sequence if the tip of the needle is located on the crack or the needle passed through the crack}.

It should be emphasized that in the previous applications of the probe method \cite{IProbe, IProbe2, IScattering0, IScattering1}
the question raised above has not been considered and it was difficult to apply the techniques developed in those papers.
However, in \cite{IProbeCrack} we found a completely different techniques which also yields Theorem 2.1 and 
an answer to the question mentioned above.
However, by reviewing the original proof, we found that it
can not cover completely the latter case in the question,
that is, the case when the needle passed the crack.  We have to impose a restriction on the position of the needle
relative to the crack mentioned below.

{\bf\noindent Definition 2.7.}  Let $x\in\Omega$ and $\sigma\in N_x$.
We say that $\sigma$ is {\it normally passing through} $\Sigma$ if $x\in\Omega\setminus\partial\Sigma$ and
there exists a bounded open set $D=D_{\sigma}$ of $\Bbb R^m$ with Lipschitz boundary, having finitely many connected components
and satisfying ($\star$) such that
$\emptyset\not=\sigma\cap\Gamma_{\sigma}\subset\text{Int}\,\Sigma$, where $\Gamma_{\sigma}=\partial D_{\sigma}$.

\begin{figure}[htbp]
\begin{center}
\epsfxsize=4.8cm
\epsfysize=4.8cm
\epsfbox{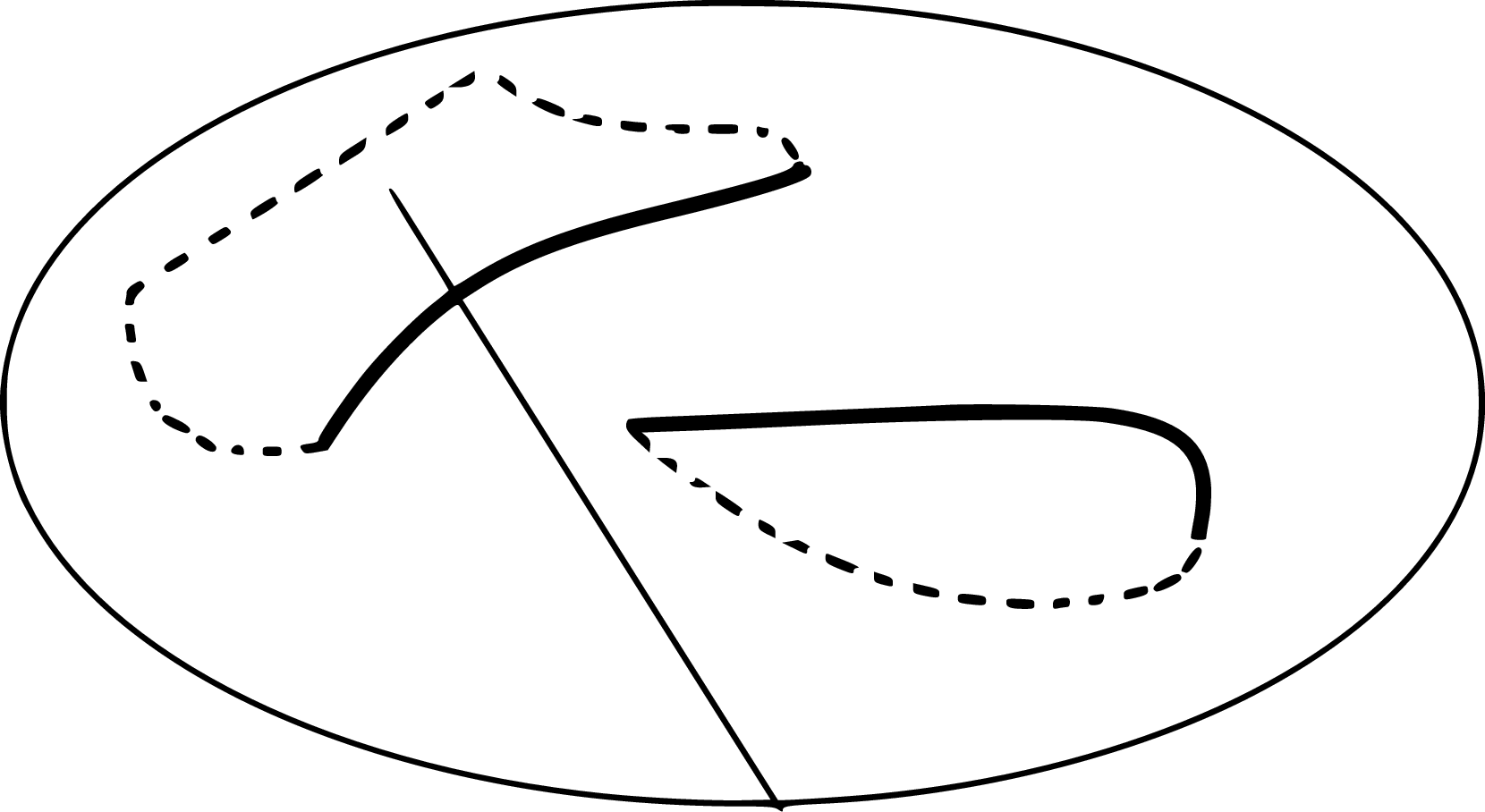}
\caption{An illustration of a needle $\sigma$ which is normally passing through $\Sigma$.
All the connected components of $D$ are indicated by the area enclosed by the dotted and thick curves.
}\label{fig1}
\end{center}
\end{figure}

\begin{figure}[htbp]
\begin{center}
\epsfxsize=4.8cm
\epsfysize=4.8cm
\epsfbox{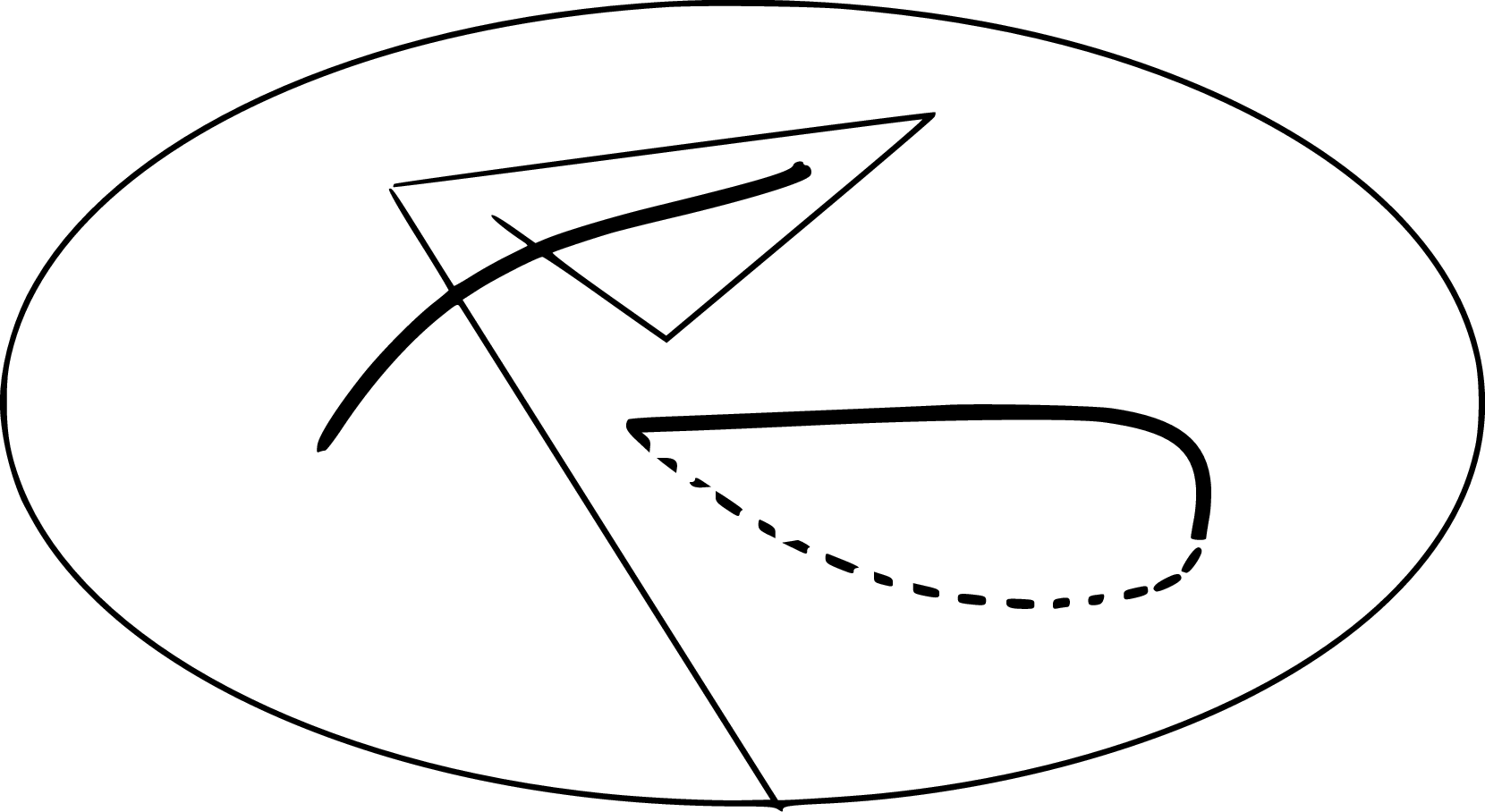}
\caption{An illustration of a needle $\sigma$ which is not normally passing through $\Sigma$.
A possible connected component of $D$ is indicated by the area enclosed by the dotted and thick curves.
}\label{fig2}
\end{center}
\end{figure}

See Figure \ref{fig1} for an illustration of $\sigma$ which is normally passing through $\Sigma$ and in Figure \ref{fig2}
not being so.  Therein $\Sigma$ is indicated by thick curves.

\vspace{0.0cm}

\vspace{0.0cm}

\newpage

The following statement is the corrected version of the original one, published in \cite{IProbeCrack} as Theorem B.
However there is no need to change the proof given there.

\proclaim{\noindent Theorem 2.2.}  
Let $x\in\Omega\setminus\partial\Sigma$ and $\sigma\in N_x$ be normally passing through $\Sigma$. 
Then, for all needle sequence $\xi=\{v_n\}$ for $(x,\sigma)$ we have $\lim_{n\rightarrow\infty}I(x,\sigma,\xi)_n=\infty$.

\endproclaim

{\it\noindent Sketch of Proof.}
By the assumption we can assume that, $D$ is given in such a way that $\sigma\cap\Gamma\subset\text{Int}\,\Sigma$.
Since $\Sigma\cap\Gamma$ is compact, there exists $\eta\in C^{\infty}_0(\Omega)$
with $0\le\eta\le 1$ such that $\eta=1$ in a neighbourhood of $\sigma\cap\Gamma$ and $\text{supp}\,(\eta\vert_{\Gamma})\subset\Sigma$.

Let $\xi=\{v_n\}$ be a needle sequence for $(x,\sigma)$ and $w_n=u_n-v_n$.
Choosing a test function $\Psi=\Psi_n$ in (2.6) which depends on $v_n$ and $\eta$, 
we obtain a lower estimate for $\Vert\nabla w_n\Vert_{L^2(\Omega_{-})}^2$ and thus, from (2.5)
one gets the lower estimate of the indicator sequence from below by using $v_n$ only.
It takes the form
$$\displaystyle
\frac{\displaystyle
\left\vert
\int_{\Omega_{-}}\vert\nabla v_n\vert^2\,dy-\int_{\Gamma}\frac{\partial v_n}{\partial\nu}(1-\eta)v_n\,dS
\right\vert^2}
{\displaystyle
M^2
\left(
\int_{\Omega_{-}}\eta\vert\nabla v_n\vert^2\,dy+\int_{\Omega_{-}}\vert\nabla\eta\vert^2\vert v_n\vert^2\,dy
\,\right)}
\le I(x,\sigma,\xi)_n,
$$
where $M$ is a positive constant independent of $v_n$.

Since $\text{supp}\,(1-\eta)\vert_{\Gamma}\subset\Omega\setminus\sigma$, from the convergence property
of $\{v_n\}$ in $\Omega\setminus\sigma$ we see that the integral
$$\displaystyle
\int_{\Gamma}\frac{\partial v_n}{\partial\nu}(1-\eta)v_n\,dS
$$
converges to
$$\displaystyle
\int_{\text{supp}\,(1-\eta)\vert_{\Gamma}}\frac{\partial G(\,\cdot-x)}{\partial\nu}(1-\eta)G(\,\cdot\,-x)\,dS.
$$
and thus is bounded.  Note also that we have
$$\displaystyle
\lim_{n\rightarrow\infty}\,\int_{\Omega_{-}}\vert\nabla v_n\vert^2\,dy=\infty.
$$
since $\sigma\cap \Omega_{-}\not=\emptyset$.

Here, using a version of the Poincar\'e inequality (e.g., \cite{Z}) and the convergence property of $\{v_n\}$ in $\Omega\setminus\sigma$, one gets
$$\displaystyle
\limsup_{n\rightarrow\infty}
\frac{\displaystyle\int_{\Omega_{-}}\vert v_n\vert^2\,dy}
{\displaystyle\int_{\Omega_{-}}\vert\nabla v_n\vert^2\,dy}<\infty.
$$
Thus one can conclude that $\lim_{n\rightarrow\infty}I(x,\sigma,\xi)_n=\infty$.

\noindent
$\Box$

Some remarks are in order.

$\bullet$  In the original version of Theorem 2.2, instead of the normality of $\sigma$ relative to $\Sigma$ it is assumed that $\emptyset\not=\sigma\cap\Sigma\subset\text{Int}\,\Sigma$.  However, the original proof
essentially makes use of the normality.  Now we know that this condition does not ensure the normality.
See again Figure \ref{fig2}.

$\bullet$  From a technical point of view the ideas of the proofs in \cite{IN} and \cite{IProbeCrack} for Theorem 2.1 are completely different.  
In \cite{IN} to show the blowing up of the indicator function when the tip of the needle approaches the crack the leading profile
of the so-called reflected solution is constructed and in contrast to this, in \cite{IProbeCrack} only integration
by parts and resulted energy estimates are employed.  Besides the result in \cite{IProbeCrack} stated as Theorem 2.2 in this section
gives also the behaviour of the indicator sequence when the needle {\it passed through} the crack.
This is not covered in \cite{IN}.

\subsection{Inverse inclusion problem}

\subsubsection{Reconstruction of inclusion}

In \cite{IProbe} and \cite{IProbe2} we considered a special version
of the Calder\'on problem \cite{Cal} which is based on the continuum model of Electrical Impedance Tomography \cite{Bo, Bo2}.  
See \cite{SCI} for more realistic model of EIT.

The problem is to extract information about the location and shape of unknown inclusions
embedded in a known background body. We model the body by $\Omega$ and denote by $D$ the set of unknown inclusions.
We denote by $\gamma$ the conductivity of $\Omega$.
We assume that $\gamma$ is given by a real valued essentially bounded function of $x\in\Omega$ having a positive lower bound.
The embedded inclusions are affected on the form of conductivity $\gamma$
$$\displaystyle
\gamma(x)=
\left\{
\begin{array}{ll}
\displaystyle
\gamma_0(x), & x\in\Omega\setminus D,\\
\\
\displaystyle
\gamma_0(x)+h(x), & x\in D,
\end{array}
\right.
$$
where $D$ is a non empty open subset of $\Omega$ with Lipschitz boundary, $\gamma_0$ denotes the background conductivity and $h$ is a deviation of the conductivity of $D$ from the background conductivity.
We assume that $\gamma_0$ is given by an essentially bounded functions of $x\in\Omega$ with a positive lower bound;
$h$ is given by an essentially bounded function of $x\in D$.

Given $f\in H^{1/2}(\partial\Omega)$ there exists a unique weak solution $u\in H^1(\Omega)$
of the the boundary value problem
$$
\left\{
\begin{array}{ll}
\nabla\cdot\gamma\nabla u=0, & x\in\Omega,\\
\\
\displaystyle
u=f, &  x\in\partial\Omega.
\end{array}
\right.
\tag {2.7}
$$
Define the Dirichlet-to-Neumann map by the formula
$$\begin{array}{ll}
\displaystyle
\Lambda_{\gamma}f=\gamma\frac{\partial u}{\partial\nu}, & x\in\partial\Omega.
\end{array}
$$

Given an arbitrary point $a\in\partial D$, we say that
$\gamma$ has a positive/negative jump from $\gamma_0$ at $a$ if
there exists a  positive number $\delta$ and an open ball
$B$ centered at $a$ such that $\gamma(x)-\gamma_0(x)\ge C$ a.e. $x\in B\cap D$
/$-(\gamma(x)-\gamma_0(x))\ge C$ a.e. $x\in B\cap D$.

$\quad$

{\bf\noindent Problem 2.2.}  Assume that the background conductivity $\gamma_0$ is known.
Extract information about the shape and location of $D$ from $\Lambda_{\gamma}$ provided
conductivity $\gamma$ has a positive or negative jump from $\gamma_0$ at each point $a\in\partial D$.

$\quad$

In \cite{IProbe} we introduced the original version of the probe method in the case when $\gamma_0$ is constant
and in \cite{IProbe2} the method has been extended to the case when $\gamma_0$ depends on $x\in\Omega$.
Here let us review the result in \cite{IProbe2} in terms of the new formulation of the probe method \cite{IProbeNew}.

We assume that $\gamma_0$ satisfies that $\gamma_0\in C^{0,1}(\overline\Omega)$ if $m=3$;
$\gamma_0\in C^{0,\theta}(\overline\Omega)$ with some $\theta\in\,]0,\,1]$ if $m=2$.

First we have to introduce a singular solution of the equation $\nabla\cdot\gamma_0\nabla v=0$ which is singular
at a given point $x\in\Omega$.  

Define
$$\displaystyle
G_0(y;x)=\frac{1}{\gamma_0(x)}G(y-x),
$$
where $G(\,\cdot\,)$ denotes the standard fundamental solution of the Laplace equation, that is,
$$\displaystyle
G(y-x)=
\left\{
\begin{array}{ll}
\displaystyle
\frac{1}{4\pi\vert y-x\vert}, & m=3,
\\
\\
\displaystyle
\frac{1}{2\pi}\log\frac{1}{\vert y-x\vert}, & m=2.
\end{array}
\right.
$$
Note that $G_0(y;x)$ satisfies, in the sense of distributions, the equation 
$$
\displaystyle
\nabla\cdot\gamma_0(x)\nabla G_0(y;x)+\delta(y-x)=0,
$$
where $\nabla=(\,\frac{\partial}{\partial y_1},\cdots,\frac{\partial}{\partial y_m}\,)^T$ and $\gamma_0(x)$
is given by freezing the original coefficient $\gamma_0$ at $y=x$.
Besides, for any finite cone $V$ with vertex at $x$ we have
$$\displaystyle
\int_V\,\vert\nabla G_0(y;x)\vert^2\,dy=\infty.
\tag {2.8}
$$
Using this singular solution, one can construct a singular solution of the equation $\nabla\cdot\gamma_0\nabla v=0$
which is an application of the Lax-Milgram theorem.

\proclaim{\noindent Lemma 2.1(\cite{IProbe2}).}  Given $x\in\Omega$
there exists a family $(G_x^0(\cdot))_{x\in\Omega}$
in $\cap_{1\le p<\frac{m}{m-2}}L^p(\Omega)$ such that
$(G_x^0(\,\cdot\,)-G_0(\,\cdot\,;x))_{x\in\Omega}$ is bounded
in $H^1_0(\Omega)$
and, for any $\varphi\in C_0^{\infty}(\Omega)$ we have
$$\displaystyle
\int_{\Omega}\gamma_0(y)\nabla G_x^0(y)\cdot\nabla\varphi(y)dy=\varphi(x).
$$
\endproclaim

{\it\noindent Sketch of Proof.}
We construct $G_x^0(\cdot)$ in the form
$$\begin{array}{ll}
\displaystyle
G_x^0(y)=G_0(y;x)+\epsilon(y;x), & y\in\Omega\setminus\{x\},
\end{array}
$$
where $\epsilon(\,\cdot\,;x)\in H^1_0(\Omega)$ is the unique weak solution of
$$\left\{
\begin{array}{ll}
\displaystyle
\nabla\cdot\gamma_0\nabla\epsilon(\cdot;x)=-f_x, & y\in\Omega,\\
\\
\displaystyle
\epsilon(\cdot;x)\vert_{\partial\Omega}=0, & y\in\partial\Omega
\end{array}
\right.
$$
and a functional $f_x$ is given by the formula
$$\begin{array}{ll}
\displaystyle
f_x(\varphi)=\int_{\Omega}
(\gamma_0(x)-\gamma_0(y))
\nabla G_0(y;x)\cdot\nabla\varphi(y)dy,
&
\displaystyle
\varphi\in C_0^{\infty}(\Omega).
\end{array}
$$

The assumption on $\gamma_0$ together with the concrete form of $G_0(y;x)$ ensures that $f_x$ is in $(H^1_0(\Omega))^*$
and bounded with respect to $x\in\Omega$.

\noindent
$\Box$

For the construction of $G_x^0$ for general $\gamma_0\in L^{\infty}(\Omega)$ one can refer \cite{GW}, however, 
for the present purpose we need only Lemma 2.1 for $\gamma_0\in C^{0,s}(\overline\Omega)$.

Now we state the existence of the needle sequence for a given needle.
\proclaim{\noindent Proposition 2.4 (\cite{IProbe2}).}  Let $B$ be an open ball with $B\cap\partial\Omega\not=\emptyset$.
Given $x\in\Omega$ and $\sigma\in N_x$
there exists a sequence $\{v_n\}$ of $H^1(\Omega)$ solutions of the equation 
$\nabla\cdot\gamma_0\nabla v=0$ in $\Omega$ with $\text{supp}\,(v_n\vert_{\partial\Omega})\,\subset
B\cap\partial\Omega$ such that $v_n\rightarrow G_x^0$
in $H^1_{\text{loc}}\,(\Omega\setminus\sigma)$.

\endproclaim

This is a direct consequence of the Runge approximation property for the equation
$\nabla\cdot\gamma_0\nabla v=0$ to the singular solution $G_x^0$
in $\Omega\setminus\sigma$.  The Runge approximation property is a consequence of the weak unique continuation
property of the equation $\nabla\cdot\gamma_0\nabla v=0$. When $m=3$, it is well known under the Lipschitz regularity
of $\gamma_0$ and in the case $m=2$ see \cite{AM}, therein the regularity of $\gamma_0$ 
is just essentially bounded.  These two properties are equivalent to each other as shown in \cite{La}.

Hereafter we denote $\Lambda_{\gamma}$ in the case when $D=\emptyset$ by $\Lambda_{\gamma_0}$.
First we introduce the indicator sequence corresponding to a given needle sequence.

{\bf\noindent Definition 2.8.}
Given $x\in\Omega$, $\sigma\in N_x$ and needle sequence $\xi=\{v_n\}$ for $(x,\sigma)$ define
the indicator sequence of this subsection by the formula
$$\begin{array}{ll}
\displaystyle
I(x,\sigma,\xi)_n=\int_{\partial\Omega}(\Lambda_{\gamma}-\Lambda_{\gamma_0})f_n\cdot f_n\,dS,
& n=1,2,\cdots,
\end{array}
$$
where $f_n(y)=v_n(y)$, $y\in\partial\Omega$.

Let $x\in\Omega\setminus\overline D$.
From Lemma 2.1 and the form of $\gamma$ we have $(\gamma-\gamma_0)\nabla G_x^0\in L^2(\Omega)$.
Thus one can find the weak solution $w=w_x$ in $H^1_0(\Omega)$ of the Dirichlet problem
$$\displaystyle
\left\{
\begin{array}{ll}
\nabla\cdot\gamma\nabla w=-\nabla\cdot(\gamma-\gamma_0)\nabla G_x^0, & y\in\Omega,\\
\\
\displaystyle
w=0, & y\in\partial\Omega.
\end{array}
\right.
\tag {2.9}
$$
We call the $w_x$ the {\it reflected solution} by $D$ and set $w=w_x$ to show the dependence on $x\in\Omega\setminus\overline D$.

By Alessandrini's identity we have
$$\displaystyle
\int_{\partial\Omega}(\Lambda_{\gamma}-\Lambda_{\gamma_0})f\cdot f\,dS
=\int_{\Omega}(\gamma-\gamma_0)\left(\vert \nabla v\vert^2+\nabla (u-v)\cdot\nabla v\right)\,dS,
\tag {2.10}
$$
where $u$ solves (2.7) and $v$ (2.7) with $\gamma=\gamma_0$.

This motivates us to introduce the indicator function of the probe method for Problem 2.2.

{\bf\noindent Definition 2.9.}
Define the indicator function $I(x;\gamma,\gamma_0)$ of independent variable $x\in\Omega\setminus\overline D$ by the formula
$$\displaystyle
I(x;\gamma,\gamma_0)=
\int_{\Omega}(\gamma-\gamma_0)\left(\vert\nabla G_x^0\vert^2+\nabla w_x\cdot\nabla G_x^0\right)\,dy.
$$

It follows from (2.9) that 
$$\displaystyle
I(x;\gamma,\gamma_0)
=\int_{\Omega}\left\{(\gamma-\gamma_0)\vert\nabla G_x^0\vert^2-\gamma\vert\nabla w_x\vert^2\right\}\,dy
\tag {2.11}
$$

Now we present the original version of the problem \cite{IProbe, IProbe2} in terms of the new formulation of the probe method
developed in \cite{IProbeNew} which is called the side A of the probe method.

\proclaim{\noindent Theorem 2.3.}

$\bullet$  (A.1)  Given $x\in\Omega\setminus\overline D$ let $\sigma\in N_x$ satisfy $\sigma\cap\overline{D}=\emptyset$.
Then, for all needle sequence $\xi=\{v_n\}$ for $(x,\sigma)$ the indicator sequence $\{I(x,\sigma,\xi)_n\}$
converges to the value of indicator function $I(\,\cdot\,;\gamma,\gamma_0)$ at $x$.

$\bullet$  (A.2)  For each $\epsilon>0$ we have $\sup_{\text{dist}\,(x,\overline D)>\epsilon}\,\vert I(x;\gamma,\gamma_0)\vert<\infty$.

$\bullet$  (A.3)  Given $a\in\partial D$  we have

$$\displaystyle
\lim_{x\rightarrow a}I(x;\gamma,\gamma_0)=
\left\{
\begin{array}{ll}
\infty, & \text{if $\gamma$ has a positive jump from $\gamma_0$ at $a$,}\\
\\
\displaystyle
-\infty, & \text{if $\gamma$ has a negative jump from $\gamma_0$ at $a$.}
\end{array}
\right.
$$

\endproclaim

{\it\noindent Sketch of Proof.}
Let $u=u_n$ solve (2.7) with $f=v_n$ on $\partial\Omega$.
The identity (2.10) yields the expression
$$\displaystyle
I(x,\sigma,\xi)_n=\int_{\Omega}\,\left\{(\gamma-\gamma_0)\vert \nabla v_n\vert^2-\gamma\vert\nabla(u_n-v_n)\vert^2\right\}\,dS.
\tag {2.12}
$$
The assertions (A.1) and (A.2) are an easy consequence of this expression and
the well-posedness of the direct problems;
the convergence property of $\nabla v_n$ to $G_x^0$ in $L^2(D)$;
the boundedness of the family $(G_x^0-G_0(\,\cdot\,;x))_{x\in\Omega}$
in $H^1_0(\Omega)$; the concrete expression of $G_0(\,\cdot\,;x)$.

The proof of assertion (A.3) is not trivial.  
From (2.11) we have
$$\displaystyle
I(x;\gamma,\gamma_0)\le \int_D h\vert\nabla G_x^0\vert^2\,dy.
\tag {2.13}
$$
On the other hand, from (2.9) we have
$$\displaystyle
\int_{\Omega}\left\{\gamma\vert\nabla w_x\vert^2+(\gamma-\gamma_0)\nabla G_x^0\cdot\nabla w_x\right\}\,dy=0.
$$
Rewriting the integrand on the left-handside, we have
$$\displaystyle
\int_{\Omega}\,\gamma\left\vert
\nabla w_x+\frac{\gamma-\gamma_0}{2\gamma}\nabla G_x^0\right\vert^2\,dy
=\frac{1}{4}\int_{\Omega}\frac{(\gamma-\gamma_0)^2}{\gamma}\,\vert\nabla G_x^0\vert^2\,dy.
$$
This yields
$$\displaystyle
\frac{1}{2}\int_{\Omega}\,\gamma\vert\nabla w_x\vert^2\,dy
\le
\frac{1}{2}\int_{\Omega}\frac{(\gamma-\gamma_0)^2}{\gamma}\,\vert\nabla G_x^0\vert^2\,dy,
$$
that is
$$\displaystyle
\int_{\Omega}\,\gamma\vert\nabla w_x\vert^2\,dy
\le
\int_{\Omega}\frac{(\gamma-\gamma_0)^2}{\gamma}\,\vert\nabla G_x^0\vert^2\,dy.
\tag {2.14}
$$
Since 
$$\displaystyle
(\gamma-\gamma_0)-\frac{(\gamma-\gamma_0)^2}{\gamma}
=\frac{\gamma_0(\gamma-\gamma_0)}{\gamma},
$$
from (2.11) and (2.14), we obtain
$$\displaystyle
I(x;\gamma,\gamma_0)
\ge \int_{D}\frac{\gamma_0h}{\gamma}\vert\nabla G_x^0\vert^2\,dy.
\tag {2.15}
$$
Thus if $\gamma$ has a negative jump from $\gamma_0$ at $a$, from (2.13) it suffices to prove
$$\displaystyle
\lim_{x\rightarrow a}\int_{D}h\vert\nabla G_x^0\vert^2\,dy=-\infty.
$$
However, by Lemma 2.1, we see that this is a consequence of (2.8) with a suitable choice of a finite cone $V$ inside $D$ with vertex at $a$.
Note that such a choice is always possible because of the Lipschitz regularity of $\partial D$.
If $\gamma$ has a positive jump from $\gamma_0$ at $a$, use (2.15).

\noindent
$\Box$

The procedure of finding the unknown inclusion based on Theorem 2.3, 
that is, the side A of the probe method is as follows.
First the needle insertion path is set in advance, and the needle is virtually inserted little by little along it.
The calculated value of the indicator function by using (A.1) will not blow uo  
unless there is an unknown inclusion on its path.  This is ensured by (A.2).
However, if there is an inclusion on the path, the value of  indicator function will blow up soon.  This is 
ensured by (A.3).  Then we can conclude that the tip of the needle at that time is on the surface of the inclusion.

The proof of Proposition 2.3 works for the needle sequence of solutions of the equation $\nabla\cdot\gamma_0\nabla v=0$ without 
any essential change provided, say $\gamma_0$ is smooth in $\Omega$.
Thus the needle sequence of this subsection has the same property as that of the Laplace equation.

Then, we see that Theorem 4.1 in \cite{IProbeNew} which is the case when $\gamma_0$ is constant 
can be extended as a variable $\gamma_0$.  More precisely we have the following result.

\proclaim{\noindent Theorem 2.4.}  Assume that $\Omega\setminus\overline D$  is connected
and that $\gamma$ satisfies the positive/negative global jump condition from $\gamma_0$, that is,
there exists a positive constant $C$ such that $\pm(\gamma(x)-\gamma_0(x))>C$ a.e. $x\in D$.

A point $x\in\Omega$ belongs to $\Omega\setminus\overline D$ if and only if
there exists a needle $\sigma\in N_x$ and needle sequence $\xi=\{v_n\}$ for $(x,\sigma)$ such that the indicator
sequence $\{I(x,\sigma,\xi)_n\}$ is bounded.  Moreover given $x\in\Omega$ and needle $\sigma\in N_x$ we have that

$\bullet$  if $x\in\Omega\setminus\overline D$ and $\sigma\cap D\not=\emptyset$ or $x\in\overline D$, then for any needle sequence
$\xi=\{v_n\}$ for $(x,\sigma)$ we have $\lim_{n\rightarrow\infty}I(x,\sigma,\xi)_n=\pm\infty$ provided
$\gamma$ satisfies the positive/negative global jump condition from $\gamma_0$

\endproclaim

{\it\noindent Sketch of Proof.}
From (2.12) we have immediately
$$\displaystyle
I(x,\sigma,\xi)_n\le\int_Dh\vert\nabla v_n\vert^2\,dy.
\tag {2.16}
$$
Applying the same argument done in deriving (2.15) to the second term of the integrand of the right-hand side on (2.10), we have 
$$\displaystyle
I(x,\sigma,\xi)_n\ge \int_D\frac{\gamma_0h}{\gamma}\vert\nabla v_n\vert^2\,dy.
\tag {2.17}
$$
Then, everything follows from inequalities (2.16) and (2.17) combining with the properties of the needle sequences same as those
of Proposition 2.3.

\noindent
$\Box$

\begin{figure}[htbp]
\begin{center}
\epsfxsize=4.8cm
\epsfysize=4.8cm
\epsfbox{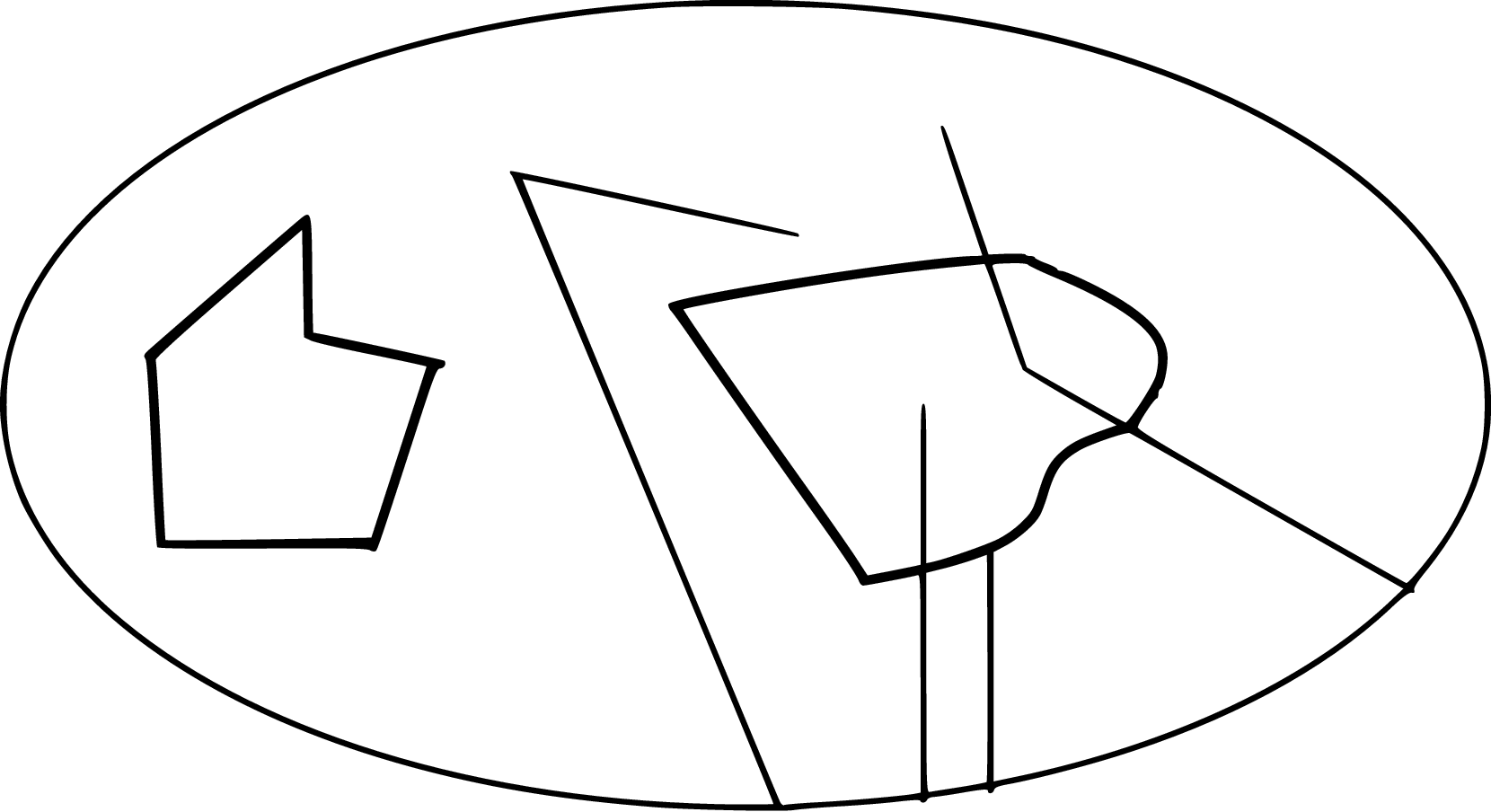}
\caption{An illustration of needles $\sigma\in N_x$ satisfying $x\in\Omega\setminus\overline D$ with $\sigma\cap\overline D=\emptyset$;
$x\in\Omega\setminus\overline D$ with $\sigma\cap D\not=\emptyset$;
$x\in\partial D$; $x\in D$.
All the connected components of $D$ are indicated by the area enclosed by the thick curves.}
\label{fig3}
\end{center}
\end{figure}

See Figure \ref{fig3} for an illustration of typical needles.  Note that the case when $x\in\Omega\setminus\overline D$
and $\sigma\in N_x$ satisfies both $\sigma\cap\partial D\not=\emptyset$ and $\sigma\cap D=\emptyset$ is excluded.
This is the {\it grazing case}.  The reason is that: the behaviour of $\Vert\nabla v_n\Vert_{L^2(V\cap\Omega)}$ for 
a given finite cone $V$ with vertex at an arbitrary point $z=\sigma(t)$ with $t\in\,]0,\,1[$ is not clear.
See Proposition 2.3.

\vspace{0cm}

As a corollary we have the complete characterization of $\overline{D}$ in terms of only the indicator sequence:
$\overline D$ is the set of all points $x$ such that for any $\sigma\in N_x$ and any needle sequence $\xi=\{v_n\}$
for $(x,\sigma)$ indicator sequence $\{I(x,\sigma,\xi)_n\}$ blows up as $n\rightarrow\infty$.

It should be pointed out that the behaviour of $\{I(x,\sigma,\xi)\}$ for general $x\in\overline D$ is not clear
if $\gamma$ does not satisfy the positive/negative global jump condition from $\gamma_0$.  However, if $\pm(\gamma-\gamma_0)$
has a positive lower bound in a neighbourhood of $\sigma\cap\overline D$, then we have
$\lim_{n\rightarrow\infty}I(x,\sigma,\xi)_n=\pm\infty$.  In other cases we do not know the behaviour.

The system (2.16) and (2.17) can be derived from a more general system described below which has been established in \cite{I0}.

\proclaim{\noindent Proposition 2.5.}
Let $\gamma_j, j =1, 2, $ be two conductivities.
Let $v_j\in H^1(\Omega)$ denote the weak solution of the Dirichlet problem
$$\left\{
\begin{array}{ll}
\displaystyle
\nabla\cdot\gamma_j\nabla v=0, &  x\in\Omega,\\
\\
\displaystyle
v=f, & x\in\partial\Omega.
\end{array}
\right.
$$
It holds that
$$
\displaystyle
\int_{\Omega}\{\gamma_1^{-1}-\gamma_2^{-1}\}\gamma_1\nabla v_1\cdot\gamma_1
\nabla v_1dx
\le <(\Lambda_{\gamma_2}-\Lambda_{\gamma_1})f\,,f>
\tag {2.18}
$$
and
$$
\displaystyle
<(\Lambda_{\gamma_2}-\Lambda_{\gamma_1})f\,,f>
\le\int_{\Omega}(\gamma_2-\gamma_1)\nabla v_1\cdot\nabla v_1 dx,
\tag {2.19}
$$
where $<\,,\,>$ denotes the dual pairing between $H^{-1/2}(\partial\Omega)$ and $H^{1/2}(\partial\Omega)$.

\endproclaim

In particular, the lower bound (2.18) is derived from a corollary of the best possible version for an early estimate of \cite{KSS}.
More precisely, in \cite{I0} the author proved also
$$
\displaystyle
\int_{\Omega}\{\gamma_1^{-1}-\gamma_2^{-1}\}\gamma_1\nabla v_1
\cdot\gamma_1\nabla v_1 dx
\le \frac{<\Lambda_{\gamma_1}f, f>}{<\Lambda_{\gamma_2}f, f>}
<(\Lambda_{\gamma_2}-\Lambda_{\gamma_1})f, f>
\tag {2.20}
$$
Since we have
$$\displaystyle
\frac{<\Lambda_{\gamma_1}f, f>}{<\Lambda_{\gamma_2}f, f>}
<(\Lambda_{\gamma_2}-\Lambda_{\gamma_1})f, f>
\le <(\Lambda_{\gamma_2}-\Lambda_{\gamma_1})f\,,f>,
$$
(2.18) follows from (2.20).
The proof of (2.20) is not trivial.  Note that if $\gamma_1$ and $ \gamma_2$ are constant,
inequalities (2.19) and (2.20) become equalities.

\subsubsection{Outside to inside argument}

In this subsection we describe a procedure introduced in \cite{IProbe2} which enables us to compute
the Dirichlet-to-Neumann map associated with the conductivity {\it inside} the inclusion from
a subset of the set 
$$\displaystyle
\left\{<(\Lambda_{\gamma}-\Lambda_{\gamma_0})f,f>\,\vert\,f\in H^{1/2}(\partial\Omega),\,\text{supp}\,f\subset\,B\cap\partial\Omega\right\}
\tag {2.21}
$$
for a fixed $B$ provided $D$ and $\gamma_0$ are known and $\gamma(x)=\gamma_0(x)$ a.e. $x\in\Omega\setminus\overline D$.
It is not assumed that $\gamma$ in $D$ is known.
This is another contribution done in \cite{IProbe2} not covered in \cite{IProbe}.

We start with two assumptions on $\gamma$ and $D$.

$\bullet$  $\gamma\in L^{\infty}(\Omega)$ satisfy $\text{ess.inf}_{x\in\Omega}\gamma(x)>0$.

$\bullet$  $D$ with a Lipschiz boundary such that $\overline{D}\subset\Omega$ and  $\Omega\setminus\overline D$ is connected.

{\bf\noindent Definition 2.10.}
Define the bounded linear operator 
$\Lambda_{\gamma\vert_D}(D):H^{1/2}(\partial D)\rightarrow H^{-1/2}(\partial D)$ by the formula
$$\displaystyle
<\Lambda_{\gamma\vert_D}(D)f,g>=\int_D\gamma\nabla u\cdot\nabla\varphi\,dx,
$$
where $g\in H^{1/2}(\partial D)$, $\varphi$ is any $H^1(D)$ function with $\varphi\vert_{\partial D}=g$, the function 
$u\in H^1(D)$ is the weak solution of the Dirichlet problem
$$
\left\{
\begin{array}{ll}
\displaystyle
\nabla\cdot\gamma\nabla u=0, & x\in D,\\
\\
\displaystyle
u=f, & x\in \partial D.
\end{array}
\right.
$$
We call $\Lambda_{\gamma\vert_D}(D)$ the Dirichlet-to-Neumann map associated with the conductivity insiside the inclusion.

{\bf\noindent Definition 2.11.}
Define the bounded linear operator 
$\Lambda_{\gamma\vert\Omega\setminus\overline D}(D):H^{1/2}(\partial D)\rightarrow H^{-1/2}(\partial D)$ by the formula
$$\displaystyle
<\Lambda_{\gamma\vert\Omega\setminus\overline D}(D)f,g>=-\int_{\Omega\setminus\overline D}\gamma\nabla u\cdot\nabla\varphi\,dx,
$$
where $g\in H^{1/2}(\partial D)$, $\varphi$ is any $H^1(\Omega\setminus\overline D)$ function with $\varphi\vert_{\partial D}=g$
and $\varphi\vert_{\partial\Omega}=0$, the function 
$u\in H^1(\Omega\setminus\overline D)$ is the weak solution of the Dirichlet problem
$$
\left\{
\begin{array}{ll}
\displaystyle
\nabla\cdot\gamma \nabla u=0, & x\in \Omega\setminus\overline D,\\
\\
\displaystyle
u=f, & x\in \partial D,\\
\\
\displaystyle
u=0, & x\in\partial\Omega.
\end{array}
\right.
$$
We call $\Lambda_{\gamma\vert\Omega\setminus\overline D}(D)$ the Dirichlet-to-Neumann map associated with the conductivity outside the inclusion.
Note that $\Lambda_{\gamma\vert\Omega\setminus\overline D}(D)$ can be calculated if both $\gamma\vert_{\Omega\setminus\overline D}$ and $D$ are known.

{\bf\noindent Definition 2.12.}
For each $f\in H^{-1/2}(\partial D)$ define
$$\begin{array}{ll}
\displaystyle
T_f(\varphi)=<f,\,\varphi\vert_{\partial D}>,
&
\displaystyle
 \varphi\in H^1_0(\Omega).
 \end{array}
$$
From the trace theorem we know $T_f\in H^{-1}(\Omega)\equiv (H^1_0(\Omega))^*$
and thus there exists a unique weak solution $u_f$ in $H^1_0(\Omega)$ of
$$
\left\{
\begin{array}{ll}
\displaystyle
\nabla\cdot\gamma\nabla u_f=-T_f, &  x\in\Omega,\\
\\
\displaystyle
u_f=0, & x\in\partial\Omega.
\end{array}
\right.
$$
Set
$$\displaystyle
G_{\gamma}(D)f=u_f\vert_{\partial D}.
$$
It is easy to see that $G_{\gamma}(D)$ is a bounded linear operator
from $H^{-1/2}(\partial D)$ to $H^{1/2}(\partial D)$.

The following fact describes a relationship between three operators.

\proclaim{\noindent Proposition 2.6.}  
$G_{\gamma}(D)$ is bijective and the formula
$$\displaystyle
\Lambda_{\gamma\vert D}(D)-\Lambda_{\gamma\vert\Omega\setminus\overline D}(D)=(G_{\gamma}(D))^{-1},
\tag {2.22}
$$
is valid.
\endproclaim
{\it\noindent Sketch of Proof.}
First we have: $\Lambda_{\gamma\vert D}(D)-\Lambda_{\gamma\vert\Omega\setminus\overline D}(D)$ is injective;
the formula $(\Lambda_{\gamma\vert D}(D)-\Lambda_{\gamma\vert\Omega\setminus\overline D}(D))G_{\gamma}(D)f=f, \forall f\in H^{-1/2}(\partial D)$,
is valid.  Thus one can conclude that
$\Lambda_{\gamma\vert D}(D)-\Lambda_{\gamma\vert\Omega\setminus\overline D}(D)$ is bijective
and thus we know $G_{\gamma}(D)=(\Lambda_{\gamma\vert D}(D)-\Lambda_{\gamma\vert\Omega\setminus\overline D}(D))^{-1}$.
Therefore $G_{\gamma}(D)$ is bijective, too and hence the desired formula
is valid.  

\noindent
$\Box$

Some remarks are in order.

$\bullet$  Any regularity of $\gamma$ on $\Omega$ is not assumed.

$\bullet$  Nachman ((6.2) in \cite{N2}) proved a fact
corresponding to formula (2.22)
in the case where $\Omega\subset\Bbb R^2$ and $\gamma\in W^{2,p}(\Omega)$
with $p>1$.  However, his formulation and proof are not based on the weak formulation.

Here we impose an additional assumption.

$\bullet$  $\gamma(x)=\gamma_0(x)$ for almost all $x\in\Omega\setminus D$ with $\gamma_0\in L^{\infty}(\Omega)$
having $\text{ess.inf}_{x\in\Omega}\,\gamma_0(x)>0$.

{\bf\noindent Definition 2.13.}
Given $F\in\,H^{-1}(\Omega)$ with $\displaystyle\text{supp}\,F\subset\Omega\setminus\overline D$,
define
$$\left\{
\begin{array}{ll}
\displaystyle
\mbox{\boldmath $G$}_{\gamma_0}F=u_0, & x\in\Omega
\\
\\
\displaystyle
\mbox{\boldmath $G$}_{\gamma}
F=\mbox{\boldmath $G$}_{\gamma_0}F+w, & x\in\Omega
\end{array}
\right.
$$
where $u_0\in H^1_0(\Omega)$ is the weak solution of
$$\left\{
\begin{array}{ll}
\displaystyle
\nabla\cdot\gamma_0\nabla u_0=-F, & x\in\Omega,\\
\\
\displaystyle
u_0=0, & x\in\partial\Omega
\end{array}
\right.
$$
and $w\in H^1_0(\Omega)$ is the weak solution of
$$\left\{
\begin{array}{ll}
\displaystyle
\nabla\cdot\gamma\nabla w=-\nabla\cdot\chi_D(\gamma-\gamma_0)\nabla u_0,
& x\in\Omega,
\\
\\
\displaystyle
w=0, & x\in\partial\Omega.
\end{array}
\right.
$$

It is a simple matter to check the following formula which connects $u_f$ and $\mbox{\boldmath $G$}_{\gamma}$.

\proclaim{\noindent Proposition 2.7.}  The formula
$$\displaystyle
F(u_f)=<f, \mbox{\boldmath $G$}_{\gamma}
F\vert_{\partial D}>,
\tag {2.23}
$$
is valid.
\endproclaim

The next formula needs the Runge approximation property of the equation $\nabla\cdot\gamma_0\nabla v=0$ in $\Omega$
and thus, at least we assume

$\bullet$  $\gamma_0\in C^{0,1}(\overline\Omega)$ ($n=3$).

Given $F, H\in\,H^{-1}(\Omega)$ with $\displaystyle\text{supp}\,F\subset\Omega\setminus\overline D$
one has two sequences $\{u_n\}$, $\{v_n\}$
of $H^1(\Omega)$ solutions of the equation $\nabla\cdot\gamma_0\nabla v=0$ in $\Omega$
with $\text{supp}\,(u_n\vert_{\partial\Omega}), \text{supp}\,(v_n\vert_{\partial\Omega})\subset\Gamma\equiv B\cap\partial\Omega$
such that $u_n\longrightarrow \mbox{\boldmath $G$}_{\gamma_0}F$ 
and $v_n\longrightarrow \mbox{\boldmath $G$}_{\gamma_0}H$ in $H^1(D)$.

Then we have the following formula.

\proclaim{\noindent Proposition 2.8.}  The formula
$$\begin{array}{l}
\displaystyle
\,\,\,\,\,\,
-H(\mbox{\boldmath $G$}_{\gamma}
F-\mbox{\boldmath $G$}_{\gamma_0}F)
\\
\\
\displaystyle
=
\frac{1}{4}\lim_{n\longrightarrow\infty}
\{<(\Lambda_{\gamma}-\Lambda_{\gamma_0})f_n^+,
f_n^+> -<(\Lambda_{\gamma}-\Lambda_{\gamma_0})f_n^-, f_n^->\},
\end{array}
\tag {2.24}
$$
is valid, where $f_n^+=u_n+v_n$ and $f_n^-=u_n-v_n$ on $\partial\Omega$.
\endproclaim

Proposition 2.8 says that, given $F$, $H$ from the set
$$\displaystyle
\left\{
<(\Lambda_{\gamma}-\Lambda_{\gamma_0})f,f>
\,\vert\, f=f_n^+, f_n^-\right\}
\tag {2.25}
$$
which is a subset of (2.21), one can recover $H(\mbox{\boldmath $G$}_{\gamma}F-\mbox{\boldmath $G$}_{\gamma_0}F)$.

Summing those up, we have the reconstruction procedure of the Dirichlet-to-Neumann map inside inclusion $\Lambda^-$
provided both $\gamma_0$ and $D$ are known and that $\gamma(x)=\gamma_0(x)$ for almost all $x\in\Omega\setminus\overline D$.

{\bf\noindent Step 1.}
Given $F, H$ being elements of $H^{-1}(\Omega)$ with $\text{supp}\,F,\,\text{supp}\,H\subset\Omega\setminus\overline D$
choose two sequences $\{u_n\}, \{v_n\}$ of $H^1(\Omega)$ solutions of equation
$\nabla\cdot\gamma_0\nabla v=0$ in $\Omega$ with
$\text{supp}\,(u_n\vert_{\partial\Omega})\subset\Gamma, \text{supp}\,(v_n\vert_{\partial\Omega})\subset\Gamma$
such that $u_n\longrightarrow \mbox{\boldmath $G$}_{\gamma_0}F, v_n\longrightarrow \mbox{\boldmath $G$}_{\gamma_0}H$ in $H^1(D)$.

{\bf\noindent Step 2.}
Calculate $f_n^+=u_n+v_n,f_n^-=u_n-v_n$ on $\partial\Omega$.

{\bf\noindent Step 3.}
Use formula (2.24) to recover $H(\mbox{\boldmath $G$}_{\gamma}F-\mbox{\boldmath $G$}_{\gamma_0}F)$ from the data (2.25).

{\bf\noindent Step 4.}
Calculate $\mbox{\boldmath $G$}_{\gamma}
F-\mbox{\boldmath $G$}_{\gamma_0}F$ in $\Omega\setminus\overline D$
from the set
$$
\{H(\mbox{\boldmath $G$}_{\gamma}
F-\mbox{\boldmath $G$}_{\gamma_0}F)\,\vert\,
\text{$H\in H^{-1}(\Omega)$ with $\text{supp}\,H\subset\Omega\setminus\overline D$}\}.
$$

{\bf\noindent Step 5.}
Calculate $\mbox{\boldmath $G$}_{\gamma}
F\vert_{\partial D}$.

{\bf\noindent Step 6.}
Given $f\in H^{-1/2}(\partial D)$ use formula (2.23) to recover $F(u_f)$ from $\mbox{\boldmath $G$}_{\gamma}
F\vert_{\partial D}$.

{\bf\noindent Step 7.}
Calculate $u_f$ in $\Omega\setminus\overline D$ from the set
$$
\{F(u_f)\,\vert\,\text{$F\in H^{-1}(\Omega)$
 with $\text{supp}\,F\subset\Omega\setminus\overline D$}\}.
$$

{\bf\noindent Step 8.}
Calculate $G_{\gamma}(D)f=u_f\vert_{\partial D}$.

{\bf\noindent Step 9.}
Calculate $G_{\gamma}(D)$ from the set
$$
\{G_{\gamma}(D)f\,\vert\, f\in H^{-1/2}(\partial D)\}.
$$

{\bf\noindent Step 10.}
Use formula (2.22) to recover $\Lambda_{\gamma\vert D}(D)$ from $G_{\gamma}(D)$.  
Note that $\Lambda_{\gamma\vert\Omega\setminus\overline D}(D)$ is known.

As a corollary we have a {\it reconstruction procedure} of $\gamma$ in $D$
together with $D$ in three dimensions under suitable regularity 
on $\gamma_0$ and $\gamma$ in $D$ and the jump condition on $\gamma$ from $\gamma_0$.  It is a consequence
of Theorem 2.3 and Nachman's reconstruction formula \cite{N1} together with outside to inside argument. 

Note that if one wishes to recover only the value of the jump of $\gamma$ from $\gamma_0$ at an arbitrary 
point $a\in\partial D$ provided $D$ is known, one can make use of a formula in \cite{IN2}
without outside to inside argument.
Therein, under the condition $\gamma_0=1$, it is shown that from the asymptotic behaviour of the indicator function
$I(x;\gamma,\gamma_0)$ as $x\rightarrow a$ one can find the value of the jump.  Its origin goes
back to a corresponding formula for the Schr\"odinger type equation established in the early stage of 
the probe method \cite{IProbe}.

Finally to clarify the power of the outside to inside argument we consider a simple situation.
Assume that $\Omega=O_0$ has an onion like layered structure, that is, there is a sequence
of open subsets $O_j$, $j=1,\cdots,l$ of $\Bbb R^2$ with smooth boundaries such that
$\overline{O_{j+1}}\subset O_j$ and $O_j\setminus\overline O_{j+1}$ is connected for $j=0,\cdots, l-1$.
Moreover, for each  $j=0,\cdots,l-1$ the conductivity of $\gamma$ on $\overline{O_j}\setminus O_{j+1}$ 
takes a constant $\gamma_j$.


\proclaim{\noindent Corollary 2.1.} Assume that for each $j=1,\cdots,l-1$ $\gamma_j\not=\gamma_{j+1}$.
Then, one can reconstruct $O_j$, $j=1,\cdots, l$ together with $\gamma$ itself from $\Lambda_{\gamma}$ acting on
a subset of $f\in H^{1/2}(\partial\Omega)$ with $\text{supp}\,f\subset\Gamma$.

\endproclaim
{\it\noindent Proof.}
Apply the outside to inside argument from $O_{j}$ to $O_{j+1}$, $j=0,\cdots,l-1$ together with
the reconstruction formula of the conductivity at the boundary \cite{Br} and the probe method.

\noindent
$\Box$

This is a kind of layer stripping approach, see \cite{S2} for the approach based on the Riccati equation.
Our approach is based on the Runge approximation property.
Note that in two dimensions the reconstruction of $\gamma$ without any regularity nor structure like above 
has been established in \cite{AP}.

The outside to inside argument has been applied also to an inverse problem 
arising from ocean acoustics as an important step, see \cite{IMN2}.

Recently in the {\it uniqueness issue} of inverse boundary value problems the outside to inside argument 
found its applications out, see \cite{CHN} and \cite{HNZ}.



\subsubsection{Extension to an elastic body: the multi probe method}

A Calder\'on type problem for the system of equations in linear elasticity \cite{Gu} was raised by the author in \cite{Ielasticity0}
motivated by understanding the data coming from the material hardness testing machine.
The problem formulated therein is as follows.

Let $\Omega$ be a bounded domain of $\Bbb R^3$ with smooth boundary.
We consider the domain $\Omega$ as an isotropic elastic body with Lam\'e parameters $\lambda$ and $\mu$.
We assume that the both parameters belong to $L^{\infty}(\Omega)$ and satisfy
$\text{ess.inf}_{x\in\Omega}\,\mu(x)>0$ and $\text{ess.inf}_{x\in\Omega}(3\lambda(x)+2\mu(x))>0$.
The governing equation of the displacement field $\mbox{\boldmath $u$}\in\,H^1(\Omega)^{n}$ is given by
$$\begin{array}{lll}
\displaystyle
\sum_{j=1}^3\,\frac{\partial}{\partial x_j}\,
\left(C_{\lambda,\mu}\,\text{Sym}\,\nabla\mbox{\boldmath $u$}\right)_{ij}=0,
& i=1,2,3, 
& x\in\Omega,
\end{array}
\tag {2.26}
$$
where the symbol $\text{Sym}\,A$ for a matrix $A$ denotes
its symmetric part and the $C_{\lambda,\mu}$ the elasticity tensor field acting on
the set of all symmetric matrix $B$ in such a way that
$$
\displaystyle
C_{\lambda,\mu}B=\lambda\,(\text{Trace}\,B)\,I_3+2\mu\,B.
$$

The associated Dirichlet-to-Neumann map $\Lambda_{\lambda,\mu}$ of $H^{1/2}(\partial\Omega)^n$
into its dual space is defined by
$$\begin{array}{ll}
\displaystyle
\Lambda_{\lambda,\mu}\mbox{\boldmath $f$}
=\left(C_{\lambda,\mu}\text{Sym}\,\nabla\mbox{\boldmath $u$}\right)\mbox{\boldmath $\nu$}=
\lambda\,(\nabla\cdot\mbox{\boldmath $u$})\mbox{\boldmath $\nu$}+2\mu\,(\text{sym}\,\nabla\mbox{\boldmath $u$})\mbox{\boldmath $\nu$},
& x\in\partial\Omega,
\end{array}
$$
where $\mbox{\boldmath $f$}\in H^{1/2}(\partial\Omega)^n$ and $\mbox{\boldmath $u$}$ solves (2.26) under the inhomogeneous 
Dirichlet boundary condition
$$\begin{array}{ll}
\displaystyle
\mbox{\boldmath $u$}=\mbox{\boldmath $f$}, & x\in\partial\Omega.
\end{array}
$$

In \cite{Ielasticity0} the author formulated the following problem.

$\quad$

{\bf\noindent Problem 2.3.}
Can one uniquely determine both $\lambda$ and $\mu$ from $\Lambda_{\lambda,\nu}$?

$\quad$

The main result in \cite{Ielasticity0} is:  the injectivity of the Frech\'et derivative of the map $(\lambda,\mu)\mapsto \Lambda_{\lambda,\mu}$
at an arbitrary fixed pair of constant Lam\'e parameters together with its {\it explicit left inverse}.

Thus the linearized problem is solvable at a pair of constant Lam\'e parameters.
This is an application of the Calder\'on method \cite{Cal}.  However,  like the Calder\'on problem where the range
of the Frech\'et derivative is not closed \cite{IkNaSurvey}, one expect that the usual inverse function theorem
in infinite dimensional spaces can not be applied to Problem 2.3 near the constant pair.
As an evidence, Problem 2.3 itself still unsolved except for some partial results \cite{ER} and \cite{NU1}.

In this section we present a result from \cite{IProbe2} for an inverse obstacle problem
which is a special, however, important version of Problem 2.3
and can be considered as an elastic body version of Problem 2.2.

We assume that domain $\Omega$ contains an inclusion denoted by $D$ which is an open subset of $\Bbb R^3$ with Lipschitz boundary
and satisfies $\overline D\subset\Omega$ and that $\Omega\setminus\overline D$ is connected.

We assume that the values of $\lambda$ and $\mu$ takes the form
$$
\displaystyle
(\lambda(x),\mu(x))=
\left\{
\begin{array}{ll}
\displaystyle
(\lambda_0(x),\mu_0(x)), & \text{if $x\in \Omega\setminus D$,}\\
\\
\displaystyle
(\lambda_0(x)+l(x), \mu_0(x)+m(x)), & \text{if $x\in D$,}
\end{array}
\right.
$$
where both $\lambda_0\in C^2(\overline\Omega)$ and $\mu_0\in C^3(\overline\Omega)$; both $l$ and $m$ belong
to $C^0(\overline D)$ (for simplicity of description).
Besides,  we assume that the Lam\'e parameters change across $\partial D$ in the following sense:
$$\displaystyle
\forall a\in\partial D\,\,(\lambda_0(a),\mu_0(a))\not=(\lambda_0(a)+l(a), \mu_0(a)+m(a)),
$$
that is,
$$
\displaystyle
\forall a\in\partial D\,\,(l(a),m(a))\not=(0,0).
\tag {2.27}
$$
Note that we never assume that $l(a)m(a)>0$ for all $a\in\partial D$ which is an easier case.
In short, the assumption (2.27) means that one of $\lambda$ and $\mu$ has a jump across $\partial D$.
This is a case not appeared in Problem 2.2.
Under assumptions listed above the author obtained the following result.

\proclaim{\noindent Theorem 2.5(\cite{IProbe2}).}
Let $B$ be an open ball with $B\cap\partial\Omega\not=\emptyset$.
There exists a subset ${\cal D}$ of $H^{1/2}(\partial D)^3$ 
such that $\partial D$ can be reconstructed from the set
$$\displaystyle
\left\{\,\int_{\partial\Omega}(\Lambda_{\lambda,\mu}-\Lambda_{\lambda_0,\mu_0}\,)\mbox{\boldmath $f$}\cdot
\mbox{\boldmath $f$}\,dS\,\vert\,\mbox{\boldmath $f$}\in\cal{D}\,\right\}.
$$
The ${\cal D}$ depends only on $\lambda_0$, $\mu_0$, $\Omega$ and $B$, and the supports of all
the members in ${\cal D}$ are contained in $B\cap\partial\Omega$.
\endproclaim

{\it\noindent Sketch of Proof.}
This should be considered as an application of the probe method of Side A.  
The singular solutions to generate needle sequences are of two types
and constructed as Lemma 2.1.
More precisely, let $x\in\Omega$.  In what follows we set $G(y-x)=\frac{1}{4\pi\vert y-x\vert}$.

(a)  One can construct $\mbox{\boldmath $u$}_x^0\in H^1_{\text{loc}}\,(\Omega\setminus\{x\})^3$ such that
$$\begin{array}{ll}
\displaystyle
\mbox{\boldmath $u$}_x^0(y)
\sim
\nabla G(y-x)
-\frac{G(y-x)}{\lambda_0(x)+2\mu_0(x)}
\,\left\{I_3-\frac{y-x}{\vert y-x\vert}\otimes\frac{y-x}{\vert y-x\vert}\,\right\}\,\nabla\mu_0(x), & y\in\Omega,
\end{array}
$$
where the relation $\sim$ means that the function on the left-hand side is given by the function on the right-hand side
plus a function in $H^1(\Omega)^3$ depending on $x\in\Omega$ and whose $H^1(\Omega)^3$-norm is bounded with respect
to $x$.

(b)  One can construct $\mbox{\boldmath $u$}_x^j\in H^1_{\text{loc}}\,(\Omega\setminus\{x\})^3$, $j=1,2,3$ 
such that
$$\begin{array}{ll}
\displaystyle
\mbox{\boldmath $u$}_x^j(y)
\sim \mbox{\boldmath $E$}_0(y-x)\mbox{\boldmath $e$}_j, & y\in\Omega,
\end{array}
$$
where the meaning of relation $\sim$ is same as (a) and $\mbox{\boldmath $E$}_0(y-x)$ is the Kelvin matrix
\cite{Gu} and its form is
$$\begin{array}{ll}
\displaystyle
\mbox{\boldmath $E$}_0(y-x)
&
\displaystyle
=\frac{1}{2}
\left(\frac{1}{\mu_0(x)}+\frac{1}{\lambda_0(x)+2\mu_0(x)}\,\right)G(y-x)\,I_3
\\
\\
\displaystyle
&
\displaystyle
\,\,\,
+\frac{1}{2}
\left(\frac{1}{\mu_0(x)}-\frac{1}{\lambda_0(x)+2\mu_0(x)}\,\right)G(y-x)\,\frac{y-x}{\vert y-x\vert}\otimes
\frac{y-x}{\vert y-x\vert}.
\end{array}
$$

Given a needle $\sigma\in N_x$, applying the Runge approximation property for the equation (2.26) which is a consequence
of the weak unique continuation property, e.g., \cite{AITY} to the singular solutions $\mbox{\boldmath $u$}_x^0$ 
and $\mbox{\boldmath $u$}_x^j$, $j=1,2,3$, we know the existence of four types of needle sequences,
$\xi^0=\{\mbox{\boldmath $v$}_m^0\}$, $\xi^j=\{\mbox{\boldmath $v$}_m^j\}$, $j=1,2,3$ in the sense:
$\mbox{\boldmath $v$}_m^j\rightarrow \mbox{\boldmath $u$}_x^j$ in $H^1_{\text{loc}}(\Omega\setminus\sigma)^3$, $j=0,1,2,3$
as $m\rightarrow\infty$.
Taking their traces onto $\partial\Omega$, we have
four members in ${\cal D}$, denoted by $\mbox{\boldmath $f$}_{m,x,\sigma}$, $\mbox{\boldmath $g$}_{m,x,\sigma}^j$, $j=1,2,3$,
respectively.  Note that their supports are in $B\cap\partial\Omega$ as well as Proposition 2.4 for the equation 
$\nabla\cdot\gamma_0\nabla v=0$.

Define the two indicator sequences
$$
\left\{
\begin{array}{l}
\displaystyle
I(x,\sigma,\xi^0)_m=\int_{\partial\Omega}\,(\Lambda_{\lambda,\mu}-\Lambda_{\lambda_0,\mu_0})\,\mbox{\boldmath $f$}_{m,x,\sigma}
\cdot\mbox{\boldmath $f$}_{m,x,\sigma}\,dS
\\
\\
\displaystyle
I(x,\sigma,\xi^1,\xi^2,\xi^3)_m=\sum_{j=1}^3\,\int_{\partial\Omega}\,(\Lambda_{\lambda,\mu}-\Lambda_{\lambda_0,\mu_0})\,\mbox{\boldmath $g$}_{m,x,\sigma}^j
\cdot\mbox{\boldmath $g$}_{m,x,\sigma}^j\,dS.
\end{array}
\right.
$$
Besides, corresponding to (2.10), we define two indicator functions
$$\left\{
\begin{array}{ll}
\displaystyle
I_0(x;(\lambda,\mu),(\lambda_0,\mu_0))
=\int_{D}
C_{l,m}\text{Sym}\,\nabla(\mbox{\boldmath $u$}_x^0+\mbox{\boldmath $w$}_x^0)\cdot\text{Sym}\,\nabla\mbox{\boldmath $u$}_x^0\,dy,
& x\in\Omega\setminus\overline D,
\\
\\
\displaystyle
I(x;(\lambda,\mu),(\lambda_0,\mu_0))
=\sum_{j=1}^3\,\int_{D}
C_{l,m}\text{Sym}\,\nabla(\mbox{\boldmath $u$}_x^j+\mbox{\boldmath $w$}_x^j)\cdot\text{Sym}\,\nabla\mbox{\boldmath $u$}_x^j\,dy,
& x\in\Omega\setminus\overline D,
\end{array}
\right.
$$
where the function $\mbox{\boldmath $w$}=\mbox{\boldmath $w$}_x^j\in H^1(\Omega)^3$ for each $j=0,1,2,3$ is the unique solution of 
the Dirichlet problem
$$
\left\{
\begin{array}{lll}
\displaystyle
\sum_{k=1}^3\,\frac{\partial}{\partial y_k}
\left(C_{\lambda,\mu}\text{Sym}\,\nabla\mbox{\boldmath $w$}\right)_{ik}
=-\sum_{k=1}^3\,\frac{\partial}{\partial y_j}\left(C_{\lambda-\lambda_0,\mu-\mu_0}\text{Sym}\,\nabla\mbox{\boldmath $u$}_x^j\right)_{ik},
& i=1,2,3, 
& y\in\Omega,
\\
\\
\displaystyle
\mbox{\boldmath $w$}=\mbox{\boldmath $0$},
& y\in\partial\Omega.
\end{array}
\right.
$$

One can easily check the convergence property of the indicator sequence and the boundedness of
indicator functions away from $\partial D$.

$\bullet$  Let $x\in\Omega\setminus\overline D$ and let $\sigma\in N_x$ satisfy $\sigma\cap\overline D=\emptyset$.
Then two indicator sequences $\{I(x,\sigma,\xi^0)_m\}$ and $\{I(x,\sigma,\xi^1,\xi^2,\xi^3)_m\}$ converge to
the value of indicator functions $I_0(\,\cdot\,;(\lambda,\mu),(\lambda_0,\mu_0))$ and $I(\,\cdot\,;(\lambda,\mu),(\lambda_0,\mu_0))$
at $x$, respectively.

$\bullet$  For each $\epsilon>0$ we have 
$$\displaystyle
\sup_{\text{dis}\,(x,\overline D)\,>\epsilon}\,\left(\vert\,I_0(x;(\lambda,\mu),(\lambda_0,\mu_0))\vert
+\vert\,I(x;(\lambda,\mu),(\lambda_0,\mu_0))\vert\right)<\infty.
$$

What really we have shown in \cite{IProbe2} in terms of Side A, that is, the essential part is as follows.

$\bullet$ Given $a\in\partial D$ we have
$$\displaystyle
\lim_{x\rightarrow a}\,\left(\vert\,I_0(x;(\lambda,\mu),(\lambda_0,\mu_0))\vert
+\vert\,I(x;(\lambda,\mu),(\lambda_0,\mu_0))\vert\right)=\infty.
$$

This is the multi probe method what we say.
It is a direct consequence of more detailed property of indicator functions listed below.

{\bf\noindent Case 1. $m(a)\not=0$.}
We have
$$\begin{array}{ll}
\displaystyle
\lim_{x\rightarrow a}\,I_0(x;(\lambda,\mu),(\lambda_0,\mu_0))=\pm\infty,
&
\text{if $\pm\,m(a)>0$.}
\end{array}
\tag {2.28}
$$

{\bf\noindent Case 2. $m(a)=0$.}
We have
$$\begin{array}{ll}
\displaystyle
\lim_{x\rightarrow a}\,I(x;(\lambda,\mu),(\lambda_0,\mu_0))=\pm\infty,
&
\text{if $\pm\,l(a)>0$.}
\end{array}
\tag {2.29}
$$

The key point in proving property (2.24) which is the reason of the choice $\mbox{\boldmath $u$}_x^0$ is
that $\nabla\cdot\mbox{\boldmath $u$}_x^0$ is weaker than $\text{Sym}\,\nabla\mbox{\boldmath $u$}_x^0$
as $x\rightarrow a$.  Thus $m(a)$ appears as a leading coefficient of 
the indicator function $I_0(x;(\lambda,\mu),(\lambda_0,\mu_0))$.

Note that we are assuming that $(l(a),m(a))\not=(0,0)$ for all $a\in\partial D$.  Thus two cases completely cover
all the possible cases.  The proof of properties (2.28) and (2.29) are done with the help of the system of integral inequalities
analogous to Proposition 2.5 established in \cite{I0}.

\noindent
$\Box$

Some remarks are in order.

$\bullet$  If one assumes the condition $m(a)\not=0$ instead of condition $(l(a),m(a))\not=(0,0)$,
then one can easily apply the usual probe method by using the single indicator function
$I_0(x;(\lambda,\mu),(\lambda_0,\mu_0))$ only.

$\bullet$  See also \cite{INT} about the uniqueness proof under a {\it monotonicity assumption} for a general anisotropic body
and its origin \cite{IEssential} which gives us a simple idea for the proof on the determination of the geometry of inclusions
in Isakov's uniqueness result \cite{Is0} for Problem 2.2.

$\bullet$  The outside to inside argument works also for the system (2.26) under the assumption that
$(\lambda(x),\mu(x))=(\lambda_0(x),\mu_0(x))$, $x\in\Omega\setminus D$ and say, 
$\lambda_0\in C^2(\overline\Omega)$ and $\mu_0\in C^3(\overline\Omega)$.  So combining with the uniqueness theorems
in \cite{ER} and \cite{NU1} and the earlier result on the boundary determination
result in \cite{ANS} one obtains the uniqueness of Lam\'e parameters inside $D$ under some additional
regularity assumptions on $\lambda_0+l$ and $\mu_0+m$ together with small ness of $(\mu_0+m)-c$ on $D$ for a positive constant $c$.
This is an extension of the result for the conductivity determination inside an inclusion which is a part of Isakov's uniqueness theorem \cite{Is0},
to an elastic body.

\subsubsection{De Giorgi-Nash-Moser theorem applied to the probe method}

In \cite{DIN} the original version of the probe method \cite{IProbe} has been applied also to the inverse inclusion problem
governed by the equation $\nabla\cdot\gamma\nabla u=F$ in $\Omega\subset\Bbb R^m$ with a known source term $F$ and two types of mixed boundary conditions.
The existence of source term $F$ causes a some technical difficulty in showing the blowing up of the indicator function.
The key point is the boundedness of $L^2$-norm over the whole domain of the reflected solution which corresponds to
the weak solution of (2.9).
However, using De Giorgi-Nash-Moser theorem (e.g., see \cite{HL}) which states the H\"older continuity of the weak solution of the governing equation
with essentially bounded coefficients, we have resolved the difficulty and established the desired result.

Let us briefly explain the idea to show $\Vert w_x\Vert_{L^2(\Omega)}$ is bounded as $x\rightarrow a\in\partial D$ in Problem 2.2.
First let $z=z_x$ be the weak solution of the Dirichlet problem
$$\left\{
\begin{array}{ll}
\displaystyle
\nabla\cdot\gamma\nabla z=-w_x & y\in\Omega,
\\
\\
\displaystyle
z=0, & y\in\partial\Omega.
\end{array}
\right.
$$
One can write
$$\begin{array}{ll}
\displaystyle
\int_{\Omega}\,\vert w_x\vert^2\,dy
&
\displaystyle
=\int_D h\nabla G_x^0\cdot\nabla z_x\,dy.
\end{array}
$$
Let $x\in\partial\Omega\setminus\overline D$ 
and assume that $h$ has a continuous extension in a neighbourhood of $\partial D$.
One can rewrite more
$$
\displaystyle
\int_{\Omega}\,\vert w_x\vert^2\,dy
=\int_D(h(y)-h(x))\nabla G_x^0\cdot\nabla z_x\,dy
+h(x)\int_D\nabla G_x^0\cdot\nabla z_x\,dy.
\tag {2.30}
$$
If we assume further that $h$ has a H\"older continuous extension with a suitable exponent $\alpha$ ($\alpha>\frac{1}{2}$ if $m=3)$)
in a neighbourhood of $\partial D$, 
we see that the first term on this right-hand side has a bound $C\Vert w_x\Vert_{L^2(\Omega)}$. 

Here we apply the H\"older continuity of $z_x$ in a neighbourhood of $\partial D$ which is a consequence
of De Giorgi-Nash-Moser theorem.  
Then, the integral of the second term in the right-hand side on (2.30) has the expression
$$
\displaystyle
\int_D\nabla G_x^0\cdot\nabla z_x\,dy
=\int_{\partial D}\frac{\partial G_x^0}{\partial\nu}\,(z_x(y)-z_x(x))\,dS(y).
$$
Note that we have used $\nabla_y\cdot\gamma^0(x)\nabla_yG_x^0(y)=0$ in $\Omega\setminus\{x\}$.
Then, we see the last term has a bounded $C\Vert w_x\Vert_{L^2(\Omega)}$ which is coming from the quantitative estimate
in De Giorgi-Nash-Moser theorem.  Combining two estimates we obtain the boundedness of $\Vert w_x\Vert_{L^2(\Omega)}$
as $x\rightarrow a\in\partial D$.

However, the author does not think that this is an only way since the coefficient is not completely general.
It is like employing a steam-hammer to crack a nut.

\subsubsection{The cutting method}

In \cite{INSlice} another method for Problem 2.2 in three dimensions with $\gamma_0=1$ has been introduced.
The method which we call the cutting method here enables us to obtain the convex hull of the section of an unknown inclusion cut by a given plane.
It is an application of the probe method in two dimensions to three dimensions.
However, instead of a needle sequence in three dimensions, we explicitly construct a sequence of harmonic functions 
in two dimensions that approximates a singular solution of the Laplace equation with an {\it isolated singularity}
outside a needle locally.
Then by extending the sequence naturally as a sequence of harmonic functions in three dimensions,
we see that the obtained sequence approximates a singular solution of the Laplace equation with a {\it line singularity}
outside a half plane locally, which should be called a {\it cutter}.  By prescribing the trace onto the surface of the body
of each member from the sequence and observing the corresponding Neumann data, one gets an indicator function
which tells us whether the cutter firstly hits an unknown inclusion or not.  This is a brief explanation of the idea of
the cutting method.
Unfortunately, unlike the probe method itself the numerical testing remains open.

\subsection{Some remarks and open problems}

\subsubsection{Revisiting inverse obstacle scattering problem}

Succeeding to \cite{IScattering1} and a numerical study of the probe method \cite{EP}(see also \cite{CLN}), 
in \cite{IProbeNew} the author reconsidered a typical inverse obstacle problem
governed by the Helmholtz equation using two sides of the probe method.  It is a reduction of the original inverse obstacle scattering problem
of acoustic wave with a fixed frequency.  Let $\Omega$ be a bounded domain of $\Bbb R^3$  with a smooth boundary.
Let $D$ be an open set of $\Bbb R^3$ with smooth boundary, having finitely many
connected components,  satisfying $\overline D\subset\Omega$ and that $\Omega\setminus\overline D$
is connected.  Let $k>0$ and assume that $k^2$ is not an eigenvalue for $\-\Delta$ in $\Omega\setminus\overline D$ with
homogeneous Dirichlet boundary value condition on $\partial\Omega$ and Neumann boundary condition on  
$\partial D$.

Given $f\in H^{1/2}(\partial\Omega)$ let $u\in H^1(\Omega\setminus\overline D)$ be the weak solution of the elliptic problem
$$\left\{
\begin{array}{ll}
\displaystyle
\Delta u+k^2 u=0, & y\in\Omega\setminus\overline D,
\\
\\
\displaystyle
\frac{\partial u}{\partial\nu}=0, & y\in\partial D,
\\
\\
\displaystyle
u=f, & y\in\partial\Omega.
\end{array}
\right.
\tag {2.31}
$$
The associated Dirichlet-to-Neumann map $\Lambda_D$ is given by
$$\displaystyle
\Lambda_D:f\longmapsto \frac{\partial u}{\partial\nu}\vert_{\partial\Omega}.
$$

The indicator function is defined as follows.  Let $G(y)$ be a solution of the Helmholtz equation in $\Bbb R^3\setminus\{(0,0,0)\}$
such that, (2.3) is satisfied with for all finite cone $V$ with vertex at $(0,0,0)$.
Hereafter we fix $G$.  

The indicator function $I(x)$, $x\in\Omega\setminus\overline D$ is defined by the formula
$$\displaystyle
I(x)=\int_D\vert\nabla G(y-x)\vert^2\,dy
-k^2\int_D \vert G(y-x)\vert^2\,dy+\int_{\Omega\setminus\overline D}\,\vert\nabla w_x\vert^2\,dy
-k^2\int_{\Omega\setminus\overline D}\,\vert w_x\vert^2\,dy,
$$
where $w_x$ is the unique weak solution of the problem
$$\left\{
\begin{array}{ll}
\displaystyle
\Delta w+k^2 w=0, & y\in\Omega\setminus\overline D,
\\
\\
\displaystyle
\frac{\partial w}{\partial\nu}=-\frac{\partial}{\partial\nu}(G(\,\cdot\,-x)),
&
y\in\partial D,\\
\\
\displaystyle
w=0, & y\in\partial\Omega.
\end{array}
\right.
$$

Given $x\in\Omega$ and $\sigma\in N_x$ we call the sequence $\{v_n\}$ of $H^1(\Omega)$ solutions of
the Helmholtz equation a {\it needle sequence} for $(x,\sigma)$ if it satisfies, for any compact set $K$ of $\Bbb R^3$ 
with $K\subset\Omega\setminus\sigma$
$$\displaystyle
\lim_{n\rightarrow\infty}\,(\Vert v_n(\,\cdot\,)-G(\,\cdot\,-x)\Vert_{L^2(K)}
+\Vert\nabla\{v_n(\,\cdot\,)-G(\,\cdot\,-x)\}\Vert_{L^2(K)}=0.
$$
that is, $\lim_{n\rightarrow\infty}v_n=G(\,\cdot\,-x)$ in $H^1_{\text{loc}}\,(\Omega\setminus\sigma)$.
Under the assumption that $k^2$ is not a Dirichlet eigenvalue for $-\Delta$ on $\Omega$, the existence is known.
See \cite{IScattering1} and Appendix in \cite{IProbeNew}.
In the following theorem we do not need this assumption.
It is the previous result \cite{IScattering1} stated under the new formulation.

\proclaim{\noindent Theorem 2.6(\cite{IProbeNew}).}  

$\bullet$  (A.1)  Given $x\in\Omega\setminus\overline D$ let $\sigma\in N_x$  satisfy $\sigma\cap\overline D=\emptyset$.
Then, for all needle sequence $\xi=\{v_n\}$ for $(x,\,\sigma)$ the indicator sequence $\{I(x,\sigma,\xi)_n\}$
defined by
$$\displaystyle
I(x,\sigma,\xi)_n
=\int_{\partial\Omega}\left(\frac{\partial v_n}{\partial\nu}\vert_{\partial\Omega}-\Lambda_D(v_n\vert_{\partial\Omega})\,\right) \overline{v_n\vert_{\partial\Omega}}\,dS,
$$
converges to the value of indicator function $I(\,\cdot\,)$ at $x$.

$\bullet$  (A.2)  For each $\epsilon>0$ we have $\sup_{\text{dist}\,(x,\overline{D})>\epsilon}\,\vert I(x)\vert<\infty$.

$\bullet$  (A.3)  Given $a\in\partial D$ we have $\lim_{x\rightarrow a}\,I(x)=\infty$.

\endproclaim

What is the new point of \cite{IProbeNew} which was not clarified nor asked in \cite{IScattering1} is 
the behaviour of indicator sequence $I(x,\sigma,\xi)_n$ as $n\rightarrow\infty$ when $x$ is just located on the boundary, 
inside or passing through the obstacles.
The following result gives an partial answer to those questions.

\proclaim{\noindent Theorem 2.7(\cite{IProbeNew}).}
Assume that $k^2$ is sufficiently small (specified below).
Let $x\in\Omega$ and $\sigma\in N_x$.  If $x\in\Omega\setminus\overline D$ and $\sigma\cap D\not=\emptyset$ 
or $x\in\overline D$, then for any needle sequence $\xi=\{v_n\}$ for $(x,\sigma)$ we have
$\lim_{n\rightarrow\infty}\,I(x,\sigma,\xi)_n=\infty$.

\endproclaim

The smallness assumption is given by the set of two inequalities
$$\left\{
\begin{array}{l}
\displaystyle
C(\Omega\setminus\overline D)^2\,k^2\le 1,\\
\\
\displaystyle
\max_{j=1,\cdots,N}\,(8C(D_j)^2)\,k^2<1,
\end{array}
\right.
\tag {2.32}
$$
where $D_1,\cdots,D_N$ denote the connected components of $D$, the constants $C_1=C(\Omega\setminus\overline D)$ and $C_2=C(D_j)$ are the Poincar\'e constants such that:
for all $w\in H^1(\Omega\setminus\overline D)$ with $w=0$ on $\partial D$
$$\displaystyle
\Vert w\Vert^2_{L^2(\Omega\setminus\overline D)}\le C_1^2\Vert\nabla w\Vert^2_{L^2(\Omega\setminus\overline D)};
$$
for all $v\in H^1(D_j)$ with $\displaystyle\int_{D_j}\,v\,dy=0$
$$\displaystyle
\Vert v\Vert^2_{L^2(D_j)}
\le C_2^2\Vert\nabla v\Vert^2_{L^2(D_j)}.
$$
Here let us explain the role of smallness assumption (2.32) in showing the blow up of the indicator sequence.
For simplicity of description let us consider only the case when $N=1$.
Using integration by parts, one gets
$$\displaystyle
I(x,\sigma,\xi)_n=
\int_D\,\vert\nabla v_n\vert^2\,dy-k^2\int_D\vert v_n\vert^2\,dy
+\int_{\Omega\setminus\overline{D}}\vert\nabla w_n\vert^2\,dy-k^2\int_{\Omega\setminus\overline{D}}\,\vert w\vert^2\,dy,
$$
where $w_n=u_n-v_n$ and $u_n$ is the solution of (2.31) with $f=v_n$ on $\partial\Omega$.
Then the first condition on (2.32) enables us to drop the third and fourth terms on this right-hand side and one gets
the lower estimate
$$\displaystyle
I(x,\sigma,\xi)_n\ge
\int_D\,\vert\nabla v_n\vert^2\,dy-k^2\int_D\vert v_n\vert^2\,dy.
$$
Let $A\subset D$ be an arbitrary measurable set of $\Bbb R^3$ with a positive measure
and $(v_n)_A$ denote the mean value of $v_n$ over $A$.
One has 
$$\displaystyle
\int_D\vert v_n\vert^2\,dy
\le 2\int_D\vert v_n-(v_n)_A\vert^2\,dy+2\vert D\vert\,\vert (v_n)_A\vert^2.
$$
Besides, using the constant $C(D)$, one gets 
$$\displaystyle
\int_D\vert v_n-(v_n)_A\vert^2\,dy
\le C(D)^2F(A,D)^2\int_D\vert\nabla v_n\vert^2\,dy,
$$
where $F(A,D)\rightarrow 2$ as the measure $\vert A\vert$ of $A$ approaches that of $D$.
This is a kind of the Poincar\'e inequality (\cite{SS, Z}).  Then the second condition on (2.32) together with a suitable choice of $A\subset D\setminus\sigma$ such that that $\overline{A}\subset\Omega\setminus\sigma$
and $\vert A\vert \approx\vert D\vert$,
enables us to show 
$$\displaystyle
I(x,\sigma,\xi)_n\ge
C_3\int_D\,\vert\nabla v_n\vert^2\,dy-2k^2\vert D\vert\,\vert(v_n)_A\vert^2,
$$
where $C_3$ is a positive number independent of $n$. 
Since $\overline{A}\cap\sigma=\emptyset$,  the second term on this right-hand side
is convergent and thus one gets the desired conclusion.

A combination of (A.1) in Theorems 2.6 and 2.7 yields the characterization of the geometry of $D$
in terms of indicator sequences.

\proclaim{\noindent Corollary 2.2(\cite{IProbeNew}.)}
Let $k$ satisfy the same conditions as Theorem 2.7.  A point $x\in\Omega$
belongs to $\Omega\setminus\overline D$ if and only if there exists a needle $\sigma\in N_x$
and a needle sequence $\xi$ for $(x,\sigma)$ such that indicator sequence $\{I(x,\sigma,\xi)_n\}$ is bounded from above.

\endproclaim

It is still an open problem whether one can remove the smallness condition (2.32) in Theorem 2.7 or not.

In \cite{IProbeImpedance}, the both sides of the probe method together with the enclosure method explained in Section 3
have been studied in the case when the boundary condition on $\partial D$ in (2.31) is replaced with the impedance boundary condition
$$\begin{array}{ll}
\displaystyle
\frac{\partial u}{\partial\nu}+\lambda(y)u=0, & y\in\partial D,
\end{array}
$$
where the coefficient $\lambda$ satisfies the condition $\text{Im}\,\lambda>0$.
See also \cite{CLN0} for the result corresponding to the side A.

\subsubsection{Yarmukhamedov's fundamental solution and a needle sequence}

In \cite{IProbeCarleman} the author pointed out an unexpected relationship between the needle sequence
for the Laplace and Helmholtz equations and a special family of fundamental solutions of the Laplace equation introduced in \cite{Y1}.

The family is parametrized by an entire function $K(w)$ satisfying the conditions:
$$\displaystyle
\left\{
\begin{array}{l}
\displaystyle
\overline{K(w)}=K(\overline w),
\\
\\
\displaystyle
K(0)=1,
\\
\\
\displaystyle
\forall R>0\,\,\sum_{m=0}^{2}\,\sup_{\vert\text{Re}\,w\vert<R}\,\vert K^{(m)}\,(w)\vert\,<\infty,
\end{array}
\right.
\tag {2.33}
$$
where $K^{(0)}=K$, $K^{(1)}=K'$ and $K^{(2)}=K''$.

Define
$$\begin{array}{ll}
\displaystyle
\Phi_K(x)=
-\frac{1}{2\pi^2}
\,\int_0^{\infty}\,\text{Im}\,\left(\frac{K(w)}{w}\right)\,\frac{du}{\sqrt{\vert x'\vert^2+u^2}}, & x'\not=(0,0),
\end{array}
$$
where $w=x_3+i\sqrt{\vert x'\vert^2+u^2}$, $x'=(x_1,x_2)$ and $x=(x_1,x_2,x_3)$.
Note that $K\equiv 1$ satisfies (2.33) and we have the expression
$$\displaystyle
\Phi_1(x)=\frac{1}{4\pi\,\vert x\vert}.
$$

In \cite{Y2} Yarmukhamedov proved the following theorem.

\proclaim{\noindent Theorem 2.8(\cite{Y2}).}
One has the expression
$$\begin{array}{ll}
\displaystyle
\Phi_K(x)=\frac{1}{4\pi\vert x\vert}+H_K(x), & x'\not=(0,0),
\end{array}
$$
where $H_K$ is $C^2$ in the whole space and satisfies the Laplace equation
$\Delta H_K=0$ in $\Bbb R^3$.

\endproclaim

The family has been applied
to the Cauchy problem for the Laplace equation \cite{Y2, Y3, Y4} as a {\it Carleman function} introduced
by Lavrentiev.  See also the books \cite{Aiz}, \cite{LRS} and the article \cite{AB} for the notion of the Carleman function and various results
on this subject itself.

Here we given an application of the family with a special choice of $K(w)$ to the probe method.
The harmonic function $H_K$ is unique.  Here we choose $K=K_{\tau,\alpha}$ depending on two parameters $\tau>0$ and $0<\alpha\le 1$
given by
$$\displaystyle
K_{\tau,\alpha}(w)=E_{\alpha}(\tau w),
$$
where the entire function $E_{\alpha}(z)$ is the Mittag-Leffler function \cite{Bateman} defined by
$$
\displaystyle
E_{\alpha}(z)=\sum_{m=0}^{\infty}\,\frac{z^m}{\Gamma(1+\alpha\,m)}.
$$
It is known that, roughly speaking, we have, as $z\rightarrow\infty$,
$E_{\alpha}(z)\sim \alpha^{-1}\,e^{z^{1/\alpha}}$ if $\vert\text{arg}\,z\vert\le \frac{\pi\alpha}{2}$
and $E_{\alpha}(z)\sim-\frac{1}{\Gamma(1-\alpha)\,z}$ if $\pi\ge \vert\text{arg}\,z\vert>\frac{\pi\alpha}{2}$.

{\bf\noindent Definition 2.14.}  Given two unit vectors $\vartheta_1$ and $\vartheta_2$ define the harmonic function $v$
given by
$$\displaystyle
v(y)=v(y;\alpha,\tau,\vartheta_1,\vartheta_2)=-H_{K}(y\cdot\vartheta,y\cdot\vartheta_2,\vartheta_1\times\vartheta_2),
$$
where $K=K_{\tau,\alpha}$.  Since $\Phi_1(x)=\frac{1}{4\pi\vert x\vert}$, we have the expression
$$\displaystyle
v(y)=-\frac{1}{2\pi^2}\,
\int_0^{\infty}\,\text{Im}\,\left(\frac{E_{\alpha}(\tau w)-1}{w}\,\right)
\frac{du}{\sqrt{\vert y\cdot\vartheta_1\vert^2+\vert y\cdot\vartheta_2\vert^2+u^2}},
\tag {2.34}
$$
where $w=y\cdot(\vartheta_1\times\vartheta_2)+i\sqrt{\vert y\cdot\vartheta_1\vert^2+\vert y\cdot\vartheta_2\vert^2+u^2}$.

The following theorem says that the regular part of Yarmukhamedov's fundamental solution generates a needle sequence
for a needle given by a line segments directed to an arbitrary direction.

\proclaim{\noindent Theorem 2.9(\cite{IProbeCarleman}).}  Let $\Omega$ be a bounded domain of $\Bbb R^3$.
Let $x\in\Omega$ and $\sigma\in N_x$ be a straight needle directed to $\omega=\vartheta_1\times\vartheta_2$,
that is, have the expression $\sigma=\{x+t\omega\,\vert\,0\le t<\infty\}\cap\Omega$.
Then, the sequence $\{v(\,\cdot\,-x;\alpha_n,\tau_n,\vartheta_1,\vartheta_2)\vert_{\Omega}\}$ is a needle sequence
for $(x,\sigma)$, where $\alpha_n$ and $\tau_n$ are suitably chosen sequences and satisfy
$0<\alpha_n<1$, $\alpha_n\rightarrow 0$, $\tau_n>0$ and $\tau_n\rightarrow\infty$.
\endproclaim



Applying the Cauchy integral theorem to representation (2.34), one gets the explicit computation formulae
of $v(x-y;\alpha,\tau,\vartheta_1,\vartheta_2)$ together with its gradient on the line $y=x+s\,\omega\,(-\infty<s<\infty)$:
$$\displaystyle
v(y-x;\alpha,\tau,\vartheta_1,\vartheta_2)
=\left\{
\begin{array}{ll}
\displaystyle
\frac{E_{\alpha}(\tau s)-1}{4\pi s}, & \text{if $y=x+s\omega$ with $s\not=0$,}
\\
\\
\displaystyle
\frac{\tau}{4\pi\,\Gamma(1+\alpha)}, & \text{if $y=x$}
\end{array}
\right.
$$
and
$$\displaystyle
\nabla v(y-x;\alpha,\tau,\vartheta_1,\vartheta_2)
=\left\{
\begin{array}{ll}
\displaystyle
\frac{d}{ds}\left\{\frac{E_{\alpha}(\tau s)-1}{4\pi s}\right\}\,\omega, & \text{if $y=x+s\omega$ with $s\not=0$,}
\\
\\
\displaystyle
\frac{\tau^2}{4\pi\,\Gamma(1+2\alpha)}\,\omega, & \text{if $y=x$.}
\end{array}
\right.
$$
From these we conclude that $v$ together with its gradient blows up on the half line
$y=x+s\,\omega\,(s\ge 0)$; the vector $\nabla v$ on the line $y=x+s\,\omega\,(-\infty<s<\infty)$ is parallel to the direction $\omega$
and its $\omega$-component is positive.

One can construct a needle sequence for the Helmholtz equation for all straight needles by using the harmonic function $v=v(y;\alpha,\tau,\vartheta_1,\vartheta_2)$.  It is an application of the Vekua transform introduced by Vekua \cite{V1,V2}.

Let $v=v(y;\alpha,\tau,\vartheta_1,\vartheta_2)$ and
define
$$\displaystyle
v^k(y;\alpha,\tau,\vartheta_1,\vartheta_2)
=v(y)-\frac{k\vert y\vert}{2}\,\int_0^1\,v(ty)\,J_1(k\vert y\vert\,\sqrt{1-t})\,\sqrt{\frac{t}{1-t}}\,dt.
$$
Note that the function on this right-hand side is nothing but the Vekua transform of harmonic function $v$.
The function $v^k(y;\alpha,\tau,\vartheta_1,\vartheta_2)$ satisfies the Helmholtz equation $\Delta u+k^2 u=0$ in the whole space.

Let $k\ge 0$ and $G_k(y)$ denote the outgoing Green function for the Helmholtz equation 
$$\displaystyle
G_k(y)=\frac{e^{ik\vert y\vert}}{4\pi \vert y\vert}.
$$

\proclaim{\noindent Theorem 2.10(\cite{IProbeCarleman}).}  Let $\Omega$ be a bounded domain of $\Bbb R^3$.
Let $x\in\Omega$ and $\sigma\in N_x$ be a straight needle directed to $\omega=\vartheta_1\times\vartheta_2$.
Then the sequence of solutions of the Helmholtz equation defined by
$$\begin{array}{ll}
\displaystyle
v^k(y-x;\alpha_n,\tau_n,\vartheta_1,\vartheta_2)+i\,\frac{\sin k\vert y-x\vert}{4\pi\vert y-x\vert}, & y\in\Omega,
\end{array}
$$
converges to $G_k(y-x)$ in $H^1_{\text{loc}}(\Omega\setminus\sigma)$ as $n\rightarrow\infty$,
where $\alpha_n$ and $\tau_n$ are suitably chosen sequences and satisfy
$0<\alpha_n<1$, $\alpha_n\rightarrow 0$, $\tau_n>0$ and $\tau_n\rightarrow\infty$.
\endproclaim

Note that the function $v^k(y-x;\alpha,\tau)$ together with its gradient at the point $y=x$ has the explicit expression
$$\left\{
\begin{array}{l}
\displaystyle
v^k(0;\alpha,\tau,\vartheta_1,\vartheta_2)=\frac{\tau}{\Gamma(1+\alpha)},
\\
\\
\displaystyle
\nabla v^k(0;\alpha,\tau,\vartheta_1,\vartheta_2)=\frac{\tau^2}{4\pi\Gamma(1+2\alpha)}\,\omega.
\end{array}
\right.
\tag {2.35}
$$
This shows that the needle sequence for the straight needle in Theorem 2.10 blows up at the tip.

Concerning with the frequency restrictions (2.32) in Theorem 2.7,
we have an open problem described below.

Let $k\ge 0$ be a constant, $x\in\Bbb R^m(m=2,3)$ and
a sequence of functions $\{v_n(x)\}$, $n=1,2,\cdots$ satisfy
the Helmholtz equation $\displaystyle\Delta v+k^2v=0$ in $\Bbb R^m$.
Let $D$ be a bounded domain of $\Bbb R^m$ and $0\in\,\overline{D}$.

Is the following statement true?

{\it If we have
$$
\displaystyle
\lim_{n\longrightarrow\infty}\int_D\vert\nabla v_n(x)\vert^2dx=\infty
$$
and
$$
\displaystyle
\lim_{n\longrightarrow\infty}\frac{v_n(0)}{\vert\nabla v_n(0)\vert}=0,
\tag {2.36}
$$
then
$$\displaystyle
\lim_{n\longrightarrow\infty}
\frac{\displaystyle\int_D\vert v_n(x)\vert^2dx}
{\displaystyle\int_D\vert\nabla v_n(x)\vert^2dx}=0.
$$
}

Note that the needle sequences constructed in Theorem 2.10 satisfies (2.36) because of (2.35).
If the statement above is right, then the second condition on  (2.32) 
can be immediately dropped in the case when $x\in\overline D$ in Theorem 2.7
and the special needle sequences constructed in Theorem 2.10 are used for the indicator sequence.

It would be interesting to do the numerical implementation of the probe method in Section 2.3.1 
using the needle sequences constructed in Theorems 2.9 and 2.10.

\subsubsection{The needle sequence and localized potential}

This subsection gives a remark on the relationship
between the basis of the probe and monotonicity method.  
The monotonicity method goes back to \cite{RTVV}, \cite{TR} and its mathematical justification
has been done in \cite{HU} by using a notion of the localized potential \cite{G}.

In this section we show that the needle sequence in the probe method generates the
localized potential.  In this sense, the base of the probe method is deeper
than that of the monotonicity method.

Let $\Omega$ be a bounded domain of $\Bbb R^m$, $m=2,3$.  
Let $y\in\Omega$ be an arbitrary point in $\Omega$ and $\sigma$ be a needle in $\Omega$ with tip $y$.
Let $\xi=\{v_n\}$ be a needle sequence for the needle $\sigma\in N_y$.

As we explained in Subsection 2.1 sequence $\{v_n\}$ satisfies:

(A)  for any finite cone $V$ with vertex $y$ we have $\lim_{n\rightarrow\infty}\Vert\nabla v_n\Vert_{L^2(V\cap\Omega)}=\infty$;

(B)  for any $z\in\sigma\cap\Omega$ and any open ball $B$ centered at $z$ we have
$\lim_{n\rightarrow\infty}\Vert\nabla v_n\Vert_{L^2(B\cap\Omega)}=\infty$.

Let $U$ be an open subset of $\Bbb R^m$ with $\sigma\subset U$.
Let $W$ be an Lebesgue measurable set of $\Bbb R^m$ such that $\overline W\subset\Omega$ and $\overline W\cap\overline U=\emptyset$.

Define
$$\begin{array}{ll}
\displaystyle
v_n^*(x)=\frac{v_n(x)}{\lambda_n}, & x\in\Omega,
\end{array}
$$
where $\lambda_n>0$.

We have
$$\displaystyle
\Vert \nabla v_n^*\Vert_{L^2(U\cap\Omega)}=\frac{\Vert \nabla v_n\Vert_{L^2(U\cap\Omega)}}{\lambda_n}
$$
Choose
$$\begin{array}{ll}
\displaystyle
\lambda_n=\Vert\nabla v_n\Vert_{L^2(U\cap\Omega)}^{\epsilon}, & 0<\epsilon<1.
\end{array}
$$
From (A) or (B) we have $\lambda_n\rightarrow\infty$ and thus
$$\displaystyle
\Vert \nabla v_n^*\Vert_{L^2(U\cap\Omega)}=\lambda_n^{\frac{1-\epsilon}{\epsilon}}\rightarrow\infty.
$$
Since $\overline {W}\subset\Omega\setminus\sigma$, we have
$$\displaystyle
\lim_{n\rightarrow\infty}
(\Vert v_n(\,\cdot\,)\Vert_{L^2(W)}
+\Vert\nabla\{v_n(\cdot)\}\Vert_{L^2(W)})=
\Vert G(\,\cdot\,-y)\Vert_{L^2(V)}
+\Vert\nabla G(\cdot-y)\Vert_{L^2(W)}<\infty.
$$
Thus one gets
$$\displaystyle
\lim_{n\rightarrow\infty}(\Vert v_n^*(\,\cdot\,)\Vert_{L^2(W)}
+\Vert\nabla\{v_n^*(\,\cdot\,)\}\Vert_{L^2(W)}=0.
$$
This means that $\{v_n^*\}$ is just the localized potential (for the Laplace equation).
This potential should be called the localized potential associated with a needle.  The potential does not depend on
$W$ itself.  Thus this is a generating function of the localized potential.

Note that if one uses the needle sequence generated by Yarmukhamedov's fundamental solution
for a straight needle, one can drop the condition $\overline W\subset\Omega$.  Instead just assume
that $W\subset\Omega$.

\subsubsection{The blow up set of the reflected solution}

In \cite{IProbeNew} we formulated a problem related to the behaviour
of the sequence of reflected solutions in the probe method applied to
a typical inverse obstacle problem.
For the purpose we introduced the notion of the blowup set for a sequence of functions.

{\bf\noindent Definition 2.15(\cite{IProbeNew}).} 
Let $U$ be an open subset of $\Bbb R^m$.
We say that a sequence
$\{g_n\}$ of $H^1(U)$ functions blows up at the point $z\in\overline U$
if for any open ball $B$ centered at $z$ it holds that
$$
\displaystyle
\lim_{n\longrightarrow\infty}
\int_{B\cap\overline U}\vert\nabla g_n(y)\vert^2dy=\infty.
$$

\noindent
We call the set of all points $z\in\overline U$ such that
$\{g_n\}$ blows up at $z$ the blowup set of $\{g_n\}$ and denote by $B(\{g_n\};U)$.

\noindent

Let $\Omega$ be a bounded domain of $\Bbb R^m$ with smooth boundary.
Given $x\in\Omega$, needle $\sigma$ with tip at $x$
and needle sequence $\xi=\{v_n\}$ for $(x,\sigma)$ let $u=u_n\in H^1(\Omega\setminus\overline D)$
solve
$$\left\{
\displaystyle
\begin{array}{ll}
\displaystyle
\Delta u=0, & x\in\Omega\setminus\overline D,\\
\\
\displaystyle
\frac{\partial u}{\partial\nu}=0, & x\in\partial D,\\
\\
\displaystyle
u=v_n, & x\in\partial\Omega,
\end{array}
\right.
$$
where $D$ is a bounded open subset of $\Bbb R^m$ with Lipschitz boundary such that $\overline D\subset\Omega$
and $\Omega\setminus\overline D$ is connected.  $\nu$ denotes the unit outward vector on $\partial\Omega$ and $\partial D$.

We call the function $u_n-v_n$ the reflected solution by the obstacle $D$.
It is easy to see that if $\sigma\cap\overline D=\emptyset$, then $\{u_n-v_n\}$
is bounded in $H^1(\Omega\setminus\overline D)$ and thus $B(\{u_n-v_n\};\Omega\setminus\overline D)=\emptyset$.

In the next proposition proved in \cite{IProbeNew} gives 
an example of $D$ in two dimensions such that $B(\{u_n-v_n\};\Omega\setminus\overline D)\not=\emptyset$.

\proclaim{\noindent Proposition 2.9.} 
Let $R>\epsilon>0$. The sets $\Omega$ and $D$ are given by the open discs centered at the
origin with radius $R$ and $\epsilon$, respectively. 
Let $x\in D$ and $\sigma\in N_x$ satisfy: $\sigma$ intersects with $\partial D$ only one time;
$\displaystyle\sigma\cap\{y\,\vert\,\vert y\vert\le\frac{\epsilon^2}{R}\}=\emptyset$.
Then we have $B(\{u_n-v_n\};\Omega\setminus\overline D)=\sigma^R$,
where the curve $\sigma^R$ is given by the formula
$$\displaystyle
\sigma^R=\left\{\frac{\epsilon^2y}{\vert y\vert^2}\,\vert\,y\in\sigma\cap\overline D\right\}.
$$
\endproclaim




The proof is based on the property of the Kelvin transform with respect to the circle centered at the origin with radius $\epsilon$.

Proposition 2.9 says that in the case when $x\in D$, the blowup set of $\{u_n-v_n\}$ is given by
a suitable curve in $\overline\Omega\setminus D$ obtained by transforming the part
of needle $\sigma\in N_x$ in $\overline D$.  This may suggest us to make use of the reflected solutions
to detect the {\it backside} of another cavity or inclusion occurred in $\Omega\setminus\overline D$.

Thus it is natural to raise

$\quad$

{\bf\noindent Problem 2.4.}   For general $\Omega$ and $D$ what can one say about $B(\{u_n-v_n\};\Omega\setminus\overline D)$
when $\sigma\cap\overline D\not=\emptyset$?  How about the case when the Laplace equation is replaced with a general elliptic equation?

$\quad$

\subsubsection{Discontinuity of anisotropic conductivity in two dimensions}

In \cite{Ideterminant} we considered the problem of identification of discontinuity in a general anisotropic inhomogeneous body
in two dimensions whose governing equation is given by the equation $\nabla\cdot A\nabla u=0$.
The coefficient $A$ of the equation takes the form
$$
A(y)=
\left\{
\begin{array}{ll}
\displaystyle
A_0(y), & y\in\Omega\setminus D,\\
\\
\displaystyle
A_0(y)+B(y), & y\in D,
\end{array}
\right.
$$
where $\Omega$ is a bounded domain of $\Bbb R^2$ with a smooth boundary, $A_0(y)$, $y\in\overline\Omega$
is a real symmetric $2\times 2$ matrix valued function and $C^2$ in a neighbourhood of $\overline\Omega$, $B(y)$ is a real symmetric $2\times 2$ matrix valued function of $y\in D$ and $C^2$ in a neighbourhood of $\overline D$.  The set $D$ is an open subset of $\Omega$
with Lipschitz boundary satisfying $\overline D\subset\Omega$ and $\Omega\setminus\overline D$ is connected.
It is assumed that both $A$ and $A_0$ are uniformly positive definite.

The problem considered in \cite{Ideterminant} is as follows.

$\quad$

{\bf\noindent Problem 2.5.}
Fix a non empty open subset $\Gamma$ of $\partial\Omega$.
Assume that both $D$ and $B$ are {\it unknown}.  Assume also that 
$$\begin{array}{ll}
\displaystyle
\text{det}(A_0(y)+B(y))\not=\text{det}\,A_0(y), & \forall y\in\partial D.
\end{array}
\tag {2.37}
$$
Identify $D$ from the knowledge of $\Lambda_Af$ for all $f\in H^{1/2}(\partial\Omega)$ with $\text{supp}\, f\subset\Gamma$.

$\quad$

In \cite{Ideterminant}, by replacing the Green's function type singular solutions used in \cite{Is0}
with other singular solutions, it is proved that the $D$ is uniquely determined by $\Lambda_Af$ for all $f\in H^{1/2}(\partial\Omega)$ with $\text{supp}\, f\subset\Gamma$ provided $B$ satisfies (2.37) and is sufficiently small on $\partial D$.
It is an application of the {\it method of singular solutions} developed in \cite{Is0}.
Unlike \cite{Is0}, we do not assume that $B$ is uniformly positive definite.
Thus one can not make use of the system in \cite{I0} like (2.18) and (2.19) directly and thus the technique based on the system developed
in \cite{IEssential} for the uniqueness can not be applied to this case.

Therefore it is interesting to consider the following problem.

$\quad$

{\bf\noindent Problem 2.6.}  Develop the probe method to reconstruct $D$ under condition (2.37).

$\quad$

If there is no discontinuity of $\text{det}\,A$ across $\partial D$, one can not uniquely determine $D$ even the case when
some of components of $B$ on $\partial D$ does not vanish (hidden singularity). There is an example pointed out by Spagnolo, S., see \cite{KV} 
and also Introduction in \cite{Ideterminant}.

Note that, in \cite{IStroh} an extension of Problem 2.5
to inverse boundary value problems for the system of equations in linearlized elasticity in two dimensions together with
a fourth order elliptic equation in the Love-Kirchhoff plate theory have been considered.
Using the relationship between them established in \cite{Ie, Ie2, Ie3} and the method of singular solutions,
the author has established a uniqueness theorems for the curve of some kind of discontinuity of piecewise constant elasticity tensor
field via the localized Dirichlet-to-Neumann map.
The smallness on $B$ in \cite{Ideterminant} may be removed by using the technique developed in \cite{IStroh}.

\section{The Enclosure Method}

The probe method developed in \cite{IProbe} is based on the Runge approximation property.
So from the beginning when the method was discovered someone had pointed out that its numerical
implementation shall be difficult.

The author seriously thought about it
and finally in \cite{E00} by considering the same problem as \cite{IProbe} the author introduced a method which is now called the enclosure method,
to extract the convex hull of unknown inclusions from the Dirichlet-to-Neumann map
acting on infinitely many input data.

The input data in \cite{E00} are the traces of the complex exponential solution used in \cite{Cal}
and later Section 4 in \cite{E01} those of the complex geometrical optics solutions constructed by Sylvester-Uhlmann \cite{SU1, SU2}.
Since those are explicitly or constructively given, the enclosure method resolves the difficulty of its numerical realization in principle.
As an evidence, soon later, in \cite{IS0} and \cite{BH0} the numerical implementation of the enclosure method was done.
See also \cite{CX} for a testing of the enclosure method in the electrical capacitance tomography.

\subsection{Two types of the original enclosure method}

Let $\Omega$ be a bounded domain in $\Bbb R^2$ with Lipschitz boundary.
Let $D$ be an open set with Lipschitz boundary satisfying
$\overline D\subset\Omega$.  $D$ may have many connected components
unless otherwise stated and we assume that $\Omega\setminus\overline D$
is connected.
Let $\nu$ denote the unit outward normal vector field to $\partial(\Omega\setminus\overline D)$.
Given $f\in H^{1/2}(\partial\Omega)$ let $u\in H^1(\Omega\setminus\overline D)$
be the weak solution of the elliptic problem
$$\left\{
\begin{array}{ll}
\displaystyle
\nabla\cdot\gamma_0\nabla u=0,
&
x\in\Omega\setminus\overline D,\\
\\
\displaystyle
\gamma_0\frac{\partial u}{\partial\nu}=0,
&
x\in\partial D,\\
\\
\displaystyle
u=f,
&
x\in\partial\Omega,
\end{array}
\right.
\tag {3.1}
$$
where, for simplicity, $\gamma_0\in C^2(\Bbb R^2)$ and $\gamma_0(x)>0$ for all $x\in\Bbb R^2$.
We consider $\Omega\setminus\overline D$ an isotropic electric conductive medium
with electrical conductivity $\gamma_0\vert_{\Omega\setminus\overline D}$; $D$ is considered
as a cavity where the electrical conductivity takes $0$.
Let $\Lambda_D$ denote the Dirichlet-to-Neumann map associated with the elliptic problem:
$$
\displaystyle\Lambda_D:f\longmapsto \gamma_0\frac{\partial u}{\partial\nu}\vert_{\partial\Omega}.
$$
$f$ is a voltage potential on $\partial\Omega$; $\Lambda_Df$
is the electric current density on $\partial\Omega$ that induces $f$.

$\quad$

{\bf\noindent Problem 3.1.}
Assume that $\gamma_0$ is known.
Find a formula that extracts information about the location
of $D$ from the data $\Lambda_D$ or its partial knowledge.

$\quad$

\subsubsection{The Dirichlet-to-Neumann map acting on CGO solutions}

In this section, we assume that $\gamma_0-1\in
C^{\infty}_0(\Bbb R^2)$ and $\gamma_0(x)>0$ for all $x\in\Bbb
R^2$. Given $\omega\in S^1$ choose $\omega^{\perp}\in S^1$
perpendicular to $\omega$. 
Set
$$\begin{array}{ll}
\displaystyle
z=\tau(\omega+i\omega^{\perp}), & \tau>0.
\end{array}
\tag {3.2}
$$
The solutions given below is called the complex geometrical optics
solutions of the equation $\nabla\cdot\gamma_0\nabla v=0$ in
$\Omega$.

To describe it recall a well known estimate in Calder\'on problem.
Let $z=(z_1,z_2)\in\text{\boldmath C}^2$ satisfy $z\cdot z=z_1^2+z_2^2=0$.
Define the tempered distribution $g_z(x)$ by the formula
$$\displaystyle
g_{z}(x)=\left(\frac{1}{2\pi}\right)^2\int_{\Bbb R^2}\frac{e^{ix\cdot\xi}}{\vert\xi\vert^2-2iz\cdot
\xi}d\xi.
\tag {3.3}
$$
$g_z(x)$ satisfies
$$\begin{array}{ll}
\displaystyle
(\Delta+2z\cdot\nabla)g_{z}(x)+\delta(x)=0, & x\in\Bbb R^2.
\end{array}
$$
Then the distribution
$$\displaystyle
G_{z}(x)=e^{x\cdot z}g_{z}(x)
$$
becomes a fundamental solution of the Laplace equation and called Faddeev's Green function.
This fundamental solution is different from the standard fundamental solution of the Laplace equation.

Given rapidly decreasing function $f$ define the tempered
distribution $g_z\ast f$ by the formula
$$\displaystyle
g_z\ast f(x)=\int_{\Bbb R^2}g_z(x-y)f(y)dy.
$$
Given $s\in\Bbb R$ let us denote the space of all tempered distributions $T$ on $\Bbb R^2$ with norm
$\Vert T\Vert_s=\Vert (1+\vert x\vert^2)^{s/2}T\Vert_{L^2(\Bbb R^2)}$ by $L^2_s(\Bbb R^2)$.

The following estimate is well known.

\proclaim{\noindent Theorem 3.1(\cite{SU1, SU2}).} Let $-1<\delta<0$ and $a>0$.
There exists a positive constant $C_{\delta,a}$ such that, for all
rapidly decreasing function $f$ on $\Bbb R^2$, multi indices
$\alpha$ with $\vert\alpha\vert\le 2$ and $z$ with $\vert
z\vert\ge a\sqrt{2}$
$$\displaystyle
\Vert D^{\alpha}g_z\ast f\Vert_{\delta}\le C_{\delta,a}\vert z\vert^{\vert\alpha\vert-1}\Vert f\Vert_{\delta+1}.
$$

\endproclaim

As a corollary we have
\proclaim{\noindent Proposition 3.1.}
There exist a positive constant $C(\gamma_0)$ and solutions $v(x;z)\in H^2(\Omega)$ with
$\tau\ge C(\gamma_0)$ of the equation $\nabla\cdot\gamma_0\nabla v=0$ in $\Omega$ such that
$$\displaystyle
v(x;z)=\frac{1}{\sqrt{\gamma_0(x)}}
e^{x\cdot z}\{1+\epsilon(x;z)\};
\tag {3.4}
$$
as $\tau\longrightarrow\infty$
$$\displaystyle
\Vert\epsilon(\,\cdot\,;z)\Vert_{L^{\infty}(\Omega)}
=\Vert\nabla\epsilon(\,\cdot\,;z)\Vert_{L^{\infty}(\Omega)}=O(\tau^{-1}).
\tag {3.5}
$$
\endproclaim

{\it\noindent Proof.}
Fix $-1<\delta<0$.  Set
$$
V=\frac{\triangle(\sqrt{\gamma_0})}{\sqrt{\gamma_0}}.
$$
Since $V$ has a compact support, the operator
$L^2_{\delta}(\Bbb R^2)\ni f\longmapsto g_z\ast(Vf)\in L^2_{\delta}(\Bbb R^2)$
is well defined and invertible provided $\tau\ge C(\gamma_0)$ and $C(\gamma_0)$ is
large enough.
Define
$$
\displaystyle
\epsilon(x;z)=-\{I+g_z\ast(V\,\cdot\,)\}^{-1}(g_z\ast V).
$$
A combination of Theorem 3.1 and the Sobolev imbedding theorem
$H^3(\Omega)\longrightarrow C^{1,\lambda}(\overline\Omega)$ for $0<\lambda<1$,
gives (3.5).  Then $v$ defined by (3.4) satisfies the equation
$\nabla\cdot\gamma_0\nabla v=0$ in $\Omega$.

\noindent
$\Box$

The good point of the construction of $v$ given by (3.4) is: it is {\it really} constructive
and does not rely on any transcendental argument, like duality combined with the method of Carleman estimate.
This is necessary for the {\it reconstruction procedure}.

Now we introduce the first indicator function in the enclosure method.

{\bf\noindent Definition 3.1.}
Define the indicator function by the formula
$$
\displaystyle
I_{\omega}(\tau)
=\int_{\partial\Omega}(\Lambda_{\emptyset}-\Lambda_D)
(v(x;z)\vert_{\partial\Omega})
\cdot\overline{v(x;z)}d\sigma(x).
$$
\noindent

The enclosure method enables us to extract the convex hull of $D$ given by
$$\displaystyle
\cap_{\omega\in S^1}\{x\in\Bbb  R^2\,\vert\,x\cdot\omega<h_D(\omega)\},
$$
where $h_D: S^1\ni\omega\mapsto\sup_{x\in\,D}x\cdot\omega$ is 
called the support function of $D$.

\proclaim{\noindent Theorem 3.2(\cite{E01}).}
The formula
$$\displaystyle
\lim_{\tau\longrightarrow\infty}\frac{\log\vert I_{\omega}(\tau)\vert}{2\tau}
=h_D(\omega),
\tag{3.6}
$$
is valid.  Moreover, we have:
$$\displaystyle
\lim_{\tau\longrightarrow\infty}e^{-2\tau t}I_{\omega}(\tau)
=
\left\{
\begin{array}{ll}
\displaystyle
0, & \text{if $t>h_D(\omega)$,}
\\
\\
\displaystyle
\infty, & \text{if $t<h_D(\omega)$,}
\end{array}
\right.
$$
and if $t=h_D(\omega)$, then
$$\displaystyle
\liminf_{\tau\longrightarrow\infty}
e^{-2\tau t}I_{\omega}(\tau)>0.
$$
\endproclaim

{\it\noindent Sketch of Proof.}   First similar to the case when $\gamma_0\equiv 1$ (see, for example, Lemma 1 in 
\cite{ISugakuEnglish}) or the Helmholtz equation case (see Subsubsection 2.3.1) at the 
{\it zero wave number}, using integration by parts, one gets
$$\displaystyle
I_{\omega}(\tau)=\int_{\Omega\setminus\overline{D}}\,\gamma_0\,\vert\nabla w\vert^2\,dx
+\int_D\,\gamma_0\,\vert\nabla v\vert^2\,dx,
$$
where $w=u-v$, $u$ is the solution of (3.1) with $f=v(\,\cdot\,\;z)$ on $\partial\Omega$
and $v=v(\,\cdot\,;z)$.  From this and the upper bound $\Vert w\Vert_{H^1(\Omega\setminus\overline{D})}
\le C\Vert\nabla v\Vert_{L^2(D)}$, we have
$$\displaystyle
C_0 I(\tau)\le e^{-2\tau h_D(\omega)}I_{\omega}(\tau)\le C_D I(\tau),
$$
where
$$\displaystyle
I(\tau)=\int_D\vert\nabla\{e^{-\tau h_D(\omega)}v(x;z)\}\vert^2dx.
$$

A combination of (3.4) and (3.5) gives, for $\tau\ge C_3>>C(\gamma_0)$
$$\begin{array}{ll}
\displaystyle
C_1\tau e^{\tau x\cdot\omega}\le\vert\nabla v(x;z)\vert\le C_2\tau e^{\tau x\cdot\omega},
& x\in\Omega.
\end{array}
$$
where $C_1, C_2$ and $C_3$ are positive constants independent of $\tau$.
From this we have, as $\tau\longrightarrow\infty$
$I(\tau)=O(\tau^2)$.

Choose a point $x_0$ on $\partial D$ such that $x_0\cdot\omega=h_D(\omega)$.
Since $\partial D$ is assumed to be Lipschitz, one can find a finite open cone $C$ with vertex at $x_0$ such that
$C\subset D$.
For a suitable $a>0$, as $s\longrightarrow 0$, we have $\vert C\cap\{x\cdot\omega=h_D(\omega)-s\}\vert\sim as$.
Then, for large $\tau_0$, small $r_0$ and all $\tau\ge\tau_0$, we have
$$\begin{array}{ll}
\displaystyle
I(\tau) & 
\displaystyle\ge C_1^2\tau^2\int_C e^{-2\tau(h_D(\omega)-x\cdot\omega)}dx\\
\\
\displaystyle
&
\displaystyle\ge C_1^2\tau^2\int_0^{r_0} e^{-2\tau s}\vert C\cap\{x\cdot\omega=h_D(\omega)-s\}\vert ds\\
\\
\displaystyle
&
\displaystyle
\ge\frac{C_1^2}{2}a\tau^2\int_0^{r_0} e^{-\tau s}sds.
\end{array}
$$
Now all the statements on the indicator function are clear.

\noindent
$\Box$

It was Siltanen \cite{SF} who gave a numerical computation procedure of Faddeev's Green function or (3.3)
and in \cite{SMI} a numerical implementation of Nachman's formula \cite{N2} in the Calder\'on problem has been done.

In \cite{E01} the case when $\gamma_0=1$, $\Omega\subset\Bbb R^2$ and the governing equation is replaced with the Helmholtz equation
$\Delta u+k^2u=0$ with a fixed $k\ge 0$ has been considered.  In that case it is assumed that
the set $\{x\in\Bbb R^3\,\vert\,x\cdot\omega=h_D(\omega)\}\cap\partial D$ consists of a single point
where the Gauss curvature of $\partial D$ is non vanishing.  Note that when $k=0$ we do not need such an assumption.
Later Sini-Yoshida \cite{SY} and \cite{SY2} removed 
the condition in the case when $k>0$.  Note that, however, if the boundary condition is replaced with $u=0$ on $\partial D$, 
then a result analogous to Theorem 3.2 has been established in \cite{E00} without the condition mentioned above.
This indicates that when $k>0$ the Neumann boundary condition case is harder than the Dirichlet one.

In \cite{IINSU} a combination of the enclosure method explained here and the complex geometrical optics solution constructed by
using the hyperbolic geometry has been introduced.

Recently the enclosure method presented here has been extended to an inverse obstacle problem
governed by the $p$-Laplace equation, see \cite{BHKS}.

Finally we mention \cite{SI} in which the enclosure method applied to Problem 2.2, that is the result in \cite{E00} and the same type as this section,
has been combined with the {\it machine learning}.  It is reported that,
compared with the original algorithm based on (3.6)-type formula in \cite{IS0} 
the accuracy of the enclosure method was significantly improved.

\subsubsection{Using a singe input}

In this subsection we always assume that $\gamma_0$ is given by a known constant.

{\bf\noindent Definition 3.2.}
Given $\omega=(\omega_1,\omega_2)\in S^1$
set $\omega^{\perp}=(\omega_2,-\omega_1)\in S^1$.
Let $\tau>0$.
Define the indicator function
$$\displaystyle
J_{\omega}(\tau)
=\int_{\partial\Omega}(\Lambda_{\emptyset}-\Lambda_D)f(x)
\cdot e^{\tau x\cdot(\omega+i\omega^{\perp})} d\sigma(x).
\tag{3.7}
$$

\noindent
In contrast to the indicator function in Definition 3.1 this is complex valued.

{\bf\noindent Definition 3.3.}
We say that $\omega\in S^1$ is regular with respect to $D$ if the set
$$\displaystyle
\partial D\cap\{x\in\Bbb R^2\,\vert\,x\cdot\omega=h_D(\omega)\}
$$
consists of only one point.

{\bf\noindent Definition 3.4.}
We say that $D$ is a polygonal cavity if
$D$ takes the form $D_1\cup\cdots\cup D_m$ with $1\le m<\infty$
where each $D_j$ is open and a polygon; $\overline D_j\cap\overline D_{j'}=\emptyset$ if $j\not=j'$.

In this subsection we  describe the following result.
\proclaim{\noindent Theorem 3.3(\cite{I1}).}
Let $f$ be a non constant function.  
Assume that $D$ is a polygonal cavity and satisfies
$$\displaystyle
\text{diam}\, D<\text{dist}\,(D,\partial\Omega).
\tag{3.8}
$$
Let $\omega$ be regular with respect to $D$.
Then the formula
$$\displaystyle
\lim_{\tau\longrightarrow\infty}\frac{\log\vert J_{\omega}(\tau)\vert}{\tau}
=h_D(\omega),
\tag{3.9}
$$
is valid.  Moreover, we have:
$$\displaystyle
\lim_{\tau\longrightarrow\infty}e^{-\tau t}\vert J_{\omega}(\tau)
\vert=
\left\{
\begin{array}{ll}
\displaystyle
0 & \text{if $t\ge h_D(\omega)$,}\\
\\
\displaystyle
\displaystyle
\infty & \text{if $t< h_D(\omega)$.}
\end{array}
\right.
$$

\endproclaim

{\it\noindent Skectch of Proof.}
Let $v=e^{\tau\,x\cdot(\omega+i\omega^{\perp})}$ and $u$ solve (3.1).
From the expression
$$\displaystyle
J_{\omega}(\tau)=\gamma_0\int_{\partial D}\,u\frac{\partial v}{\partial\nu}\,d\sigma,
$$
we see that it suffices to study the asymptotic behaviour of the function
$$\displaystyle
\tau\mapsto e^{-\tau h_D(\omega)}\int_{\partial D}\,u\frac{\partial v}{\partial\nu}\,d\sigma.
\tag {3.10}
$$
Using the exponential decaying property of the function $e^{-\tau h_D(\omega)}v$ in $x\cdot\omega<h_D(\omega)-\delta$ with a small $\delta>0$,
one can localize the integral around a single corner point $x_0$ on the convex hull of $D$ by virtue of the assumption that $\omega$ is regular.
Then using the {\it convergent} series expansion of $u$ which is a consequence of the eigenfunction expansion of the ordinary differential equation
with the homogeneous Neumann boundary condition (see \cite{Gr}), one gets the complete asymptotic expansion of function (3.10)
as $\tau\rightarrow\infty$ having the form
$$\displaystyle
e^{-\tau h_D(\omega)}\int_{\partial D}\,u\frac{\partial v}{\partial\nu}\,d\sigma
\sim
e^{i\tau\,x_0\cdot\omega^{\perp}}\sum_{m=1}^{\infty} \frac{A_m}{\tau^{\mu_m}},
$$
where $0<\mu_1<\mu_2<\cdots\,\rightarrow\infty$ and $A_m$ are constants.

What we have to do here is: there exists a $m\ge 1$ such that $A_m\not=0$.
This is a non trivial statement since coefficients $A_m$ contain those of series expansion of $u$
and the geometry of $\partial D$ around $x_0$.
This is proved by a contradiction argument with the help of expansion of $u$ around $x_0$,
the precise structure of $A_m$ and Friedman-Isakov's extension argument \cite{FI}.

\noindent
$\Box$

Some remarks are in order.

$\bullet$  Even now, removing condition (3.8) remains open.  However, the proof in \cite{I1}
tells us that condition (3.8) is redundant if $\frac{\Upsilon}{\pi}$ is an {\it irrational} number, where
$\Upsilon$ denotes the {\it interior angle} of $D$ at the vertex $x_0$ in Sketch of Proof.  
However, such an information is {\it unstable}, so we don not think 
imposing such an unstable restriction is a right way. Mathematically, the interesting and valuable thing
is finding an argument to prove $A_m\not=0$ for a $m\ge 1$ without making use of condition (3.8)
in the case when $\frac{\Upsilon}{\pi}$ is a {\it rational number}.

$\bullet$  The set of all directions that are not regular with respect to given polygonal cavity $D$ is a finite set.
Therefore we do not have to worry about the choice of directions that are not regular.

$\bullet$  In \cite{IConductivity} the idea of the proof of Theorem 3.3 has been applied to Problem 2.2 provided $\gamma_0$ is constant,
$\gamma(x)$ on $D$ is also an unknown constant different from $\gamma_0$, $D$ is polygonal.  
See also \cite{Ipartial} for an analogous result using the data $u(P)-u(Q)$ with fixed two points $P$ and $Q$ on $\partial\Omega$, where $u$ solves
$\nabla\cdot\gamma_0\nabla u=0$ in $\Omega$ with the Neumann boundary condition
$$\begin{array}{ll}
\displaystyle
\gamma_0\frac{\partial u}{\partial\nu}=\gamma_0\frac{\partial}{\partial\nu}e^{\tau\,x\cdot(\omega+i\omega^{\perp})}, & x\in\partial\Omega.
\end{array}
$$

$\bullet$  An application to a similar inverse boundary value problem for an elastic body in two dimensions also has been done in a crack case 
\cite{IICrackIso, IICrackAniso} and cavity case \cite{IICavity}.  The main task is to establish a convergent series expansion
of the solution of the governing equation around a corner which is not trivial.  It was done by using an analytic continuation of 
a stress function \cite{Mu} and the method of the Mellin transform, respectively.

$\bullet$  Numerical implementation of an algorithm based on formula (3.9) has been done in \cite{IO0}.  The basic idea shares with that
of \cite{IS0}:  finding a least square fitting line that is approximately passing the calculated points
$(\tau_j,\log\vert J_{\omega}(\tau_j)\vert)$ with $0<\tau_1<\cdots<\tau_m$.  Then consider its slope as $h_D(\omega)$.
See also \cite{IO1} in which the same idea has been tested for an algorithm based on the formula in \cite{IConductivity}.

\subsubsection{Unknown constant background case}

In this section we consider the case when $\gamma_0$ is given by an {\it unknown} constant.
We have only the {\it qualitative information} that the background body is uniformly conductive.
How does the enclosure method work for this case ?

We can not make use of the full form of indicator function (3.7) since $\Lambda_{\emptyset}f$ is unknown.
So, instead of it, here we use only the term 
$$\displaystyle
\int_{\partial\Omega}\,\Lambda_Df\cdot e^{\tau x\cdot(\omega+i\omega^{\perp})}\,d\sigma(x).
$$

In \cite{IEnclosureCubo} the author found that: under some assumptions on $f$ and the geometry of $\Omega$
one can establish a formula for some information about the geometry of $D$.

The main result is as follows.

Let $\Omega$ be the domain enclosed by an ellipse.
By choosing a suitable system of orthogonal coordinates one can write
$$\displaystyle
\Omega=\left\{(x_1,x_2)\,\vert\,\left(\frac{x_1}{a}\right)^2+\left(\frac{x_2}{b}\right)^2<1\right\},
$$
where $a\ge b>0$.  In what follows we always use this coordinates system.

We denote by $E(\Omega)$ the set of all points on the segment that connects the focal points
$(-\sqrt{a^2-b^2},\, 0)$ and $(\sqrt{a^2-b^2},\,0)$ of $\Omega$.

Let $h_{E(\Omega)}$ denote the support function of set $E(\Omega)$.
We have $h_{E(\Omega)}(\omega)=\sqrt{a^2-b^2}\,\vert\omega_1\vert$, where
$\omega_1=\omega\cdot\mbox{\boldmath $e$}_1$.

Given a function $f$ on $\partial\Omega$ 
write
$$\begin{array}{ll}
\displaystyle
f(\theta)=f(a\cos\theta,b\sin\theta), \theta\in\Bbb R.
\end{array}
$$
Its Fourier expansion takes the form
$$\displaystyle
f(\theta)=\frac{1}{2}\alpha_0+\sum_{m=1}^{\infty}(\alpha_m\,\cos\,m\theta+\beta_m\sin m\theta),
$$
where
$$\left\{
\begin{array}{l}
\displaystyle
\alpha_m=\frac{1}{\pi}\,\int_0^{2\pi}\,f(\theta)\cos\,m\theta\,d\theta,
\\
\\
\displaystyle
\beta_m=\frac{1}{\pi}\,\int_0^{2\pi}\,f(\theta)\sin\,m\theta\,d\theta.
\end{array}
\right.
$$
Set
$$
\begin{array}{llll}
\displaystyle
c_0=\frac{\alpha_0}{2},
&
\displaystyle
c_m=\frac{\alpha_m-i\beta_m}{2},
&
\displaystyle
c_{-m}=\overline{c_m},
&
\displaystyle
m\ge 1,
\end{array}
$$
$$
\displaystyle
A_{\pm}=\frac{1}{2}
\left(\frac{1}{a}\pm\frac{1}{b}\right)
$$
and
$$
\displaystyle
C_m(f)=
\left\{
\begin{array}{ll}
\displaystyle
A_{-}\overline{c_1}+A_+c_1, & m=0,\\
\\
\displaystyle
(A_{-}c_{m-1}+A_+c_{m+1})\left(\frac{a+b}{a-b}\right)^{\frac{m}{2}}
+(A_{-}\overline{c_{m+1}}+A_+\overline{c_{m-1}}\,)\left(\frac{a-b}{a+b}\right)^{\frac{m}{2}},
&
m\ge 1.
\end{array}
\right.
$$

We say that $f$ is {\it band limited} if
there exists a natural number $N\ge 1$ such that, for all $m\ge N+1$ 
$\alpha_m=\beta_m=0$.  Thus we have $C_m(f)=0$ for all $m\ge N+2$.

\proclaim{\noindent Theorem 3.4(\cite{IEnclosureCubo}).}  Assume that $D$ is polygonal and satisfies (3.8).  
Let $\omega$ be regular with respect to $D$.  Let $f$ be band limited.

\noindent
(1)  Let $a>b$.  Let $\omega$ satisfy $\omega_1\not=0$.
Let $f$ satisfy
$$\displaystyle
\sum_{m=1}^{\infty}\,(\text{sgn}\,\omega_1)^m\,m^2\,C_m(f)\not=0.
\tag {3.11}
$$
Then, we have
$$\displaystyle
\lim_{\tau\rightarrow\infty}
\frac{1}{\tau}
\,
\log\left\vert\int_{\partial\Omega}\,\Lambda_Df\cdot e^{\tau x\cdot(\omega+i\omega^{\perp})}\,d\sigma(x)\right\vert
=\max\,(h_D(\omega),\,h_{E(\Omega)}(\omega)).
\tag {3.12}
$$

\noindent
(2)  Let $a=b$.  Let $f$ satisfy: for some $N\ge 1$ $\alpha_m=\beta_m=0$ for all $m$ with $m\ge N+1$
and $\alpha_N^2+\beta_N^2\not=0$.  Then, we have (3.12).

\endproclaim

{\it\noindent Sketch of proof.}  We only mention (1).
Let $v=e^{\tau\,x\cdot(\omega+i\omega^{\perp})\,}$ and $u$ solve (3.1).
Integration by parts yields
$$\displaystyle
\int_{\partial\Omega}\Lambda_Df\cdot v\,d\sigma=\gamma_0\int_{\partial\Omega}\,f\frac{\partial v}{\partial\nu}\,d\sigma-\gamma_0\int_{\partial D}\,u\frac{\partial v}{\partial\nu}\,d\sigma.
\tag {3.13}
$$
In the proof of Theorem 3.3 we have already known the asymptotic behaviour of the second term on this right-hand side
provided $f$ is not a constant function.  Thus the main task is to study the first term, that is, the integral
$$\displaystyle
\int_{\partial\Omega}\,u\frac{\partial v}{\partial\nu}\,d\sigma
=\tau(\omega_1+i\omega_2)\int_{\partial\Omega}\,fv(\nu_1-i\nu_2)\,d\sigma,
$$
where $\omega=(\omega_1,\omega_2)$ and $\nu=(\nu_1,\nu_2)$.

We can explicitly compute the integral on this right-hand side ant the result is
$$\displaystyle
\int_{\partial\Omega}\,fv(\nu_1-i\nu_2)\,d\sigma
=2\pi ab\sum_{m=0}^{\infty}\,i^m
J_m(-i\sqrt{a^2-b^2}\,\tau(\omega_1+i\omega_2)\,)\,C_m(f).
$$
Then, applying the {\it compound asymptotic expansion} of the Bessel function due to Hankel
(see page 133 in  \cite{Ol}) to the term on this right-hand side and noting the equation
$$\displaystyle
\sum_{m=0}^{N+1}\,C_m(f)=0,
$$
one gets, as $\tau\rightarrow\infty$
$$\displaystyle
\int_{\partial\Omega}\,fv(\nu_1-i\nu_2)\,d\sigma
=i\pi ab z^{-1}
\left(\frac{1}{2\pi z}\right)^{1/2}
e^{iz}e^{-i\frac{\pi}{4}}
\left(
\sum_{m=1}^{N+1}m^2C_m(f)+O(\tau^{-1})\right),
\tag {3.14}
$$
where $z=-i\sqrt{a^2-b^2}\,\tau(\omega_1+i\omega_2)$.  Note that $e^{iz}=e^{\tau h_{E(\Omega)}(\omega)}e^{i\tau\sqrt{a^2-b^2}\,\omega_2}$.
Hereafter carefully comparing the asymptotic behaviour of the right-hand side on (3.14) with the second term on (3.13), one concludes
the quantity
$$\displaystyle
e^{-\tau\max\,(h_D(\omega),h_{E(\Omega)}(\omega))}\,\left\vert\int_{\partial\Omega}\,\Lambda_Df\cdot v\,d\sigma\right\vert
$$
is truly algebraic decaying as $\tau\rightarrow\infty$.

\noindent
$\Box$

A typical example of a band limited $f$ that satisfies (3.11) for all $\omega$ with $\omega_1\not=0$ is the $f$ 
given by
$$\displaystyle
f(\theta)=A\cos\,N\theta+B\sin\,N\theta,
$$
where $N\ge 1$ and $A^2+B^2\not=0$.  Note that this $f$ is {\it independent} of $\text{sgn}\,\omega_1$.  This means that,
using only the single $f$, we have the formula (3.12) for all regular $\omega$ with $\omega_1\not=0$.
In general the condition (3.11) depends on $\text{sgn}\,\omega_1$, so we need two $f$s.

Since the function $\omega\mapsto \max\,(h_D(\omega), h_{E(\Omega)}(\omega))$ is continuous and 
the set of all directions that are not regular with respect to polygonal $D$ is a finite set, we have the following uniqueness result.

\proclaim{\noindent Corollary 3.1 (\cite{IEnclosureCubo}).}
Assume that $D$ is polygonal and satisfies (3.8).  

(1)  Let $a>b$.  Let $f_+$ and $f_{-}$ be band limited and satisfy
$$\displaystyle
\sum_{m=1}^{\infty}\,(\pm)^m\,m^2\,C_m(f_{\pm})\not=0.
$$
Then the data $\Lambda_Df_{+}$ and $\Lambda_{D}f_{-}$ uniquely determine the convex hull
of the set $D\cup\,E(\Omega)$.

(2)  Let $a=b$.  Let $f$ be band limited and non constant.  Then the data $\Lambda_Df$ uniquely
determines the convex hull of the set $D\cup\,\{(0,0)\}$.

\endproclaim

The author thinks that this is the first uniqueness result using at most two Neumann data
in the case when the background conductivity is unknown.  And one can easily understand of
the difficulty of applying usual {\it traditional} step, that is,
the uniqueness of the Cauchy problem, in the proof of uniqueness
in inverse obstacle problems with a single input, like those of \cite{FI, FV}, since the background conductivity
is unknown.

Theorem 3.4 together with Corollary 3.1 has been extended also to the case when
the governing equation of $u$ takes the form
$\nabla\cdot\gamma\nabla u=0$ in $\Omega$, where $\gamma(x)=\gamma_0$ if $x\in\Omega\setminus D$
and $\gamma(x)$, $x\in D$ given by an unknown positive constant different from $\gamma_0$.
See \cite{IEnclosureCubo}.

\subsubsection{Using a Carleman function in an unbounded domain}

As a typical simplification of the geometry of domain we consider a domain between two parallel planes.
Let $\Omega=\Bbb R^2\times\,]0,\,h[$ with a fixed $h>0$. 
This is an infinite slab.
We assume that the conductivity of $\Omega$ has the form
$$\displaystyle
\gamma(x)
=
\left\{
\begin{array}{ll}
\displaystyle
I_3, & x\in\Omega\setminus D,\\
\\
\displaystyle
I_3+h(x), & x\in D,
\end{array}
\right.
$$
where $D$ is an open subset of $\Omega$ with Lipschitz boundary and $h$ is an essentially bounded real symmetric matrix-valued
function on $D$; $I_3$ is the $3\times 3$-identity matrix; $\gamma$ is uniformly positive definite over $\Omega$.

We assume that $D$ is {\it bounded}.  This means that $D$ is a model of inclusion existing only a small part of $\Omega$. 
And the assumption that $\Omega$ is an infinite slab means: we only pay attention to a some restricted part of a given large conductive body 
and the effect from the other part can be ignored.  So, this is a kind of a localization of the original situation.

In \cite{Islab} we considered the problem of extracting information about the geometry of $D$ from the well-defined
Dirichlet-to-Neumann map $\Lambda_{\gamma}$ on $\partial\Omega$, where $\Lambda_{\gamma}: H^{1/2}(\partial\Omega)\rightarrow H^{-1/2}(\partial\Omega)$.
In the case when $\Omega$ is a bounded domain, as before like Theorem 3.2,
under a suitable jump condition for $\gamma$ on $\partial D$ depending on a given direction $\omega\in S^2$
one can extract the value of the support function for $D$ at direction $\omega$ from
the indicator function \cite{E00, IS0}
$$\displaystyle
I_{\omega,\omega^{\perp}}(\tau)=\int_{\partial\Omega}\,(\Lambda_{\gamma}-\Lambda_1)\,
(v(x;z)\vert_{\partial\Omega})\cdot
\overline{v(x;z)}\,dS(x),
$$
where $v(x;z)=e^{x\cdot z}$ and $z=\tau(\omega+i\omega^{\perp})$, $\tau>0$ and $\omega^{\perp}\in S^2$ with $\omega\cdot\omega^{\perp}=0$.

However, in the present situation that $\Omega$ is an infinite slab, not only the trace of $v(\,\cdot\,;z)$ onto $\partial\Omega$
does not belong to $H^{1/2}(\partial\Omega)$ but also $v(\,\cdot\,;z)$ itself not to $H^1(\Omega)$ where the direct problem
is well-posed.  So it is quite natural to propose the following problem.

$\quad$

{\bf\noindent Problem 3.2.}
What function does play a role similar to $v(\,\cdot\,;z)$ ?

$\quad$

In \cite{Islab} the author considered this and gave a solution.
The point is: change the role of the Carleman function.

The findings therein are as follows.

$\bullet$  Originally the Carleman function gives a direct computation formula of a solution of the Laplace equation (in the present situation)
in a domain using the Cauchy data given on a part of the boundary.
Instead use the formula as an approximation with a good decaying property for a given solution of the Laplace equation.

$\bullet$ For the purpose a special member of Yarmukhamedov's fundamental solution \cite{Y2} denoted by $\Phi_K$ in Theorem 2.8 is suitable.

Let us explain more.  We choose $K(w)=e^{m\,w^2}$ with a parameter $m>0$ which fulfiles (2.33).
Then $\Phi_m\equiv \Phi_K$ takes the form
$$\displaystyle
\Phi_m(x)
=-\frac{e^{m(x_3^2-\vert x'\vert^2)}}{2\pi^2}\,\int_0^{\infty} e^{-m\,u^2}\,\frac{w_m(x,u)\,du}{\vert x\vert^2+u^2},
$$
where
$$\displaystyle
w_m(x,u)
=\frac{x_3\,\sin\,2m\,x_3\sqrt{\vert x'\vert^2+u^2}}{\sqrt{\vert x'\vert^2+u^2}}
-\cos\,2m\,x_3\sqrt{\vert x'\vert^2+u^2}.
$$

Recall $\Omega=\Bbb R^2\times\,]0,\,h[$.
Given $\epsilon>0$ and $R>0$ define
$$\left\{
\begin{array}{l}
\displaystyle
\Omega_0(R)=\left\{(x',x_3)\,\vert\,\vert x'\vert<R, 0<x_3<h\,\right\},\\
\\
\displaystyle
\Omega_{\epsilon}(R)=\left\{(x',x_3)\,\vert\,\vert x'\vert<R, -\epsilon<x_3<h\,\right\}.
\end{array}
\right.
$$ 
We have the decomposition of the boundary 
$$\displaystyle
\partial\Omega_{\epsilon}(R)=\Gamma_{\epsilon}^+(R)\cup\,\Gamma_{\epsilon}^-(R)\cup\,\Gamma_{\epsilon}'(R),
$$
where 
$$\left\{
\begin{array}{l}
\displaystyle
\Gamma_{\epsilon}^+(R)=\left\{(x',x_3)\,\vert\,\vert x'\vert<R,\,x_3=h+\epsilon\,\right\}\\
\\
\displaystyle
\Gamma_{\epsilon}^-(R)=\left\{(x',x_3)\,\vert\,\vert x'\vert<R,\,x_3=-\epsilon\,\right\}\\
\\
\displaystyle
\Gamma_{\epsilon}'(R)=\left\{(x',x_3)\,\vert\,\vert x'\vert=R,\,-\epsilon\le x_3\le h+\epsilon\,\right\}.
\end{array}
\right.
$$
Given $0<\delta<1$ define
$$\left\{
\begin{array}{l}
\displaystyle
C_y(\delta)
=\text{int}\,\left(C_y(\mbox{\boldmath $e$}_3, \alpha)
\cup C_y(-\mbox{\boldmath $e$}_3, \alpha)\right),\\
\\
\displaystyle
C_{\epsilon}^{\delta}(R)=\cup_{y\in\,\Gamma_{\epsilon}'(R)}\,C_y(\delta),
\end{array}
\right.
$$
where $C_y(\pm\mbox{\boldmath $e$}_3,\alpha)=C_y(\omega,\alpha)$ with $\omega=\pm\mbox{\boldmath $e$}_3=\pm(0,0,1)^T$
is the three-dimensional cone about  $\omega$ of opening angle $\frac{\pi\alpha}{2}$ with vertex at $y$
and
$$\displaystyle
\alpha=1-\frac{2}{\pi}\,\tan^{-1}\delta.
$$

The idea of this section is based on the formula of Carleman type stated below.

\proclaim{\noindent Theorem 3.5 (\cite{Islab}).}
Let $v$ be a continuously differentiable function in a neighbourhood of $\Omega_{\epsilon}(R)$
and satisfy the Laplace equation in $\Omega_{\epsilon}(R)$.
The formula
$$
\displaystyle
v(x)
=\lim_{m\rightarrow\infty}\,
\int_{\Gamma_{\epsilon}^+(R)\cup\,\Gamma_{\epsilon}^-(R)}
\left(\Phi_m(y-x)\frac{\partial v}{\partial\nu(y)}
-v\,\frac{\partial}{\partial\nu(y)}\,\Phi_m(y-x)\right)\,dS(y),
\tag {3.15}
$$
is valid for all $x\in\Omega_0(R)\setminus C_{\epsilon}^{\delta}(R)$ for each fixed $\delta>0$.
The convergence is uniform on $\Omega_0(R)\setminus C_{\epsilon}^{\delta}(R)$ together
with its first-order derivatives.

\endproclaim

{\it\noindent Skectch of Proof.}
Since $\Phi_m$ is a fundamental solution of the Laplace equation,
we have
$$
\displaystyle
v(x)
=
\int_{\Gamma_{\epsilon}^+(R)\cup\,\Gamma_{\epsilon}^-(R)\cup\Gamma_{\epsilon}'(R)}
\left(\Phi_m(y-x)\frac{\partial v}{\partial\nu(y)}
-v\,\frac{\partial}{\partial\nu(y)}\,\Phi_m(y-x)\right)\,dS(y).
$$
The role of $\delta$ is clearly indicated by the decaying property of $\Phi_m$ as $m\rightarrow\infty$:
if $x\not\in C_{(0,0,0)}(\delta)$ and $\vert x'\vert\ge c_1$ and $\vert x_3\vert\le c_2$
with some $c_1,c_2>0$, then $\Phi_m(x)=O(m e^{-m(1-\delta^2)\,\vert x'\vert^2})$ (Lemma 2.1 in \cite{Islab}).
Roughly speaking, using this property we can omit the contribution of the Cauchy data on $\Gamma_{\epsilon}'(R)$
in the formula above.

\noindent
$\Box$

Theorem 3.5 means that the values of the harmonic function in $\Omega_0(R)\setminus C_{\epsilon}^{\delta}(R)$
can be uniquely and constructively determined by its Cauchy data on $\Gamma_{\epsilon}^+(R)\cup\Gamma_{\epsilon}^-(R)$.
However, the values on the set $\Omega_0(R)\cap C_{\epsilon}^{\delta}(R)$ are {\it out of control} by the same Cauchy data.

So based on the formula (3.15) we introduce the harmonic function in a neighbourhood of $\overline\Omega$ by the formula
$$\displaystyle
e_m(x;\tau,\omega,\omega^{\perp})
=\int_{\Gamma_{\epsilon}^+(R)\cup\,\Gamma_{\epsilon}^-(R)}
\left(\Phi_m(y-x)\frac{\partial v}{\partial\nu(y)}
-v\,\frac{\partial}{\partial\nu(y)}\,\Phi_m(y-x)\right)\,dS(y),
$$
where $v(y)=e^{y\cdot z}$.

By Theorem 3.5, we see that the function $e_m$ approximates $e^{x\cdot z}$ in $\Omega_0(R)\setminus C_{\epsilon}^{\delta}(R)$ as $m\rightarrow\infty$.
Besides we see that, if $x\in\Omega$, then function $e_m$ is rapidly decaying as $\vert x'\vert\rightarrow\infty$ and thus
$e_m\in H^1(\Omega)$.  So we can consider $\Lambda_{\gamma}$ acting on the trace of $e_m$ onto $\partial\Omega$
and define the indicator function 
$$\displaystyle
I_m(\tau;\omega,\omega^{\perp})=\int_{\partial\Omega}\,(\Lambda_{\gamma}-\Lambda_1)\,
(e_m\vert_{\partial\Omega})\cdot
\overline{\varphi(x')e_m}\,dS(x),
$$
where $\varphi\in C_0^{\infty}(\Bbb R^2)$ satisfies $\varphi(x')=1$ in a some large bounded domain which contains
$\vert x'\vert\le R$.  Under the a-priori assumption $D\subset \Omega_0(R)\setminus C_{\epsilon}^{\delta}(R)$,
by taking limit $m\rightarrow\infty$ and next limit $\tau\rightarrow\infty$ one can extract the
value of $h_D(\omega)$ from the indicator function above.

In two dimensions we have two results based on the formula of the Carleman type \cite{Ytwo}:
one is an extension of Theorem 3.2 to the case when $\Omega$ is an infinite strip; another is that of Theorem 3.3
to the same case.  See \cite{ITwomethods} for the precise statements together with their proofs.

The method presented here provides us how to use some special Carleman function {\it itself} in inverse obstacle problems.
It is different from another one \cite{IProbeCarleman} as a {\it generator} of an explicit needle sequence in the probe method
since therein only its {\it regular part} is used.

\subsection{Reflection, refraction and Kelvin transform}

\subsubsection {Discontinuity in a layered medium}

Let $\Omega\subset\Bbb R^n$, $n=2,3$ be a bounded domain with Lipschitz boundary.
We assume that the conductivity of $\Omega$ is given by $\gamma$ having the form
$$\displaystyle
\gamma(x)
=
\left\{
\begin{array}{ll}
\displaystyle
\gamma_0(x)\,I_n, & x\in\Omega\setminus D,\\
\\
\displaystyle
\gamma_0(x)\,I_n+h(x), & x\in D,
\end{array}
\right.
\tag {3.16}
$$
where $D$ is an open subset of $\Omega$ with Lipschitz boundary and $h$ is an essentially bounded real symmetric matrix-valued
function on $D$; $I_n$ is the $n\times n$-identity matrix;
$\gamma\in L^{\infty}(\Omega)$ and both $\gamma$ and $\gamma_0$ are uniformly positive definite over $\Omega$.

\begin{figure}[htbp]
\begin{center}
\epsfxsize=8cm
\epsfysize=4cm
\epsfbox{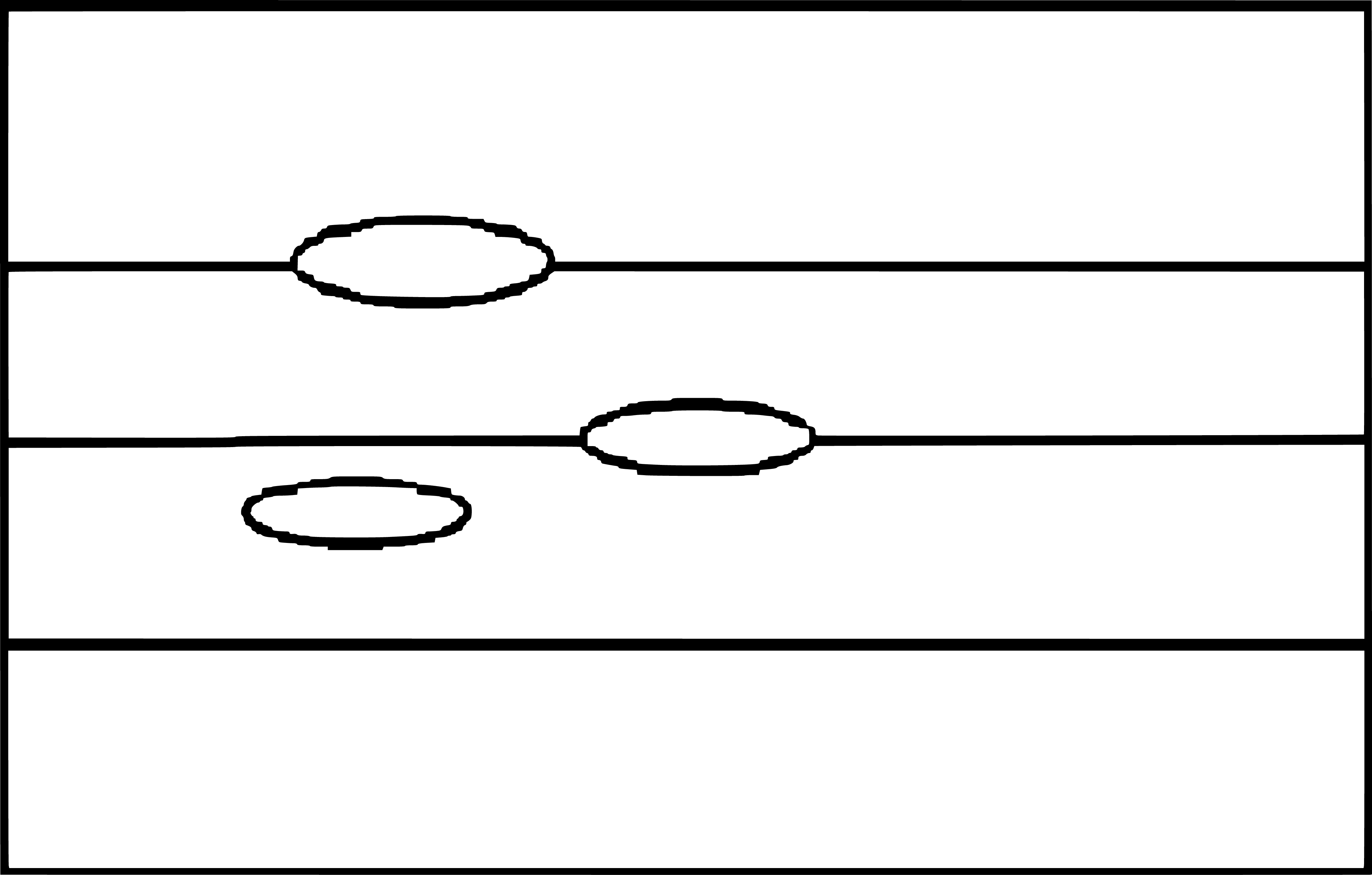}
\caption{An illustration of inclusions embedded in $\Omega$.
All the inclusions  are indicated by the area enclosed by the curves.
}\label{fig4}
\end{center}
\end{figure}

In this section we consider Problem 2.2 in the case when $\gamma_0$ has a $m$-layered structure.
More precisely, we assume that $\gamma_0$ is given by the restriction
to $\Omega$ of the function $\tilde{\gamma_0}(x)$:
$$\begin{array}{lll}
\displaystyle
\tilde{\gamma_0}(x)=\gamma_j, & x\in\,\Bbb R^{n-1}\times\,I_j, & j=1,\cdots, m.
\end{array}
\tag {3.17}
$$
Here $\gamma_1,\gamma_2,\cdots,\gamma_m$, $m\ge 2$ are positive constants, $I_j=]c_j,\,c_{j-1}[$,
$\infty=c_0>c_1>c_2>\cdots>c_{m-1}>c_m=-\infty$ and they are all {\it known}.
See Figure \ref{fig4} for an illustration of inclusions embedded in $\Omega$ which are
denoted by three areas surrounded by closed curves.

\vspace{0.0cm}

In what follows, for simplicity, we use the same symbol $\gamma_0$ for $\tilde{\gamma_0}$.

Given the direction $\omega\in S^{n-1}$ and $\delta>0$ set
$$\displaystyle
D_{\omega}(\delta)=\{x\in D\,\vert\,h_D(\omega)-\delta<x\cdot\omega<h_D(\omega)\}.
$$

We say that $\gamma$ has a jump from $\gamma_0$ on $\partial D$ from direction $\omega$ if
the one of the following two conditions is satisfied.

$\bullet$  There exist constants $C=C(\omega)>0$ and $\delta=\delta(\omega)$ such that, for almost all $x\in D_{\omega}(\delta)$
the lowest eigen value of $h(x)$ is greater than $C$.

$\bullet$  There exist constants $C=C(\omega)>0$ and $\delta=\delta(\omega)$ such that, for almost all $x\in D_{\omega}(\delta)$
the lowest eigen value of $-h(x)$ is greater than $C$.

In \cite{ILayered}, using the enclosure method, the author established the following theorem.

\proclaim{\noindent Theorem 3.6(\cite{ILayered}).}  Let $\omega\in S^{n-1}$ satisfy $\omega\cdot\mbox{\boldmath $e$}_n\not=0$.
Assume that $\gamma$ has a jump from $\gamma_0$ on $\partial D$ from the direction $\omega$ and that $\partial D$ is $C^2$.
Then one can extract $h_D(\omega)$ from $\Lambda_{\gamma}$.

\endproclaim

The key point is to: construct a solution of the equation $\nabla\cdot\gamma_0\nabla v=0$ in $\Bbb R^n$
which plays the same role as the complex geometrical optics solution in the case when $\gamma_0$ is smooth
as done in Theorem 3.2.

Their construction is analogous to that of the reflected/refracted plane wave solutions caused by the presence of multilayers
having different propagation speeds of sound wave.

The construction is as follows.  Given $\omega\in S^{n-1}$ let $\omega^{\perp}\in S^{n-1}$ be a vector perpendicular to $\omega$.
Let $z$ be the complex vector exactly same as (5.0).
Here we only describe the case when $\omega\cdot\mbox{\boldmath $e$}_n>0$, that is, $\text{Re}\,z\cdot\mbox{\boldmath $e$}_n>0$.

Let $e^+(x)=e^+(x;\gamma_0,z)$ be a function having the form
$$\displaystyle
e^+(x;\gamma_0,z)
=
\left\{
\begin{array}{ll}
\displaystyle
e^{x\cdot z}+B_1e^{K_{c_1}(x)\cdot z}, & x\in\Bbb R^{n-1}\times I_1,\\
\\
\displaystyle
A_2e^{x\cdot z}+B_2e^{K_{c_2}(x)\cdot z}, & x\in\Bbb R^{n-1}\times I_2,
\\
\\
\displaystyle
\vdots &
\\
\\
\displaystyle
A_{m-1}e^{x\cdot z}+B_{m-1}e^{K_{c_{m-1}}(x)\cdot z}, & x\in\Bbb R^{n-1}\times I_{m-1},
\\
\\
\displaystyle
A_me^{x\cdot z}, & x\in\Bbb R^{n-1}\times I_m,
\end{array}
\right.
\tag {3.18}
$$
where $K_{c_j}$ is the {\it reflection} across the plane $x_n=c_j$, that is the map
$$\displaystyle
K_{c_j}:(x_1,\cdots,x_{n-1},x_n)\mapsto\,(x_1,\cdots,x_{n-1}, 2c_j-x_n).
$$
The coefficients $A_2,\cdots, A_m$ and $B_1,\cdots, B_{m-1}$ are unknown constants to be determined in such a way that
$e^+$ belongs to $H^1_{\text{loc}}\,(\Bbb R^n)$ and satisfies the equation
$\nabla\cdot\gamma_0\nabla e^+=0$ in $\Bbb R^3$ in the weak sense, that is, for all $\varphi\in C^{\infty}_0(\Bbb R^3)$
$$\displaystyle
\int_{\Bbb R^n}\,\gamma_0\nabla e^+\cdot\nabla\varphi\,dx=0.
$$
Since in each layer function $e^{+}$ satisfies the equation, those conditions are satified if and only if, for $j=1,\cdots,m-1$
the equations
$$\left\{
\begin{array}{l}
\displaystyle
\lim_{x_\downarrow c_j}\,e^+(x;\gamma_0,z)=\lim_{x_n\uparrow c_j}\,e^+(x;\gamma_0,z),\\
\\
\displaystyle
\lim_{x_\downarrow c_j}\,\gamma_j\frac{\partial}{\partial x_n}\,e^+(x;\gamma_0,z)=
\lim_{x_n\uparrow c_j}\,\gamma_{j+1}\frac{\partial}{\partial x_n}e^+(x;\gamma_0,z),
\end{array}
\right.
\tag {3.19}
$$
are satisfied.  Writing the equations on (3.19) as a linear system for unknown coefficients
$A_2,\cdots, A_m$ and $B_1,\cdots, B_{m-1}$ and carefully analyzing the resulted equations,  we obtain the unique solvability.

\proclaim{\noindent Theorem 3.7(\cite{ILayered}).}
Let $m\ge 3$.  There exists a positive number $R_m=R_m(\gamma_0)$ such that,
if $\text{Re}\,z\cdot\mbox{\boldmath $e$}_n>R_m$, then there exists the unique
$U^+=(B_1,A_2,B_2\cdots,A_{m-1},B_{m-1},A_m)$ depending on $z\cdot\mbox{\boldmath $e$}_n$ such that function $e^+$ given by (3.18)
belongs to $H^1_{\text{loc}}\,(\Bbb R^n)$ and satisfies $\nabla\cdot\gamma_0\nabla e^+=0$ in $\Bbb R^n$ in the weak sense.
Besides, as $\text{Re}\,z\cdot\mbox{\boldmath $e$}_n\rightarrow\infty$ we have, uniformly with $\text{Im}\,z\cdot\mbox{\boldmath $e$}_n$
$$\left\{
\begin{array}{l}
\displaystyle B_1\rightarrow R_{12}\\
\\
\displaystyle
A_2\rightarrow T_{12},\,\,B_2\rightarrow T_{12}R_{23}
\\
\\
\displaystyle
A_3\rightarrow T_{12}\,T_{23},\,\,B_3\rightarrow T_{12}\,T_{23}\,R_{34}
\\
\\
\displaystyle
\vdots
\\
\\
\displaystyle
A_{m-1}\rightarrow T_{12}\,T_{23}\cdots\,T_{m-2,m-1},\,\,B_{m-1}\rightarrow T_{12}\,T_{23}\cdots\,T_{m-2,m-1}\,R_{m-1,m}
\\
\\
\displaystyle
A_m\rightarrow T_{12}\,T_{23}\cdots\,T_{m-2,m-1}\,T_{m-1,m},
\end{array}
\right.
\tag {3.20}
$$
where the constants $T_{kl}$ and $R_{kl}$ are given by
$$\displaystyle
T_{kl}=\frac{2\gamma_k}{\gamma_k+\gamma_l},\,\,R_{kl}=\frac{\gamma_k-\gamma_l}{\gamma_k+\gamma_l}.
$$
\endproclaim

Using the asymptotic property on (3.20) and applying an argument similar to that of the proof of Theorem 3.2 with
the help of the system of inequalities on (2.18) and (2.19) of Proposition 2.5, 
one gets the analogous formula to (3.6) for the direction $\omega$ with $\omega\cdot\mbox{\boldmath $e$}_n>0$.

Note that we can also construct a solution named $e^-=e^-(x;\gamma_0,z)$ of $\nabla\cdot\gamma_0\nabla u=0$ in $\Bbb R^n$ corresponding
to the case $\text{Re}\,z\cdot\mbox{\boldmath $e$}_n<0$ by a reflection which has a similar asymptotic behaviour to listed on (3.20).
From $e^-$ instead of $e^+$ one gets also the value of the support function at the direction $\omega$ with $\omega\cdot\mbox{\boldmath $e$}_n<0$.
This yields Theorem 3.6.

Remarks are in order.

$\bullet$  The case $m=2$, we can explicitly give $e^{\pm}(x;\gamma_0,z)$.

$\bullet$  An extension of Theorems 3.6 and 3.8 to the case when $\gamma_j$ on (3.17) is replaced with a homogeneous anisotropic conductivity $C_j$
are given in \cite{ILayered2}.

$\bullet$  When $n=2$, the inclusion  $D$ is replaced with a thin perfectly insulated inclusion, that is a crack located on the line $x_2=c_j$ for some $j$,
Theorem 3.7 combined with the enclosure method with a single input has been applied in \cite{ICrack0}.
Its extension to the case $n=3$ remains open.

\subsubsection{Estimating cracks and spot welded parts}

As we did in studying the blow up set of the refracted solution in the probe method, 
one can deform also the complex plane wave $e^{x\cdot z},x\in\Bbb R^n$ with $z\cdot z=0$. 
In \cite{IIS} we made use of its Kelvin transform for an inverse crack problem in two dimensions.

Given $x\in\Bbb R^2$ define 
$$\begin{array}{lll}
\displaystyle
v_{\tau}(y;x)=\exp\,\left(\displaystyle -\tau\,\frac{y-x}{\vert y-x\vert^2}\cdot z\,\right), 
& y\in\Bbb R^2\setminus\{x\}, 
& \tau>0,
\end{array}
$$
where 
$$\displaystyle
z=\mbox{\boldmath $e$}_2+i\mbox{\boldmath $e$}_1.
$$  
The function $v_{\tau}(y;x)$ of $y\in\Bbb R^2\setminus\{x\}$
is nothing but the Kelvin transform of the function $e^{-\tau\,y\cdot\,z}$.

Since we have
$$\displaystyle
\frac{y-x}{\vert y-x\vert^2}\cdot z
=\frac{i}{(y_1-x_1)-i(y_2-x_2)},
$$
the function $v_{\tau}(y;x)$ satisfies the Laplace equation in $\Bbb R^2\setminus\{x\}$.

Let $s>0$.  One can write
$$\displaystyle
\frac{1}{2s}+\frac{y-x}{\vert y-x\vert^2}\cdot z
=\frac{\vert y-(x-s\mbox{\boldmath $e$}_2)\,\vert^2-s^2}{2s\vert y-x\vert^2}.
$$
Thus we have
$$\displaystyle
\lim_{\tau\rightarrow\infty} e^{-\frac{\tau}{2s}}\vert v_{\tau}(y;x)\vert=
\left\{
\begin{array}{ll}
\displaystyle
0 & \text{if $\vert y-(x-s\mbox{\boldmath $e$}_2)\vert>s$,}
\\
\\
\displaystyle
\infty &  \text{if $\vert y-(x-s\mbox{\boldmath $e$}_2)\vert<s$,}
\end{array}
\right.
$$
and if $\vert y-(x-s\mbox{\boldmath $e$}_2)\vert=s$, then $e^{-\frac{\tau}{2s}}v_{\tau}(y;x)$
is highly oscillating as $\tau\rightarrow\infty$.

In \cite{ICrack0} we considered an inverse crack problem in a two dimensional $m$-layered medium whose governing equation is the 
equation $\nabla\cdot\tilde{\gamma_0}\nabla u=0$ in a domain $\Omega=]a, b[\,\times\,]c, d[$ with $\tilde{\gamma_0}$ given by
(3.16), $a<b$ and $c<c_{m-1}<\cdots<c_1<d$.
It is assumed that the unknown cracks are completely inside the medium and lying on some layer joining lines.  Using the enclosure method
with a single input combined with the complex plane wave $e^{\pm}$ mentioned in the previous section, 
we gave an extraction formula of the values of the support function of the union of all unknown cracks at regular directions $\omega$
with $\omega\cdot\mbox{\boldmath $e$}_2\not=0$ with respect to the union.

However, the obtained information is not enough to distinguish, for example, two cracks located on a single layer joining line.
To get more information, in \cite{IIS} using function $v_{\tau}(y;x)$, we applied the enclosure with a single input
to  a two layered medium version of the inverse crack problem.
This is the case: $\Omega=]0, a[\times\,]0, b[$ with $0<a, b$ and 
$$\begin{array}{ll}
\displaystyle
\Sigma\subset [0, a]\times\{c\}, & 0<c<b;
\end{array}
$$
the observation data is the pair $u\vert_{\partial\Omega}$ with a {\it fixed} $g\not=0$ satisfying $\int_{\partial\Omega}g\,d\sigma=0$,
where $u=u(y)$ solves
$$
\left\{
\begin{array}{ll}
\displaystyle
\Delta u=0, & y\in\Omega\setminus\Sigma,
\\
\\
\displaystyle
\frac{\partial u}{\partial\nu}=0, & y\in\Sigma,
\\
\\
\displaystyle
\frac{\partial u}{\partial\nu}=g, & y\in\partial\Omega.
\end{array}
\right.
\tag {3.21}
$$
For explanation, here we consider only the case when $\Sigma=([0, c_1]\times\{c\})\cup\,([c_2, a]\times\{c\})$ with unknown 
numbers $c_1<c_2$ in the open interval $]0,a[$.  Note that the segment 
$$\displaystyle
W=]c_1, c_2[\times\{c\},
$$ 
is considered as a model of an unknown single {\it spot welded part} of two layers   $\Omega^{+}=]0,\,a[\times\,]0,\,c[$ and
$\Omega^{-}=]0,\,a[\times\,]c,\,b[$ whose conductivities are same, homogeneous and isotropic.

Define the indicator function
$$\begin{array}{ll}
\displaystyle
I(\tau;x)=e^{-\frac{\tau}{2s_0}}\int_{\partial\Omega}
\left(g\,v_{\tau}(y;x)-\frac{\partial v_{\tau}}{\partial\nu}(y;x)\,u\right)\,d\sigma(y), & \tau>0,
\end{array}
\tag {3.22}
$$
where $x$ runs on the segment $]0,\,a[\times\{b+\epsilon\}$ with a {\it fixed} $\epsilon>0$ and
the positive number $s_0$ is chosen
in such a way that the point $x-2s_0\mbox{\boldmath $e$}_2$ is located on the line $y_2=c$ for all $x\in\,]0,\,a[\times\{b+\epsilon\}$,
that is $s_0=\frac{b+\epsilon-c}{2}$.

The following result gives us a characterization of $c_1$ and $c_2$ in terms of a difference of the asymptotic behaviour
of indicator function $I(\tau;x)$ as $\tau\rightarrow\infty$ for each $x\in\,]0, a[\times\,\{b+\epsilon\}$.

\proclaim{\noindent Theorem 3.8(\cite{IIS}).}
Let $g$ satisfy $\text{supp}\,g\subset\,]0,\,a[\times\{b\}$.
We have

$\bullet$ if $x\in \{c_1,c_2\}\times\{b+\epsilon\}$, then there exists an integer $N\ge 1$ and a complex number $A\not=0$
such that $\lim_{\tau\rightarrow\infty}\,\tau^{\frac{2N-1}{2}}\,I(\tau;x)=A$.

$\bullet$  if $x\not\in\{c_1,c_2\}\times\{b+\epsilon\}$, then $I(\tau;x)$ is exponentially decaying as $\tau\rightarrow\infty$.

\endproclaim

{\it\noindent Sketch of proof.} The proof consists of three parts.

$\bullet$  The convergent series expansion of solution $u$ of (3.21) around the point $(c_j,c)$ and the complete 
asymptotic expansion of the indicator function at $x\in \{c_1,c_2\}\times\{b+\epsilon\}$.
The proof basically follows \cite{ICrack0}, however, function $v_{\tau}(y;x)$ is more complicated than the complex plane wave,
we need the help of the {\it method of the steepest descent}.

$\bullet$  An extension argument of \cite{AD} to show nonvanishing of a coefficient in the expansion of the indicator function above.
For this the condition $\text{sup}\,g\subset\partial\Omega^+\setminus([0,\,a]\times\{c\})$ is essential.

$\bullet$  A {\it harmonic continuation} of $u$ from $\Omega^+$ into $\Omega^-$ locally at a given point, especially,
on the set $]0,\,a[\times\{c\}\setminus\overline W$ (this is a nontrivial case) to show the exponential decay of the indicator function
at $x\not\in\{c_1,c_2\}\times\{b+\epsilon\}$.  This can be done by a reflection.

\noindent
$\Box$

Since $e^{-\frac{\tau}{2s_0}}v_{\tau}(y;x)$ is exponentially decaying for $y$ in $\vert y-s_0\mbox{\boldmath $e$}_2\vert>s_0+\delta$ with a small $\delta>0$,
one can reduce $\partial\Omega$ on (3.22) with the smaller pat $\partial\Omega\cap\{y\,\vert\,y_2>c-\delta\}$ and gets the same result.
Besides, under some assumption on the location of $W$ from {\it above}, the domain of integral on (3.22) can be reduced to a part on the up side of $\Omega$.

The $s_0$ on (3.22) is always fixed in Theorem 3.8.  There is another approach which uses $s_0$ as another variable.
Let $x\in\,]0,\,a[\times\{b+\epsilon\}$.
The approach enables us to find the largest circle centered at $x-s\mbox{\boldmath $e$}_2$ with radius $s$ that firstly touches $W$, that is,
the quantity defined by
$$\displaystyle
s_{\Sigma}(x)=\sup\left\{s>0\,\vert\,B_s(x-s\mbox{\boldmath $e$}_2)\subset\Bbb R^2\setminus\Sigma\right\},
$$
where $B_s(\xi)=\{y\in\Bbb R^2\,\vert\,\vert y-\xi\vert<s\}$.  Note that if $s$ is sufficiently small, 
then $B_s(x-s\mbox{\boldmath $e$}_2)\subset\Omega^+$.

In \cite{HIIS}, roughly speaking, in the present case,  
we established that, for all $x\in\,]0,\,a[\times\{b+\epsilon\}$ except for the points on the projection
of $W$ onto the set $]0,\,a[\,\times\{b+\epsilon\}$
the function
$$\displaystyle
e^{\frac{\tau}{2s_0}}e^{-\frac{\tau}{2s_{\Sigma}(x)}}\,I(\tau;x)
$$
is truly algebraic decaying as $\tau\rightarrow\infty$.  Thus we have the formula
$$\displaystyle
\lim_{\tau\rightarrow\infty}\,\frac{1}{\tau}\log\vert I(\tau;x)\vert=\frac{1}{2s_{\Sigma(x)}}-\frac{1}{2s_0}.
$$
Note that this right-hand side has a peak at the projection of two endpoints of $W$ onto set $]0,\,a[\,\times\{b+\epsilon\}$.
Paying attention to this property, we proposed an algorithm to find the points $(c_1,c)$ and $(c_2,c)$ between them the welded
part $W$ lies. We did its numerical experiment and the result suggests a possibility of monitoring a spot welding process.

\subsection{The Mittag-Leffler function and the enclosure method}

The original enclosure method yields only the convex hull of unknown discontinuity, however, it is a constructive method 
compared with the probe method.  It is natural to consider whether one can obtain more information
than the convex hull by modifying the enclosure method.
In \cite{IMittag} the author introduced a method of using the Mittag-Leffler function instead of the complex exponential function.
This is natural since the Mittag-Leffler function is an extension of the complex exponential function.
The method has been applied to Problem 2.2 in the case when the background conductivity $\gamma_0$ is given by a constant.
It is based on a notion of a generalized support function.  
A numerical implementation together with a introduction of an {\it alternative version} which is based on a method of searching
by a cone has been done in \cite{IS1}. 
 
In what follows, $\Omega$ denotes a bounded domain of $\Bbb R^2$ with smooth boundary
and $\gamma$ be the conductivity of $\Omega$ given by (3.16) with $\gamma_0(x)=1$.
We assume that $\gamma$ has the global jump condition: there exists a positive constant $C>0$
such that $h(x)\xi\cdot\xi\ge C\vert\xi\vert^2$ for all $\xi\in\Bbb R^2$ and a.e, $x\in D$
or $-h(x)\xi\cdot\xi\ge C\vert\xi\vert^2$ for all $\xi\in\Bbb R^2$ and a.e, $x\in D$;
$\partial D$ is Lipschitz.

Let $0<\alpha\le 1$.  Given $y\in\Bbb R^2$, $\omega\in S^1$ and $t\in\,]-\infty, 0]$ define the harmonic function $v_{\tau}^{\alpha}(x;y,\omega,t), x\in\Bbb R^2$ by the formula
$$\displaystyle
v_{\tau}^{\alpha}(x;y,\omega,t)
=E_{\alpha}(\tau\{(x-y-t\omega)\cdot\omega+i(x-y-t\omega)\cdot\omega^{\perp}\}),
$$
where $\omega^{\perp}\in S^1$ and satisfy $\omega\cdot\omega^{\perp}=0$ and $\text{det}\,(\omega^{\perp}\,\omega)>0$.
Note that $\omega^{\perp}$ is uniquely determined by $\omega$.

From the asymptotic behaviour of the Mittag-Leffler function, we have:
$$\displaystyle
\lim_{\tau\rightarrow\infty}\,\vert v_{\tau}^{\alpha}(x;y,\omega,t)\vert
=
\left\{
\begin{array}{ll}
\displaystyle
\infty & \text{if $x\in C_{y+t\omega}(\omega,\alpha)$,}
\\
\\
\displaystyle
0 & \text{if $x\not\in C_{y+t\omega}(\omega,\alpha)$,}
\end{array}
\right.
$$
where $C_{y+t\omega}(\omega,\alpha)$ denotes the two-dimensional cone
about $\omega$ of opening angle $\frac{\pi\,\alpha}{2}$ with vertex at $y+t\omega$.
Note that $C_{y+t\omega}(\omega,\alpha)$ is closed.

{\bf\noindent Definition 3.5.}
Define the indicator function in this section by the formula
$$\begin{array}{ll}
\displaystyle
I^{\alpha}_{(y,\omega)}(\tau,t)
=\int_{\partial\Omega}\,(\Lambda_{\gamma}-\Lambda_1)f^{\alpha}\cdot\overline{f^{\alpha}}\,d\sigma,
& \tau>0,
\end{array}
$$
where $f^{\alpha}(x)=v_{\tau}^{\alpha}(x;y,\omega,t), x\in\partial\Omega$.

From this one we extract information about the geometry of an unknown inclusion $D$.
For this we have two approaches as described above.

\subsubsection{Generalized support function}

In \cite{IMittag} the author introduced an extension of the support function.

{\bf\noindent Definition 3.6.}
Given $(y,\omega)\in\,(\Bbb R^2\setminus\overline\Omega)\times S^1$ with $C_y(\omega,\alpha)\cap\overline{\Omega}=\emptyset$,
define the generalized support function of $D$ by the formula
$$\displaystyle
h^{\alpha}_D(y,\omega)=
\inf\left\{t\in\,]-\infty,\,0[\,\vert\,
\forall s\in\,]t,\,0[\, C_{y+s\omega}(\omega,\alpha)\cap\overline D=\emptyset\right\}.
$$
Note that we have $D\subset \Bbb R^2\setminus C_{y+h^{\alpha}_D(y,\omega)\,\omega}\,(\omega,\alpha)$.
Thus knowing the generalized support function, one gets an upper bound of the geometry of $D$ in the sense above.

The following is the first approach which has been announced in \cite{IHokkaido}.

\proclaim{\noindent Theorem 3.9(\cite{IMittag}).}  Let $\alpha<1$.  Let $(y,\omega)\in(\Bbb R^2\setminus\overline\Omega)\times S^1$
satisfy $C_y(\omega,\alpha)\cap\overline{\Omega}=\emptyset$.
We have
$$\displaystyle
\lim_{\tau\rightarrow\infty}\,\vert I^{\alpha}_{(y,\omega)}(\tau,t)\vert =
\left\{
\begin{array}{ll}
0 & \text{if $t>h^{\alpha}_D(y,\omega)$,}
\\
\\
\displaystyle
\infty & \text{if $t<h^{\alpha}_D(y,\omega)$}
\end{array}
\right.
$$
and if $t=h^{\alpha}_D(y,\omega)$, then 
$$\displaystyle
\liminf_{\tau\rightarrow\infty}\,\vert I^{\alpha}_{(y,\omega)}(\tau,t)\vert>0.
$$
\endproclaim

{\it\noindent Sketch of Proof.}
Proposition 2.5 together with the global jump condition for $\gamma$ everything reduced to
check the corresponding asymptotic behaviour of the integral
$$\displaystyle
\int_D\vert\nabla v^{\alpha}_{\tau}(x;y,\omega,t)\,\vert^2\,dx.
$$
The most delicate part is the case when $t=h^{\alpha}_D(y,\omega)$, that is to show that
$$\displaystyle
\liminf_{\tau\rightarrow\infty}\,
\int_D\vert\nabla v^{\alpha}_{\tau}(x;y,\omega,t)\,\vert^2\,dx>0.
$$
Unlike before which is the case $\alpha=1$, we have to divide the case into two subcases:
(a) $y+t\omega\in\partial D$; (b) $y+t\omega\not\in\partial D$.
For both we use an asymptotic expansion of the Mittag-Leffler function stated as Lemma 2.1 in \cite{IMittag}.
The proof of Lemma 2.1 which follows the line described in \cite{Bateman} and is given as Appendix in \cite{IMittag}.

Note that the case (b) means that there is a point $x_0$ on $\partial C_{y+t\omega}(\omega,\alpha)\cap\partial D\setminus\{y+t\omega\}$.
Since $\partial D$ is Lipschitz, one can replace $D$ at $x_0$ locally from below the interior of an triangle with vertex $x_0$,
that is, one may assume that $D$ is the interior of a triangle in advance.  In the analysis of case (b)
we fully making use of this special shape.

\noindent
$\Box$

An extension of Theorem 3.9 to the inverse inclusion problem
governed by the equation $\nabla\cdot(\sigma-i\omega\epsilon)\nabla u=0$
has been established in \cite{IComplexConductivity}.

\subsubsection{Searching by a cone and regularization}

The indicator function of Theorem 3.9 has the large parameter $\tau>0$ and another parameter $-\infty<t<0$.
The following theorem is an alternative version of Theorem 3.9.  The idea is: instead of searching parameter $t$ 
use direction $\omega$ as another searching parameter.

\proclaim{\noindent Theorem 3.10(\cite{IS1}).}
Given $(y,\omega)\in\Omega\times S^1$ we have
$$\displaystyle
\lim_{\tau\rightarrow\infty}\,\vert I^{\alpha}_{(y,\omega)}(\tau,0)\vert =
\left\{
\begin{array}{ll}
0 & \text{if $C_y(\omega,\alpha)\cap\overline D=\emptyset$,}
\\
\\
\displaystyle
\infty & \text{if $\text{int}\,(C_y(\omega,\alpha))\cap D\not=\emptyset$}
\end{array}
\right.
$$
and if $C_y(\omega,\alpha)\cap\overline D\not=\emptyset$ and $\text{int}\,(C_y(\omega,\alpha))\cap D=\emptyset$, then 
$$\displaystyle
\liminf_{\tau\rightarrow\infty}\,\vert I^{\alpha}_{(y,\omega)}(\tau,0)\vert>0.
$$

\endproclaim

The essential part of the proof is exactly the same as in Theorem 3.9.

From Theorem 3.10 we have the formula:
$$\displaystyle
V(D)=\cup_{0<\alpha<1}\,\cup_{\omega\in S^1}\,
\left\{y\in\Omega\,\vert\,\lim_{\tau\rightarrow\infty}\,I^{\alpha}_{(y,\omega)}(\tau,0)=0\right\},
\tag {3.23}
$$
where $V(D)$ is the set of all points $y\in\Omega$ such that: there exists a direction $\omega$ 
with $\left\{y+s\omega\,\vert\,0\le s<\infty\right\}\cap \overline D=\emptyset$.  Clearly we have
$D\subset\Omega\setminus V(D)$.

In \cite{IS1} an algorithm based on (3.23) has been introduced.  Besides as a main result we gave a modification
of Theorem 3.10 to cover the case of {\it finitely many noisy} measurements.  
The idea is making use of a special truncation of power series expansion of the Mittag-Leffler function
with $\alpha=\frac{1}{n}$, $n\ge 2$:
$$\displaystyle
E_{\frac{1}{n}}(z)\sim \sum_{m=0}^{nN}\,\frac{z^m}{\Gamma(\frac{m}{n}+1)},
$$
where $N\ge 1$.  By using the closed form of $E_{\frac{1}{n}}(z)$ in \cite{Bateman} one can give a precise error estimate
for the error term caused by this truncation.

This form suggests another indicator function
$$\begin{array}{ll}
\displaystyle
I^{\frac{1}{n},nN}{(y,\omega)}(\tau;\{\delta_j^i\})
&
\displaystyle
=\sum_{0\le m,\,l\le nN}
\,\frac{\tau^{m+l}}{\Gamma(\frac{m}{n}+1)\,\Gamma(\frac{l}{n}+1)}
\\
\\
\displaystyle
&
\displaystyle
\,\,\,
\times
\int_{\partial\Omega}\left\{(\Lambda_{\gamma}-\Lambda_1)f_m(x;y,\omega)+\delta_m^1(x)\right\}
\cdot
\overline{\left(f_l(x;y,\omega)+\delta_l^2(x)\right)}\,d\sigma(x).
\end{array}
$$
where for $j=0,1,\cdots, nN$,
$$\begin{array}{ll}
\displaystyle
f_j(x;y,\omega)=\{(x-y)\cdot(\omega+i\omega^{\perp})\}^j, & x\in\partial\Omega.
\end{array}
$$
Note that the functions $\delta^i_j\in L^2(\partial\Omega)$, $j=0,1,\cdots, nN$ and $i=1,2$
are of noise and error caused by the measurement process.  

Define
$$\displaystyle
\Vert \delta^i_j\Vert\equiv\max_{i=1,2}\max_{0\le j\le nN}\,\Vert \delta^i_j\Vert_{L^2(\partial\Omega)}.
$$

Let $\beta_0$ denote the unique positive solution of the equation
$$\displaystyle
\frac{2}{e}\,\beta+\log\beta=0.
$$
The number $\beta_0$ has the bound $\frac{1}{\sqrt{e}}<\beta_0<1$.  

Given $y\in\Omega$ define
$$\displaystyle
C(y)=\frac{\beta}{\theta e}\,\left(1+\frac{1}{c(y)^n}\right),
\tag {3.24}
$$
with arbitrary {\it fixed} $\beta\in\,]0,\,\beta_0[$ and $\theta\in\,]0,\,1[$ and $c(y)$ is an arbitrary function of $y\in\Omega$
such that 
$$\displaystyle
c(y)\ge\sup_{x\in\Omega}\,\vert x-y\vert.
$$
For example, one can choose $\beta=\theta\beta_0$ and $c(y)=\sup_{x\in\Omega}\,\vert x-y\vert$.
Then (3.24) becomes
$$\displaystyle
C(y)=\frac{\beta_0}{e}\,\left(1+\frac{1}{(\sup_{x\in\Omega}\,\vert x-y\vert)^n}\right).
$$
Note that $c(y)$ becomes large when $y\in\Omega$ is near $\partial\Omega$ and thus $C(y)$ becomes small.

Define
$$\begin{array}{ll}
\displaystyle
N(\delta;y)=\left[-\frac{\log\delta}{C(y)}\right], & 0<\delta\le e^{-C(y)}.
\end{array}
$$
and
$$\displaystyle
\tau_y(N)=\left(\frac{N+\frac{1}{n}}{e}\,\right)^{\frac{1}{n}}\,\frac{\beta^{1/n}}{c(y)}.
$$
Here the symbol $[*]$ denotes the largest integer that does not exceed $*$.

\proclaim{\noindent Theorem 3.11(\cite{IS1}).}  Given $y\in\Omega$ and $\delta\in\,]0,\,e^{-C(y)}]$
let $N=N(\delta;y)$.  We have
$$
\left\{\displaystyle
\begin{array}{ll}
\displaystyle
\lim_{\delta\rightarrow 0}\,\left(\sup_{\Vert\delta_j^i\Vert\le\delta}
\,I^{\frac{1}{n},nN}{(y,\omega)}(\tau_y(N);\{\delta_j^i\})\right)
=0 & \text{if $C_y(\omega,\alpha)\cap\overline D=\emptyset$,}
\\
\\
\displaystyle
\lim_{\delta\rightarrow 0}\,\left(\inf_{\Vert\delta_j^i\Vert\le\delta}
\,I^{\frac{1}{n},nN}{(y,\omega)}(\tau_y(N);\{\delta_j^i\})\right)
=\infty & \text{if $\text{int}\,(C_y(\omega,\alpha))\cap D\not=\emptyset$}
\end{array}
\right.
$$
and if $C_y(\omega,\alpha)\cap\overline D\not=\emptyset$ and $\text{int}\,(C_y(\omega,\alpha))\cap D=\emptyset$, then 
$$\displaystyle
\liminf_{\delta\rightarrow 0}\,\left(\inf_{\Vert\delta_j^i\Vert\le\delta}
\,I^{\frac{1}{n},nN}{(y,\omega)}(\tau_y(N);\{\delta_j^i\})\right)>0.
$$

\endproclaim

In Theorem 3.11 it is assumed that $\delta$ which is an upper bound of $\Vert\delta_i^j\Vert$, is {\it known}.
Roughly speaking, the result states: if $\delta$ is sufficiently small and choose $N=N(\delta;y)\sim\vert\log\,\delta\vert$ and 
$\tau=\tau_y(N)\sim N^{\frac{1}{n}}$, then the ``indicator function'' $I^{\frac{1}{n},nN}{(y,\omega)}(\tau_y(N);\{\delta_j^i\})$ behaves
like the original indicator function $I^{\alpha}_{(y,\omega)}(\tau,0)$ even the case when the data contain error and noise whose size 
is dominated by $\delta$.  This is a {\it regularized version} of the enclosure method whose idea goes back to the corresponding formula
for the enclosure method with a single input in \cite{IRegularized} and also see \cite{ITwomethods}.

The idea of combining the enclosure method with the Mittag-leffler function has been applied also to Problem 2.1 in two dimensions
which is an inverse crack problem.  See \cite{IO2}.

\subsection{Inverse obstacle scattering and the enclosure method}

The ideas of the enclosure method \cite{I1, E00} have or still now may have many applications to inverse obstacle problems
governed by various partial differential equations which describe the wave phenomena in the frequency domain.  
In this section we present three of them.

\subsubsection{Logarithmic differential of the indicator function}

In this section we describe an application of the enclosure method \cite{I1}
to the prototype inverse obstacle scattering problem.

Let $d\in S^1$ and $k>0$.
Let $B_R$ be an open disc with radius $R$ centered at a fixed point.
Assume that an unknown obstacle $D$ is given by an open set with Lipschitz boundary of $\Bbb R^2$ with $\overline D\subset B_R$;
$B_R\setminus\overline D$ is connected.

Let $u=u(x)$ be the solution of scattering problem governed by the Helmholtz equation
$$\left\{
\begin{array}{ll}
\displaystyle
\Delta u+k^2u=0, & x\in\Bbb R^2\setminus\overline D,
\\
\\
\displaystyle
\frac{\partial u}{\partial\nu}=0, & x\in\partial D,
\\
\\
\displaystyle
\sqrt{r}\left(\frac{\partial w}{\partial r}-ikw\right)\longrightarrow 0, & r=\vert x\vert\longrightarrow\infty,
\end{array}
\right.
$$
where $w=u-e^{ikx\cdot d}$ is the scattered wave and the last condition is called the Sommerfeld outgoing radiation condition.

Note that $\nu$ stands for the unit outward normal vector field to
$\partial(B_R\setminus\overline D)$.
The boundary condition for $\partial u/\partial\nu$ on $\partial D$
means that $D$ is a {\it sound-hard} obstacle and
should be considered as a weak sense.

$\quad$

{\bf\noindent Problem 3.2.}
Fix $k$ and $d$.  Find a formula that extracts information about the location
of $D$ from the Cauchy data of $w$ on $\partial B_R$.

$\quad$

In this section we always assume that $D$ is polygonal, that is,
$D$ takes the form $D_1\cup\cdots\cup D_m$ with
$1\le m<\infty$ where each $D_j$ is a connected component of $D$ and a polygon; $\overline
D_j\cap\overline D_{j'}=\emptyset$ if $j\not=j'$.

Let $\omega$ and $\omega^{\perp}$ be two unit
vectors perpendicular to each other.  We always assume that the orientation of $\omega^{\perp}$ and $\omega$ coincides with
$\mbox{\boldmath $e$}_1$ and $\mbox{\boldmath $e$}_2$ and thus $\omega^{\perp}$ is unique.

Define the indicator function
$$\begin{array}{ll}
\displaystyle
I(\tau;\omega, d, k)
=\int_{\partial B_R}\left(\frac{\partial u}{\partial\nu}v_{\tau}
-\frac{\partial v_{\tau}}{\partial\nu}u\right)
d\sigma, & \tau>0,
\end{array}
\tag {3.25}
$$
where
$$\displaystyle
v_{\tau}(x)=e^{x\cdot(\tau\,\omega+i\sqrt{\tau^2+k^2}\,\omega^{\perp})}.
$$
The function $v_{\tau}$ satisfies the Helmholtz equation $(\Delta+k^2)v=0$
in the whole plane.

We say that $\omega$ is regular with respect to $D$ if the set
$\displaystyle
\partial D\cap\{x\in\Bbb R^2\,\vert\,x\cdot\omega=h_D(\omega)\}$
consists of only one point. 
In other words, when $t$ moves from 
$t=\infty$ to $-\infty$, the line $x\cdot\omega=t$ which is
perpendicular to direction $\omega$ descends from infinity and firstly touches
a single point on $\partial D$ at $t=h_D(\omega)$.
Since $D$ is assumed to be polygonal,
the point should be a vertex of the convex hull $[D]$.

In \cite{IScattering2}, the author has established the formula
\proclaim{\noindent Theorem 3.12(\cite{IScattering2}).}
Let $\omega$ is regular with respect to $D$.
We have
$$\displaystyle
\lim_{\tau\longrightarrow\infty}\frac{1}{\tau}
\log
\left\vert I(\tau;\omega, d, k)
\right\vert
=h_D(\omega),
\tag {3.26}
$$
Moreover, we have the following:
$$
\displaystyle
\lim_{\tau\longrightarrow\infty}e^{-\tau t}\left\vert
I(\tau; \omega, d, k)
\right\vert
=
\left\{
\begin{array}{ll}
0
&
\text{if $t\ge h_D(\omega)$,}\\
\\
\displaystyle
\infty
&
\text{if $t<h_D(\omega)$.}
\end{array}
\right.
$$

\endproclaim

In this formula one makes use of only the absolute
value of indicator function $I(\tau;\omega,d,k)$.
Thus, one needs two regular directions
$\omega$ for determining a single vertex of the convex hull of $D$
since formula (3.26) gives only a single line on which the vertex
lies.

To overcome this redundancy, in \cite{IEnclosure5}
the author has introduced
a method for the use of the complex values
of the indicator function which directly yields the coordinates of a
vertex of the convex hull of $D$ with indicator functions for a
single regular direction $\omega$.

The method makes use of not only original indicator function (3.25) but also its derivative
with respect to $\tau$:
$$\displaystyle
I'(\tau;\omega,d,k)
=\int_{\partial B_R}
\left(\frac{\partial u}{\partial\nu}\cdot\partial_{\tau}v_{\tau}-
\frac{\partial}{\partial\nu}(\partial_{\tau}v_{\tau})u\right)d\sigma.
\tag {3.27}
$$

Th main result is the following theorem.

\proclaim{\noindent Theorem 3.13(\cite{IEnclosure5}).}
Let $\omega$ be regular with respect to $D$.
Let $x_0\in\partial D$ be the point with $x_0\cdot\omega=h_D(\omega)$.
There exists a $\tau_0>0$ such that, for all $\tau\ge\tau_0$
$\vert I(\tau;\omega,d,k)\vert>0$
and the formula
$$\displaystyle
\lim_{\tau\longrightarrow\infty}
\frac{\displaystyle I'(\tau;\omega,d,k)}
{I(\tau;\omega,d,k)}
=h_D(\omega)+ix_0\cdot\omega^{\perp}
\tag {3.28}
$$
is valid.

\endproclaim

Note that it is also shown that the
convergence rate of (3.28) is better than that of (3.26).

Here we mention here two implications of formula (3.28).

(1).  One can estimate the set of all directions that are not regular with respect to $D$. 
Let us denote the set by $I(D)$.  Since $I(D)$ is finite, one can write
$$\displaystyle
I(D)=\{(\cos\,\theta_j,\sin\,\theta_j)\,\vert\,0\le \theta_1<\cdots<\theta_N<2\pi\}.
$$
Let $\omega\in S^1\setminus I(D)$.
We denote by $x(\omega)=(x(\omega)_1,x(\omega)_2)$ the single point
in $\displaystyle
\partial D\cap\{x\in\Bbb R^2\,\vert\,x\cdot\omega=h_D(\omega)\}$.
Formula (3.28) is equivalent to the formulae
$$\left\{
\begin{array}{l}
\displaystyle
x(\omega)_1=\text{Re}\,\left\{(\omega_1+i\omega_2)\lim_{\tau\longrightarrow\infty}
\overline{\left(\frac{\displaystyle I'(\tau;\omega,d,k)}
{I(\tau;\omega,d,k)}\right)}\right\},
\\
\\
\displaystyle
x(\omega)_2=\text{Im}\,\left\{(\omega_1+i\omega_2)\lim_{\tau\longrightarrow\infty}
\overline{\left(\frac{\displaystyle I'(\tau;\omega,d,k)}
{I(\tau;\omega,d,k)}\right)}\right\}.
\end{array}
\right.
\tag {3.29}
$$
Since $D$ is polygonal, we see that $x(\omega)$ is constant on each open arc $\{(\cos\theta,\sin\theta)\,\vert\,
\theta_j<\theta<\theta_{j+1}\}$, $j=1,\cdots, N$ with $\theta_{N+1}=\theta_1$
and one of $x(\omega)_1$ or $x(\omega)_2$ has a first kind of discontinuity at each $\omega\in I(D)$.
Therefore one can expect
that computing both $x(\omega)_1$ and $x(\omega)_2$ for
sufficiently many $\omega$ via formulae (3.29), one can
estimate $I(D)$.

(2).   As a corollary of (3.28) we have 
$$\displaystyle
\lim_{\tau\longrightarrow\infty}
\left\vert\frac{I'(\tau;y,\omega,d,k)}
{I(\tau;y,\omega,d,k)}\right\vert
=\vert x(\omega)-y\vert,
$$
where $\omega\in S^1\setminus I(D)$, $y$ be an arbitrary point of $\Bbb R^2$ and
$$\displaystyle
I(\tau;y,\omega,d,k)
=e^{-y\cdot(\tau\omega+i\sqrt{\tau^2+k^2}\omega^{\perp})}I(\tau;\omega,d,k).
$$
Since $x(\omega)$ is the unique zero point of the function $y\mapsto \vert x(\omega)-y\vert$, one possible alternative idea to find
$x(\omega)$ is to consider the minimization problem of the
following function for a suitable $\tau$:
$$\displaystyle
y\mapsto \left\vert\frac{I'(\tau;y,\omega,d,k)}
{I(\tau;y,\omega,d,k)}\right\vert.
$$

In \cite{IEnclosure5}, Theorem 3.13 has been extended also to a thin obstacle case and
the case when the observation data is given by the far-field pattern of 
the scattered wave.

Let us describe this latter case.
It is well known that scattered wave $w$ has the asymptotic expansion
as $r\longrightarrow\infty$ uniformly with respect to $\varphi\in\,S^1$:
$$\displaystyle
w(r\varphi)=\frac{e^{ikr}}{\sqrt{r}}F(\varphi;d,k)+O\left(\frac{1}{r^{3/2}}\right).
$$
The coefficient $F(\varphi;d,k)$ is called the far-field pattern of
the scattered wave $w$ at direction $\varphi$, see \cite{CK}.

We identify $\varphi=(\varphi_1,\varphi_2)$ with the complex number
given by $\varphi_1+i\varphi_2$ and denote it by the same symbol $\varphi$.
Given $N=1,\cdots$, $\tau>0$, $\omega\in S^1$ and $k>0$ define the function $g_N(\,\cdot\,;\tau,k,\omega)$ on $S^1$
by the formula
$$
\displaystyle
g_{N}(\varphi;\tau,k,\omega)
=\frac{1}{2\pi}
\sum_{\vert m\vert\le N}
\left\{\frac{ik\varphi}{(\tau+\sqrt{\tau^2+k^2})\omega}\right\}^m.
\tag {3.30}
$$
Let the center point of $B_R$ be the origin of the coordinates and assume that
$\overline D\subset B_R$.

\proclaim{\noindent Theorem 3.14(\cite{IEnclosure5}).}  Let $\omega$ be regular with
respect to $D$.  Let $\beta_0$ be the unique positive solution of the equation
$$\displaystyle
\frac{2}{e}s+\log s=0.
$$
Let $\beta$ satisfy $0<\beta<\beta_0$.  Let $\{\tau(N)\}_{N=1,\cdots}$ be an arbitrary
sequence of positive numbers satisfying, as $N\longrightarrow\infty$
$$\displaystyle
\tau(N)=\frac{\beta N}{eR}+O(1).
$$
Then the formula
$$\displaystyle
\lim_{N\longrightarrow\infty}
\frac
{\displaystyle\int_{S^1}F(-\varphi;d,k)\partial_{\tau}g_N(\varphi;\tau(N),k,\omega)d\sigma(\varphi)}
{\displaystyle\int_{S^1}F(-\varphi;d,k)g_N(\varphi;\tau(N),k,\omega)d\sigma(\varphi)}
=h_D(\omega)+ix_0\cdot\omega^{\perp}
\tag {3.31}
$$
is valid.
\endproclaim

The origin of formula (3.31) is the formula established in \cite{IHerglotz}:
$$\displaystyle
\lim_{N\longrightarrow\infty}
\frac{1}{\tau(N)}
\log\left\vert\int_{S^1}F(-\varphi;d,k)g_N(\varphi;\tau(N),k,\omega)d\sigma(\varphi)\right\vert
=h_D(\omega),
\tag {3.32}
$$
where the coice of $\tau(N)$ is same as Theorem 3.14 and $\omega$ is regular with respect to $D$.

The $g=g_N$ given by (3.30) with $\tau=\tau(N)$ is chosen in such a way that, as $N\rightarrow\infty$ 
$$\begin{array}{ll}
\displaystyle
\int_{S^1} e^{iky\cdot\varphi}\,g_N(\varphi;\tau,k,\omega)\,d\sigma(\varphi)
\approx e^{y\cdot(\tau\,\omega+i\sqrt{\tau^2+k^2}\,\omega^{\perp})},
& y\in\overline{B_R}.
\end{array}
$$
Note that the function on this left-hand side is the Herglotz wave function with the density $g$ \cite{CK} which is an
entire solution of the Helmholtz equation.

The idea of taking the logarithmic differential of the indicator function can be applied
also to the previously published works \cite{ISource0, I1, IConductivity, ICrack0, Ipartial, IEnclosure4, IScattering3}.

Note also that a numerical algorithm based on formula (3.32) together with its limited aperture version \cite{IHerglotz}
has been proposed and tested in \cite{INS}.

\subsubsection{Inverse obstacle scattering with a single point source}

The enclosure method \cite{I1} covers also another inverse obstacle scattering problem in which
the scattered wave is generated by a single point source located within a finite distance from an unknown obstacle 
and observed on a known circle surrounding the obstacle.  
This type of problem can be seen in, e.g., a mathematical formulation of microwave tomography \cite{Se}, subsurface radar \cite{DGS}, etc..

In this section we assume that $D$ is polygonal.
Let $y\in\Bbb R^2\setminus\overline D$.   Let $E=E_D(x,y)$ be the unique solution of the scattering problem:
$$\left\{
\begin{array}{ll}
\displaystyle
\Delta E+k^2 E=0, & x\in\Bbb R^2\setminus\overline D,
\\
\\
\displaystyle
\frac{\partial}{\partial\nu}E=-\frac{\partial}{\partial\nu}\Phi_0(\,\cdot\,,y),
& x\in\partial D,
\\
\\
\displaystyle
\sqrt{r}\left(\frac{\partial E}{\partial r}-ikE\right)\rightarrow 0, & r=\vert x\vert\rightarrow\infty,
\end{array}
\right.
$$
where
$$\displaystyle
\Phi_0(x,y)=\frac{i}{4}H^{(1)}_0(k\vert x-y\vert)
$$
and $H^{(1)}_0$ denotes the Hankel function of the first kind \cite{Ol}.

The total wave outside $D$ exerted by the point source located at $y$ is given by the formula:
$$\displaystyle
\Phi_D(x,y)=\Phi_0(x,y)+E_D(x,y),\,x\in\Bbb R^2\setminus\overline D.
$$

$\quad$

{\bf\noindent Problem 3.3.}  Let $R_1>R$ and denote by $B_{R_1}$ and $B_R$ the open balls centered at a common arbitrary fixed point 
with radii $R_1$ and $R$, respectively.  Assume that $D$ satisfies $\overline D\subset B_R$.
Fix $k>0$ and $y\in\partial B_{R_1}$.
Extract information about the location and shape
of $D$ from $\Phi_D(x,y)$ given at all $x\in\partial B_{R}$.

$\quad$

Define
$$\displaystyle
J(\tau; \omega, y, k)
=
\int_{\partial B_{R}}
\left(
\frac{\partial}{\partial\nu}\Phi_D(x,y)\cdot v_{\tau}(x;\omega)
-\frac{\partial}{\partial\nu}v_{\tau}(x;\omega)\cdot\Phi_D(x,y)\right)d\sigma(x).
$$

It should be pointed out that $(\partial/\partial\nu)\Phi_D(x,y)$ for $x\in\partial B_{R}$ can be computed
from $\Phi_D(x,y)$ for $x\in\partial B_{R}$ by the formula
$$\begin{array}{ll}
\displaystyle
\frac{\partial}{\partial\nu}\Phi_D(x,y)=\frac{\partial}{\partial\nu}\Phi_0(x,y)+
\frac{\partial}{\partial\nu}\tilde{E}(x),
&
x\in\partial B_{R},
\end{array}
$$
where $\tilde{E}(x)$ solves the exterior Dirichlet problem for the Helmholtz equation:
$$\left\{
\begin{array}{ll}
\displaystyle
\Delta\tilde{E}+k^2\tilde{E}=0, & x\in\Bbb R^2\setminus\overline{B_{R}},\\
\\
\displaystyle
\tilde{E}=\Phi_D(\,\cdot\,,y)-\Phi_0(\,\cdot\,,y), x\in\partial B_{R},\\
\\
\displaystyle
\sqrt{r}\left(\frac{\partial \tilde{E}}{\partial r}-ik\tilde{E}\right)\rightarrow 0,
& r=\vert x\vert\rightarrow\infty.
\end{array}
\right.
$$

The first result of \cite{IEnclosure4} is as follows.

\proclaim{\noindent Theorem 3.15(\cite{IEnclosure4}).}
Assume that $\omega$ is regular with respect to $D$ and that
$$\displaystyle
\text{diam}\,D<\text{dist}\,(D,\partial B_{R_1}).
\tag {3.33}
$$
It holds that
$$\displaystyle
\lim_{\tau\longrightarrow\infty}
\frac{1}{\tau}\log
\left\vert J(\tau; \omega, y, k)
\right\vert
=h_D(\omega).
\tag {3.34}
$$
Moreover, we have the following:
$$
\displaystyle
\lim_{\tau\longrightarrow\infty}e^{-\tau t}\left\vert
J(\tau; \omega, y, k)
\right\vert
=
\left\{
\begin{array}{ll}
0
&
\text{if $t\ge h_D(\omega)$,}\\
\\
\displaystyle
\infty
&
\text{if $t<h_D(\omega)$.}
\end{array}
\right.
$$

\endproclaim

Some remarks on condition (3.33) are in order.

$\bullet$  Since $D$ satisfies $\text{diam}\,D<2R$ and 
$$
\displaystyle
\text{dist}\,(D,\partial B_{R_1})>\text{dist}\,(B_R,\partial B_{R_1})=R_1-R,
$$
condition (3.33) is satisfied if $R_1\ge 3R$.

$\bullet$ In Theorem 3.12 there is no restriction like (3.33).  This is due to a difference character of incident fields.

$\bullet$  It is an open problem whether one can drop (3.33) completely or not.

In \cite{IEnclosure4} applications of the enclosure method to thin obstacles, obstacles in a two layered medium
and an unexpected application to the linear sampling method \cite{CKi}, that is, the unsolvability of the far-field equation
for polygonal obstacles are given.

Note that, in \cite{IEnclosure5} the logarithmic differential version of formula (3.34) has been established.

\subsubsection{The enclosure method for the stationary Schr\"odinger equation}

In this section we present an unpublished application of the enclosure method \cite{E00}
to an inverse obstacle scattering problem governed by the the stationary Schr\"odinger equation\footnote{The contents of this section
are taken from the unpublished paper \cite{Idis}.}.

Let $\Omega$ be a bounded domain of $\Bbb R^3$ with smooth boundary.
We consider an inverse problem for the stationary Schr\"odinger equation
$$\begin{array}{ll}
\displaystyle
\Delta u+k^2u+Vu=0,
&
x\in\Omega.
\end{array}
\tag {3.35}
$$
Here $k$ is a positive constant and $V$ is a real valued $L^{\infty}(\Omega)$ function.
We assume that $0$ is not an eigenvalue for the Schr\"odinger operator $\Delta+k^2+V$ 
in $\Omega$ with the Dirichlet condition.  Then, given $f\in H^{1/2}(\partial\Omega)$ there exists a unique weak 
solution $u\in H^1(\Omega)$ of equation (3.35) with $u=f$ on $\partial\Omega$.

The map
$$\displaystyle
\Lambda_V:f\longmapsto \frac{\partial u}{\partial\nu}\vert_{\partial\Omega}
$$
is called the Dirichlet-to-Neumann map.  Here $\nu$ denotes the unit outward normal vector field to $\partial\Omega$.

We consider the following.

$\quad$

{\bf\noindent Problem 3.4.}
Extract information about the discontinuity of $V$ from $\Lambda_V$.

$\quad$

We have already known that there is a reconstruction formula of $V$ itself \cite{N1, N3, No} in three dimensions.
Key steps of the formula consist of two steps.

$\bullet$  Find a unique solution $\Psi(\,\cdot\,;\zeta,\Lambda_V)$ defined on $\partial\Omega$ of an integral equation depending on
a complex vector $\zeta$ with $\zeta\cdot\zeta=0$ and $\vert\zeta\vert>>1$ and $\Lambda_V$.

$\bullet$ Calculate $\Lambda_V\Psi(\,\cdot\,;\zeta,\Lambda_V)\vert_{\partial\Omega}$.

The point is: one has to know full knowledge of $\Lambda_V$ for obtaining the trace $\Psi(\,\cdot\,;\Lambda_V)\vert_{\partial\Omega}$ which is 
used in the second step.

Here we consider the case when the support of $V$ is contained in the closure of $D$
which is an open subset of $\Omega$ with $C^2$ boundary and $V$ satisfies a kind of a jump condition across $\partial D$ specified later.  
In this case we show that the information about the convex hull of $D$ can be obtained by knowing $\Lambda_V(v\vert_{\partial\Omega})$ 
for explicitly given solutions of the Helmholtz equation $\Delta v+k^2v=0$.  
Those solutions are different from solutions $\Psi(\,\cdot\,;\zeta,\Lambda_V)$
and much simpler.

Given $\tau>0$, $\omega\in S^2$ and $\omega^{\perp}\in S^2$ let
$v=e^{x\cdot(\tau\,\omega+i\sqrt{\tau^2+k^2}\,\omega^{\perp}\,)}$.
Define
$$\displaystyle
I_{\omega,\omega^{\perp}}\,(\tau)
=\int_{\partial\Omega}\,
\left\{
\Lambda_V(v\vert_{\partial\Omega})-\frac{\partial v}{\partial\nu}\right\}\,\overline{v}\,dS.
$$
We say that $V$ has a jump on $\partial D$ from the direction $\omega$ if
there exist $C=C(\omega)>0$ and $\delta=\delta(\omega)>0$ such that $V(x)\ge C$ for almost all $x\in D_{\omega}(\delta)$
or $-V(x)\ge C$ for almost all $x\in D_{\omega}(\delta)$.

\proclaim{\noindent Theorem 3.16.}
Let $\omega$ and $D$ satisfy $x\cdot\omega<h_D(\omega)$ for almost all $x\in\partial D$ with respect to the standard surface measure
induced on $\partial D$ 
from the Euclidean metric of $\Bbb R^3$.
If $V$ has a jump on $\partial D$ from direction $\omega$, then we have
$$\displaystyle
\lim_{\tau\rightarrow\infty}\,\frac{1}{2\tau}\,\log\vert I_{\omega,\omega^{\perp}}(\tau)\vert=h_D(\omega).
$$
Besides, we have
$$\displaystyle
\lim_{\tau\rightarrow\infty}\,e^{-2\tau t}\tau^2\vert I_{\omega,\omega^{\perp}}(\tau)\vert
=
\left\{
\begin{array}{ll}
\displaystyle
0 & \text{if $t>h_D(\omega)$,}
\\
\\
\displaystyle
\infty & \text{if $t<h_D(\omega)$}
\end{array}
\right.
$$
and
$$\begin{array}{ll}
\displaystyle
\liminf_{\tau\rightarrow\infty}e^{-2\tau t}\tau^2\vert I_{\omega,\omega^{\perp}}(\tau)\vert>0 & \text{if $t=h_D(\omega)$.}
\end{array}
$$

\endproclaim

Before describing the proof of Theorem 3.16 we prepare two lemmas. 

\proclaim{\noindent Lemma 3.1.}
Let $v\in H^2(\Omega)$ be a solution of the Helmholtz equation
$$\begin{array}{ll}
\displaystyle
\Delta v+k^2v=0, & x\in\Omega
\end{array}
$$
and let $u\in H^1(\Omega)$ solve
$$\left\{
\begin{array}{ll}
\displaystyle
\Delta u+k^2u+Vu=0, & x\in\Omega,
\\
\\
\displaystyle
u=v, & x\in\partial\Omega.
\end{array}
\right.
$$
Then, there exists $C=C(V)>0$ such that 
$$\displaystyle
\Vert u-v\Vert_{L^2(\Omega)}\le C\Vert Vv\Vert_{L^1(D)}.
$$

\endproclaim

{\it\noindent Proof.}
Set $w=u-v$.
One can find the unique $p\in H^1(\Omega)$ such that
$$\left\{
\begin{array}{ll}
\displaystyle
\Delta p+k^2\,p+V\,p=\overline {w}, & x\in\Omega,
\\
\\
\displaystyle
p=0, & x\in\partial\Omega.
\end{array}
\right.
\tag {3.36}
$$
By the elliptic regularity, we have $p\in H^2(\Omega)$ and $\Vert p\Vert_{H^2(\Omega)}\le C_1(V)\,\Vert w\Vert_{L^2(\Omega)}$.
By the Sobolev imbedding theorem, we have $\Vert p\Vert_{L^{\infty}(\Omega)}\le C_2(V)\Vert p\Vert_{H^2(\Omega)}$.
Thus one gets 
$$\displaystyle
\Vert p\Vert_{L^{\infty}(\Omega)}\le C_1(V)C_2(V)\Vert w\Vert_{L^2(\Omega)}.
\tag {3.37}
$$

The function $w$ satisfies
$$\left\{
\begin{array}{ll}
\displaystyle
\Delta w+k^2\,w+V\,w=-V\,v, & x\in\Omega,
\\
\\
\displaystyle
w=0, & x\in\partial\Omega.
\end{array}
\right.
$$
Using this and (3.36) together with integration by parts, we have the expression
$$\begin{array}{ll}
\displaystyle
\Vert w\Vert_{L^2(\Omega)}^2
&
\displaystyle
=-\int_{\Omega}\nabla p\cdot\nabla w\,dx+\int_{\Omega}\,p(k^2+V)w\,dx
\\
\\
\displaystyle
&
\displaystyle
=-\int_{D}\,p V\,v\,dx
\end{array}
$$
and thus one
gets
$$\displaystyle
\Vert w\Vert_{L^2(\Omega)}^2
\le
\Vert p\Vert_{L^{\infty}(\Omega)}\,\Vert V\,v\Vert_{L^1(D)}.
$$
Now from (3.37) we obtain the desired estimate.

\noindent
$\Box$

\proclaim{\noindent Lemma 3.2.}
Let $v=e^{x\cdot(\tau\omega+i\sqrt{\tau^2+k^2}\,\omega^{\perp})\,}$. 
Let $\omega$ and $D$ satisfy $x\cdot\omega<h_D(\omega)$ for almost all $x\in\partial D$ with respect to the standard surface measure on $\partial D$.
We have
$$\displaystyle
\lim_{\tau\rightarrow\infty}\,\frac{\Vert v\Vert_{L^1(D)}}{\Vert v\Vert_{L^2(D)}}=0.
\tag {3.38}
$$
\endproclaim

{\it\noindent Proof.}
In \cite{E00}, using the interior sphere condition for $\partial D$ we have already proved:
there exist $C>0$ and $\tau_0>0$ such that, for all $\tau\ge\tau_0$
$$\displaystyle
e^{-2\tau h_D(\omega)}\,\Vert v\Vert_{L^2(D)}^2\ge C\tau^{-2}.
\tag {3.39}
$$
The divergence formula yields
$$\begin{array}{ll}
\displaystyle
e^{-\tau h_D(\omega)}\Vert v\Vert_{L^1(D)}
&
\displaystyle
=\frac{1}{\tau}\,\int_D\,\nabla\cdot\{e^{\tau(x\cdot\omega-h_D(\omega)}\,\omega\}\,dx\\
\\
\displaystyle
&
\displaystyle
=\frac{1}{\tau}
\int_{\partial D}\,e^{\tau(x\cdot\omega-h_D(\omega))}\,\nu\cdot\omega\,dS(x).
\end{array}
$$
Then from this and (3.39) we obtain
$$
\displaystyle
\frac{\Vert v\Vert_{L^1(D)}}{\Vert v\Vert_{L^2(D)}}
\le
\frac{1}{\sqrt{C}}\int_{\partial D}\,e^{\tau(x\cdot\omega-h_D(\omega))}\,\nu\cdot\omega\,dS(x).
$$
Now under the assumption on $\omega$ and $D$ one can apply Lebesgue's dominated convergence theorem
to this right-hand side and conclude that (3.38) is valid.

\noindent
$\Box$

The proof of Theorem 3.16 is as follows.  We start with the well known identity
$$
\displaystyle
-\int_{\partial\Omega}\,\left\{\Lambda_V(v\vert_{\partial\Omega})-\frac{\partial v}{\partial\nu}\right\}\,\overline{v}\,dS
=\int_D\,V\vert v\vert^2\,dx+\int_D\,V\,(u-v)\overline {v}\,dx.
\tag {3.40}
$$
This together with Lemma 3.1 yields
$$\begin{array}{l}
\displaystyle
\,\,\,\,\,\,
e^{-2\tau h_D(\omega)}\vert I_{\omega,\omega^{\perp}}(\tau)\vert
\\
\\
\displaystyle
\ge
e^{-2\tau h_D(\omega)}
\,
\Vert v\Vert_{L^2(D)}^2
\,\left\{
\frac{\displaystyle \left\vert
\int_D\,V\vert e^{-\tau h_D(\omega)}\,v\vert^2\,dx
\right\vert}
{\displaystyle
e^{-2\tau h_D(\omega)}\,\Vert v\Vert_{L^2(D)}^2}
-C(V)\Vert V\Vert_{L^{\infty}(D)}^2\,\frac{\Vert v\Vert_{L^1(D)}}{\Vert v\Vert_{L^2(D)}}\,
\right\}.
\end{array}
\tag {3.41}
$$

First consider the case when $V(x)\ge C$ for almost all $x\in D_{\omega}(\delta)$.
We see that
$$\begin{array}{ll}
\displaystyle
\int_D\,V\vert e^{-\tau h_D(\omega)}\,v\vert^2\,dx
&
\displaystyle
=\int_{D_{\omega}(\delta)}\,V\vert e^{-\tau h_D(\omega)}\,v\vert^2\,dx
+\int_{D\setminus D_{\omega}(\delta)}\,V\vert e^{-\tau h_D(\omega)}\,v\vert^2\,dx
\\
\\
\displaystyle
&
\displaystyle
\ge
C\int_{D_{\omega}(\delta)}\,\vert e^{-\tau h_D(\omega)}\,v\vert^2\,dx
+O(e^{-2\tau\delta})
\\
\\
\displaystyle
&
\displaystyle
=Ce^{-2\tau h_D(\omega)}\,\Vert v\Vert_{L^2(D)}^2+O(e^{-2\tau\delta}).
\end{array}
$$
Thus a combination of this and (3.39) yields:
there exist $C'>0$ and $\tau_0'>0$ such that, for all $\tau\ge\tau_0'$
$$\displaystyle
\frac{\displaystyle
\int_D\,V\vert e^{-\tau h_D(\omega)}\,v\vert^2\,dx}
{e^{-2\tau h_D(\omega)}\,\Vert v\Vert_{L^2(D)}}
\ge C'.
$$
Now from this together with (3.41) and Lemma 3.2 we conclude that 
there exist $C">0$ and $\tau_0">0$ such that, for all $\tau\ge\tau_0"$ 
$$\displaystyle
e^{-2\tau h_D(\omega)}\tau^2\vert I_{\omega,\omega^{\perp}}(\tau)\vert\ge C".
$$ 
Note that this is valid also for the case when $-V(x)\ge C$ for almost all
$x\in D_{\omega}(\delta)$.  

On the other hand, from Lemma 3.1 and (3.40) we see that $e^{-2\tau h_D(\omega)}\,\vert I_{\omega,\omega^{\perp}}(\tau)\vert$ is bounded 
from above with respect to $\tau>0$.  Now a standard argument yields the desired conclusions.

\subsection{Two applications more}

\subsubsection{Extracting geometry of polygonal sources}

Historically, the author submitted the first paper \cite{E00} to J. Inverse Ill-Posed Probl. on 1998
in which the enclosure method using the Dirichlet-to-Neumann map acting on the complex geometrical optics solutions
was introduced.  At the same time in \cite{ISource0} 
he found an idea that would serve as the basis for another version of the method in \cite{I1}.

The one of problems considered in \cite{ISource0} is as follows.

Let $\Omega$ be a bounded domain of $\Bbb R^2$ with Lipschitz boundary and let $F\in L^2(\Omega)$ and
satisfy $\text{supp}\,F\subset\Omega$.
Let $u\in H^1(\Omega)$ satisfy the Helmholtz equation
$$\begin{array}{ll}
\displaystyle
\Delta u+k^2u=F, & x\in\Omega.
\end{array}
\tag {3.42}
$$
Here $k\ge 0$ is fixed.

$\quad$

{\bf\noindent Problem 3.5.}
Assume that $k$ is known and the source term $F$ is unknown.
Extract information about the geometry of $\text{supp}\,F$ from the Cauchy data $u\vert_{\partial\Omega}$
and $\frac{\partial u}{\partial\nu}\vert_{\partial\Omega}$ for a fixed solution $u$ of (3.42).

$\quad$

When $k=0$ this is a reduction to bounded domain $\Omega$ of the two dimensional version of inverse source problem
in geophysics.  When $k>0$, $F$ is a mathematical model of a radiation source.  Note that we never specify
the boundary condition of the field $u$ on $\partial\Omega$.

It is well known that there is a non uniqueness result for the problem without restricting the class of unknown source $F$,
that is a-priori assumption of the form of $F$.
In \cite{ISource0} it is assumed that $F$ takes the form
$$\displaystyle
F(x)=
\left\{\begin{array}{ll}
\displaystyle
\rho(x), & x\in D,\\
\\
\displaystyle
0, & x\in\Omega\setminus D,
\end{array}
\right.
$$
where

$\bullet$  $D$ is an open set of $\Bbb R^2$ with $\overline D\subset\Omega$ and given by an inside of a polygon.

$\bullet$  $\rho\in L^2(D)$ and, for each vertex $p$ on $\partial D$ there exists an open disk $B_{\eta}(p)$ centered at $p$
with radius $\eta>0$, number $\alpha\in\,]0,\,1]$ and a function $\tilde{\rho}\,\in C^{0,\alpha}(\overline{B_{\eta}(p)})$ such that
$$\left\{
\begin{array}{ll}
\displaystyle
\rho=\tilde{\rho} & \text{on $B_{\eta}(p)\cap \overline D$,}
\\
\\
\displaystyle
\tilde{\rho}(p)\not=0. &
\end{array}
\right.
$$

Define the indicator function
$$\begin{array}{ll}
\displaystyle
I_{\omega}(\tau)=\int_{\partial\Omega}
\left(\frac{\partial u}{\partial\nu}\,v-\frac{\partial v}{\partial\nu}\,u\right)\,d\sigma, & \tau>0,
\end{array}
$$
where $\omega,\omega^{\perp}\in S^1$ with $\omega\cdot\omega^{\perp}=0$ and
$$\begin{array}{ll}
\displaystyle
v=e^{x\cdot\,(\tau\omega+i\sqrt{\tau^+k^2}\,\omega^{\perp})}, & x\in\Bbb R^2.
\end{array}
$$

\proclaim{\noindent Theorem 3.17(\cite{ISource0}).}  Let $\omega$ be regular with respect to $D$.  We have
$$\displaystyle
\lim_{\tau\rightarrow\infty}\,
\frac{1}{\tau}\log\,\vert I_{\omega}(\tau)\vert=h_D(\omega).
$$
\endproclaim

{\it\noindent Sketch of Proof.}
Integration by partys yields that we have
$$\displaystyle
e^{-\tau h_D(\omega)}I_{\omega}(\tau)=e^{-\tau h_D(\omega)}\int_D\rho(x)v(x)\,dx.
$$
Thus it suffices to study the asymptotic behaviour of the integral on this right-hand side.
Let $p$ be the point on $\partial D$ with $p\cdot\omega=h_D(\omega)$.
By localizing around $p$, one can reduce the domain $D$ to a triangle with vertex $p$ which
lies in the domain $h_D(\omega)-\delta<x\cdot\omega<h_D(\omega)$ with a small $\delta>0$.
Then Lemma 3.3 mentioned below which is nothing but Lemma 2.1 in \cite{ISource0}, provides us the necessary information to complete the proof.

\noindent
$\Box$

\proclaim{\noindent Lemma 3.3.} Let $D$ be given by the inside of a triangle and $\omega$ be regular with respect to $D$.
Let $p$ denote the point on $\partial D$ with $p\cdot\omega=h_D(\omega)$.
There exists a non zero complex
number $C_D$ depending only on $D$ and $\omega$ such that, for any $\rho\in
C^{0,\theta}(\overline D)$ as $\tau\longrightarrow\infty$ the
formula
$$\displaystyle
e^{-\tau p\cdot\omega}\int_D \rho(x)v(x)dx=C_D\frac{e^{i\sqrt{\tau^2+k^2}\,p\cdot\omega^{\perp}}}{2\tau^2}
\left(\rho(p)+O(\tau^{-\theta})\right),
$$
is valid.
\endproclaim

See Lemma 3.1 in \cite{ESurvey} for the three dimensional version of Lemma 3.3.

Note that, from the proof it is easy to see that one can replace $D$ with the inside of a finitely many union of polygons.
The vertex $p\in\partial D$ can be reduced to only the vertex of the convex hull of $D$ and the $\tilde{\rho}$ should be
taken only for such $p$.

The method works also for the inverse obstacle problem governed by the equation
$$\begin{array}{ll}
\displaystyle
\Delta u+k^2\,n(x)u=0, & x\in\Omega,
\end{array}
$$ 
where $k>0$, $n(x)=1+F(x)$ with $F(x)$ being the same as above.  We assume that $u\in H^2(\Omega)(\subset C^{0,\beta}(\overline\Omega))$, however,
no boundary conditions on $\partial\Omega$ are specified.

\proclaim{\noindent Theorem 3.18(\cite{ISource0}).}
Let $\omega\in S^1$ be regular with respect to $D$ appearing in $F$.
Let $p\in\partial D$ be the unique point on the line $x\cdot\omega=h_D(\omega)$.
If $u(p)\not=0$, then one can calculate $h_D(\omega)$ from the Cauchy data of $u$ on $\partial\Omega$.
\endproclaim

The idea presented here of connecting the {\it corner singularity} of the source domain
with the complex geometrical optical solution to extract the information of the vertices of the convex hull has long been
almost ignored or forgotten in the area of inverse source problems.  Recently in \cite{Bl} the idea has been used in showing
nonexistence of nonradiating sources having a corner.  See also \cite{BL} for an application to an inverse
source problem for the Navier equations.

Let us explain briefly the problem and the idea {\it in terms of the enclosure method} in two dimensions.
Fix $k\ge 0$.
Given $F\in L^2(\Bbb R^2)$ with compact support, let $u$ be the radiating solution of the Helmholtz equation
$$\begin{array}{ll}
\displaystyle
\Delta u+k^2 u=F, & x\in\Bbb R^2,
\end{array}
$$
where the radiating means $u$ satisfies the outgoing Sommerfeld radiation condition
$$\begin{array}{ll}
\displaystyle
\sqrt{r}\,\left(\frac{\partial u}{\partial r}-ik u\,\right)\longrightarrow 0, & r=\vert x\vert\longrightarrow 0.
\end{array}
$$
Then $u$ has the asymptotic expansion as $r\rightarrow\infty$ (cf. \cite{CK}):
$$\begin{array}{ll}
\displaystyle
u(r\vartheta)\sim\frac{e^{ikr} f(\vartheta)}{\sqrt{r}}, & \vartheta\in S^1.
\end{array}
$$
If $f(\vartheta)\equiv 0$, then the source $F$ is called  nonradiating.
The problem is: what is a nonradiating source?

Here we present an argument that ensures:
the source having the form $F(x)=\chi_D(x)\,\rho(x)$ can not be nonradiating for $D$ and $\rho$ the same as Theorem 3.17.
It is a contraction argument and shows a closed relationship between the probelm and the enclosure method.

The argument is as follows.
Choose a large open disc $B$ that contains $\overline{D}$.  
Since the restriction of $u$ to $B$ satisfies (3.42) with $\Omega=B$, one can compute the indicator function
$I_{\omega}(\tau)$ for an arbitrary $\omega\in S^1$ and $\tau>0$. 
Choose a regulr $\omega\in S^1$ with respect to $D$.  Then Lemma 3.3 yields
the indicator function $I_{\omega}(\tau)\not=0$ for a large $\tau>0$.
Now assume that $f(\vartheta)\equiv 0$. 
Then by the Rellich lemma (cf. \cite{CK}),
we have $u(x)\equiv 0$ for $x\in\Bbb R^2\setminus B$.  Thus the Cauchy data of $u$ on $\partial B$ vanish.
Therefore $I_{\omega}(\tau)\equiv 0$.  This is a contradiction.   The essence of this argument is:
the non trivialness of the indicator function $I_{\omega}(\tau)$ at a direction $\omega$ yields
that of $f(\theta)$.  Note also that, instead of the Rellich lemma, one can compute
the Cauchy data of $u$ on $\partial B$ from $f(\theta)$ by using an explicit formula (cf. \cite{CK}), 
and thus the indicator function too.

\subsubsection{The Cauchy problem for the stationary Schr\"odinger equation}

The stationary Schr\"odinger equation is a typical second order elliptic equation and the Cauchy problem
has various applications.  In \cite{ICauchy} the author considered a solution of the stationary Schr\"odinger equation
in a domain that is given by the intersection of a bounded convex set with a half-space and found a computation formula
of the value of the solution at an arbitrary point in the domain from the Cauchy data on a part of the boundary
lying in the half-space.

The key idea is as follows.

Modify the equation by using a {\it tetrahydron} having the given point as one of vertexes and construct the complex geometrical optics solution of \cite{SU1, SU2} for the modified equation.
Then, the asymptotic behaviour of the complex geometrical optics solution with respect to a parameter combined with
integration by parts yields the desired formula.

The result in \cite{ICauchy} covers three and two dimensional cases, however, 
let us explain the idea in the two-dimensional case in a simpler situation. See also \cite{ESurvey}
for an explanation in the three dimensional case.

We denote by $B$ the open disc centered at $0$ with radius $1$.
Let $\vert t\vert<1$ and set
$$\displaystyle
\Omega=B\cap\{x_2>t\}.
$$
Let $u\in H^2(\Omega)$ be an arbitrary solution of the stationary Schr\"odinger equation
$$\begin{array}{ll}
\displaystyle
-\Delta u+Vu=0, & x\in\Omega.
\end{array}
\tag {3.43}
$$
Here $V=V(x)$ is a known,  essentially bounded and complex-valued function on $\Omega$.
Set $\Gamma=\partial B\cap\{x\,\in\Bbb R^2\,\vert\,\,x_2>t\}$.

In this subsection we give a formula that yields an analytical way of calculating $u$ in $\Omega$
by using the Cauchy data of $u$ on $\Gamma$.
Note that since $\partial\Omega$  is Lipschitz, from the trace theorem we have
$u\vert_{\partial\Omega}\in L^2(\partial\Omega)$ and $\partial u/\partial\nu\vert_{\partial\Omega}\in L^2(\partial\Omega)$.

\noindent
The formula obtained is based on the property of Faddeev's Green function and an {\it exotic} representation
of the Dirac delta function.

Let $y$ be an arbitrary point in $\Omega$.
Choose an open set $D$ such that the shape of $D$ is triangle with vertex at $y$
and other vertexes are in the half plane $x_2<y_2$.
Let $\tau>0$.
We denote by $\chi_D$ the characteristic function of $D$.
Our formula is based on the existence of a special solution of the equation:
$$\begin{array}{ll}
\displaystyle
-\Delta v+Vv=\chi_D e^{\tau(x_2-y_2)}e^{i\tau x_1},
& x\in\Omega.
\end{array}
\tag {3.44}
$$
\proclaim{\noindent Proposition 3.2.} There exist a positive number $\tau_0$
and solutions $v_{\tau}\in H^2(\Omega)$ with $\tau\ge\tau_0$ of
(3.44) and depending on $y$ such that $w_{\tau}\equiv
e^{-\tau(x_2-y_2)}e^{-i\tau x_1}v_{\tau}=O(\tau)$ in $H^2(\Omega)$
as $\tau\longrightarrow\infty$.
\endproclaim
{\it\noindent Proof.} $\tilde{V}$ stands for the zero extension of $V$
outside $\Omega$. Set $z=\tau((0,1)^T+i(1,0)^T)$.
From Theorem 3.1 one knows that the operator
$L^2_{\delta}(\Bbb R^2)\ni f\longmapsto f+g_{z}\ast(\tilde{V}f)\in L^2_{\delta}(\Bbb R^2)$
is invertible for fixed $\delta$ with $-1<\delta<0$ provided
$\tau$ is large enough.
Define
$$\displaystyle
w'=\{I+g_{z}\ast(\tilde{V}\,\cdot\,)\}^{-1}(g_{z}\ast\chi_D).
$$
Then
$$\displaystyle
w'+g_z\ast(\tilde{V}w')
=g_{z}\ast\chi_D
\tag {3.45}
$$
and this thus gives
$$\begin{array}{ll}
\displaystyle
(-\Delta+2z\cdot\nabla) w'+\tilde{V}w'=\chi_D, & x\in\Bbb R^2.
\end{array}
\tag {3.46}
$$
Define $v_{\tau}\equiv e^{x\cdot z}e^{-\tau y_2}w'\vert_{\Omega}$.
Then (3.46) yields (3.44) and other desired properties of $v_{\tau}$
is a consequence of Theorem 3.1 and (3.45).

\noindent
$\Box$

Now we state the desired formula.

\proclaim{\noindent Theorem 3.19.}
Let $u$ be an arbitrary $H^2(\Omega)$ solution of (3.43).
Then the formula
$$\displaystyle
u(y)=\frac{2}{C_D}\lim_{\tau\longrightarrow\infty}
\tau^2 e^{-i\tau y_1}\int_{\partial B\cap\{x_2>t\}}
(\frac{\partial u}{\partial\nu}v_{\tau}-\frac{\partial v_{\tau}}{\partial\nu}u)d\sigma(x),
\tag {3.47}
$$
is valid.
\endproclaim
{\it\noindent Proof.}
Using the solutions $v_{\tau}$, from (3.43) and (3.44) we obtain
$$\begin{array}{c}
\displaystyle
\int_D e^{\tau(x_2-y_2)}e^{i\tau x_1}u(x)dx\\
\\
\displaystyle
=\int_{\partial B\cap\{x_2>t\}}(\frac{\partial u}{\partial\nu}v_{\tau}-\frac{\partial v_{\tau}}{\partial\nu}u)d\sigma(x)
+\int_{B\cap\{x_2=t\}}(\frac{\partial u}{\partial \nu}v_{\tau}-\frac{\partial v_{\tau}}{\partial\nu}u)d\sigma(x).
\end{array}
\tag {3.48}
$$
Choose $\epsilon>0$ in such a way that $0<\epsilon<y_2-t$.
Applying the trace theorem to the convex domain $B\cap\{t<x_2<t+\epsilon/2\}$,
from Proposition 3.2, we see that both  $\Vert v_{\tau}\Vert_{L^2(B\cap\{x_2=t\})}$ and
$\Vert\nabla v_{\tau}\Vert_{L^2(B\cap\{x_2=t\})}$ are exponentially decaying as $\tau\longrightarrow\infty$.
Since $u\in H^2(\Omega)$, by the Sobolev imbedding theorem,
$u$ is H\"older continuous.  Thus one can apply Lemma 3.3 in the case when $\rho=u$, $\omega=\mbox{\boldmath $e$}_2$ to the left-hand side on (3.48).  From these one obtains (3.47).

\noindent
$\Box$

Note that it is a routine to show that from the formula (3.47) in the special domain $\Omega$
one can deduce the weak unique continuation theorem for the Schr\"odinger equation in a general connected domain.

Since the complex geometrical optics solution is constructed by using Faddeev's Green function, this is considered as an application 
of his Green function to the Cauchy problem for an elliptic equation with a variable coefficient.
Note that, in \cite{IS2} a numerical implementation of the formula has been done.

\section{The Time Domain Enclosure Method}

All of the applications of the enclosure method presented in the previous sections are devoted to
inverse problems governed by elliptic equations.  This means that the observation data are independent 
of time variable.  However, there are a lot of important and interesting inverse problems whose governing
equations are not elliptic and involve time variable.
As a typical example, we can mention the inverse problem governed by the heat equation, the wave equation
or the Maxwell system.

\subsection{Inverse source problem for the heat equation}

Let $\Omega$ be a bounded domain of $\Bbb R^n$ ($n=2,3$) with smooth boundary.  Let $T$ be an arbitrary positive number.
Let $u=u(x,t)$ satisfy
$$\left\{
\begin{array}{ll}
u_t=\Delta u+f(x,t), & (x,t)\in\Omega\times\,]0,\,T[,
\\
\\
\displaystyle
u(x,0)=0, & x\in\Omega.
\end{array}
\right.
$$
A typical inverse source problem for the heat equation is: to extract information about the unknown source $f(x,t)$
from the data $u\vert_{\partial\Omega\times\,]0,\,T[}$ and $\frac{\partial u}{\partial\nu}\vert_{\partial\Omega\times\,]0,\,T[}$.

It is well known that one can not uniquely determine general source $f(x,t)$ from the data.
However, it is also known that under some a priori assumption on the form of $f(x,t)$ one can 
extract the full or partial information about the source, e. g., see \cite{EH, YO}.

In \cite{ISource} the author considered the following inverse problem.

$\quad$

{\bf\noindent Problem 4.1.}  Assume that there exist a positive number $T_0$ less than $T$ and point $x_0\in\Omega$
such that $f(x_0,T_0)\not=0$ and $f(x,t)=0$ for all $(x,t)\in\Omega\times\,]0,\,,T_0[$.
Extract $T_0$ and information about the location of the set $\{x\in\Omega\,\vert\,f(x,T_0)\not=0\}$ from
the data $u\vert_{\partial\Omega\times\,]0,\,T[}$ and $\frac{\partial u}{\partial\nu}\vert_{\partial\Omega\times\,]0,\,T[}$.

$\quad$

The time $T_0$ is just the time when the unknown source becomes {\it active}.  It seems that there is no results for determining $T_0$
together with the place where the source firstly active.

In \cite{ISource} using the idea of enclosure method, the author gave two explicit extraction formulae.
For this purpose we see exponential solutions for the homogeneous backward heat equation
$$\displaystyle
\partial_tv+\Delta v=0.
$$

The method of separation of variables tells us that there are two types of exponential solutions.  The first is
the real-valued function given by
$$\displaystyle
v=v_{r}(x,t)=e^{\sqrt{\tau}\,x\cdot\omega-\tau t},
$$
where $\omega$ is a given unit vector and $\tau>0$ a parameter. The asymptotic behaviour of $v_r(x,t)$ as $\tau\rightarrow\infty$
divides the space-time into two parts:
$$\displaystyle
\lim_{\tau\rightarrow\infty}v_r(x,t)
=
\left\{
\begin{array}{ll}
\displaystyle \infty & \text{if $t<0$,}
\\
\\
\displaystyle 0 & \text{if $t>0$.}
\end{array}
\right.
$$
In this section we call this function a {\it real} exponential solution.

The second is the complex-valued function given by
$$\displaystyle
v=v_c(x,t)=e^{x\cdot z-(z\cdot z)\,t},
$$
where the complex vector $z$ has the form
$$\displaystyle
z=c\tau\left(\omega+i\sqrt{1-\frac{1}{c^2\tau}}\,\omega^{\perp}\,\right),
\tag {4.1}
$$
the unit vectors $\omega$ and $\omega^{\perp}$ are perpendicular each other
and $\tau$ a parameter satisfying $\tau>c^{-2}$  with a fixed $c>0$.  Note that complex vector
$z$ satisfies the equation $z\cdot z=\tau$.
We call function $v_c$ a complex exponential solution.  We have
$$\displaystyle
\vert v_c(x,t)\vert=e^{\tau(c\,x\cdot\omega-t)}.
$$
Thus the asymptotic behaviour of $\vert v_c(x,t)\vert$ as $\tau\rightarrow\infty$
divides the space-time into two parts:
$$\displaystyle
\lim_{\tau\rightarrow\infty}\vert v_c(x,t)\vert
=
\left\{
\begin{array}{ll}
\displaystyle \infty & \text{if $c\,x\cdot\omega>t$,}
\\
\\
\displaystyle 0 & \text{if $c\,x\cdot\omega<t$.}
\end{array}
\right.
$$

In this section we present applications of solutions $v_r$ and $v_c$ to Problem 4.1.

\subsubsection{Application of complex exponential solution}
Let $n=2$.  We assume that the unknown source takes the form
$$\displaystyle
f(x,t)=\sum_{j=1}^N\chi_{P_j\times\,[T_j,T]}(x,t)\,\rho_j(x,t),
$$
where: each $P_j\subset\Omega$ is given by the interior of a polygon;
if $j\not=j'$, then $\overline{P_j}\cap\overline{P_{j'}}=\emptyset$;
$T_j$ satisfies $0\le T_j<T$; $\rho_j$ is H\"older continuous on $\overline{P_j}\times\,[T_j,\,T]$;
$\rho_j(p,T_j)\not=0$ at all vertices of the convex hull of $P_j$.

{\bf\noindent Definition 4.1.}  Define the indicator function of independent variable $\tau\in\,]c^{-2},\,\infty[$ by the formula
$$\displaystyle
I_{\omega,c}(\tau)=\int_0^T\int_{\partial\Omega}\left(\frac{\partial v_c}{\partial\nu}\,u-\frac{\partial u}{\partial\nu}\,v_c\,\right)\,dsdt.
$$

Before describing the result concerning the asymptotic behaviour of the indicator function we introduce some notation and notion.
Let $\omega(c)$ denote the unit vector in the space-time equipped with the Euclidean norm given by
$$\displaystyle
\omega(c)=\frac{1}{\sqrt{c^2+1}}\left(\begin{array}{c}
\displaystyle c\omega\\
\displaystyle
-1
\end{array}
\right).
\tag {4.2}
$$

We set
$$\displaystyle
D=\cup_{j=1}^N(P_j\times\,]T_j,\,T[)\subset \Omega\times\,]0,\,T[.
$$
Let $\vartheta$ be an arbitrary unit vector in the space-time.  Define
$$\displaystyle
h_D(\vartheta)=\sup_{(x,t)\in D}\,
\left(\begin{array}{c}
\displaystyle
x\\
\displaystyle
t
\end{array}
\right)
\cdot
\vartheta.
$$
The function $h_D:\vartheta\longmapsto h_D(\vartheta)$ is called the support function for $D$.

We say that the direction $\omega(c)$ in the space-time is {\it regular} with respect to $D$ if
the set
$$\displaystyle
\left\{(x,t)\in\Bbb R^2\times\Bbb R\,\left\vert\,
\right.
\left(\begin{array}{c}
\displaystyle
x\\
\displaystyle
t
\end{array}
\right)
\cdot
\omega(c)=h_D(\omega(c))
\right\}
\cap
\overline D
$$
consists of a single point in the space-time.

Now we state a result in \cite{ISource} in the case $n=2$.

\proclaim{\noindent Theorem 4.1.}  
Assume that direction $\omega(c)$ is regular with respect to $D$ and the observation time $T$ and $c$ 
satisfy
$$\displaystyle
\sup_{x\in\Omega}\left(\begin{array}{c}
\displaystyle
x
\\
\displaystyle
T
\end{array}
\right)
\cdot\omega(c)<h_D(\omega(c)).
\tag {4.3}
$$
Then, the formula
$$\displaystyle
\lim_{\tau\rightarrow\infty}\,
\frac{1}{\tau}\log\vert I_{\omega,c}(\tau)\vert
=\sqrt{c^2+1}\,h_D(\omega(c)),
\tag {4.4}
$$
is valid.  Moreover we have
$$\displaystyle
\lim_{\tau\rightarrow\infty}\,e^{\tau s}\vert I_{\omega,c}(\tau)\vert=
\left\{
\begin{array}{ll}
\displaystyle
0, & \text{if $s\le-\sqrt{c^2+1}\,h_D(\omega(c))$,}
\\
\\
\displaystyle
\infty, &
\text{if $s>-\sqrt{c^2+1}\,h_D(\omega(c))$.}
\end{array}
\right.
$$

\endproclaim

Note that the assumption (4.3) means that the set $\Omega\times\{T\}$ is
contained in the half-space time 
$\left(\begin{array}{c}
\displaystyle
x\\
\displaystyle
t
\end{array}
\right)
\cdot\omega<h_D(\omega(c)).
$

\subsubsection{Application of real exponential solution}

In this section we consider the case when $n=2, 3$.
In \cite{ISource}, using a real exponential solution we gave an extraction formula of $T_0$ in Problem 4.1
from the data.  Here we present a further application of the real exponential solution which is not presented in
\cite{ISource}.  The contents presented here are taken from the unpublished manuscript \cite{IDPE}.

{\bf\noindent Definition 4.2.}
Define the indicator function of independent variable $\tau\in\,]0,\,\infty[$ given by the formula
$$\displaystyle
I_{\omega,r}(\tau)=\int_0^T\int_{\partial\Omega}\left(\frac{\partial v_r}{\partial\nu}\,u-\frac{\partial u}{\partial\nu}\,v_r\,\right)\,dSdt.
$$

Here is the list of the assumtions on the unknown source $f(x,t)$.

$\bullet$  Source $f(x,t)$ takes the form $f(x,t)=\chi_D(x,t)\rho(x,t)$,
where $D\subset\overline{\Omega}\times\,[T_0,\,T]$ is a Lebesgue measurable set with $
D\cap\{(x,t)\,\vert\,t=T_0\}=\{x_1,\cdots,x_N\}\subset\Omega$; function $\rho(x,t)$ is essentially bounded on $D$ and coincides with a uniformly
H\"older continuous function with exponent $\theta_j\in\,]0,\,1]$ on the intersection of $D$ with an open neighbourhood
of each point $x_j$ of $x_1,\cdots,x_N$.

\noindent
Thus the source may appear at the time $T_0$ firstly and points $x_1,\cdots, x_N$.

$\bullet$  There exist a positive number $\delta<T-T_0$ and mutually disjoint sets $W_1,\cdots,W_N$ of
$\overline\Omega$ having positive $n$-dimensional Lebesgue measures such that the set $D\cap\,(\Bbb R^n\times\,[T_0,\,T_0+\delta[)$
is given by the disjoint union of the cone $V_j$ with the vertex at $x_j$ and the base $W_j$, that is
$$\displaystyle
V_j=\left\{(x_j+\frac{s}{\delta}\,(y-x_j),\,T_0+s)\,\vert\,y\in W_j,\,0\le s<\delta\right\}.
$$
\noindent
This assumption describes the standing behaviour of the source at $t=T_0$.  In particular, 
it requires that the section  $D\cap\,(\,\Bbb R^n\times\{T_0+s\})$
with $0<s<<1$ of source domain $D$ evaluates at each point $x_j$ of $x_1,\cdots, x_N$ in the future direction like a cone with vertex at the point.

$\bullet$  We assume that 
$$\displaystyle
\sum_{x_i\cdot\omega=\max_{j=1,\cdots,N}\,x_j\cdot\omega}\,\rho(x_i,T_0)\,\vert W_i\vert\not=0.
\tag {4.5}
$$

This last assumption ensures, in some sense, that the total source viewed from the direction $\omega$ is really active at $t=T_0$ at all the points $x_1,\cdots,x_N$ 
on the plane $x\cdot\omega=h_{D_0}(\omega)$, where $D_0=\{x_1,\cdots,x_N\}$ 
and $h_{D_0}(\omega)$ the value of the support function for $D_0$ at direction $\omega$.
For example, if $\rho(x_i,T_0)$ for all $i$ with $x_i\cdot\omega=h_{D_0}(\omega)$ have the same sign, then the condition (4.5) is satisfied.
In this case the set of all $x\in\Omega$ with $f(x,T_0)\not=0$ is contained in set $D_0$.

Under the assumptions listed above, we obtain the following result which is an extension of Theorem 3.2 in \cite{ISource}.

\proclaim{\noindent Theorem 4.2.}  As $\tau\rightarrow\infty$ we have
$$\begin{array}{ll}
\displaystyle
\log\vert I_{\omega,r}(\tau)\vert
&
\displaystyle
=-\tau\,T_0+\sqrt{\tau}\,h_{D_0}(\omega)-(n+1)\log\tau
\\
\\
\displaystyle
&
\displaystyle
\,\,\,
+\log\left\vert\sum_{x_i\cdot\omega=h_{D_0}(\omega)}\,\rho(x_i,T_0)\,\frac{\vert W_i\vert}{\delta^n}\,n!\right\vert
+O(\tau^{-\min\,(\frac{1}{2},\theta_0(\omega))}),
\end{array}
\tag {4.6}
$$
where 
$$\displaystyle
\theta_0(\omega)=\min\{\theta_i\,\vert\,x_i\cdot\omega=h_{D_0}(\omega)\}.
\tag {4.7}
$$
Moreover, we have
$$\displaystyle
\lim_{\tau\rightarrow\infty}\,e^{\tau s}\vert I_{\omega,r}(\tau)\vert=
\left\{
\begin{array}{ll}
0 & \text{if $s<T_0$,}\\
\\
\displaystyle
\infty
&
\text{if $s>T_0$.}
\end{array}
\right.
$$

\endproclaim

{\it\noindent Proof.}
It suffice to show that, as $\tau\rightarrow\infty$
$$
\displaystyle
e^{-\sqrt{\tau}h_{D_0}(\omega)}\,\tau^{n+1}\,e^{\tau T_0}I_{\omega,r}(\tau)
=\sum_{x_i\cdot\omega=h_{D_0}(\omega)}\,\rho(x_i,T_0)\,\frac{\vert W_i\vert}{\delta^n}\,n!+O(\tau^{-\min\,(\frac{1}{2},\theta_0(\omega))}).
\tag {4.8}
$$
Integration by parts gives us the representation formula of the indicator function
$$\begin{array}{ll}
\displaystyle
e^{\tau T_0}I_{\omega,r}(\tau)
&
\displaystyle
=e^{\tau T_0}\int_0^T\int_{\Omega}f(x,t)v_r(x,t)\,dxdt
-e^{\tau T_0}\int_{\Omega}u(x,T)v_r(x,T)\,dx
\\
\\
\displaystyle
&
\displaystyle
\equiv
I+II.
\end{array}
\tag {4.9}
$$
Since $T>T_0$, we see that, as $\tau\rightarrow\infty$
$$\displaystyle
II=O(e^{-\tau(T-T_0)/2}).
\tag {4.10}
$$
From the first two assumptions listed above on $f$ we have
$$
\displaystyle
I
=
\sum_{j=1}^N\,e^{\tau T_0}\,\int_{V_j}\rho(x,t)e^{\sqrt{\tau}\,x\cdot\omega-\tau t}\,dxdt
+O(e^{-\tau\delta/2}).
\tag {4.11}
$$
Here we claim that, as $\tau\rightarrow\infty$
$$\displaystyle
\tau^{n+1}\,e^{\tau T_0}\,\int_{V_j}\rho(x,t)e^{\sqrt{\tau}\,x\cdot\omega-\tau t}\,dxdt
=e^{\sqrt{\tau}\,x_j\cdot\omega}
\left(\rho(x_j,T_0)\frac{\vert W_j\vert}{\delta^n}\,n!+O(\tau^{-\min\,(\frac{1}{2},\,\theta_j\,)})\right).
\tag {4.12}
$$
This is proved as follows.  Using change of variables, one can write
$$\begin{array}{l}
\displaystyle
\,\,\,\,\,\,
e^{\tau T_0}\int_{V_j}\rho(x,t)e^{\sqrt{\tau}\,x\cdot\omega-\tau t}\,dxdt
\\
\\
\displaystyle
=e^{\sqrt{\tau}\,x_j\cdot\omega}\,
\int_0^{\delta}
\left(\frac{s}{\delta}\right)^n\,
e^{-\tau s}\,ds\,\int_{W_j}\,\rho(x_j+\frac{s}{\delta}\,(y-x_j),T_0+s)\,
e^{\frac{\sqrt{\tau}\,s}{\delta}\,(y-x_j)\cdot\omega}\,dy
\\
\\
\displaystyle
\equiv
III+IV,
\end{array}
\tag {4.13}
$$
where
$$\displaystyle
III=e^{\sqrt{\tau}\,x_j\cdot\omega}\,
\rho(x_j,T_0)
\int_0^{\delta}
\left(\frac{s}{\delta}\right)^n\,
e^{-\tau s}\,ds\,\int_{W_j}\,
e^{\frac{\sqrt{\tau}\,s}{\delta}\,(y-x_j)\cdot\omega}\,dy
$$
and
$$\displaystyle
IV
=e^{\sqrt{\tau}\,x_j\cdot\omega}\,
\int_0^{\delta}
\left(\frac{s}{\delta}\right)^n\,
e^{-\tau s}\,ds\,\int_{W_j}\,\left(\rho(x_j+\frac{s}{\delta}\,(y-x_j),T_0+s)-\rho(x_j,T_0)\right)\,
e^{\frac{\sqrt{\tau}\,s}{\delta}\,(y-x_j)\cdot\omega}\,dy.
$$
Set $c_j(\tau,y)=1-\frac{(y-x_j)\cdot\omega}{\delta\sqrt{\tau}}$ for $y\in W_j$.
One can write
$$\begin{array}{ll}
\displaystyle
e^{-\sqrt{\tau}\,x_j\cdot\omega}\,\tau^{n+1}\,III
&
\displaystyle
=\rho(x_j,T_0)
\int_{W_j}\left(\int_0^{\infty}\left(\frac{\xi}{\delta}\right)^n
\,e^{-c_j(\tau,y)\xi}\,d\xi\right)\,dy
\\
\\
\displaystyle
&
\displaystyle
\,\,\,
+\rho(x_j,T_0)
\int_{W_j}\left(\int_{\tau\delta}^{\infty}\left(\frac{\xi}{\delta}\right)^n
\,e^{-c_j(\tau,y)\xi}\,d\xi\right)\,dy.
\end{array}
$$
Since we have $\inf_{y\in W_j}c_j(\tau,y)>0$ for all $\tau\ge\tau_0$ with a positive number $\tau_0$ independent of $y\in W_j$,
it is easy to see that the second term on this right-hand side is dominated by
$O(e^{-C\tau})$ with a positive number $C$.

On the other hand we have
$$\begin{array}{ll}
\displaystyle
\int_{W_j}\left(\int_0^{\infty}\left(\frac{\xi}{\delta}\right)^n
\,e^{-c_j(\tau,y)\xi}\,d\xi\right)\,dy
&
\displaystyle
=\int_{W_j}\frac{dy}{c_j(\tau,y)^{n+1}\delta^n}\int_0^{\infty}s^n e^{-s}\,ds
\\
\\
\displaystyle
&
\displaystyle
=\int_{W_j}\frac{dy}{c_j(\tau,y)^{n+1}}\frac{n!}{\delta^n}\\
\\
\displaystyle
&
\displaystyle
=\frac{\vert W_j\vert}{\delta^n}\,n!+O(\tau^{-1/2}).
\end{array}
$$
Thus we obtain
$$\displaystyle
\tau^{n+1}\,III
=e^{\sqrt{\tau}\,x_j\cdot\omega}
\left(\rho(x_j,T_0)\,\frac{\vert W_j\vert}{\delta^n}\,n!+O(\tau^{-1/2})\right).
\tag {4.14}
$$
Besides we have
$$\begin{array}{ll}
\displaystyle
e^{-\sqrt{\tau}\,x_j\cdot\omega}\,\tau^{n+1}\,\vert IV\vert
&
\displaystyle
\le C\tau^{n+1}\,\int_0^{\delta}\left(\frac{s}{\delta}\right)^n\,e^{-\tau s}\,ds
\int_{W_j}\left(\left\vert \frac{s}{\delta}(y-x_j)\right\vert^2+s^2\right)^{\theta_j/2}
\,e^{\frac{\sqrt{\tau}\,s}{\delta}\,(y-x_j)\cdot\omega}\,dy\\
\\
\displaystyle
&
\displaystyle
\le
C'\tau^{n+1}\,\int_0^{\delta}s^{n+\theta_j}\,e^{-\tau s}\,ds\,\int_{W_j}\,e^{\frac{\sqrt{\tau}\,s}{\delta}\,(y-x_j)\cdot\omega}\,dy
\\
\\
\displaystyle
&
\displaystyle
=C'\tau^{-\theta_j}\,\int_{W_j}\,\left(\int_0^{\tau\delta}\xi^{n+\theta_j}\,e^{-c_j(\tau,y)\xi}\,d\xi\right)\,dy
\\
\\
\displaystyle
&
\displaystyle
=O(\tau^{-\theta_j}).
\end{array}
$$
Thus, from this together with (4.13) and (4.14) we obtain (4.12).

Now substituting (4.12) into (4.11), we obtain
$$\begin{array}{ll}
\displaystyle
e^{-\sqrt{\tau}\,h_{D_0}(\omega)}\,\tau^{n+1}\,I
&
\displaystyle
=\sum_{x_i\cdot\omega=h_{D_0}(\omega)}\,\left(\frac{\vert W_i\vert}{\delta^n}\,n!+O(\tau^{-\min\,(\frac{1}{2},\theta_i)})\right)\\
\\
\displaystyle
&
\displaystyle
\,\,\,
+\sum_{x_j\cdot\omega<h_{D_0}(\omega)}\,e^{-\sqrt{\tau}\,(h_{D_0}(\omega)-x_j\cdot\omega)}
\left(\frac{\vert W_i\vert}{\delta^n}\,n!+O(\tau^{-\min\,(\frac{1}{2},\theta_i)})\right)
+O(e^{-\tau\delta/4})
\\
\\
\displaystyle
&
\displaystyle
=\sum_{x_i\cdot\omega=h_{D_0}(\omega)}\,\frac{\vert W_i\vert}{\delta^n}\,n!
+O(\tau^{-\min\,(\frac{1}{2},\theta(\omega))}).
\end{array}
$$
From this, (4.10) and (4.9) yield (4.8).

\noindent
$\Box$

The formula (4.6) tells us that one can extract the occurring time $T_0$ and $h_{D_0}(\omega)$ for a given $\omega$
from the data in Problem 4.1.  As a corollary we obtain an estimation of the convex hull of the set $D_0=\{x_1,\cdots,x_N\}$.
Note that we do not assume that $N$ is known.

It should be pointed out that the space information of unknown source in the formula (4.6) is behind the time information.
Especially information about the source strength at the occurring time is very weak.
However, formula (4.4) tells us that we can obtain the space and time information at the same time.  
This shows an advantage of using complex exponential solutions in theoretical sense.

\subsection{Analytical method for an inverse heat conduction problem}

In this subsection we consider the following problem.

$\quad$

{\noindent\bf Problem 4.2.}  Give a point inside a heat conductive body.
Extract information about the time evolution of the temperature at the point
from the pair of the temperature and heat flux observed on a part of the surface of the body
over a finite time interval.

$\quad$

This is a typical and important ill-posed problem which has various applications \cite{A, BBC}.
When the situation is reduced to one space dimensional case and the governing equation of temperature field
is given by the heat equation, the problem is referred to as sideways heat equation and there are extensive mathematical
studies, e.g., \cite{Cara, Eld, Le, T}.
See also \cite{Jo} for general heat operators and references therein.

In \cite{IHC}, the author found two types of extraction formulae for the problem
when the governing equation of the temperature field is given by the heat equation in mutli space dimensions.

More precisely let $\Omega\subset\Bbb R^n$ ($n=1,2,3$) be a bounded domain with a smooth boundary and $0<T<\infty$.
We assume that the temperature field $u=u(x,t)$ satisfies the heat equation
$$\begin{array}{ll}
\displaystyle
u_t=\Delta u, & (x,t)\in\Omega\times\,]0,\,T[.
\end{array}
\tag {4.15}
$$
Let $\nu$ denote the unit outward normal vector field to $\partial\Omega$.

In short, the author found the following:  given an arbitrary point $(x_0,t_0)\in\,\Omega\times\,]0,\,T[$ where and when
one wants to know the value of $u$, there are non-empty open subsets $\Gamma\subset\partial\Omega\times\,]0,\,T[$ and $U\subset\Omega$ 
properly chosen for $(x_0,t_0)$ such that from the data $u(x,0)$ given at all $x\in U$
and the pair $(u(x,t), \frac{\partial u}{\partial\nu}(x,t))$ given at all $(x,t)\in\Gamma$
one can calculate the value of $u$ at $(x_0,t_0)$.

The two approaches that yield the value $u(x_0,t_0)$ are: one is an application of the enclosure method \cite{ESurvey};
another can be considered as an extension to the heat equation of Yarmukhamedov's formula \cite{Y2} which yields a computation formula
of the solution to the Cauchy problem for the Laplace equation.

Both approaches are based on the construction of an exotic fundamental solution for the backward heat equation
and its decaying property in a half spacetime.

\subsubsection{A special fundamental solution for the backward heat equation}

Given a complex vector $z=(z_1,\cdots,z_n)^T$ we know that
the function $(x,t)\longmapsto e^{x\cdot z-t(z\cdot z)}$ solves the backward heat equation.

Given a function $g=g(x,t)$ we construct a special solution of the equation
$$\begin{array}{ll}
-v_t-\Delta v=e^{x\cdot z-t(z\cdot z)}\,g, & (x,t)\in\Bbb R^n\times\Bbb R\equiv \Bbb R^{n+1},
\end{array}
\tag {4.16}
$$
that takes he form
$$\displaystyle
v(x,t)=e^{x\cdot z-t(z\cdot z)}\,w(x,t).
$$
Thus the construction is reduced to solving the equation
$$\begin{array}{ll}
\displaystyle
-w_t-2\,z\cdot\nabla w-\Delta w=g, & (x,t)\in\Bbb R^{n+1}.
\end{array}
$$
Taking the Fourier transform of both sides, we obtain
$$\displaystyle
P_z(\xi,\eta)\,\hat{w}(\xi,\eta)=\hat{g}(\xi,\eta),
$$
where
$$\displaystyle
P_z(\xi,\eta)=-i\eta-2i\,z\cdot\xi+\vert\xi\vert^2.
$$
We see that the set $\{(\xi,\eta)\in\Bbb R^{n+1}\,\vert\,P_z(\xi,\eta)=0\}$ is compact
and forms a submanifold of $\Bbb R^{n+1}$ with codimension $2$ provided $\text{Im}\,z\not=\mbox{\boldmath $0$}$.
Thus, for complex vector $z$ with $\text{Im}\,z\not=\mbox{\boldmath $0$}$ the function $\frac{1}{P_z(\xi,\eta)}$
defines a tempered distribution
on $\Bbb R^{n+1}$.

$\quad$

{\bf\noindent Definition 4.3.}  Given complex vector $z$ with $\text{Im}\,z\not=\mbox{\boldmath $0$}$
define the tempered distribution $G_z(x,t)$ as the inverse Fourier transform of $\frac{1}{P_z(\xi,\eta)}$,
that is
$$
\displaystyle
G_z(x,t)=\frac{1}{(2\pi)^{n+1}}
\int\,e^{i(x\cdot\xi)+it\eta}\frac{d\xi\,d\eta}{P_z(\xi,\eta)}
$$
and the distribution $K_z(x,t)$ by the formula
$$\displaystyle
K_z(x,t)=e^{x\cdot z-t(z,z)}G_z(x,t).
\tag {4.17}
$$
Note that $K_z$ is a solution of the equation
$$\begin{array}{ll}
\displaystyle
v_t+\Delta v+\delta(x,t)=0, & (x,t)\in\Bbb R^{n+1},
\end{array}
\tag {4.18}
$$
in the sense of distribution.

Given $s\in\Bbb R$ we denote by $L^2_s(\Bbb R^{n+1})$ the set of all tempered distributions
$g=g(x,t)$ such that $(1+\vert x\vert^2+t^2)^{\frac{s}{2}}\,g\in L^2(\Bbb R^{n+1})$
and set
$$\displaystyle
\Vert g\Vert_s
=\left(
\int \vert g(x,t)\vert^2(1+\vert x\vert^2+t^2)^s\,dx\,dt\right)^{\frac{1}{2}}.
$$
Let $G_z* g$ denote the convolution of tempered distribution $G_z$ with an arbitrary
rapidly decreasing function $g$ on $\Bbb R^{n+1}$.
Note that the set of all rapidly decreasing functions on $\Bbb R^{n+1}$ 
is dense in $L^2_s(\Bbb R^{n+1})$.

\proclaim{\noindent Theorem 4.3\,(\cite{IHC}).}  Let $-1<\delta<0$ and $R>0$.
Let $z=\mbox{\boldmath $a$}+i\,\mbox{\boldmath $b$}$ with $\vert\mbox{\boldmath $b$}\vert\ge R$.
Then, for all rapidly decreasing functions $g$ on $\Bbb R^{n+1}$ we have
$$
\left\{
\begin{array}{ll}
\Vert D_x^{\alpha}\,G_z*g\Vert_{\delta}\le C(R,\delta)\,(\sqrt{1+\vert\mbox{\boldmath $a$}\vert^2}+\vert\mbox{\boldmath $a$}\vert)
\,\vert\mbox{\boldmath $b$}\vert^{\vert\alpha\vert}\,\Vert g\Vert_{\delta+1}, & \vert\alpha\vert\le 2,
\\
\\
\displaystyle
\Vert D_t\,G_z*g\Vert_{\delta}
\le
C(R,\delta)\,(\sqrt{1+\vert\mbox{\boldmath $a$}\vert^2}+\vert\mbox{\boldmath $a$})
\,(2\vert\mbox{\boldmath $a$}\vert\,\vert\mbox{\boldmath $b$}\vert+\vert\mbox{\boldmath $b$}\vert^2)\,
\Vert g\Vert_{\delta+1}. &
\end{array}
\right.
$$

\endproclaim
{\it\noindent Sketch of Proof.}
Since we have
$$\begin{array}{ll}
\displaystyle
G_z(x,t)=G_{i\,\text{Im}\,z}\,(x-2t\text{Re}\,z,t), & \text{Im}\,z\not=\mbox{\boldmath $0$},
\end{array}
\tag {4.19}
$$
everything is reduced to the case when $\text{Re}\,z=\mbox{\boldmath $0$}$.
So we introduce
$$\displaystyle
F_{\tau}(x,t)=G_{i\,\tau\omega}(x,t),
$$
where $\omega\in S^{n-1}$ and $\tau>0$.
To study the property of $F_{\tau}*\,$ acting on the set of all rapidly decreasing functions on $\Bbb R^{n+1}$ we employ
the argument done in \cite{IB}.

The argument consists of three steps.  

$\bullet$  Step 1.   A relationship between the operator $F_{\tau}*\,$ and $F_1*\,$.  Given a distribution $g(x,t)$
define 
$$\begin{array}{ll}
g_{\lambda}(x,t)=g(\lambda x, \lambda^2 t), & \lambda>0.
\end{array}
$$
We have the relationship
$$\displaystyle
D_x^{\alpha}D_t^{\beta}\,F_{\tau}*g=\tau^{-2+\vert\alpha\vert+2\beta}\,\left\{D_x^{\alpha}D_t^{\beta}\,F_1*(g_{\tau^{-1}})\right\}.
\tag {4.20}
$$

$\bullet$  Step 2.  An estimation of a scaling effect on weighted $L^2$-norms.  
Given $R>0$ let $\tau\ge R$.  
We have
$$\left\{
\begin{array}{ll}
\displaystyle
\Vert g_{\tau}\Vert_s\le \frac{(C(R)^{s/2}}{\tau^{2s+\frac{n+2}{2}}}\,\Vert g\Vert_{s}, & s<0,\\
\\
\displaystyle
\Vert g_{\tau^{-1}}\Vert_{s'}
\le
\frac{\tau^{2s'+\frac{n+2}{2}}}{C(R)^{s'/2}}
\,\Vert g\Vert_{s'}, & s'>0,
\end{array}
\right.
\tag {4.21}
$$
where $C(R)=\min\,\{R^4, R^2, 1\}$.

$\bullet$  Step 3.  An weighted $L^2$-estimate for the operator $F_1*$.  We have the lower estimate
$$\begin{array}{ll}
\displaystyle
\vert P_{\,i\omega}(\xi,\eta)\vert
\ge \frac{1}{2}(\vert\xi\vert^2+\vert\eta\vert), &
\vert (\xi,\eta)\vert\ge 8\sqrt{1+8^2}.
\end{array}
$$
Let $-1<\delta<0$.  Using this inequality, a local representation of $\frac{1}{P_{i\,\omega}\,(\xi,\eta)}$ in each neighbourhood of some zero points
of $P_{i\,\omega}(\xi,\eta)$ and Lemma 3.1 in \cite{SU2}, we have
$$\begin{array}{ll}
\displaystyle
\Vert D_x^{\alpha}D_t^{\beta}F_1*g\Vert_{\delta}
\le C_{\delta}\Vert g\Vert_{1+\delta},
&
\vert\alpha\vert+2\beta\le 2.
\end{array}
\tag {4.22}
$$

Given $R>0$ let $\tau\ge R$.  From (4.20), (4.21) and (4.22) we obtain
$$\begin{array}{ll}
\displaystyle
\Vert D_x^{\alpha}D_t^{\beta}G_{i\tau\omega}*g\Vert_{\delta}
\le \frac{C_{\delta}\tau^{\vert\alpha\vert+2\beta}}{C(R)^{1/2}}\,\Vert g\Vert_{1+\delta},
&
\vert\alpha\vert+2\beta\le 2.
\end{array}
$$
A combination of this and (4.19) yields the desired estimate.

\noindent
$\Box$

In what follows we describe only the cases $n=2, 3$.

\subsubsection{First approach}

The first approach basically follows  the idea of the enclosure method originally applied to the Cauchy problem
for the stationary Schr\"odinger equation in \cite{ICauchy}.

Given $c>0$ and $\omega$, $\omega^{\perp}$ two real unit vectors of $\Bbb R^n$  with $\omega\cdot\omega^{\perp}=0$.
Let $z$ be the complex vector given by (4.1).  

Let $(x_0,t_0)\in\Bbb R^{n+1}$ and $D$ be a bounded open subset of $\Bbb R^{n+1}$ with $\overline D\subset\Omega\times\,]0,\,T[$
and $(x_0,t_0)\in\partial D$.

We assume that that $D$ is {\it visible} at $(x_0,t_0)$ as $\tau\rightarrow\infty$ from the complex direction $z$
as defined in \cite{IHC}, that is,
there exist $\mu>0$ and 
constant $C_D\not=0$ such that, for all $\rho\in C^{\infty}(\overline D)$ 
$$\displaystyle
\lim_{\tau\rightarrow\infty}
e^{-x_0\cdot z+t_0(z\cdot z)}\tau^{\mu}
\int_D\,e^{x\cdot z-t(z\cdot z)}\rho(x,t)\,dx\,dt=C_D\,\rho(x_0,t_0).
\tag{4.23}
$$
Note that $C_D$ may depend on $\omega$, $\omega^{\perp}$ and $c$ in (4.1) and $(x_0,t_0)$.

Let $\chi_D$ denote the characteristic function of $D$.  
From Theorem 4.3 we see that the map $g\longmapsto G_z*g\in L^2_{\delta}(\Bbb R^{n+1})$ has the unique extension
as a bounded linear operator of $L_{\delta+1}(\Bbb R^{n+1})$ into $L^2_{\delta}(\Bbb R^{n+1})$.  
We denote it by the same symbol.  
the function
$$\begin{array}{ll}
\displaystyle
v(x,t)=e^{x\cdot z-t(z\cdot z)}\,(G_z*g)(x,t), & g=\chi_D
\end{array}
\tag {4.24}
$$
satisfies the backward heat equation (4.16) in the sense of distribution
\footnote{
In addition, we have $e^{-x\cdot z+t(z\cdot z)}v(x,t)\in L^2_{\delta}(\Bbb R^{n+1})$
and such $v$ is unique.  This fact is a consequence of Theorem 7.1.27 of \cite{H} as done
in \cite{SU2}.  However, it should be pointed out that, for our purpose this information is unnecessary since
the $v$ is explicitly constructed as (4.24).}.

Let $u$ be a solution of (4.15).
Given $(x_0,t_0)\in \Omega\times\,]0,\,T[$ one wishes to
extract the value  $u(x_0,t_0)$ from the data $u\vert_{\Gamma}$, $\frac{\partial u}{\partial\nu}\vert_{\Gamma}$
and $u(\,\cdot\,,0)\vert_U$.

\proclaim{\noindent Theorem 4.4 (\cite{IHC}).}  Let $\omega(c)$ be the unit vector in the spacetime given by (4.2).
Assume that $T$, $\Gamma$ and $U$ satisfy the following conditions:
$$\displaystyle
\sup_{x\in\Omega}\,
\left(\begin{array}{c}
\displaystyle
x\\
\displaystyle
T
\end{array}
\right)
\cdot\omega(c)<
\left(\begin{array}{c}
\displaystyle
x_0\\
\displaystyle
t_0
\end{array}
\right)
\cdot\omega(c);
\tag {4.25}
$$
$$\displaystyle
\sup_{x\in\Omega\setminus U}\,
\left(\begin{array}{c}
\displaystyle
x\\
\displaystyle
0
\end{array}
\right)
\cdot\omega(c)<
\left(\begin{array}{c}
\displaystyle
x_0\\
\displaystyle
t_0
\end{array}
\right)
\cdot\omega(c);
\tag {4.26}
$$
$$\displaystyle
\sup_{(x,t)\in\,(\partial\Omega\times\,]0,\,T[)\setminus\Gamma}\,
\left(\begin{array}{c}
\displaystyle
x\\
\displaystyle
t
\end{array}
\right)
\cdot\omega(c)<
\left(\begin{array}{c}
\displaystyle
x_0\\
\displaystyle
t_0
\end{array}
\right)
\cdot\omega(c).
\tag {4.27}
$$
Let $v$ be the function given by (4.24).    Then, we have
$$\begin{array}{l}
\displaystyle
\,\,\,\,\,\,
C_D\,u(x_0,t_0)
\\
\\
\displaystyle
=-\,\lim_{\tau\rightarrow\infty}\,\tau^{\mu}\,e^{-x_0\cdot z+t_0(z\cdot z)}\,
\left\{\int_{\Gamma}\,
\left(\frac{\partial v}{\partial\nu}\,u-\frac{\partial u}{\partial\nu}\,v\right)\,dSdt
-\int_U\,v(x,0)u(x,0)\,dx\right\}.
\end{array}
\tag {4.28}
$$

\endproclaim

The conditions (4.25), (4.26) and (4.27) are the restriction on the size of $T$,  positioning of $\Gamma$ and $U$ depending 
on $(x_0,t_0)$ and $\omega(c)$.

For example, (4.25) means that the set $\overline\Omega\times\{T\}$
is contained in the half-spacetime
$$
\left\{(x,t)\,\left\vert\right.\,
\left(\begin{array}{c}
\displaystyle
x\\
\displaystyle
t
\end{array}
\right)
\cdot\omega(c)
<
\left(\begin{array}{c}
\displaystyle
x_0\\
\displaystyle
t_0
\end{array}
\right)
\cdot\omega(c)
\right\}.
$$
Since the time component of $\omega(c)$ is negative, if $\omega(c)$ is given, then by choosing large $T$ the condition (4.25) is satisfied.
Or if $T$ is given, then by choosing small $c$, we have the same conclusion.  Note that making $c$ in $\omega(c)$ small means that
the ``propagation speed'' $c^{-1}$ of $\vert v_c(x,t)\vert=e^{\tau(c\,x\cdot\omega-t)}$ large.  So we called the $c$ the {\it virtual slowness}
in \cite{IBrazil}.

Conditions (4.25), (4.26) and (4.27) together with Theorem 4.3  ensure
that $e^{-x_0\cdot z+t_0(z\cdot z)}v$ together with derivatives
is exponentially decaying on $\Omega\times\{T\}$, $\Omega\setminus U$ and $(\partial\Omega\times\,]0,\,T[)\setminus\Gamma$ as $\tau\rightarrow\infty$.
Then, integration by parts together with (4.23) with $\rho=u$ yields the desired formula.

So the next problem is to find a suitable $D$ which is visible at $(x_0,t_0)$ as $\tau\rightarrow\infty$
from the direction $z$ given by (4.1).  For this we have

\proclaim{\noindent Proposition 4.1 (\cite{IHC}).}  Let $\delta>0$.
Let $D\subset\Bbb R^{n+1}$ be a finite cone with a vertex at $P=(x_0,t_0)$ and a bottom face $Q\not=\emptyset$ that is a bounded open subset of
$n$-dimensional hyper plane
$$
\displaystyle
\left(\begin{array}{c}
\displaystyle
x\\
\displaystyle
t
\end{array}
\right)
\cdot\omega(c)=
\left(\begin{array}{c}
\displaystyle
x_0\\
\displaystyle
t_0
\end{array}
\right)
\cdot\omega(c)-\delta.
\tag {4.29}
$$
If $\rho\in C^{0,\theta}(\overline D)$ with $0<\theta\le 1$, then we have
$$\displaystyle
\lim_{\tau\rightarrow\infty}\,\frac{2}{n!}\,(c\tau)^{n+1}\,e^{-x_0\cdot z+t\,(z\cdot z)}
\,\int_D\,e^{x\cdot z-t\,(z\cdot z)}\,\rho(x,t)\,dxdt
=K_D\,\rho(P),
$$
where
$$\displaystyle
K_D=2\delta\,\int_Q\,\frac{dS(y)}
{\displaystyle
\left(\frac{\delta\,\sqrt{c^2+1}}{c}-i\,(y-P)\cdot
\left(\begin{array}{c}
\displaystyle
\omega^{\perp}\\
\displaystyle
0
\end{array}
\right)
\right)^{n+1}
}.
\tag {4.30}
$$

\endproclaim

The key point of the proof is the representation of $D$
$$\displaystyle
D=\cup_{0<s<\delta}\,\left\{P+\frac{s}{\delta}\,(y-P)\,\vert y\in Q\right\}
$$
and a reduction to the case when $\rho(x,t)\equiv \rho(P)$.  See the proof of Theorem 2.2 and Lemma 4.1 in \cite{ISource}.

Thus, if $K_D\not=0$, then the $\mu$ and $C_D$ in (4.23) are given by $\mu=n+1$ and
$$\displaystyle
C_D=\frac{n!\,K_D}{2\,c^{n+1}}.
$$

It seems that it is not easy to show that $K_D\not=0$ for $D$ with general $Q$.  In \cite{IHC} we showed that
it is true for a special choice of $Q$.

The concrete choice of $Q$ in the case when $n=2$ is as follows.
Given $\delta>0$ choose arbitrary two points $x_1$, $x_2$ on the line $x\cdot\omega=x_0\cdot\omega-\delta\,\frac{\sqrt{1+c^2}}{c}$
in such a way that the orientation of the two vectors $\omega$ and $x_1-x_2$ coincides with that of the standard bases $\mbox{\boldmath $e$}_1$
and $\mbox{\boldmath $e$}_2$ of $\Bbb R^2$.  Let $Q$ be the inside of the triangle in $\Bbb R^{2+1}$ with the vertices $(x_1,t_0)$,
$(x_2,t_0)$ and $(x_0,t_0+\delta\sqrt{1+c^2}\,)$.  We see that the $Q$ is on the plane (4.29).
In this case the $K_D$ given by (4.30) satisfies $K_D\not=0$ and $C_D$ takes the form
$$\displaystyle
C_D=
\frac{\displaystyle
\vert(\mbox{\boldmath $\nu$}_3\times\mbox{\boldmath $\nu$}_2)
\times(\mbox{\boldmath $\nu$}_1\times\mbox{\boldmath $\nu$}_3)\vert}
{\displaystyle
\{(\mbox{\boldmath $\nu$}_3\times\mbox{\boldmath $\nu$}_2)\cdot\vartheta\}
\{(\mbox{\boldmath $\nu$}_1\times\mbox{\boldmath $\nu$}_3)\cdot\vartheta\}
},
$$
where $\mbox{\boldmath $\nu$}_1$, $\mbox{\boldmath $\nu$}_2$, $\mbox{\boldmath $\nu$}_3$ are the unit outward normal vectors to all faces
of $D$ except $Q$ which are numbered in a suitable order and
$$\displaystyle
\vartheta=\left(\begin{array}{c}
\displaystyle
c(\omega+i\omega^{\perp})\\
\displaystyle
-1
\end{array}
\right).
$$ 
For more explanation of the derivation together with the case when $n=3$ we refer the reader to \cite{IHC}.

Since $v$ in Theorem 4.4 is given by (4.24), it is important to compute $G_z$.  In \cite{IHC} we have given
its explicit form, that is the formula
$$\begin{array}{ll}
\displaystyle
G_z(x,t)
&
\displaystyle
=e^{-i(x-2t\,\mbox{\boldmath $a$})\cdot\mbox{\boldmath $b$}-\vert\mbox{\boldmath $b$}\vert^2\,t}\\
\\
\displaystyle
&
\displaystyle
\,\,\,
\times
\left\{
-\left(\frac{\vert\mbox{\boldmath $b$}\vert}{2\pi}\right)^n
\,\int_{\vert\xi\vert<1}
\,
e^{i\vert\mbox{\boldmath $b$}\vert\,(x-2t\,\mbox{\boldmath $a$})\cdot\xi}
\,e^{\vert\xi\vert^2\,\vert\mbox{\boldmath $b$}\vert^2\,t}\,d\xi
+H(-t)
\left(\frac{1}{2\sqrt{\pi\,\vert t\vert}}\right)^n\,e^{\frac{\vert x-2t\,\mbox{\boldmath $a$}\vert^2}{4t}}\,
\right\},
\end{array}
$$
where $H(t)$ is the Heaviside function, $\text{Re}\,z=\mbox{\boldmath $a$}$ and $\text{Im}\,z=\mbox{\boldmath $b$}$.

It will be interested to do a numerical implementation of formula (4.28) as done in \cite{IS2} for the Laplace equation.

\subsubsection{Second approach}

The second approach is to make use of $K_z$ given by (4.17).  Since $K_z$ satisfies (4.18), we have an integral representation of the solution $u$ at $(x_0,t_0)$ in terms of the surface integrals on $\partial(\Omega\times\,]0,\,T[)$.
Then vanishing of some of the surface integrals  as $\tau\rightarrow\infty$ yields the computation formula of $u(x_0,t_0)$ from the data $u\vert_{\Gamma}$, $\frac{\partial u}{\partial\nu}\vert_{\Gamma}$
and $u(\,\cdot\,,0)\vert_U$.  This can be considered as an extension of Yarmukhamedov's formula
for the Laplace equation to the heat equation.

The key point is the decaying property of $K_z$ as $\tau\rightarrow\infty$ in a half-spacetime.

\proclaim{\noindent Proposition 4.2 (\cite{IHC}).}
For each $\alpha$ and $\beta$ the quantity
$$\displaystyle
e^{\tau\,\sqrt{1+\tau^2}\,\delta}\,
\sup\,
\left\{\vert\partial_x^{\alpha}\,\partial_t^{\beta}\,K_z(x,t)\vert\,\vert\,
\left(\begin{array}{c}
\displaystyle
x\\
\displaystyle
t
\end{array}
\right)
\cdot\omega(c)
<-\delta
\right\},
$$
is at most algebraically growing as $\tau\rightarrow\infty$.

\endproclaim

As a corollary under the same assumptions as Theorem 4.4  we obtain the formula
$$
\displaystyle
u(x_0,t_0)
=-\,\lim_{\tau\rightarrow\infty}\,
\left\{\int_{\Gamma}\,
\left(\frac{\partial v}{\partial\nu}\,u-\frac{\partial u}{\partial\nu}\,v\right)\,dSdt
-\int_U\,v(x,0)u(x,0)\,dx\right\},
$$
where $v(x,t)=K_z(x-x_0,t-t_0)$ with $z$ given by (4.1).

\subsection{Extracting discontinuity in a heat conductive body}

In this subsection we consider an inverse initial boundary value problem for parabolic equation with 
discontinuous coefficients.  The problem is to extract information about the location and shape 
of unknown inclusions embedded in a known isotropic heat conductive body from a set of input heat flux across 
the surface of the body and output temperature observed on the same surface.

Developing the enclosure method to the problem was initiated by the author in \cite{I4}.
Therein one-space dimensional case was considered and various formulations of the method
proposed.  Later, in higher-space dimensional case some of the formulations are realized in the articles \cite{IK1, IK2, IFR, IK4} together with a technical work \cite{IK3}.  We have also an application \cite{IIVisco} to an inverse cavity problem in a visco elastic body.
In \cite{IKW} numerical implementation
of a method in \cite{I4} has been done.  And there is an application of the time domain enclosure method to an inverse obstacle problem
governed by the Stokes system, see \cite{MST}.  However, the author thinks there should be something to do more based on 
the recent development of the method presented in the last section of this review article.

In this subsection we consider only two or three-space dimensional case.

Let $\Omega$ be a bounded domain of $\Bbb R^n$, $n=2,3$, with a smooth boundary.
Let $T>0$.  Given $f$ let $u=u_f$ be the solution of the initial boundary value problem for the parabolic
equation:
$$
\left\{
\begin{array}{ll}
\displaystyle
u_t-\nabla\cdot\gamma\nabla u=0, & (x,t)\in\Omega\times\,]0,\,T[,\\
\\
\displaystyle
\gamma\nabla u\cdot\nu=f, & (x,t)\in\partial\Omega\times\,]0,\,T[,\\
\\
\displaystyle
u(x,0)=0, & x\in\Omega,
\end{array}
\right.
\tag {4.31}
$$
where $\gamma=\gamma(x)=(\gamma_{ij}(x))$ satisfies the following conditions.

$\bullet$  For each $i,j=1,\cdots,n$, $\gamma_{ij}(x)$ is real, belongs to $L^{\infty}(\Omega)$ and satisfies
$\gamma_{ij}(x)=\gamma_{ji}(x)$.

$\bullet$  There exists a positive constant $C$ such that $\gamma(x)\xi\cdot\xi\ge C\vert\xi\vert^2$ for all
$\xi\in\Bbb R^n$ and a.e. $x\in\Omega$.

We assume that there exists an ope set $D$ of $\Bbb R^n$ with a smooth boundary such that $\overline D\subset\Omega$
and $\gamma(x)$ a.e. $x\in\Omega\setminus D$ coincides with the $n\times n$ identity matrix $I_n$ multiplied by a
smooth positive function $\gamma_0(x)$ of $x\in\overline\Omega$ and satisfies one of the following
two conditions.

(A1)  There exists a positive constant $C'$ such that $-(\gamma(x)-\gamma_0(x)I_n)\xi\cdot\xi
\ge C'\vert\xi\vert^2$ for all $\xi\in\Bbb R^n$ and $a.e.x\in D$.

(A2)  There exists a positive constant $C'$ such that $(\gamma(x)-\gamma_0(x)I_n)\xi\cdot\xi
\ge C'\vert\xi\vert^2$ for all $\xi\in\Bbb R^n$ and $a.e.x\in D$.

Thus $\gamma$ takes the form
$$\displaystyle
\gamma(x)=
\left\{\begin{array}{ll}
\displaystyle
\gamma_0(x)\,I_n, & x\in\Omega\setminus D,\\
\\
\displaystyle
\gamma_0(x)\,I_n+h(x), & x\in D,
\end{array}
\right.
$$
where $h(x)=\gamma(x)-\gamma_0(x)\,I_n$ a.e. $x\in D$.  The set $D$ is a model of an inclusion embedded in $\Omega$
and (A1)/(A2) means that $D$ has a lower/higher conductivity from reference conductivity $\gamma_0\,I_n$.

$\quad$

{\bf\noindent Problem 4.3.} Fix $T$.  Assume that $\gamma_0$ is known and that both $D$ and $h$ are unknown.
Extract information about the location and shape of $D$ from a set of the pair
of temperature $u_f(x,t)$ and heat flux $f(x,t)$ given at all $x\in\partial\Omega\times\,]0,\,T[$.

$\quad$

\subsubsection{A reduction to a general problem depending on a large parameter}

The time domain enclosure method developed for Problem 4.3 starts with making the transform
$$\begin{array}{lll}
\displaystyle
w_f(x,\tau)=\int_0^T e^{-\tau t}\,u_f(x,t)\,dt, & x\in\Omega, & \tau>0.
\end{array}
$$
We see that the function $w=w_f$ satisfies
$$
\left\{\begin{array}{ll}
\displaystyle
\nabla\cdot\gamma\nabla w-\tau w=e^{-\tau T}\,u_f(x,T), & x\in\Omega,\\
\\
\displaystyle
\gamma\nabla w\cdot\gamma=g(x,\tau), & x\in\partial\Omega,
\end{array}
\right.
\tag {4.32}
$$
where
$$\begin{array}{ll}
\displaystyle
g(x,\tau)=\int_0^Te^{-\tau t}f(x,t)\,dt, & x\in\partial\Omega.
\end{array}
\tag {4.33}
$$

We generalize this as below.

Given $F(\,\cdot\,,\tau)\in L^2(\Omega)$ and $g(\,\cdot\,,\tau)\in H^{-1/2}(\partial\Omega)$ with $\tau>0$,
let $w=w(\,\cdot\,\tau)\in H^1(\Omega)$ be the weak solution of
$$\left\{
\begin{array}{ll}
\displaystyle
\nabla\cdot\gamma\nabla w-\tau w=e^{-\tau T}\,F(x,\tau), & x\in\Omega,\\
\\
\displaystyle
\gamma\nabla w\cdot\nu=g(x,\tau), & x\in\partial\Omega.
\end{array}
\right.
\tag {4.34}
$$

And consider the following problem.

$\quad$

{\bf\noindent Problem 4.4.}  Fix $T$.
Assume that $\gamma_0$ is known and $F(x,\tau)$ , $D$ and $h$ are all unknown.
Extract information about the geometry of $D$ from a set of the pair of $w(\,\cdot\,,\tau)\vert_{\partial\Omega}$
and $g(\,\cdot\,,\tau)$ with $\tau>>1$.

$\quad$

As a solution of Problem 4.4, in \cite{IFR} we obtained the following theorem.

\proclaim{\noindent Theorem 4.5}
Let $\sigma\in\{\frac{1}{2}, 1\}$.
Assume that there exist constants $C_1$ and $\kappa_1$ such that, as $\tau\rightarrow\infty$
$$\displaystyle
\Vert F(\,\cdot\,,\tau)\Vert_{L^2(\Omega)}=O(\tau^{\kappa_1}\exp{(C_1\tau^{\sigma})}).
\tag {4.35}
$$
Let $(v_{\tau})_{\tau\ge\tau_0}$ be a family of $H^1(\Omega)$ solutions of the equation
$$\begin{array}{ll}
\displaystyle
\nabla\cdot\gamma_0\nabla v-\tau v=0, & x\in\Omega
\end{array}
\tag {4.36}
$$
and satisfy the conditions, for some constants $\kappa_2$, $\kappa_3$, $\kappa_4$,
$C_2$, $C_3$ and $C_4>0$
$$\displaystyle
\Vert\nabla v_{\tau}\Vert_{L^2(D)}=O(\tau^{\kappa_2}\exp{(C_2\tau^{\sigma})}),
\tag {4.37}
$$
$$\displaystyle
\Vert\nabla v_{\tau}\Vert_{L^2(D)}\ge C_4\,\tau^{\kappa_3}\exp{(C_2\tau^{\sigma})},
\tag {4.38}
$$
$$\displaystyle
\Vert v_{\tau}\Vert_{H^1(\Omega)}=O(\tau^{\kappa_4}\exp{(C_3\tau^{\sigma})}).
\tag {4.39}
$$
Let $g=g(\,\cdot\,,\tau)$ have the form
$$\begin{array}{ll}
\displaystyle
g(x,\tau)=\Psi(\tau)\,\gamma_0\frac{\partial v_{\tau}}{\partial\nu}, & x\in\partial\Omega,
\end{array}
\tag {4.40}
$$
where $\Psi$ satisfies the conditions, for constants $\mu$ and $\mu'$
$$\left\{
\begin{array}{l}
\displaystyle
\liminf_{\tau\rightarrow\infty}\,\tau^{\mu}\,\vert\Psi(\tau)\vert>0,
\\
\\
\displaystyle
\vert\Psi(\tau\vert=O(\tau^{\mu'}).
\end{array}
\right.
\tag {4.41}
$$
If $\sigma=1$ and $T$ satisfies
$$\displaystyle
T>C_1+C_3-2C_2,
\tag {4.42}
$$
or $\sigma=\frac{1}{2}$ and $T$ is an arbitrary positive number,
then we have
$$\displaystyle
\lim_{\rightarrow\infty}\frac{1}{2\tau^{\sigma}}
\,\log\left\vert\,\int_{\partial\Omega}
\,\left(g\overline{v_{\tau}}-w\,\gamma_0\frac{\partial\overline v_{\tau}}{\partial\nu}\,\right)\,dS\,\right\vert
=C_2.
\tag {4.43}
$$

\endproclaim

Some remarks are in order.

$\bullet$  The concrete values of the constants $\kappa_1$, $\kappa_2$, $\kappa_3$, $\kappa_4$, $C_4$, $\mu$ and $\mu'$
do not affect the conclusion.

$\bullet$  From (4.38) and (4.39) we have $C_3\ge C_2$.

$\bullet$  It follows from (4.37) and (4.38) that the constant $C_2$ is uniquely determined by the formula
$$\displaystyle
\lim_{\tau\rightarrow\infty}\,\frac{1}{\tau^{\sigma}}\,\log\,\Vert\nabla v_{\tau}\Vert_{L^2(D)}=C_2.
$$
Thus $C_2$ should have some information about the geometry of $D$.

The point in the application of Theorem 4.5  to Problem 4.3 is to prescribe the special heat flux $f$ with parameter $\tau$
given by
$$\begin{array}{ll}
\displaystyle
f=f_{\tau}(x,t)=\varphi(t)\,\gamma_0\frac{\partial v_{\tau}}{\partial\nu}(x), & (x,t)\in\partial\Omega\times\,]0,\,T[,
\end{array}
\tag {4.44}
$$
where $\varphi=\varphi(t)$ is a function of $t\in\,]0,\,T[$.  
Then $g$ given by (4.33) takes the form (4.40) with $\Psi(\tau)$ given by
$$\displaystyle
\Psi(\tau)=\int_0^Te^{-\tau t}\,\varphi(t)\,dt.
$$
Under a suitable condition on the behaviour of $\varphi(t)$ as $t\downarrow 0$ together with $\varphi\in L^2(0,\,T)$,
we see that conditions on (4.41) are satisfied.

Besides $F(\,\cdot\,,\tau)$ in (4.34) becomes $F(\,\cdot\,,\tau)=u_{f_{\tau}}(\,\cdot\,,T)$.
We always consider the solution class of $u_f$ for general $f$ such that
$u_f(\,\cdot\,,T)$ has a meaning as an element
in $L^2(\Omega)$ and has the estimate
$$\displaystyle
\Vert u_f(\,\cdot\,,T)\Vert_{L^2(\Omega)}
\le C\Vert f\Vert _{L^2(0,\,T;H^{-1/2}(\partial\Omega))},
$$
see \cite{IFR}.  This is coming from the well-posedness of the direct problem (4.31) in suitable function spaces.
Thus it follows from this and (4.39) that this $F(\,\cdot\,,\tau)$ satisfies (4.35) with $C_1=C_3$.

Therefore as a corollary of Theorem 4.5 we conclude that: if $\sigma=1$ and $T>2(C_3-C_2)$
or $\sigma=\frac{1}{2}$ and $T$ is an arbitrary positive constant, then one has the formula (4.43).

Thus the main task is to find a suitable family $(v_{\tau})$ in Theorem 4.5 such that
constant $C_2$ yields some concrete information about the geometry of $D$.

\subsubsection{Using complex geometrical optics solutions}

One of possible choice is to make use of the complex geometrical optics solutions of the equation (4.36).

The construction is as follows.  We assume that $\gamma_0-1\in C^{\infty}_0(\Bbb R^n)$ and $\gamma_0(x)>0$ for all $x\in\Bbb R^n$.

Let $z$ be the complex vector given by (4.1), that is,
$$\displaystyle
z=c\tau\left(\omega+i\sqrt{1-\frac{1}{c^2\tau}}\,\omega^{\perp}\,\right),
$$
the unit vectors $\omega$ and $\omega^{\perp}$ are perpendicular each other
and $\tau$ a parameter satisfying $\tau>c^{-2}$  with a fixed $c>0$.

We construct a special solution of equation (4.36) having the form
$$\displaystyle
v_{\tau}(x)\sim\frac{e^{x\cdot z}}{\sqrt{\gamma_0(x)}}
$$
as $\tau\rightarrow\infty$, in an appropriate sense, and apply Theorem 4.5 in the case when $\sigma=1$.

For the purpose we follow the approach in \cite{SU2}. 
It is based on the change of dependent variable
$$\displaystyle
\frac{1}{\sqrt{\gamma_0}}\,\nabla\cdot\gamma_0\nabla\left(\frac{1}{\sqrt{\gamma_0}}\,\cdot\,\right)
=\Delta-V,
$$
where
$$\displaystyle
V=\frac{\Delta\sqrt{\gamma_0}}{\sqrt{\gamma_0}}.
$$

Set 
$$\displaystyle
v_{\tau}=\frac{e^{x\cdot z}}{\sqrt{\gamma_0}}\,(1+\epsilon_{\tau}),
\tag {4.45}
$$
where $\epsilon_{\tau}$ is a new unknown function of independent variables $x\in\Bbb R^n$.  
We see that the function $v_{\tau}$ satisfies equation (4.36) if
$\epsilon_{\tau}$ satisfies the equation
$$\begin{array}{ll}
\displaystyle
\left\{
\Delta+2z\cdot\nabla-\tau\left(\frac{1}{\gamma_0}-1\right)-V\right\}\epsilon_z
=\tau\left(\frac{1}{\gamma_0}-1\right)+V, & x\in\Bbb R^n.
\end{array}
\tag {4.46}
$$
In \cite{IFR} we have proved the unique existence of the solution of equation (4.46) in a weighted $L^2$-space.

More precisely, it is an application of the following result.

\proclaim{\noindent Theorem 4.6(\cite{IFR}).}
Let $-1<\delta<0$ and $a, b\in C^{\infty}_0(\Bbb R^n)$.  Given $\eta>0$,
there exist positive constants $C_j=C_j(a,b,\Omega,\delta,\eta)$. $j=5,6$ such that
if $c\ge C_5$ and $\tau\ge C_6$, then $c^2\tau>1$
and there exists a unique $\epsilon_z\in L^2_{\delta}(\Bbb R^n)$ with $z$ given by (4.1) such that
$$\begin{array}{ll}
\displaystyle
(\Delta+2z\cdot\nabla-\tau a-b)\epsilon_z=\tau a+b, & x\in\Bbb R^n.
\end{array}
$$
Moreover, $\epsilon_z\vert_{\Omega}$ can be identified with a function in $C^1(\overline\Omega)$
and 
$$\displaystyle
\Vert\epsilon_z\Vert_{L^{\infty}(\Omega)}+
\Vert\nabla\epsilon_z\Vert_{L^{\infty}(\Omega)}\le\eta.
\tag {4.47}
$$

\endproclaim

{\it\noindent Sketch of Proof.}
One can construct the special solution
of the equation $(-\Delta-2\,z\cdot\nabla)g+\delta(x)=0$ in $\Bbb R^n$ by 
the formula
$$\displaystyle
g=g_z(x)=-\frac{1}{(2\pi)^n}\,\int_{\Bbb R^n}\,\frac{e^{ix\cdot\xi}\,d\xi}{\vert\xi\vert^2-2i\,z\cdot\xi}.
$$
By \cite{SU2} we know that, given $f\in L^2_{\delta+1}(\Bbb R^n)$ the solution of the equation
$(-\Delta-2\,z\cdot\nabla)\Psi+f=0$ in $\Bbb R^n$ is unique in $L^2_{\delta}(\Bbb R^n)$
and has the expression $\Psi=g_z*f$, where $g_z*\,\cdot$ denotes the unique extension
of the original convolution operator $g_z*\,\cdot$ acting on all rapidly decreasing function on $\Bbb R^n$
as the bounded linear operator of $L^2_{\delta+1}(\Bbb R^n)$ into $L_{\delta}(\Bbb R^n)$.

We construct $\epsilon_z\in L^2_{\delta}(\Bbb R^n)$ as the solution of the integral equation
$$\displaystyle
\epsilon_z=-\tau\,g_z*(a\epsilon_z)-g_z*(b\epsilon_z)-\tau\,g_z*a-g_z*b,
$$
that is,
$$\displaystyle
(I-A_z)\epsilon_z=-\tau\,g_z*a-g_z*b,
\tag {4.48}
$$
where
$$\begin{array}{ll}
\displaystyle
A_zh=-\tau\,g_z*(a h)-g_z*(b h), & h\in L^2_{\delta}(\Bbb R^n).
\end{array}
$$
The key point is the estimate concerning with the operator $g_z*\,\cdot$:
$$\displaystyle
\Vert D^{\alpha}\,g_z*f\Vert_{L^2_{\delta}(\Bbb R^n)}
\le (c\tau\lambda)^{\vert\alpha\vert-1}\,C_{\delta,R}\,\Vert f\Vert_{L^2_{\delta+1}(\Bbb R^n)},
\tag {4.49}
$$
where $\vert\alpha\vert\le 2$, $c\tau\lambda\ge R$ and
$$\displaystyle
\lambda=\lambda(c,\tau)=\sqrt{1-\frac{1}{c^2\tau}}.
$$
The proof of (4.49) is a combination of Sylvester-Uhlmann's argument \cite{SU2} and
a scaling argument \cite{IB}.  See Section 3.2 in \cite{IFR} for the detail.  
From this we know that
$$\displaystyle
\Vert A_z\Vert_{B(L^2_{\delta}(\Bbb R^n)}
\le (c\lambda)^{-1}
C(\delta, R)\,(\Vert<x>\,a\Vert_{L^{\infty}(\Bbb R^n)}
+\tau^{1}\Vert<x>\,b\Vert_{L^{\infty}(\Bbb R^n)}),
$$
where $<x>=(1+\vert x\vert^2)^{1/2}$ and $c\tau\lambda\ge R$.   Roughly speaking, from this, we see that:
if $c$ is sufficiently large independent of $\tau\ge 1$, then the operator norm of $A_z$ is less than $1$.
Thus using the Neumann series, one gets the unique solution of (4.48).  The size control (4.47) can be also done
by rechoosing a sufficiently large $c$ if necessary and noting that the right-hand side on (4.48) also
has a bound $O((c\lambda)^{-1})$.

\noindent
$\Box$

Note that the constants $C_5$ and $C_6$ are independent of unit vectors $\omega$ and $\omega^{\perp}$
in (4.1).

Now let
$$\begin{array}{ll}
\displaystyle
a=\frac{1}{\gamma_0}-1, & b=V.
\end{array}
$$
Then the family of $v_{\tau}$ given by (4.45) satisfies
(4.37) and (4.39) with $\sigma=1$ and
$$\left\{
\displaystyle
\begin{array}{l}
C_2=c\,h_D(\omega),
\\
\\
\displaystyle
C_3=c\,h_{\Omega}(\omega),
\end{array}
\right.
$$
where $h_D(\omega)=\sup_{x\in D}\,x\cdot\omega$ and $h_{\Omega}(\omega)=\sup_{x\in\Omega}\,x\cdot\omega$.
Besides, Choosing a small $\eta$ in Theorem 4.6 and using a standard technique previously used in the enclosure method,
we see that the lower estimate (4.38) with $\sigma=1$ is also satisfied with $C_2$ given above.

Therefore we obtain the following corollary.

\proclaim{\noindent Corollary 4.1(\cite{IFR}).}
Assume that $\gamma_0-1\in
C_0^{\infty}(\Bbb R^n)$ and $\gamma_0(x)>0$ for all $x\in\Bbb R^3$.  Choose $\eta$ in Theorem 4.6 small enough. 
Fix the virtual slowness $c$ in (4.1) as $c=C_5$,
where $C_5$ is just the same as Theorem 4.6. Let $f=f_{\tau}$ be the
function of $(x,t)\in\partial\Omega\times]0,\,T[$, also
depending on a parameter $\tau\ge C_6$, defined by the equation (4.44),
where $v_{\tau}$ is given by (4.45) and a real-valued function
$\varphi\in L^2(0,\,T)$
satisfying the condition 
$$\displaystyle
\liminf_{\tau\rightarrow\infty}\tau^{\mu}\left\vert\int_0^Te^{-\tau t}\varphi(t)\,dt\right\vert>0
$$
for a $\mu\in\Bbb R$.

If $T$ satisfies
$$\displaystyle
T>2c(h_{\Omega}(\omega)-h_D(\omega)),
\tag {4.50}
$$
then we have
$$\displaystyle
\lim_{\tau\longrightarrow\infty}
\frac{1}{2\tau}
\log\left\vert
\int_{\partial\Omega}
\int_0^Te^{-\tau t}\left(-\overline{v_{\tau}(x)}f_{\tau}(x,t)+u_{f_{\tau}}(x,t)\gamma_0\frac{\partial\overline{v_{\tau}}}{\partial\nu}(x)\right)dtdS\right\vert
=ch_D(\omega).
$$

\endproclaim

Some remarks are in order.

$\bullet$  Theorem 4.5 is valid also for general inhomogeneous anisotropic background conductivity instead of $\gamma_0\,I_n$.

Thus to obtain a result like Corollary 4.1, we have to construct the complex geometrical optics solution
for general inhomogeneous anisotropic conductivity which plays the role of $v_{\tau}$ in Corollary 4.1.
In two dimensions, it may be possible to do it
by using the idea of isothermal coordinates as used in \cite{S} for the Calder\'on problem together with that of
choosing large virtual slowness $c$ in (4.1) in a reduced equation.
However, in three dimensions, we have no idea.

\subsubsection{Using real exponential solutions}

Theorem 4.5 with $\sigma=\frac{1}{2}$ covers some of previous results \cite{IK1, IK2} in the case when $\gamma_0$ is a constant, say $\gamma_0=1$.

Let $n=3$.  The following is a list of explicit examples of the family $(v_{\tau})$ of Theorem 4.5 with $\sigma=\frac{1}{2}$:

$$\left\{
\begin{array}{lll}
\displaystyle
v_{\tau}(x;\omega)=e^{\sqrt{\tau}\,x\cdot\omega}, & x\in\Bbb R^3, & \omega\in S^2,\\
\\
\displaystyle
v_{\tau}(x;p)=\frac{e^{-\sqrt{\tau}\,\vert x-p\vert}}{\vert x-p\vert}, & x\in\Bbb R^3\setminus\{p\}, & p\in\Bbb R^3\setminus\overline\Omega
\end{array}
\right.
$$
and, for an arbitrary point $y\in\Bbb R^3$
$$\displaystyle
v_{\tau}(x;y)=
\left\{
\begin{array}{ll}
\displaystyle
\frac{e^{\sqrt{\tau}\,\vert x-y\vert}-e^{-\sqrt{\tau}\,\vert x-y\vert}}{\vert x-y\vert},
& x\in\Bbb R^3\setminus\{y\},\\
\\
\displaystyle
2\tau, & x=y.
\end{array}
\right.
$$
Then, the constants $C_1=c(\Omega)$ and $C_2=c(D)$ of Theorem 4.5 are given by
$$\displaystyle
c(A)=\left\{
\begin{array}{ll}
\displaystyle
h_{A}(\omega), & \text{if $v_{\tau}=v_{\tau}(x;\omega)$,}
\\
\\
\displaystyle
-d_{A}(p), & \text{if $v_{\tau}=v_{\tau}(x;p)$,}
\\
\\
\displaystyle
R_{A}(y), & \text{if $v_{\tau}=v_{\tau}(x;y)$},
\end{array}
\right.
$$
where for $A=\Omega, D$ $d_A(p)$ and $R_A(y)$ denote the distance of $p$ to set $A$ and minimum radius of the open ball
that contains $A$ and centered at $y$.

Therefore formula (4.43) of Theorem 4.5 with $\sigma=\frac{1}{2}$ under the assumption that $\gamma_0=1$ reproduces Theorems 1.2 and 1.4 in \cite{IK2}.

In the case when $\gamma_0$ is not necessary a constant, it will be possible to construct
of a {\it geometrical optics solution} of equation (4.36) having the form
$$\displaystyle
v_{\tau}\sim e^{\sqrt{\tau}\,\varphi}\left(a_0+\frac{a_1}{\sqrt{\tau}}+\cdots\right)
$$
as $\tau\rightarrow\infty$ provided the eikonal equation
$$\begin{array}{ll}
\displaystyle
\gamma_0\nabla\varphi\cdot\nabla\varphi=1, & x\in\Omega,
\end{array}
$$
is solvable.  See Theorem 3.1 in \cite{ISource} and its application to the inverse source problem for the heat equation
explained in the previous subsection 4.1.

However, the solvability of the eikonal equation is not a simple matter unlike construction of the complex geometrical optics solutions.
See Section 3 in \cite{U} for an assumption for $\gamma_0$ on the solvability of the eikonal equation globally in a neighbourhood of 
$\overline\Omega$.

\subsection{Inverse obstacle problems using the wave propagation}

This section is concerned with the time domain enclosure method
to some of typical inverse obstacle problems which employ the dynamical scattering data of acoustic or electro magnetic waves
over a finite time interval

\subsubsection{Extracting from bistatic data}

In \cite{IW00} by considering two inverse obstacle scattering problems in an exterior domain or the whole space,
the time domain enclosure method for the wave equations has been introduced.  It is a three-dimensional realization 
in wave phenomena of the idea introduced in \cite{I4} (see also \cite{IBrazil}).
The idea has been developed in the succeeding series \cite{IEO2, IEIII}.

Here we present some of the results in \cite{IEIII}.
It is assumed that the unknown obstacle is sound-soft one.
The governing equation of the wave is given by the classical wave equation.
The wave is generated by the initial data localized outside the obstacle
and observed over a finite time interval at a different place from the support of the initial data,
that is, the observed data are the so-called bistatic data.
This is a simple mathematical model of the data collection process using an acoustic wave/electromagnetic wave
such as, bistatic active sonar, radar, etc.

Let $D$ be a nonempty bounded open subset of $\Bbb R^3$ with $C^2$-boundary
such that $\Bbb R^3\setminus\overline D$ is connected.
Let $0<T<\infty$.
Let $f\in\,L^2(\Bbb R^3)$ satisfy $\text{supp}\,f\cap\overline D=\emptyset$.
Let $u=u_f(x,t)$ denote the weak solution (\cite{DL})
of the following initial boundary value problem for the classical wave equation:
$$\left\{
\begin{array}{ll}
\displaystyle
u_{tt}-\Delta u=0, & (x,t)\in (\Bbb R^3\setminus\overline D)\times\,]0,\,T[,
\\
\\
\displaystyle
u(x,0)=0, & x\in\Bbb R^3\setminus\overline D,
\\
\\
\displaystyle
u_t(x,0)=f(x) & x\in\Bbb R^3\setminus\overline D,\\
\\
\displaystyle
u=0, & (x,t)\in\partial D\times\,]0,\,T[.
\end{array}
\right.
\tag {4.51}
$$
The boundary condition for $u$ in (4.51) means that $D$ is a sound-soft obstacle.

The problem to be considered is as follows.

$\quad$

{\bf\noindent Problem 4.5.} 
Let $B$ and $B'$ be two {\it known} open balls centered at $p\in\Bbb R^3$ and $p'\in\Bbb R^3$ with radius $\eta$ 
and $\eta'$, respectively such that $(\overline B\cup\overline{B'})\cap\overline D=\emptyset$.
Let $\chi_B$ denote the characteristic function of $B$ and set $f=\chi_B$.
Assume that $D$ is unknown.
Extract information about the geometry of $D$ from the data $u_f(x,t)$ given at all
$x\in B'$ and $t\in\,]0,\,T[$.

$\quad$

We consider only the pair $B$ and $B'$ satisfying $\overline{B'}\cap\overline B=\emptyset$.
Then the data in Problem 4.5 are called the bistatic data.

Why use the wave generated by the incident wave propagating to the whole directions (non directive incident wave)?
A typical example of the directive incident wave is Gaussian jet, e.g., \cite{KZ} and references therein.
The primary part of the reflected wave produced by a directive wave 
depends on the incident direction relative to the normal at a hitting point on the surface
of unknown obstacle (Snell's law).  
Thus the position of the receiver/observation point is important!
Gaussian jet contains a large parameter $p>0$ (carrying frequency) to asymptotically
localize as $p\longrightarrow\infty$ in the space-time.
In contrast to this our incident wave is free from any parameter.  Only in the data analysis we introduce a large
parameter.

{\bf\noindent Definition 4.4.}
Define the function of $\tau>0$ by the formula
$$\displaystyle
I(\tau;B,B')=\int_{\Bbb R^3\setminus\overline D}(fv_g-w_fg)dx,
$$
where $f=\chi_B$, $g=\chi_{B'}$,
$$\displaystyle
v_g(x)=v_g(x,\tau)=\frac{1}{4\pi}
\int_{\Bbb R^3}\frac{e^{-\tau\vert x-y\vert}}{\vert x-y\vert}g(y)dy
$$
and
$$\displaystyle
w_f(x)=w_f(x,\tau)
=\int_0^Te^{-\tau t}u_f(x,t)dt,\,\,x\in\Bbb R^3\setminus\overline D.
$$
In what follows we simply write $I(\tau;B,B')=I(\tau)$.

The indicator function can be computed from $u_f$ on $B'\times\,]0,\,T[$
since we have
$$\displaystyle
I(\tau)=\int_Bv_gdx-\int_{B'}w_fdx.
$$

Define
$$\displaystyle
\phi(x;y,y')
=\vert y-x\vert+\vert x-y'\vert,\,\,(x,y,y')\in\Bbb R^3\times\Bbb R^3\times\Bbb R^3.
$$
This is the length of the broken path connecting $y$ to $x$ and $x$ to $y'$.

We denote the {\it convex hull} of the set $F\subset\Bbb R^3$ by $[F]$.

\proclaim{\noindent Theorem 4.7 (\cite{IEIII}).}
Let $[\overline B\cup\overline B']\cap\partial D=\emptyset$
and $T$ satisfy
$$\displaystyle
T>\min_{x\in\partial D,\,y\in\partial B,\, y'\in\partial B'}\phi(x;y,y').
\tag {4.52}
$$
Then, there exists a $\tau_0>0$ such that, for all $\tau\ge\tau_0$,
$I(\tau)>0$
and the formula
$$\begin{array}{c}
\displaystyle
\lim_{\tau\longrightarrow\infty}
\frac{1}{\tau}
\log I(\tau)
=-\min_{x\in\partial D,\,y\in\partial B,\, y'\in\partial B'}\phi(x;y,y')
\end{array}
\tag {4.53}
$$
is valid.
\endproclaim

{\it\noindent Sketch of Proof.}
It suffices to prove the following two estimates:
there exist $\mu_j\in\Bbb R$, $C_j>0$ with $j=1,2$ and $\tau_0>0$ which are independent of $\tau$
such that, for all $\tau\ge\tau_0$,
$$\displaystyle
e^{\tau\min_{x\in\partial D,\,y\in\partial B,\,y'\in\partial B'}\phi(x;y,y')}\int_{\Bbb R^3\setminus\overline D}(fv_g-w_fg)dx
\le C_1\tau^{\mu_1}
\tag {4.54}
$$
and
$$\displaystyle
C_2\tau^{\mu_2}
\le e^{\tau\min_{x\in\partial D,\,y\in\partial B,\,y'\in\partial B'}\phi(x;y,y')}\int_{\Bbb R^3\setminus\overline D}(fv_g-w_fg)dx.
\tag {4.55}
$$
The proof of (4.54) is not so difficult.  So we describe only the proof of (4.55) which is the central part.

Given $q\in\partial D$ satisfying $\phi(q;p,p')=\min_{x\in\partial D}\phi(x;p,p')$,
there exists an open ball $\tilde{D}$ contained in $D$ such that
$\min_{x\in\partial\tilde{D}}\phi(x;p,p')=\phi(q;p,p')$.

Let $\tilde{u}=\tilde{u}_f$ denote the weak solution of the following initial boundary value problem:
$$\left\{
\begin{array}{ll}
\displaystyle
\tilde{u}_{tt}-\Delta\tilde{u}=0, & (x,t)\in(\Bbb R^3\setminus\overline{\tilde{D}})\times\,]0,\,T[,\\
\\
\displaystyle
\tilde{u}(x,0)=0, & x\in\Bbb R^3\setminus\overline{\tilde{D}},\\
\\
\displaystyle
\tilde{u}_t(x,0)=f(x), & x\in\Bbb R^3\setminus\overline{\tilde{D}},\\
\\
\displaystyle
\tilde{u}=0, & (x,t)\in\partial\tilde{D}\times\,]0,\,T[.
\end{array}
\right.
$$

Define
$$\begin{array}{ll}
\displaystyle
\tilde{w}_f(x,\tau)
=\int_0^Te^{-\tau T}\tilde{u}(x,t)dt,\,\,x\in\Bbb R^3\setminus\overline{\tilde{D}} & \tau>0.
\end{array}
$$
We have
$$\displaystyle
\int_{\Bbb R^3\setminus\overline D}(fv_g-w_fg)dx
\ge\int_{\Bbb R^3\setminus\overline{\tilde{D}}}(fv_g-\tilde{w}_fg)dx+O(\tau^{-1}e^{-\tau T}).
\tag {4.56}
$$
This is because of $(\text{supp}\,f\cup\text{supp}\,g)\cap\overline D=\emptyset$
and the fact that, for some function $Z=Z_{\tau}\in L^2(\Bbb R^3\setminus\overline D)$ satisfying 
$\Vert Z\Vert_{L^2(\Bbb R^3\setminus\overline D)}=O(\tau^{-1}e^{-\tau T})$
we have
$$\begin{array}{ll}
\displaystyle
\tilde{w_f}-w_f\ge Z, &
x\in\Bbb R^3\setminus\overline D.
\end{array}
$$
The last inequality is a consequence of the maximum principle for the modified Helmholtz equation.

From (4.56) we see that now the original obstacle $D$ is replaced with $\tilde{D}$ which is convex.
Then a combination of an argument in \cite{LP} which makes use of the maximum principle and the convexity of $\tilde{D}$ yields
$$\displaystyle
\int_{\Bbb R^3\setminus\overline{\tilde{D}}}(fv_g-\tilde{w}_fg)dx
\ge \int_{\tilde{D}}(\nabla v_f\cdot\nabla v_g+\tau^2 v_fv_g)\,dx+O(\tau^{-1}e^{-\tau T}).
$$
Thus everything is reduced to give a lower estimate for the integral on this right-hand side as $\tau\rightarrow\infty$.

An explicit computation yields
$$\begin{array}{l}
\,\,\,\,\,\,
\displaystyle
\tau^2e^{-\tau(\eta+\eta')}e^{\tau\min_{x\in\partial D}\,\phi(x;p,p')}
\int_{\tilde{D}}(\nabla v_f\cdot\nabla v_g+\tau^2 v_fv_g)\,dx
\\
\\
\displaystyle
\ge C_{\tilde{D}}(p,p')e^{\tau\min_{x\in\partial D}\,\phi(x;p,p')}\int_{\tilde{D}}e^{-\tau\phi(x;p,p')}\,dx+O(\tau^{-\infty}),
\end{array}
$$
where
$$\displaystyle
C_{\tilde{D}}(p,p')\sim
\inf_{x\in\tilde{D}}\,\left\{1+\frac{(p-x)\cdot(p'-x)}{\vert p-x\vert\,\vert p-x'\vert}\right\}.
$$
Note that this right-hand side is positive by virtue of $[\{p,p'\}]\cap\overline{\tilde{D}}=\emptyset$
and that
$$\displaystyle
\min_{x\in\partial D,\, y\in \partial B,\,y'\in \partial B'}\phi(x;y,y')
=\min_{x\in\partial D}\phi(x;p,p')-(\eta+\eta').
$$
Finally the lower estimate 
$$\displaystyle
\liminf_{\tau\rightarrow\infty}\tau^3e^{\tau\min_{x\in\partial D}\,\phi(x;p,p')}\int_{\tilde{D}}e^{-\tau\phi(x;p,p')}\,dx>0
$$
yields (4.55).

\noindent
$\Box$

The quantity $\min_{x\in\partial D}\phi(x;p,p')$ coincides with the shortest length
of the broken paths connecting $p$ to a point $q$ on $\partial D$ and $q$ to $p'$,
that is, the {\it first arrival time} of the signal with speed $1$ in the
geometrical optics sense.

{\bf\noindent Definition 4.5.} Let $p$ and $p'$ satisfy
$[\{p,p'\}]\cap\partial D=\emptyset$. Define
$$
\displaystyle
\Lambda_{\partial D}(p,p')
=\{q\in\partial D\,\vert\,
\phi(q;p,p')=\min_{x\in\partial D}\phi(x;p,p')\}.
$$
We call this the {\it first reflector} between $p$ and $p'$.
The points in the first reflector are called the {\it first reflection points} between $p$ and $p'$.

Given $c>\vert p-p'\vert$ define
$$\displaystyle
E_c(p,p')=\{x\in\Bbb R^3\,\vert\,\phi(x;p,p')=c\}.
$$
This is a {\it spheroid} with focal points $p$ and $p'$.
We have $\Lambda_{\partial D}(p,p')=\partial D\cap E_c(p,p')$ with 
$c=\min_{x\in\partial D}\phi(x;p,p')$ and
the two tangent planes at $q\in\Lambda_{\partial D}(p,p')$ of $\partial D$ and $E_c(p,p')$ coincide.
We denote by $S_q(\partial D)$ and $S_q(E_c(p,p'))$ the shape operators (or the Weingarten maps)
at $q$ with respect to $\nu_q$ which is the outward normal to $\partial D$ and inward normal to $E_c(p,p')$.
Those are symmetric linear operators on the common tangent space at $q$ of $\partial D$ and $E_c(p,p')$.
It is easy to see that $S_q(E_c(p,p'))-S_q(\partial D)\ge 0$ as the quadratic form
on the same tangent space at $q$.

The following theorem describes the leading term of the asymptotic expansion
of $I(\tau)$ as $\tau\longrightarrow\infty$.

\proclaim{\noindent Theorem 4.8.}
Let $B$ and $B'$ satisfy $[\overline B\cup\overline{B'}]\cap\overline D=\emptyset$.
Let $T$ satisfy (4.52).  Assume that $\Lambda_{\partial D}(p,p')$ is finite
and for all $q\in\Lambda_{\partial D}(p,p')$,
$$\displaystyle
\text{det}\,(S_q(E_c(p,p'))-S_q(\partial D))>0.
$$
If $D$ is convex and $\partial D$ is $C^3$, then we have
$$
\displaystyle
\lim_{\tau\longrightarrow\infty}
\tau^4e^{\tau\min_{x\in\partial D,\,y\in\partial B,\, y'\in\partial B'}\phi(x;y,y')}
I(\tau)
=\frac{\pi}{2}\sum_{q\in\Lambda_{\partial D}(p,p')}
\frac{R(q;B,B')}
{\sqrt{\text{det}\,(S_q(E_c(p,p'))-S_q(\partial D))}},
$$
where
$$
\displaystyle
R(q;B,B')
=
\left(\frac{\text{diam}\,B}{2\vert q-p\vert}\right)\cdot
\left(\frac{\text{diam}\,B'}{2\vert q-p'\vert}\right).
$$

\endproclaim

We omit to describe the sketch of the proof.
The central part of the proof uses the reflection argument made in \cite{LP}, and the convexity assumption of $D$ is essential.

Using Theorems 4.7 and 4.8, we can obtain the information listed below:

$\bullet$  $E_c(p,p')$ with $c=\min_{x\in\partial D}\phi(x;p,p')$, that is the largest spheroid with focal points $p$ and $p'$
whose exterior encloses $D$.

$\bullet$ all the first reflection points between $p$ and $p'$ together with the normal of $\partial D$ at the points.

$\bullet$ the Gauss curvature and a modification of the mean curvature of $\partial D$
at an arbitrary first reflection point between $p$ and $p'$ provided $D$ is convex.

Finally we note that as a corollary of Theorems 4.7 and 4.8 we have a reconstruction formula of a ball $D$ from the bistatic data
$u_f(x,t)$, $x\in B'$, $0<t<T$  with $f=\chi_B$ and a large fixed $T$ under the assumption 
$(\overline {B}\cup\overline{B'})\cap\overline D=\emptyset$.  See Corollary 1.2 and Theorem 1.5 in \cite{IEIII}.

Further explanation of the results in \cite{IEIII} are given in \cite{IRIMS2013}, however, it is written in Japanese.

\subsubsection{Extracting from monostatic data}

It is interesting to extract information about not only the geometry of
an unknown obstacle but also the qualitative or quantitative state of the obstacle surface from dynamical scattering data.
Here we consider an inverse obstacle-scattering problem which is described by the classical wave equation outside an obstacle
with a dissipative boundary condition.

First, we describe the problem.
Let $D$ be a nonempty bounded open subset of $\Bbb R^3$ with $C^2$-boundary such that
$\Bbb R^3\setminus\overline D$
is connected. Let $\gamma$ be a function belonging to $L^{\infty}(\partial D)$ and satisfy $\gamma\ge 0$.
Let $0<T<\infty$.  Let $B$ be an open ball satisfying $\overline B\cap \overline D=\emptyset$.
We denote by $\chi_B$ the characteristic function of $B$; $p$ and $\eta$ the center and (very small) radius of $B$, respectively.

Let $u=u_B(x,t)$ denote the weak solution of the following initial boundary value problem for the classical wave equation:
$$\displaystyle
\left\{
\begin{array}{ll}
\displaystyle
u_{tt}-\Delta u=0, & (x,t)\in\,(\Bbb R^3\setminus\overline D)\times\,]0,\,T[,\\
\\
\displaystyle
u(x,0)=0, & x\in\Bbb R^3\setminus\overline D,\\
\\
\displaystyle
\partial_tu(x,0)=\chi_B(x), & x\in\Bbb R^3\setminus\overline D,\\
\\
\displaystyle
\frac{\partial u}{\partial\nu}-\gamma(x)\partial_tu=0, & (x,t)\in\partial D\times\,]0,\,T[.
\end{array}
\right.
\tag {4.57}
$$
Here, $\nu$ denotes the unit normal to $\partial D$, oriented towards the exterior of $D$.

The solution of (4.57) is a model of the wave that loses the energy on the surface of the obstacle
since  formal computation yields
$$\displaystyle
{\cal E}'(t)=-\int_{\partial D}\gamma(x)\vert\partial_t u\vert^2dS\le 0,
$$
where
$$\displaystyle
{\cal E}(t)=\frac{1}{2}\int_{\Bbb R^3\setminus\overline D}(\vert\partial_t u\vert^2+\vert\nabla u\vert^2)dx,\,t\in [0,\,T].
$$
The distribution of $\gamma$ represents the state of the surface of the obstacle.

In this section we consider the following problem.

$\quad$

{\bf\noindent Problem 4.6.}  Fix a large $T$ (to be determined later).
Assume that both $D$ and $\gamma$
are unknown.  Extract information about the
location and shape of $D$ together with the values of $\gamma$ from the wave field $u_B(x,t)$ given at
all $x\in B$ and $t\in\,]0,\,T[$.

$\quad$

A solution to Problem 4.6 by using the time domain enclosure method has been obtained in \cite{IEO2}.

Define
$$
\begin{array}{ll}
\displaystyle
I_B(\tau)=\int_B(w-v)dx, & \tau>0,
\end{array}
$$
where
$$\begin{array}{lll}
\displaystyle
w(x)=w_B(x,\tau)=\int_0^Te^{-\tau t}u_B(x,t)dt,
&
x\in\Bbb R^3\setminus\overline D,
&
\tau>0.
\end{array}
\tag {4.58}
$$
and  $v=v_B(\,\cdot\,,\tau)\in H^1(\Bbb R^3)$ denote the weak solution of the modified Helmholtz equation
$$
\begin{array}{ll}
\displaystyle
(\Delta-\tau^2)v+\chi_B=0,
&
x\in\Bbb R^3.
\end{array}
\tag {4.59}
$$
The function $\tau\longmapsto I_B(\tau)$ is the indicator function in the enclosure method for Problem 4.6.

We write $\gamma>>1$ if $\exists C'>0\,\,\gamma(x)-1\ge C'$ a.e. $x\in\partial D$;
$\gamma<<1$ if $\exists C'>0\,\,-(\gamma(x)-1)\ge C'$ a.e. $x\in\partial D$.

\proclaim{\noindent Theorem 4.9(\cite{IEO2}).}
Let $T$ satisfies
$$\displaystyle
T>2\text{dist}\,(D,B).
\tag {4.60}
$$
We have:

if $\gamma<<1$ is satisfied, then there exists a $\tau_0>0$ such that, for all $\tau\ge\tau_0$
$$\displaystyle
I_B(\tau)>0;
\tag {4.61}
$$

if $\gamma>>1$ is satisfied, then there exists a $\tau_0>0$ such that, for all $\tau\ge\tau_0$
$$\displaystyle
I_B(\tau)<0.
\tag {4.62}
$$

In both cases, the formula
$$
\lim_{\tau\longrightarrow\infty}\frac{1}{2\tau}\log\left\vert I_B(\tau)\right\vert
=-\text{dist}\,(D,B)
\tag {4.63}
$$
is valid.

\endproclaim

Some remarks are in order.

$\bullet$  Roughly speaking, from (4.61) and (4.62) one can distinguish $\gamma>1$ or $\gamma<1$.  This shows also $\gamma\equiv 1$ is very special.
In fact, if $\gamma\equiv 1$, the space dimension is one and  $D=]a,\,\infty[$ with $a\in\Bbb R$,
then one can not determine $a$ itself from the corresponding data.  See p. 27 in \cite{IEO2}.

$\bullet$  Since $\text{dist}\,(D,B)=d_{\partial D}(p)-\eta$ with $d_{\partial D}(p)=\inf_{y\in\partial D}\vert y-p\vert$,
formula (4.63) gives us
the largest sphere $B_{d_{\partial D}(p)}(p)$ whose exterior encloses $D$.

$\bullet$  After having $d_{\partial D}(p)$, using formula (4.63) for {\it infinitely many} $B$,
one can detect all the points on the first reflector
$$\displaystyle
\Lambda_{\partial D}(p)=\{q\in\partial D\,\vert\,\vert q-p\vert
=d_{\partial D}(p)\},
$$
whose members are on the sphere $B_{d_{\partial D}(p)}(p)$.
The procedure is as follows.
Given $\omega\in S^2$ take a small open ball $B'$ centered at
$p+s\omega$ with a small $s\in]0,\,d_{\partial D}(p)[$, produce the data
$u_{B'}\equiv u_B\vert_{B\rightarrow B'}$ on $B'$ over a large but finite time interval.
Then one gets $d_{\partial D}(p+s\omega)$ via formula (4.63) with $I_{B'}\equiv I_B\vert_{B=B'}$.
By checking $d_{\partial D}(p+s\omega)-(d_{\partial D}(p)-s)>0$ or not one can make a decision whether the point $p+d_{\partial D}(p)\omega$
belongs to $\partial D$, that is, to $\Lambda_{\partial D}(p)$ or not.
It is ensured by the following simple fact.

\proclaim{\noindent Proposition 4.3(\cite{ICA}).}
If the point $p+d_{\partial D}(p)\,\omega$ belongs to $\partial D$, then $d_{\partial D}(p+s\omega)=d_{\partial D}(p)-s$.
If $p+d_{\partial D}(p)\,\omega$ does not, then $d_{\partial D}(p+s\omega)>d_{\partial D}(p)-s$.
\endproclaim

So the next problem is to extract the quantitative state of obstacle surface.

To describe the formula, we introduce some notion
in differential geometry same as the previous section.
Let $q\in\Lambda_{\partial D}(p)$.
Let $S_q(\partial D)$ and $S_q(\partial
B_{d_{\partial D}(p)}(p))$ denote the {\it shape operators} (or {\it Weingarten maps}) at $q$
of $\partial D$ and $\partial B_{d_{\partial D}(p)}(p)$ with
respect to $\nu_q$ which is the outward normal to $\partial D$ and inward normal to $\partial B_{d_{\partial D}}(p)$.
Because $q$ attains the minimum of the function: $\partial D\ni y\longmapsto
\vert y-p\vert$, we have always $S_q(\partial B_{d_{\partial D}(p)}(p))-S_{q}(\partial D)\ge 0$
as the quadratic form on the common tangent space at $q$.

Note also that we have
$$\displaystyle
\text{det}\,(S_q(\partial B_{d_{\partial D}(p)}(p))-S_{q}(\partial D))
=\lambda^2-2H_{\partial D}(q)\lambda+K_{\partial D}(q),
\tag {4.64}
$$
where $\lambda=1/d_{\partial D}(p)$, $H_{\partial D}(q)$ and $K_{\partial D}(q)$ are
the mean and Gauss curvatures at $q$ with respect to $\nu_q$.
Thus if the curvatures and $d_{\partial D}(p)$ are known, then the left-hand side on (4.64) is known.

Here we write $0<<\gamma<<1$ if $\exists C'>0\,\,\exists C''>0\,\,\,C''\le \gamma(x)\le 1-C'$ a.e. $x\in\partial D$.
This is stronger than $\gamma<<1$.

The second formula is the following.

\proclaim{\noindent Theorem 4.10(\cite{IMP}).}
Assume that $\partial D$ is $C^3$ and $\gamma\in C^2(\partial D)$.
Assume that $\gamma$ satisfies $0<<\gamma<<1$ or $\gamma>>1$.
Let $T$ satisfy (4.60).
Assume that the set $\Lambda_{\partial D}(p)$ consists of finite points
and
$$\displaystyle
\text{det}\,(S_q(\partial B_{d_{\partial D}(p)}(p))-S_{q}(\partial D))>0\,\,\,
\forall q\in\Lambda_{\partial D}(p).
\tag {4.65}
$$
Then, we have
$$
\displaystyle
\,\,\,\,\,\,
e^{2\tau\text{dist}\,(D,B)}I_B(\tau)
=\frac{1}{\tau^4}
\cdot
\frac{\pi}{2}
\left(\frac{\eta}{d_{\partial D}(p)}\right)^2
\sum_{q\in\Lambda_{\partial D}(p)}
k_q(p)\frac{1-\gamma(q)}{1+\gamma(q)}
+o(\tau^{-4}),
\tag {4.66}
$$
where
$$\displaystyle
k_q(p)=\frac{1}{\displaystyle
\sqrt{\text{det}\,(S_q(\partial B_{d_{\partial D}(p)}(p))-S_{q}(\partial D))}}.
$$

\endproclaim

Note that if $D$ is {\it convex}, then $\Lambda_{\partial D}(p)$
consists of a single point and (4.65) is satisfied.

Here we describe two applications of Theorem 4.10 in remote sensing.
Fix a point $p\in\Bbb R^3\setminus\overline D$ and we have a known point $q\in\Lambda_{\partial D}(p)$.

$\bullet$ Finding the value of $\gamma$ at $q$ provided the Gauss curvature $K_{\partial D}(q)$ 
and mean curvature $H_{\partial D}(q)$ of $\partial D$
with respect to $\nu_q$ are known.

The procedure is as follows.

{\bf\noindent Step 1.}  From point $p$, let us go a little bit forward to $q$.
We denote by $p'$ the point. We see that the set $\Lambda_{\partial D}(p')$ 
consists of only the single point $q$
and (4.65) for $p$ replaced with $p'$ is satisfied.

{\bf\noindent Step 2.}  Choose a small open ball $B'$ centered at $p'$ and generate the wave $u_{B'}$.

{\bf\noindent Step 3.}  Observe the data $u_{B'}$ on $B'$ over time interval $]0,\,T'[$ for a large  $T'$.

{\bf\noindent Step 4.}  Compute the indicator function $I_{B'}(\tau)$ and from (4.66)
find the value
$$\displaystyle
{\cal F}_{B'}(q)\equiv\frac{\pi}{2}
\left(\frac{\eta'}{d_{\partial D}(p')}\right)^2k_q(p')
\,\frac{1-\gamma(q)}{1+\gamma(q)},
$$
where $\eta'$ denotes the radius of $B'$ and we have $d_{\partial D}(p')=d_{\partial D}(p)-\vert p-p'\vert$.

{\bf\noindent Step 5.}  By solving the equation in Step 4, find $\gamma(q)$.

Note that if $D$ is {\it convex}, then $\Lambda_{\partial D}(p)$
consists of a single point and (4.65) is satisfied.  Thus from formula (4.66) one can directly extract
the value $\gamma(q)$.

$\bullet$  Finding the curvatures and the value of $\gamma$ at $q$.

The procedure is as follows.

{\bf\noindent Step 1.}  Choose three points  $p_1$, $p_2$, and $p_3$ on the segment $p\rightarrow q$ different from its end points.

{\bf\noindent Step 2.}  Choose three small balls $B_1$, $B_2$, and $B_3$ with a common radius $\eta'$ centered at $p_1$, $p_2$, and $p_3$, respectively
and generate three waves $u_{B_j}$, $j=1,2,3$.

{\bf\noindent Step 3.} Apply formula (4.66) to $I_{B'}(\tau)$ with $B'=B_j$, $j=1,2,3$ and find three quantities
$$\begin{array}{ll}
\displaystyle
{\cal F}_j\equiv
\frac{2}{\pi}
\left(\frac{d_{\partial D}(p_j)}{\eta'}\right)^2
{\cal F}_{B'}(q)\vert_{B'=B_j}=\frac{A}{\sqrt{\lambda_j^2-2H\lambda_j+K}}, &
\displaystyle
j=1,2,3,
\end{array}
\tag {4.67}
$$
where $A=(1-\gamma(q))/(1+\gamma(q))$, $H=H_{\partial D}(q)$ and $K=K_{\partial D}(q)$;
$\lambda_j=1/d_{\partial D}(p_j)$, $j=1,2,3$.

{\bf\noindent Step 4.}  Find $H$ and $K$ by solving the linear system numerically
$$
\displaystyle
\left(
\begin{array}{cc}
\displaystyle
-(\lambda_1{\cal F}_1^2-\lambda_2{\cal F}_2^2)
&
\displaystyle
{\cal F}_1^2-{\cal F}_2^2\\
\\
\displaystyle
-(\lambda_2{\cal F}_2^2-\lambda_3{\cal F}_3^2)
&
\displaystyle
{\cal F}_2^2-{\cal F}_3^2
\end{array}
\right)
\left(\begin{array}{c}
\displaystyle
2 H
\\
\displaystyle
K
\end{array}
\right)
=
\left(\begin{array}{c}
\displaystyle
{\cal F}_2^2\lambda_2^2-{\cal F}_1^2\lambda_1^2
\\
\\
\displaystyle
{\cal F}_3^2\lambda_3^2-
{\cal F}_2^2\lambda_2^2
\end{array}
\right).
\tag {4.68}
$$

{\bf\noindent Step 5.}  Find the value $A^2$ by the formula
$$\displaystyle
A^2=\frac{1}{3}
\sum_{j=1}^3{\cal F}_j^2(\lambda_j^2-2H\lambda_j+K).
\tag {4.69}
$$

{\bf\noindent Step 6.}  Determine whether $\gamma(q)>1$ or $\gamma(q)<1$ from one of three
equations (4.67) and find the value $(1-\gamma(q))/(1+\gamma(q))$ and hence $\gamma(q)$ itself
by taking the square root of the both sides on (4.69).

There is another application of Theorem 4.10 which aims at monitoring the state of obstacle surface,
that is, a change of the value of $\gamma$ at a monitoring point.
To describe the application, we denote the indicator function $I_B(\tau)$ by $I_B(\tau;\gamma)$.

\proclaim{\noindent Corollary 4.2.}
Assume that $\partial D$ is $C^3$.
Let $\gamma_0$ and $\gamma_1$ belong to $C^2(\partial D)$ and satisfy $0<<\gamma<<1$ or $\gamma>>1$.
Let $T$ satisfy (4.60).
Assume that the set $\Lambda_{\partial D}(p)$ consists of finite points
and satisfies (4.65).
Then, the limit $\lim_{\tau\longrightarrow\infty} I_B(\tau;\gamma_1)/I_B(\tau;\gamma_0)$ exists,
and we have
$$\displaystyle
\min_{q\in\Lambda_{\partial D}(p)}
\frac{
\displaystyle\frac{1-\gamma_1(q)}{1+\gamma_1(q)}
}
{
\displaystyle\frac{1-\gamma_0(q)}{1+\gamma_0(q)}
}
\le
\lim_{\tau\longrightarrow\infty}
\frac{I_B(\tau;\gamma_1)}
{I_B(\tau;\gamma_0)}
\le
\max_{q\in\Lambda_{\partial D}(p)}
\frac{
\displaystyle\frac{1-\gamma_1(q)}{1+\gamma_1(q)}
}
{
\displaystyle\frac{1-\gamma_0(q)}{1+\gamma_0(q)}
}.
\tag {4.70}
$$
In particulr, if $\Lambda_{\partial D}(p)$ consists of a single point $q\in\partial D$, we have
$$\displaystyle
\lim_{\tau\longrightarrow\infty}
\frac{I_B(\tau;\gamma_1)}
{I_B(\tau;\gamma_0)}
=\frac{
\displaystyle\frac{1-\gamma_1(q)}{1+\gamma_1(q)}
}
{
\displaystyle\frac{1-\gamma_0(q)}{1+\gamma_0(q)}
}.
\tag {4.71}
$$

\endproclaim

Estimates on (4.70) give us some global information about the values of $\gamma_1$ relative to $\gamma_0$ 
at all the points on $\Lambda_{\partial D}(p)$ without knowing the curvatures.

Formula (4.71) tells us whether there is a deviation of the value of $\gamma_1$ from $\gamma_0$
at a fixed monitoring point on the obstacle surface or not.  
Note that given an arbitrary $q\in\partial D$ if $p=q+s\nu_q$ and $s$ is a 
sufficiently small positive number, then $\Lambda_{\partial D}(p)=\{q\}$, and (4.65) is satisfied.
We do not need any information about the geometry of $\partial D$ around $q$.  This is remarkable.
However we need just two observed waves which were
generated from the same initial data for case $\gamma=\gamma_0$ on one day and 
case $\gamma=\gamma_1$ on another day.

Some further remarks are in order.

$\bullet$  There is a well-known and classical result due to Majda \cite{Mo} in
the context of the Lax-Phillips scattering theory for the wave equation with the same boundary condition.  
Under the assumption that the obstacle is strictly convex, 
he clarified that the leading term of the  high frequency asymptotics for the {\it scattering amplitude}
which is observed at infinity, contains information about only the Gauss curvature.
This is a big difference from Theorem 4.10 and one can easily imagine its reason by taking the limit $d_{\partial D}(p)\rightarrow\infty$
in the formula (4.64):
$$\displaystyle
\text{det}\,(S_q(\partial B_{d_{\partial D}(p)}(p))-S_{q}(\partial D))
\rightarrow K_{\partial D}(q).
$$
From a technical point of view, it should be emphasized that the approach in \cite{Mo} is a geometrical optics based method which employs a parametrix
in time domain by Lax, Majda and Taylor and calculus of Fourier integral and pseudo differential operators.
Thus his approach is completely different from our integration by parts based one.

$\bullet$  In \cite{IEE} the case when the boundary condition in (4.57)
is given by the Robin boundary condition
$$\begin{array}{ll}
\displaystyle
\frac{\partial u}{\partial\nu}-\beta(x)u=0, & (x,t)\in\partial D\times\,]0,\,T[,
\end{array}
$$
where $\beta$ is a real-valued function on $\partial D$ and belongs to $C^2(\partial D)$, has been considered.
Comparing the result in \cite{IEE} with (4.66), we see that information about the values of $\gamma$ is visible
rather than those of $\beta$.

$\bullet$  Extending the results in this subsection to the bistatic data case is an intersting open problem.

$\bullet$  Recently the time domain enclosure method had an application to a problem on a superconducting radio-frequency cavity
           \cite{LS}.

\subsubsection{Inclusion in a rough background medium}

The time domain enclosure method enables us to extract information about the geometry
of the set of unknown inclusions embedded in a rough background medium.
Such situation appears as the problem of through-wall imaging \cite{B} by using electromagnetic waves.
When the wall is electromagnetically penetrable, a typical problem is as follows.

Generate the electromagnetic wave at the place where the observer is.  And observe the wave at the same place
over a finite time interval.  The observed wave should include information about the unknown object behind the several walls 
or embedded in a complicated medium.  Clearly, 
free-space assumptions no longer apply after the electromagnetic waves propagate through the walls.
How can one extract the information from the observed wave?

In \cite{Iwall} as a first step the author considered the case when the governing equation is given by a single equation for
a scalar wave in the time domain and gave a mathematical method based on the governing equation, that should be called
a PDE approach, as an application of the time domain enclosure method.

Let $\alpha\in L^{\infty}(\Bbb R^3)$ and satisfy $\text{ess}.\inf_{x\in\Bbb R^3}\alpha(x)>0$.
Let $0<T<\infty$.  Given $f\in L^2(\Bbb R^3)$, let $u=u_f(x,t)$ be the weak solution
of
$$\displaystyle
\left\{
\begin{array}{ll}
\displaystyle
\alpha(x)u_{tt}-\Delta u=0, & (x,t)\in\Bbb R^3\times\,]0,\,T[,
\\
\\
\displaystyle
u(x,0)=0, & x\in\Bbb R^3,
\\
\\
\displaystyle
u_t(x,0)=f(x), & x\in\Bbb R^3.
\end{array}
\right.
$$
We assume that $\alpha$ takes the form
$$\displaystyle
\alpha(x)=\left\{
\begin{array}{ll}
\displaystyle
\alpha_0(x), & \mbox{if $x\in \Bbb R^3\setminus D$},
\\
\\
\displaystyle
\alpha_0(x)+h(x), & \mbox{if $x\in D$},
\end{array}
\right.
$$
where $D\subset\Bbb R^3$ is a bounded open subset with Lipschitz boundary.
The function $\alpha_0$ belongs to $L^{\infty}(\Bbb R^3)$
and satisfies $m_0^2\le\alpha_0(x)\le M_0^2$ a.e.$x\in\Bbb R^3$ for positive constants $m_0$ and $M_0$;
$h$ belongs to $L^{\infty}(D)$.

We write $\alpha>>\alpha_0$ if $\exists C>0\,\,h(x)\ge C$ a.e. $x\in D$;
$\alpha<<\alpha_0$ if $\exists C>0\,\, -h(x)\ge C$ a.e. $x\in D$.

Note that $\alpha_0$ and $h$ are just essentially bounded on $\Bbb R^3$ and $D$, respectively.
In particular, $\alpha_0$ can be a model for various background media such as a multilayered medium with complicated 
interfaces or unions of various domains with different refractive indexes.

Let $B$ be an open ball satisfying $\overline B\cap\overline
D=\emptyset$.  We choose, for simplicity, $f=\chi_B$,
where $\chi_B$ is the characteristic function of ball $B$.

Generate $u_f$ by this initial data $f$. In this paper, we consider
the following inverse problems under the assumption that
$\alpha_0$ is known and that both $D$ and $h$ are
unknown.

$\quad$

{\bf\noindent Problem 4.7.}
Find a criterion whether $D=\emptyset$ or not in terms of only $u_f$ on $B$ over
time interval $]0,\,T[$.

$\quad$

{\bf\noindent Problem 4.8.}
Assume that $D\not=\emptyset$.
Extract information about $D$  together with property of $h$ from $u_f$ on $B$ over time interval $]0,\,T[$ for a fixed large $T$.

$\quad$

The results to Problems 4.7 and 4.8 in \cite{Iwall} are as follows.

Let $\tau>0$ and define
$$\begin{array}{ll}
\displaystyle
w(x,\tau)=\int_0^T e^{-\tau t}u(x,t)dt, & x\in\Bbb R^3.
\end{array}
$$
Let $v=v(x,\tau)\in H^1(\Bbb R^3)$ be the weak solution of
$$\begin{array}{ll}
\displaystyle
\Delta v-\alpha_0\tau^2 v+\alpha_0 f=0, & x\in\Bbb R^3.
\end{array}
$$
Define
$$\displaystyle
I_f(\tau,T)
=\int_B\alpha_0f(w-v)dx.
$$

The following theorem gives us solutions to Problems 4.7 and 4.8.

\proclaim{\noindent Theorem 4.11(\cite{Iwall}).}
We have:
$$\displaystyle
\lim_{\tau\longrightarrow\infty}e^{\tau T}I_f(\tau,T)
=
\left\{
\begin{array}{ll}
\displaystyle
0 & \text{if $D=\emptyset$ and $T>0$,}
\\
\\
\displaystyle
-\infty & \text{if $D\not=\emptyset$, $\alpha>>\alpha_0$ and  $T>2M_0\text{dist}(D,B)$,}
\\
\\
\displaystyle
\infty & \text{if $D\not=\emptyset$, $\alpha<<\alpha_0$ and $T>2M_0\text{dist}(D,B)$.}
\end{array}
\right.
$$
Besides, in case of $D\not=\emptyset$, if $\alpha>>\alpha_0$ or $\alpha<<\alpha_0$, then we have, for all $T>2M_0\text{dist}(D,B)$
$$
\left\{
\begin{array}{l}
\displaystyle
-M_0\,\text{dist}\,(D,B)\le\liminf_{\tau\longrightarrow\infty}
\frac{1}{2\tau}
\log\left\vert I_f(\tau,T)\right\vert,
\\
\\
\displaystyle
\limsup_{\tau\longrightarrow\infty}
\frac{1}{2\tau}
\log\left\vert I_f(\tau,T)\right\vert
\le -m_0\,\text{dist}\,(D,B).
\end{array}
\right.
$$

\endproclaim

{\it\noindent Sketch of proof.}
First we have, as $\tau\rightarrow\infty$,
$$\left\{
\begin{array}{l}
\displaystyle
I_{f}(\tau,T)
\le
\tau^2\int_{\Bbb R^3}\,\frac{\alpha_0}{\alpha}\,(\alpha_0-\alpha)\,v^2\,dx+O(\tau^{-1}e^{-\tau T}),
\\
\\
\displaystyle
I_{f}(\tau,T)\ge\tau^2\int_{\Bbb R^3}\,(\alpha_0-\alpha)\,v^2\,dx+O(\tau^{-1}\,e^{-\tau T}).
\end{array}
\right.
$$
The proof of the upper estimate is more technical than that of the lower estimate.
Anyway everything is reduced to studying the leading profile of $v$, in the case when $D\not=\emptyset$.

For the purpose we have
$$
\left\{
\begin{array}{l}
\displaystyle
v(x)\ge\frac{1}{4\pi}\,\int_B\,\alpha_0(y)\,\frac{e^{-M_0\,\tau\vert x-y\vert}}{\vert x-y\vert}\,dy,
\\
\\
\displaystyle
v(x)\le\frac{1}{4\pi}\,\int_B\,\alpha_0(y)\,\frac{e^{-m_0\,\tau\vert x-y\vert}}{\vert x-y\vert}\,dy.
\end{array}
\right.
$$
Thus by virtue of assumptions $\alpha>>\alpha_0$ or $\alpha<<\alpha_0$,
we have reduced the problem to the asymptoic behaviour of the following integral as $\tau\rightarrow\infty$:
$$\begin{array}{ll}
\displaystyle
\int_B\frac{e^{-\gamma\,\tau\vert x-y\vert}}{\vert x-y\vert}\,dy &  x\in D,
\end{array}
$$
where $\gamma$ is a positive constant.

\noindent
$\Box$

Theorem 4.11 suggests a new direction of the enclosure method in the case when the background medium is inhomogeneous
and quite complicated.

Some remarks are in order.

$\bullet$  If we have the a-priori information $\text{dist}\,(D,B)<M$ with a known constant $M$,
then we can specify an explicit waiting time  for correcting the observation data in advance such as
$T\ge 2M_0M$.

$\bullet$  Theorem 4.11 remains valid if $I_f(\tau,T)$ is replaced with $I'_f(\tau,T)$
given by
$$\displaystyle
I'_f(\tau,T)
=\int_B\alpha_0(w-v')\,dx,
$$
where 
$$\displaystyle
v'(x,\tau)=\int_0^T\,e^{-\tau t}\,V(x,t)\,dt
$$
and $V$ solves
$$\left\{
\begin{array}{ll}
\displaystyle
\alpha_0(x)V_{tt}-\Delta V=0, & (x,t)\in\Bbb R^3\times\,]0,\,T[,\\
\\
\displaystyle
V(x,0)=0, & x\in\Bbb R^3,\\
\\
\displaystyle
\partial_tV(x,0)=\chi_B(x), & x\in\Bbb R^3.
\end{array}
\right.
$$
This indicator function can be used for detecting a change appearing as the distribution of $\alpha$
in some area without full knowledge of the background medium.  
For example, let $\alpha_0$ and $\alpha$ 
be the $n$-th day's and $n+1$-th day's $\alpha$.  
We assume that $\alpha_0=\alpha_1=1$ a.e. $B$.
Let $T$ be a large positive number.
By sending a wave in media $\alpha_n$ and $\alpha_{n+1}$ we have $V$ and $u$ on $B$, respectively.
Then the asymptotic behaviour of the indicator function $I'_f(\tau,T)$ may tell us a deviation or change of $\alpha_{n+1}$
from $\alpha_n$.

$\bullet$  We have the conjecture that under the assumption $\alpha_0\in L^{\infty}(\Bbb R^3)$,
it holds that the value
$$\displaystyle
\lim_{\tau\rightarrow\infty}\,\frac{1}{2\pi}\,\log\vert I_f(\tau,T)\vert
$$
exists and coincides with the first arrival time of the signal governed by the equation
$$\displaystyle
\alpha_0(x)u_{tt}-\Delta u=0,
$$
which started from $\partial B$ at $t=0$ and arrived at $\partial D$. 
We have two results that support this conjecture.

(1)  If $\alpha_0$ is a piecewise constant function and the space dimension is one, 
then we have a positive result \cite{Iwall}.

(2)  Recently, in \cite{IK5, IKK} the time domain enclosure method
has been applied to the problem of finding an unknown inclusion embedded in
a two layered medium, in the sense that the coefficient of the governing equation for the background medium
takes a constant value in each layered medium.  The result also supports the conjecture.

As a closely related problem, it would be interested to apply the enclosure method with the monostatic data
to the problem of finding information about the geometry
of an unknown impenetrable obstacle behind a known one.  
In \cite{IShadow} an asymptotic relation between an indicator function and some information about the geometry of the unknown obstacle
has been established.

\subsubsection{Inverse obstacle problems governed by the Maxwell system}

In this section we consider an inverse obstacle scattering problem using the electromagnetic wave.
The governing equation is the Maxwell system over a finite time interval.
The wave is generated at the initial time by a volumetric current source supported on a very small ball placed outside the obstacle.
Only the electric component of the wave is observed on the same ball (monostatic data) over a finite time interval.

It is assumed that the wave satisfies the Leontovich boundary condition 
on the surface of an unknown obstacle (causes the dissipation of the energy on the surface).
The condition is described by using an unknown positive function on the surface of the obstacle which is called the surface admittance.
It is shown that from the observed data one can extract information about the value of the surface admittance 
and the curvatures at the points on the surface nearest to the center of the ball.
This shows that a single shot over a finite time interval contains a meaningful information about 
the quantitative state of the surface of the obstacle.
This may have an application to the surface monitoring to check the state of the wall surface.

Let us formulate the problem more precisely.
We denote by $D$ the unknown obstacle.

It is assumed that $D$ is a non empty bounded open set of $\Bbb R^3$ with $C^2$-boundary such that $\Bbb R^3\setminus\overline D$
is connected The obstacle is embedded in
a medium like air (free space) which has constant electric permittivity $\epsilon(\,>0)$ and 
magnetic permeability $\mu(>0)$.
We assume that the electric field $\mbox{\boldmath $E$}=\mbox{\boldmath $E$}(x,t)$ and magnetic field $\mbox{\boldmath $H$}
=\mbox{\boldmath $H$}(x,t)$
are generated only by the current density $\mbox{\boldmath $J$}=\mbox{\boldmath $J$}(x,t)$ at the initial time
located not far away from the unknown obstacle.

The wave $(\mbox{\boldmath $E$}, \mbox{\boldmath $H$})$ satisfies
$$\displaystyle
\left\{
\begin{array}{ll}
\displaystyle
\frac{\partial}{\partial t}\mbox{\boldmath $E$}
-\epsilon^{-1}\nabla\times\mbox{\boldmath $H$}
=
\epsilon^{-1}\mbox{\boldmath $J$}, & (x,t)\in\Bbb R^3\setminus\overline D\times\,]0,\,T[,
\\
\\
\displaystyle
\frac{\partial}{\partial t}\mbox{\boldmath $H$}
+\mu^{-1}\nabla\times\mbox{\boldmath $E$}
=\mbox{\boldmath $0$}, & (x,t)\in\Bbb R^3\setminus\overline D\times\,]0,\,T[,
\\
\\
\displaystyle
\mbox{\boldmath $E$}(x,0)=\mbox{\boldmath $0$}, & x\in\Bbb R^3\setminus\overline D,
\\
\\
\displaystyle
\mbox{\boldmath $H$}(x,0)=\mbox{\boldmath $0$}, & x\in\Bbb R^3\setminus\overline D;
\end{array}
\right.
$$
and the boundary condition
$$
\displaystyle
\begin{array}{ll}\mbox{\boldmath $\nu$}\times\mbox{\boldmath $H$}(t)
-\lambda\,
\mbox{\boldmath $\nu$}\times(\mbox{\boldmath $E$}(t)\times\mbox{\boldmath $\nu$})=\mbox{\boldmath $0$},
& (x,t)\in\,\partial D\times\,]0,\,T[,
\end{array}
\tag {4.72}
$$
where $\lambda\in C^1(\partial D)$ and satisfies $\inf_{x\in\partial D}\lambda(x)>0$;
$\mbox{\boldmath $\nu$}$ denotes the unit outward normal to $\partial D$.

There should be several choices of current density $\mbox{\boldmath $J$}$ as a model of antenna, see \cite{Ba, CB}.

Here we assume that $\mbox{\boldmath
$J$}$ takes the form
$$
\displaystyle
\mbox{\boldmath $J$}(x,t)=f(t)\chi_B(x)\mbox{\boldmath $a$},
\tag {4.73}
$$
where $\mbox{\boldmath $a$}$ is an arbitrary unit vector, $B$ is a (very small)
open ball satisfying $\overline B\cap\overline D=\emptyset$,
$\chi_B$ denotes the characteristic function of $B$ and $f\in C^1[0,\,T]$ with $f(0)=0$.

$\bullet$  The boundary condition (4.72) is called the Leontovich boundary condition \cite{AT, CK, KH, N} and,
in particular,  \cite{KK} in the case when $\lambda$ is constant.

$\bullet$  The quantity $1/\lambda$ is called the surface impedance, see \cite{AT} and thus $\lambda$ is called the admittance.

$\bullet$  The existence of the admittance $\lambda$ causes the loss of the energy of the solution 
on the surface of the obstacle after stopping of the source supply.

See \cite{Kapitonov} for the study of the energy decay of the solution as $t\longrightarrow\infty$
Note that we always consider the case $T<\infty$.  Thus our result is independent of such study.

We consider the following problem.

$\quad$

{\bf\noindent Problem 4.9.}  Fix a large (to be determined later) $T<\infty$.  Observe $\mbox{\boldmath $E$}(t)$ on $B$
over the time interval $]0,\,T[$.
Extract information about the geometry of $D$ and the values of $\lambda$ on $\partial D$ from the observed data.

$\quad$

Denote the solution of the system mentioned above in the case when $D=\emptyset$ by
$(\mbox{\boldmath $E$}_0(t), \mbox{\boldmath $H$}_0(t))$ with $\mbox{\boldmath $J$}$ given by (4.73).

This means that

\noindent
(i)  the pair $(\mbox{\boldmath $E$}_0(t), \mbox{\boldmath $H$}_0(t))
\equiv (\mbox{\boldmath $E$}_0(\,\cdot\,,t), \mbox{\boldmath $H$}_0(\,\cdot\,,t))$ 
belongs to $C^1([0,\,T];L^2(\Bbb R^3)^3
\times L^2(\Bbb R^3)$ as a function of $t$;

\noindent
(ii)  for each $t\in[0,\,T]$, the pair $(\nabla\times\mbox{\boldmath $E$}_0(t),
\nabla\times\mbox{\boldmath $H$}_0(t))$ belongs to $L^2(\Bbb R^3)^3
\times L^2(\Bbb R^3)^3$;

\noindent
(iii) it holds that
$$\displaystyle
\left\{
\begin{array}{ll}
\displaystyle
\frac{d}{dt}\mbox{\boldmath $E$}_0
-\epsilon^{-1}\nabla\times\mbox{\boldmath $H$}_0
=
\epsilon^{-1}\mbox{\boldmath $J$},
&
\displaystyle
\frac{d}{dt}\mbox{\boldmath $H$}_0
+\mu^{-1}\nabla\times\mbox{\boldmath $E$}_0
=\mbox{\boldmath $0$},
\\
\\
\displaystyle
\mbox{\boldmath $E$}_0(0)=\mbox{\boldmath $0$},
&
\displaystyle
\mbox{\boldmath $H$}_0(0)=\mbox{\boldmath $0$}.
\end{array}
\right.
$$
Note that in this case, the solvability has been ensured by applying theory of $C_0$
contraction semigroups \cite{Y}.

{\bf\noindent Assumption 1.}
We assume also that the pair $(\mbox{\boldmath $E$}(t), \mbox{\boldmath $H$}(t))$ with $\mbox{\boldmath $J$}$ given by ({\bf C})
satisfies

\noindent
(i)'  the pair $(\mbox{\boldmath $E$}(t), \mbox{\boldmath $H$}(t))
\equiv (\mbox{\boldmath $E$}(\,\cdot\,,t), \mbox{\boldmath $H$}(\,\cdot\,,t))$ 
belongs to $C^1([0,\,T];L^2(\Bbb R^3\setminus\overline D)^3
\times L^2(\Bbb R^3\setminus\overline D)^3)$ as a function of $t$;

\noindent
(ii)'  for each $t\in[0,\,T]$, the pair $(\nabla\times\mbox{\boldmath $E$}(t),
\nabla\times\mbox{\boldmath $H$}(t))$ belongs to $L^2(\Bbb R^3\setminus\overline D)^3
\times L^2(\Bbb R^3\setminus\overline D)^3$;

\noindent
(iii)' it holds that
$$\displaystyle
\left\{
\begin{array}{ll}
\displaystyle
\frac{d}{dt}\mbox{\boldmath $E$}
-\epsilon^{-1}\nabla\times\mbox{\boldmath $H$}
=
\epsilon^{-1}\mbox{\boldmath $J$},
&
\displaystyle
\frac{d}{dt}\mbox{\boldmath $H$}
+\mu^{-1}\nabla\times\mbox{\boldmath $E$}
=\mbox{\boldmath $0$},
\\
\\
\displaystyle
\mbox{\boldmath $E$}(0)=\mbox{\boldmath $0$},
&
\displaystyle
\mbox{\boldmath $H$}(0)=\mbox{\boldmath $0$};
\end{array}
\right.
$$

\noindent
(iv) for each $t\in\,[0,\,T]$, $(\mbox{\boldmath $E$}(t), \mbox{\boldmath $H$}(t))$ satisfies,
in the sense of trace
$\displaystyle
\mbox{\boldmath $\nu$}\times\mbox{\boldmath $H$}(t)
-\lambda\,\mbox{\boldmath $\nu$}\times(\mbox{\boldmath $E$}(t)\times\mbox{\boldmath $\nu$})=\mbox{\boldmath $0$}$
on $\partial D$ as described in \cite{Kapitonov}.

Note that, at this stage, each term in the boundary condition does not have a point-wise meaning.  
What we know is: the left-hand side $\displaystyle
\mbox{\boldmath $\nu$}\times\mbox{\boldmath $H$}(t)
-\lambda\,\mbox{\boldmath $\nu$}\times(\mbox{\boldmath $E$}(t)\times\mbox{\boldmath $\nu$})$
just belongs to the dual space of $H^{1/2}(\partial D)^3$.

{\bf\noindent Definition 4.6.}
Define the indicator function
$$\begin{array}{ll}
\displaystyle
I_{\mbox{\boldmath $J$}}(\tau,T)
=\int_B
\mbox{\boldmath $f$}(x,\tau)\cdot(\mbox{\boldmath $W$}_e-\mbox{\boldmath $V$}_e)dx,
&
\tau>0
\end{array}
\tag {4.74}
$$
where
$$
\begin{array}{ll}
\displaystyle
\mbox{\boldmath $W$}_e(x,\tau)
=\int_0^T e^{-\tau t}\mbox{\boldmath $E$}(x,t)dt,
&
\displaystyle
\mbox{\boldmath $V$}_e(x,\tau)
=\int_0^Te^{-\tau t}\mbox{\boldmath $E$}_0(x,t)dt
\end{array}
$$
and
$$\begin{array}{ll}
\displaystyle
\mbox{\boldmath $f$}(x,\tau) & 
\displaystyle
=-\frac{\tau}{\epsilon}\int_0^Te^{-\tau t}\mbox{\boldmath $J$}(x,t)\,dt
\\
\\
\displaystyle
&
\displaystyle
=-\frac{\tau}{\epsilon}\int_0^Te^{-\tau t}f(t)dt\chi_B(x)\mbox{\boldmath $a$}.
\end{array}
$$

Set
$$
\displaystyle
\lambda_0=\sqrt{\frac{\epsilon}{\mu}}.
$$
We write $\lambda>>\lambda_0$ and $\lambda<<\lambda_0$ if
$\exists C>0$\,$\lambda(x)\ge\lambda_0+C$ for all $x\in\partial D$
and $\exists C>0$\,
$\lambda(x)\le\lambda_0-C$ for all $x\in\partial D$, respectively.

Roughly speaking, we can say that: the condition $\lambda>>\lambda_0$/$\lambda<<\lambda_0$
means that the admittance $\lambda$ is greater/less than the special value 
$\lambda_0$ which is the admittance of free space \cite{AT}.

Here we introduce a technical assumption which should be removed.

Define also
$$\displaystyle
\mbox{\boldmath $W$}_m(x,\tau)
=\int_0^T e^{-\tau t}\mbox{\boldmath $H$}(x,t)dt.
$$
We have
$$
\left\{
\begin{array}{ll}
\displaystyle
\nabla\times\mbox{\boldmath $W$}_e+\tau\mu \mbox{\boldmath $W$}_m=-e^{-\tau T}\mu\mbox{\boldmath $H$}(x,T),
& x\in\,\Bbb R^3\setminus\overline D,\\
\\
\displaystyle
\nabla\times\mbox{\boldmath $W$}_m-\tau\epsilon \mbox{\boldmath $W$}_e
-\frac{\epsilon}{\tau}\mbox{\boldmath $f$}(x,\tau)
=e^{-\tau T}\epsilon \mbox{\boldmath $E$}(x,T),
& x\in\,\Bbb R^3\setminus\overline D
\end{array}
\right.
$$
and 
$$\begin{array}{ll}
\displaystyle
\mbox{\boldmath $\nu$}\times\mbox{\boldmath $W$}_m
-\lambda\,\mbox{\boldmath $\nu$}\times(\mbox{\boldmath $W$}_e\times\mbox{\boldmath $\nu$})=\mbox{\boldmath $0$},
& x\in\partial D.
\end{array}
\tag {4.75}
$$

Here we introduce a technical assumption concerning the regularity:

$\quad$

{\bf\noindent Assumption 2.}  The functions $\mbox{\boldmath $W$}_e$ and $\mbox{\boldmath $W$}_m$
belong to $H^1$ on the intersection of an open neighbourhhod of $\partial D$ in $\Bbb R^3$ with $\Bbb R^3\setminus\overline D$.

This assumption makes us possible to treat the each term in the Leontovich boundary condition for $(\mbox{\boldmath $W$}_e, \mbox{\boldmath $W$}_m)$
pointwise.  More precisely,
from this assumption we see that
the boundary condition (4.75) is satisfied in the sense of the usual trace in $H^{1/2}(\partial D)^3$.
See \cite{Kapitonov} and \cite{IMax2,IMax3} for discussion related to this assumption.

In \cite{IMax2} we have obtained the following result.

\proclaim{\noindent Theorem 4.12(\cite{IMax2}).} 
Let $\mbox{\boldmath $a$}_j$,
$j=1,2$ be linearly independent unit vectors. Let
$\mbox{\boldmath $J$}_j(x,t)=f(t)\chi_B(x)\mbox{\boldmath $a$}_j$
and $f$ satisfy
$$
\displaystyle
\exists\gamma\in\Bbb R\,\,\liminf_{\tau\rightarrow\infty}\tau^{\gamma}
\left\vert\int_0^Te^{-\tau t}f(t)\,dt\right\vert>0.
\tag {4.76}
$$
Then, we have:
$$
\displaystyle
\lim_{\tau\rightarrow\infty}e^{\tau T}\sum_{j=1}^2I_{\mbox{\boldmath $J$}_j}(\tau,T)
=
\left\{
\begin{array}{ll}
\displaystyle
0,
& \text{if $T\le2\sqrt{\mu\epsilon}\text{dist}\,(D,B)$,}\\
\\
\displaystyle
\infty,
&
\text{if $T>2\sqrt{\mu\epsilon}\text{dist}\,(D,B)$ and $\lambda>>\lambda_0$,}\\
\\
\displaystyle
-\infty,
&
\text{if $T>2\sqrt{\mu\epsilon}\text{dist}\,(D,B)$ and $\lambda<<\lambda_0$.}
\end{array}
\right.
$$
Moreover, if $\lambda$ satisfies $\lambda>>\lambda_0$ or $\lambda<<\lambda_0$, then
for all $T>2\sqrt{\mu\epsilon}\,\text{dist}\,(D,B)$
$$
\displaystyle
\lim_{\tau\rightarrow\infty}\frac{1}{\tau}
\log\left\vert
\sum_{j=1}^2I_{\mbox{\boldmath $J$}_j}(\tau,T)\right\vert
=-2\sqrt{\mu\epsilon}\,\text{dist}\,(D,B).
$$

\endproclaim

Why two?

$\bullet$  Any $\mbox{\boldmath $\nu$}_q$ with $q\in\Lambda_{\partial D}(p)$ never satisfies
$\mbox{\boldmath $a$}_1\times\mbox{\boldmath $\nu$}_q=\mbox{\boldmath $a$}_2\times\mbox{\boldmath $\nu$}_q=\mbox{\boldmath $0$}$
if $\mbox{\boldmath $a$}_1$ and $\mbox{\boldmath $a$}_2$ are linearly independent.

$\bullet$  Indicator function (4.74)
becomes weak if $\mbox{\boldmath $a$}\times\mbox{\boldmath $\nu$}_q\approx \mbox{\boldmath $0$}$
with $q\in\Lambda_{\partial D}(p)$.  See also Theorem 4.13 described later.
This is consistent with a property of the dipole antenna.  We should not direct $\mbox{\boldmath $a$}$ to the obstacle.

Some remarks are in order.

All the statements of Theorem 4.12 are valid if $\mbox{\boldmath $V$}_e$ in $I_{\mbox{\boldmath $J$}}(\tau,T)$
is replaced with the unique weak solution $\mbox{\boldmath $V$}\in L^2(\Bbb R^3)^3$ with $\nabla\times
\mbox{\boldmath $V$}\in L^2(\Bbb R^3)^3$
of
$$\begin{array}{ll}
\displaystyle
\frac{1}{\mu\epsilon}\nabla\times\nabla\times\mbox{\boldmath $V$}
+\tau^2\mbox{\boldmath $V$}
+\mbox{\boldmath $f$}(x,\tau)=\mbox{\boldmath $0$},
& x\in\Bbb R^3.
\end{array}
$$ 
In what follows, we denote by $\mbox{\boldmath $V$}_e^0$ the weak solution.

The reason why such a replacement is possible is the following.
Introduce another indicator function by the formula
$$
\displaystyle
\tilde{I}_{\mbox{\boldmath $f$}}(\tau,T)
=\int_B\mbox{\boldmath $f$}(x,\tau)\cdot(\mbox{\boldmath $W$}_e-\mbox{\boldmath $V$}_e^0)\,dx.
\tag {4.77}
$$
Using the simple facts
$$\displaystyle
\Vert\mbox{\boldmath $V$}_e-\mbox{\boldmath $V$}_e^0\Vert_{L^2(\Bbb R^3\setminus\overline D)}
=O(\tau^{-1}e^{-\tau T})
$$
and
$$
\displaystyle
\Vert\mbox{\boldmath $f$}\Vert_{L^2(B)}=O(\tau^{-1/2}),
$$
one has
$$\displaystyle
I_{\mbox{\boldmath $J$}}(\tau,T)
=\tilde{I}_{\mbox{\boldmath $J$}}(\tau,T)
+O(\tau^{-3/2}e^{-\tau T}).
$$
Thus, one can transplant all the results in Theorem 4.12 into the case when the indicator function is given by (4.77).

This version's advantage is: no need of time domain computation to compute $\mbox{\boldmath $V$}_e^0$
unlike $\mbox{\boldmath $V$}_e$.

Succeeding to \cite{IMax} which treats the case when $\lambda=\infty$, in \cite{IMax3} the author clarified
the leading profile of the indicator functions (4.74) or (4.77) as $\tau\longrightarrow\infty$
which yields a result beyond Theorem 4.12 in the following sense: one can extract quantitative information about the state of the surface of an unknown
obstacle using the time domain enclosure method.

\proclaim{\noindent Theorem 4.13(\cite{IMax3}).}
Assume that $\partial D$ is $C^4$ and $\lambda\in C^1(\partial D)$.
Let $f$ satisfy (4.76) and $T>2\sqrt{\mu\epsilon}\,\text{dist}\,(D,B)$.
Assume that the first reflector $\Lambda_{\partial D}(p)$ consists of finite points
and
$$\displaystyle
\text{det}\,(S_q(\partial B_{d_{\partial D}(p)}(p))-S_{q}(\partial D))>0\,\,\,
\forall q\in\Lambda_{\partial D}(p).
\tag {4.78}
$$
Assume also that there exists a point $q\in\Lambda_{\partial D}(p)$ such that $\lambda(q)\not=\lambda_0$
and $\mbox{\boldmath $\nu$}_q\times\mbox{\boldmath $a$}\not=\mbox{\boldmath $0$}$.
Then, we have
$$
\displaystyle
\lim_{\tau\longrightarrow\infty}\tau^2e^{2\tau\,\sqrt{\mu\epsilon}\,\text{dist}\,(D,B)}\frac{\tilde{I}_{\mbox{\boldmath $J$}}(\tau,T)}
{\tilde{f}(\tau)^2}
=
\frac{\pi}{2}
\left(\frac{\eta}{d_{\partial D}(p)}\right)^2\frac{\lambda_0^2}{\epsilon^4}
\,
\sum_{q\in\Lambda_{\partial D}(p)}
k_q(p)
\frac{\displaystyle
\lambda(q)-\lambda_0}
{\displaystyle
\lambda(q)+\lambda_0}
\vert\mbox{\boldmath $\nu$}_q\times\mbox{\boldmath $a$}\vert^2,
$$
where
$$\displaystyle
\tilde{f}(\tau)=\int_0^T e^{-\tau t}f(t)\,dt,\,\,
k_q(p)=\frac{1}{\displaystyle
\sqrt{\text{det}\,(S_q(\partial B_{d_{\partial D}(p)}(p))-S_{q}(\partial D))}}.
$$

\endproclaim

{\it\noindent Sketch of Proof.}
Define
$$\displaystyle
\mbox{\boldmath $V$}_m^0=-\frac{1}{\tau\mu}\,\nabla\times\mbox{\boldmath $V$}_e^0.
$$
We have
$$
\left\{
\begin{array}{ll}
\displaystyle
\nabla\times\mbox{\boldmath $V$}_e^0+\tau\mu \mbox{\boldmath $V$}_m^0=\mbox{\boldmath $0$},
& x\in\,\Bbb R^3,\\
\\
\displaystyle
\nabla\times\mbox{\boldmath $V$}_m^0-\tau\epsilon \mbox{\boldmath $V$}_e^0
-\frac{\epsilon}{\tau}\mbox{\boldmath $f$}(x,\tau)
=\mbox{\boldmath $0$},
& x\in\,\Bbb R^3.
\end{array}
\right.
$$
Define
$$\displaystyle
\left\{
\begin{array}{l}
\displaystyle
\mbox{\boldmath $R$}_e=\mbox{\boldmath $W$}_e-\mbox{\boldmath $V$}_e^0,
\\
\\
\displaystyle
\mbox{\boldmath $R$}_m=\mbox{\boldmath $W$}_m-\mbox{\boldmath $V$}_m^0.
\end{array}
\right.
$$
This is the reflected solution.

The proof starts with obtaining the expression
$$
\displaystyle
\tilde{I}_{\mbox{\boldmath $f$}}(\tau,T)
=J(\tau)+E(\tau)+O(e^{-\tau T}\tau^{-1}),
$$
where
$$
\left\{
\begin{array}{l}
\displaystyle
J(\tau)
=\frac{1}{\mu\epsilon}\int_{\partial D}(\mbox{\boldmath $\nu$}\times\mbox{\boldmath $V$}_e^0)
\cdot\nabla\times\mbox{\boldmath $V$}_e^0\,dS-
\frac{\tau}{\epsilon}\int_{\partial D}\frac{1}{\lambda}\vert\mbox{\boldmath $V$}_m^0\times\mbox{\boldmath $\nu$}\vert^2dS,
\\
\\
\displaystyle
E(\tau)
=\frac{\tau}{\epsilon}\left\{\int_{\Bbb R^3\setminus\overline D}(\tau\mu\vert\mbox{\boldmath $R$}_m\vert^2
+\tau\epsilon\vert\mbox{\boldmath $R$}_e\vert^2)\,dx
+\int_{\partial D}\frac{1}{\lambda}\vert\mbox{\boldmath $R$}_m\times\mbox{\boldmath $\nu$}\vert^2\,dS\right\}.
\end{array}
\right.
$$

Thus, the essential part of the proof of Theorem 4.13 should be the study of the asymptotic behaviour of
$J(\tau)$ and $E(\tau)$ as
$\tau\longrightarrow\infty$.  
The asymptotic behaviour of $J(\tau)$ can be reduced to that of a Laplace-type integral \cite{IMax2}.
For $E(\tau)$ we have Theorem 4.14 mentioned below.  It reduces also to a Laplace-type integral.

\noindent
$\Box$

\proclaim{\noindent Theorem 4.14(\cite{IMax3}).}
Assume that $\partial D$ is $C^4$ an $\lambda\in C^2(\partial D)$.
Assume that first reflector $\Lambda_{\partial D}(p)$ consists of finite points
and (4.78) is satisfied; there exists a point $q\in\Lambda_{\partial D}(p)$ such that $\lambda(q)\not=\lambda_0$
and that $\displaystyle
\mbox{\boldmath $\nu$}_q\times(\mbox{\boldmath $a$}\times\mbox{\boldmath $\nu$}_q)\not=\mbox{\boldmath $0$}$.
Let $f$ satisfy (4.76) and $\displaystyle
T>\sqrt{\mu\epsilon}\,\text{dist}\,(D,B)$.
Then, we have $J^*(\tau)>0$ for all $\tau>>1$ and
$$
\displaystyle
\lim_{\tau\longrightarrow\infty}\frac{E(\tau)}
{\displaystyle
J^*(\tau)
}=1,
\tag {4.79}
$$
where
$$\left\{
\begin{array}{l}
\displaystyle
J^*(\tau)
=\frac{\tau}{\epsilon}
\int_{\partial D}
\frac{\displaystyle
\lambda-\lambda_0}
{\displaystyle
\lambda+\lambda_0}\,
(\mbox{\boldmath $\nu$}\times\mbox{\boldmath $V$}_m^0)\cdot
\mbox{\boldmath $V$}_{em}^0\,dS,
\\
\\
\displaystyle
\mbox{\boldmath $V$}_{em}^0=\mbox{\boldmath $\nu$}\times(\mbox{\boldmath $V$}_e^0\times\mbox{\boldmath $\nu$})
-\frac{1}{\lambda}\,\mbox{\boldmath $\nu$}\times\mbox{\boldmath $V$}_m^0.
\end{array}
\right.
$$
\endproclaim

{\it\noindent Sketch of Proof.}
First recall some preliminary knowledge about the reflection \cite{GT}.
Let $x^r$ denote the reflection across $\partial D$ of the point $x\in\Bbb R^3\setminus D$ with
$d_{\partial D}(x)<2\delta_0$ for a sufficiently small $\delta_0>0$.
It is given by $x^r=2q(x)-x$, where $q(x)$ denotes the unique point on $\partial D$ such that
$d_{\partial D}(x)=\vert x-q(x)\vert$. 
Note that $q(x)$ is $C^2$ for $x\in\Bbb R^3\setminus D$ with $d_{\partial D}(x)<2\delta_0$ if $\partial D$ is $C^3$.

In Proposition 3 of \cite{IMax} we have already established the {\it reflection principle} for the Maxwell system across the {\it curved surface}.
Applying the principle to the present case, one obtains the vector field for $x\in\,\Bbb R^3\setminus\overline D$ with $d_{\partial D}(x)<2\delta_0$
denoted by $(\mbox{\boldmath $V$}_e^0)^*$ from vector field $\mbox{\boldmath $V$}_e^0$ in $D$ in such a way that
$$\left\{
\begin{array}{ll}
\displaystyle
(\mbox{\boldmath $V$}_e^0)^*=-\mbox{\boldmath $V$}_e^0, & x\in\partial D,
\\
\\
\displaystyle
\mbox{\boldmath $\nu$}\times(\mbox{\boldmath $V$}_m^0)^*=
\mbox{\boldmath $\nu$}\times\mbox{\boldmath $V$}_m^0, & x\in\partial D,
\end{array}
\right.
$$
where
$$\begin{array}{lll}
\displaystyle (\mbox{\boldmath $V$}_m^0)^*\equiv -\frac{1}{\tau\mu}\nabla\times\{(\mbox{\boldmath $V$}_e^0)^*\},
& x\in\Bbb R^3\setminus\overline{D}, & d_{\partial D}(x)<2\delta_0.
\end{array}
$$
For $x\in\Bbb R^3\setminus\overline D$ with $d_{\partial D}(x)<2\delta_0$ 
we have
$$\begin{array}{ll}
\displaystyle
\frac{1}{\mu\epsilon}\nabla\times\nabla\times(\mbox{\boldmath $V$}_e^0)^*+\tau^2(\mbox{\boldmath $V$}_e^0)^*
&
\displaystyle
=\text{terms from $\mbox{\boldmath $V$}_e^0(x^r)$ and $(\mbox{\boldmath $V$}_e^0)'(x^r)$}
\\
\\
\displaystyle
&
\displaystyle
\,\,\,
+2d_{\partial D}(x)
\times\text{terms from $(\nabla^2\mbox{\boldmath $V$}_e^0)(x^r)$}\\
\\
\displaystyle
&
\displaystyle
\equiv \mbox{\boldmath $Q$}(x).
\end{array}
$$
Note that: all the coefficients of $\mbox{\boldmath $V$}_e^0(x^r)$, $(\mbox{\boldmath $V$}_e^0)'(x^r)$
and $(\nabla^2\mbox{\boldmath $V$}_e^0)(x^r)$ in $\mbox{\boldmath $Q$}(x)$ are 
independent of $\tau$ and continuous in a
tubular neighbourhood of $\partial D$;
the coefficients of $(\nabla^2\mbox{\boldmath $V$}_e^0)(x^r)$ in $\mbox{\boldmath $Q$}(x)$ are $C^1$ therein;
the construction of $(\mbox{\boldmath $V$}_e^0)^*$ involves the {\it shape operator} of $\partial D$ (affecting the curvatures).

Define
$$\left\{
\begin{array}{l}
\displaystyle
\mbox{\boldmath $R$}_e^0
=\frac{\displaystyle\tilde{\lambda}-\lambda_0}
{\displaystyle\tilde{\lambda}+\lambda_0}
\phi_{\delta}
\,(\mbox{\boldmath $V$}_e^0)^*,
\\
\\
\displaystyle
\mbox{\boldmath $R$}_m^0
=\frac{\displaystyle\tilde{\lambda}-\lambda_0}
{\displaystyle\tilde{\lambda}+\lambda_0}
\,\phi_{\delta}\,(\mbox{\boldmath $V$}_m^0)^*,
\end{array}
\right.
$$
where

$\bullet$  $\tilde{\lambda}(x)=\lambda(q(x))$ for $x\in\Bbb R^3\setminus D$ with $d_{\partial D}(x)<2\delta_0$.
The function $\tilde{\lambda}$ is $C^2$ therein and coincides with $\lambda(x)$ for $x\in\partial D$.

$\bullet$  a cutoff function $\phi_{\delta}\in C^{2}(\Bbb R^3)$ with
$0<\delta<\delta_0$ which satisfies
$0\le\phi_{\delta}(x)\le 1$; $\phi_{\delta}(x)=1$ if $d_{\partial D}(x)<\delta$;
$\phi_{\delta}(x)=0$ if $d_{\partial D}(x)>2\delta$;
$\vert\nabla\phi_{\delta}(x)\vert\le C\delta^{-1}$;
$\vert\nabla^2\phi_{\delta}(x)\vert\le C\delta^{-2}$.

The pair $(\mbox{\boldmath $R$}_e^0,\mbox{\boldmath $R$}_m^0)$ depends on $\delta$.

Define
$$
\left\{
\begin{array}{l}
\displaystyle
\mbox{\boldmath $R$}_e^1=\mbox{\boldmath $R$}_e-\mbox{\boldmath $R$}_e^0,\\
\\
\displaystyle
\mbox{\boldmath $R$}_m^1=\mbox{\boldmath $R$}_m-\mbox{\boldmath $R$}_m^0.
\end{array}
\right.
$$
We see that
$$
\left\{
\begin{array}{ll}
\displaystyle
\nabla\times\mbox{\boldmath $R$}_e^1+\tau\mu\mbox{\boldmath $R$}_m^1
=-(\nabla\times\mbox{\boldmath $R$}_e^0+\tau\mu\mbox{\boldmath $R$}_m^0)
-e^{-\tau T}\mu\mbox{\boldmath $H$}(x,T),
& x\in\Bbb R^3\setminus\overline D,\\
\\
\displaystyle
\nabla\times\mbox{\boldmath $R$}_m^1-\tau\epsilon\mbox{\boldmath $R$}_e^1
=-(\nabla\times\mbox{\boldmath $R$}_m^0-\tau\epsilon\mbox{\boldmath $R$}_e^0)
+e^{-\tau T}\epsilon\mbox{\boldmath $E$}(x,T),
& x\in\Bbb R^3\setminus\overline D
\end{array}
\right.
$$
and
$$
\begin{array}{ll}
\displaystyle
\mbox{\boldmath $\nu$}\times\mbox{\boldmath $R$}_m^1
-\lambda\,\mbox{\boldmath $\nu$}\times(\mbox{\boldmath $R$}_e^1\times\mbox{\boldmath $\nu$})
=-\mbox{\boldmath $V$}_1, & x\in\partial D,
\end{array}
$$
where
$$
\begin{array}{ll}
\displaystyle
\mbox{\boldmath $V$}_1
=-\frac{2\lambda\lambda_0}{\lambda+\lambda_0}\mbox{\boldmath $V$}_{em}^0\vert_{\lambda=\lambda_0},
&
x\in\partial D.
\end{array}
$$

Then, roughly speaking, one gets an asymptotic profile of $E(\tau)-J^*(\tau)$ as $\tau\longrightarrow\infty$
which extracts the main term involving $\mbox{\boldmath $\nu$}\times\mbox{\boldmath $R$}_m^1$
on $\partial D$:
$$\begin{array}{ll}
\displaystyle
E(\tau)
&
\displaystyle
=J^*(\tau)
+\frac{\tau}{\epsilon}\int_{\partial D}
\mbox{\boldmath $V$}_{em}^0\cdot(\mbox{\boldmath $\nu$}\times\mbox{\boldmath $R$}_m^1)\,dS\\
\\
\displaystyle
&
\displaystyle
\,\,\,
+O(e^{-\tau T}(\tau^{-2}e^{-\tau\sqrt{\mu\epsilon}\,\text{dist}\,(D,B)}\vert\tilde{f}(\tau)\vert+\tau^{-1}e^{-\tau T})).
\end{array}
\tag {4.80}
$$
Note that the term $O(e^{-\tau T}(\tau^{-2}e^{-\tau\sqrt{\mu\epsilon}\,\text{dist}\,(D,B)}\vert\tilde{f}(\tau)\vert+\tau^{-1}e^{-\tau T}))$
is uniform with respect to $\delta\in\,]0,\,\delta_0[$.

Here we see that, by choosing $\delta=\tau^{-1/2}$,
$$\displaystyle
\lim_{\tau\longrightarrow\infty}\tau^3\,e^{2\tau\,\sqrt{\mu\epsilon}\,\text{dist}\,(D,B)}
\frac{\displaystyle
\Vert\mbox{\boldmath $\nu$}\times\mbox{\boldmath $R$}_m^1\Vert_{L^2(\partial D)}^2}
{\tilde{f}(\tau)^2}=0.
\tag {4.81}
$$
This enables us to ignore the second term on (4.80) relative to $J^*(\tau)$ and yields (4.79).
The key point in showing the validity of (4.81) is:
as $\tau\rightarrow\infty$, $\mbox{\boldmath $V$}_{em}^0\vert_{\lambda=\lambda_0}$ in $\mbox{\boldmath $V$}_1$ on $\partial D$ becomes {\it small}.
This is a mathematical meaning of the admittance of free space.

\noindent
$\Box$

The proof of (4.81) is a modification of the Lax-Phillips reflection argument \cite{LP} originally developed for the study of the right-end point
of the support of the scattering kernel for the scalar wave equation in the context of the Lax-Phillips scattering theory.
However, it should be emphasized that our version of their argument (not given here in detail, see \cite{IMax3}) is rather straightforward.

Some additional remarks are in order.

$\bullet$  The factor $2$ in the restriction $T>2\sqrt{\mu\epsilon}\,\text{dist}\,(D,B)$ in
Theorem 4.13 is dropped in Theorem 4.14.  The quantity $\sqrt{\mu\epsilon}\,\text{dist}\,(D,B)$
corresponds to the first arrival time of the wave generated at $t=0$ on $B$ and reached at $\partial D$ firstly.

$\bullet$  The asymptotic formula (4.79)
clarifies the effect on the leading profile of the energy of the reflected solutions
$\mbox{\boldmath $R$}_e$ and $\mbox{\boldmath $R$}_m$ in terms of the deviation of th surface admittance
from that of free-space admittance and the energy density of the incident wave.

$\bullet$  In Theorem 4.13 neither jump condition $\lambda>>\lambda_0$ nor $\lambda<<\lambda_0$ is assumed.
Instead it is assumed that there exists a point $q\in\Lambda_{\partial D}(p)$ such that $\lambda(q)\not=\lambda_0$
and $\mbox{\boldmath $\nu$}_q\times\mbox{\boldmath $a$}\not=\mbox{\boldmath $0$}$.

$\bullet$  Theorem 4.13 remain valid if $\tilde{I}_{\mbox{\boldmath $J$}}$ is replaced with
$I_{\mbox{\boldmath $J$}}$.

As explained in Section 4.4.2 for the scalar wave equation case \cite{IMP}, Theorem 4.13 suggests us a procedure for:
finding Gauss and mean curvatures and $\lambda$ at an arbitrary point $q$ on $\Lambda_{\partial D}(p)$.
In particular, one can know the mean curvature at a known point on $\partial D$ nearest to $p$!  
This is the advantage of near field measurements.
It is a translated one of that described in \cite{IMP} in which the scalar wave equation
is considered.  Thus we do not mention it's details here.

A problem to be solved is the following.

$\bullet$ Extend the results in \cite{Iwall} for scalar wave case (obstacle detection, together with a rough distinction)
to a penetrable or impenetrable obstacle embedded in a general inhomogeneous medium
(imagine the case $\mu$, $\epsilon$ and $\sigma$ are inhomogeneous
in the whole space and have a jump from the background medium on a bounded domain).

\subsubsection{Enclosing an unknown obstacle using a single response}

To pursuit the possibility of the time domain enclosure method,  in \cite{EIV} and \cite{EV} the author
considered the problem of detecting an unknown obstacle embedded in
a bounded domain by generating a single wave at the boundary and receiving only the wave propagating inside.

The problem considered therein is as follows.
Let $\Omega$ be a bounded domain with $C^2$-boundary.
Le $D$ be a nonempty bounded open subset of $\Omega$ with $C^2$-boundary such that
$\Omega\setminus\overline D$ is connected.
Let $0<T<\infty$.
Given $f=f(x,t), (x,t)\in\partial\Omega\times\,]0,\,T[$ let $u=u_f(x,t), (x,t)\in\,(\Omega\setminus\overline D)\times\,]0,\,T[$ 
denote the solution of the following initial boundary value problem for the classical wave equation:
$$\displaystyle
\left\{
\begin{array}{ll}
\displaystyle
u_{tt}-\Delta u=0, & (x,t)\in\,(\Omega\setminus\overline D)\times\,]0,\,T[,\\
\\
\displaystyle
u(x,0)=0, & x\in\Omega\setminus\overline D,
\\
\\
\displaystyle
u_t(x,0)=0, & x\in\Omega\setminus\overline D,\\
\\
\displaystyle
\frac{\partial u}{\partial\nu}=0, & (x,t)\in\partial D\times\,]0,\,T[,\\
\\
\displaystyle
\frac{\partial u}{\partial\nu}=f(x,t), & (x,t)\in\partial\Omega\times\,]0,\,T[.
\end{array}
\right.
\tag {4.82}
$$
We use the same symbol $\nu$ to denote both the outer unit normal vectors of both $\partial D$
and $\partial\Omega$.

$\quad$

{\bf\noindent Problem 4.10.}  Fix a large $T$ (to be determined later) and a single $f$ (to be specified later).
Assume that set $D$ is  unknown.  Extract information about the
location and shape of $D$ from the wave field $u_f(x,t)$, which is given for
all $x\in\partial\Omega$ and $t\in\,]0,\,T[$.

$\quad$

The correspondence
$$\begin{array}{lll}
\displaystyle
f\longmapsto 
u_f(x,t), & x\in\partial\Omega, & 0<t<T,
\end{array}
$$
is called the response operator \cite{Be}.
Problem 4.10 asks us, in other words, to extract information
about the geometry of $D$ from a {\it single point} on the graph of the response operator.
Note that the BC method \cite{Be, Be2} gives us a reconstruction formula of an unknown
variable coefficient of the governing equation of wave propagation, however, it makes use of full knowledge of the response operator.

The time domain enclosure method presented here employs the data over the {\it finite time interval} $]0,\,T[$.
Because of the finite propagation property of the wave equation any exact method to extract information about unknown discontinuity being not near
the surface needs a restriction on $T$ from below.  It should be truly used in the proof to justify the method.
The time domain enclosure method meets the requirement.

Here we describe two types of input Neumann data.

(i)  First let $B$ be an open ball centered at $p$ and radius $\eta$ satisfying $\overline B\cap\overline\Omega=\emptyset$.
Define
$$\begin{array}{ll}
\displaystyle
f_B=\frac{\partial}{\partial\nu}\,v_B(x,t), & (x,t)\in\partial\Omega\times\,]0,\,T[,
\end{array}
$$
where $v=v_B$ solves
$$\displaystyle
\left\{
\begin{array}{ll}
\displaystyle
v_{tt}-\Delta v=0, & (x,t)\in\Bbb R^3\times\,]0,\,T[,\\
\\
\displaystyle
v(x,0)=0, & x\in\Bbb R^3,\\
\\
\displaystyle
v_t(x,0)=g(x), & x\in\Bbb R^3,
\end{array}
\right.
\tag {4.83}
$$
under the special choice of the initial velocity field $g\in H^1(\Bbb R^3)$
given by
$$\begin{array}{ll}
\displaystyle
g(x)=\Psi_B(x)\equiv (\eta-\vert x-p\vert)\chi_B(x), & x\in\Bbb R^3.
\end{array}
\tag {4.84}
$$

(ii)  Second let $B$ be an arbitrary open ball centered at $p$ and radius $\eta$.
In contrast to (i) there is no restriction on the location of $B$ such as $\overline B\cap\overline\Omega=\emptyset$.
Define
$$\begin{array}{ll}
\displaystyle
f_{B}^*=\frac{\partial}{\partial\nu}\,v_B(x,T-t), & (x,t)\in\partial\Omega\times\,]0,\,T[,
\end{array}
$$
where $v=v_B$ solves (4.83).  Note that in the right-hand side the time reversal operation appears.  
Thus $f_{B}^*$ plays a role of time reversal mirror \cite{F} when obstacle $D$ is absent.

Both the functions $f_B$ and $f_B^*$ do not contain any large parameter.
Let $u=u_B, u_B^*$ be the solutions of (4.82) with $f=f_B, f_B^*$, respectively.

Define two indicator functions $I_{T}(\tau;B)$ and $I^*_{T}(\tau;B)$ of independent variable $\tau>0$
by the formulae
$$\left\{
\begin{array}{l}
\displaystyle
I_{T}(\tau;B)=\int_{\partial\Omega}\,(w_B-w_B^0)\,\frac{\partial w_B^0}{\partial\nu}\,dS,
\\
\\
\displaystyle
I_{T}^*(\tau;B)=\int_{\partial\Omega}\,(w_B^*-w_B^{0,*})\,\frac{\partial w_B^{0,*}}{\partial\nu}\,dS,
\end{array}
\right.
$$
where
$$\left\{
\begin{array}{ll}
\displaystyle
w_B(x,\tau)=\int_0^Te^{-\tau t}u_B(x,t)dt,
&
x\in\Omega\setminus\overline D,
\\
\\
\displaystyle
w_B^*(x,\tau)=\int_0^Te^{-\tau t}u_B^*(x,t)dt,
&
x\in\Omega\setminus\overline D,
\\
\\
\displaystyle
w_B^0(x,\tau)=\int_0^Te^{-\tau t}v_B(x,t)dt,
&
x\in\Bbb R^3,
\\
\\
\displaystyle
w_B^{0,*}(x,\tau)=\int_0^Te^{-\tau t}v_B(x,T-t)dt,
&
x\in\Bbb R^3.
\end{array}
\right.
\tag {4.85}
$$

These indicator functions can be computed from the response $u_B$ and $u_B^*$ on $\partial\Omega$ over
the time interval $]0,\,T[$ which is the solution of (4.82) with $f=f_B, f_{B}^*$.

First we describe the asymptotic behaviour of the indicator function $I_T(\tau;B)$.

\proclaim{\noindent Theorem 4.15(\cite{EIV}).}

(i)  If $T$ satisfies
$$\displaystyle
T>2\text{dist}\,(D,B)-\text{dist}\,(\Omega,B),
\tag {4.86}
$$
then, there exists a positive number $\tau_0$ such that
$I_{T}(\tau;B)>0$ for all $\tau\ge\tau_0$, and we have
$$\displaystyle
\lim_{\tau\longrightarrow\infty}
\frac{1}{\tau}
\log I_{T}(\tau;B)=-2\text{dist}\,(D,B).
\tag {4.87}
$$

(ii)  We have
$$\displaystyle
\lim_{\tau\longrightarrow\infty}e^{\tau T}I_{T}(\tau;B)=
\left\{
\begin{array}{ll}
\displaystyle
\infty & \text{if}\,\,T>2\text{dist}\,(D,B),\\
\\
\displaystyle
0     & \text{if}\,\,T<2\text{dist}\,(D,B).
\end{array}
\right.
\tag {4.88}
$$

\endproclaim

Since we have $\text{dist}\,(D,B)=d_{\partial D}(p)-\eta$, from Theorem 4.15 we see that the indicator function
for each $B$ uniquely determines $d_{\partial D}(p)$ and hence the sphere $\vert x-p\vert=d_{\partial D}(p)$
on which there exists a point on $\partial D$.  Moving $p$ outside $\Omega$, say $p=p_1,\cdots, p_m$, $m<\infty$,
one can obtain an estimation of the geometry of $D$ by using the intersection of the exteriors of the spheres
$\vert x-p_j\vert=d_{\partial D}(p_j)$, $j=1,\cdots,m$.  When data becomes a lot, obtained information
becomes detailed.  This is not like a method that one input yields everything.

Some further remarks are in order.

$\bullet$  If $T=2\,\text{dist}\,(D,B)$, we have $e^{\tau T}\,I_T(\tau;B)=O(\tau^4)$ as $\tau\rightarrow\infty$.

$\bullet$  It follows from (ii) in Theorem 4.15 that
$$\displaystyle
2\,\text{dist}\,(D,B)
=
\sup\left\{T\in\,]0,\,\infty[\,\vert\,
\lim_{\tau\rightarrow\infty}\,e^{\tau T}\,I_T(\tau;B)=0
\right\}.
$$
This formula has a similarity with the original version of the enclosure method \cite{E00}.

$\bullet$  We have, see \cite{IW00}, 
$$\displaystyle
2\,\text{dist}\,(D,B)-\text{dist}\,(\Omega,B)
\ge
\inf\,\left\{
\vert x-y\vert+\vert y-z\vert\,\vert\,x\in\partial B,\,y\in\partial D,\,z\in
\partial\Omega\,\right\}.
$$
Thus, the restriction (4.86) on $T$ from below does not against the geometrical optics.

$\bullet$  We observe the reflected wave on $\partial\Omega$.  So we have no need of waiting until
the time $T=2\,\text{dist}\,(D,B)$ which is greater than $2\,\text{dist}\,(D,B)-\text{dist}\,(\Omega, B)$, to get information about the geometry of $D$.

$\bullet$  There is no assumption like $p$ being outside the convex hull of $\overline\Omega$.

It should be pointed out that, in \cite{Ithermo} the idea developed in Theorem 4.5 has been applied to an inverse cavity problem
governed by a linear system in thermoelasticity \cite{Car}, which is a coupled system of the heat and wave equations.

Second we describe the asymptotic behaviour of another indicator function $I_T^*(\tau;B)$.
In what follows we always choose $T$ in such a way that $\Omega\subset B_{T-\eta}(p)$, that is,
$$\displaystyle
T-\eta\ge R_{\Omega}(p),
\tag {4.89}
$$
where $R_{\Omega}(p)=\sup_{x\in\Omega}\,\vert x-p\vert$.
Recall also the following quantity that is our main concern:
$$\displaystyle
R_D(p)=\sup_{x\in D}\,\vert x-p\vert.
$$
The quantity $R_U(p)$, $U=\omega, D$ gives us the radius of the minimum sphere centered at $p$ that encloses $U$.

\proclaim{\noindent Theorem 4.16(\cite{EV}).}
(i)  Let $\eta$ satisfy 
$$\displaystyle
\eta+2R_D(p)>R_{\Omega}(p).
\tag {4.90}
$$
Then, there exists a positive number $\tau_0$ such that $I_{T}^*(\tau;B)>0$
for all $\tau\ge\tau_0$ and
we have
$$\displaystyle
\lim_{\tau\rightarrow\infty}
\frac{1}{\tau}\log I_{T}^*(\tau;B)
=-2\left\{(T-\eta)-R_D(p)\right\}.
\tag {4.91}
$$

(ii) If  $T>2\{(T-\eta)-R_D(p)\}$, then
$$\displaystyle
\lim_{\tau\rightarrow\infty}
e^{\tau T}I_{T}^*(\tau;B)
=\infty.
$$

(iii)  Assume instead of (4.89) the stronger condition
$$\displaystyle
T-\eta>R_{\Omega}(p).
\tag {4.92}
$$
If $T<2\{(T-\eta)-R_D(p)\}$, then
$$\displaystyle
\lim_{\tau\rightarrow\infty}
e^{\tau T}I_{T}^*(\tau;B)=0.
$$

\endproclaim

Some remarks are in order.

$\bullet$  The restriction (4.90) on $\eta$ says that ball $B$ can not be arbitrary small.
However, if $2R_D(p)\ge R_{\Omega}(p)$, then $\eta$ can be arbitrary small.
Besides, if $\eta$ is large in the sense that $\eta\ge R_{\Omega}(p)$, then (4.90) is satisfied.
Hence as a corollary of Theorem 4.16 we have
$$\displaystyle
\lim_{\tau\rightarrow\infty}\,e^{\tau T}\,I_T^*(\tau;B)
=
\left\{
\begin{array}{ll}
\displaystyle
0 & \text{if $T>2(\eta+R_D(p))$,}\\
\\
\displaystyle
\infty & \text{if $\eta+R_{\Omega}(p)\le T<2(\eta+R_D(p))$.}
\end{array}
\right.
$$
provided the condition $\eta\ge R_{\Omega}(p)$ which is independent of $D$.

$\bullet$  The condition  (4.90) is equivalent to more intuitive condition
$$\displaystyle
T>\{(T-\eta)-R_D(p)\}+(R_{\Omega}(p)-R_D(p)\,).
\tag {4.93}
$$
Note that, see \cite{EV}, the quantity on this right-hand side hsas the lower bound
$$\begin{array}{l}
\displaystyle
\,\,\,\,\,\,
\{(T-\eta)-R_D(p)\}+(R_{\Omega}(p)-R_D(p)\,)
\\
\\
\displaystyle
\ge
\inf\left\{
\vert x-y\vert+\vert y-z\vert\,\vert\,x\in\partial B_{T-\eta}(p),
\,y\in\partial D,\, z\in\partial\Omega\,\right\}
\end{array}
$$
provided $T$ satisfies (4.89).

The meaning of the restriction (4.89) is as follows.  The $v_B$ has the form
$$\displaystyle
v_B(x,t)=\frac{1}{4\pi t}\,\int_{\partial B_t(x)}\,\Psi_B(y)\,dy.
$$
This is because of the Kirchhoff formula and we have the lacuna of the final data
$$\displaystyle
\text{supp}\,v_B(\,\cdot\,,T)\cup \text{supp}\,(v_B)_t(\,\cdot\,,T)
\subset \Bbb R^3\setminus B_{T-\eta}(p).
\tag {4.94}
$$
Thus condition (4.89) yields $v_B(x,T)=(v_B)_t(x,T)=0$, $x\in\Omega$,
Since the wave equation has the time reversal invariance, the function 
$v_B(x,T-t)$, $x\in\Bbb R^3$, $0<t<T$ satisfies the wave equation
and the initial data is supported on the spherical shell $\overline B_{T+\eta}(p)\setminus B_{T-\eta}(p)$
which encloses $\Omega$.
Thus, Theorem 4.16 is another version of Theorem 4.15 in which ball $B$ is replaced with 
a spherical shell.  To generate a necessary Neumann data a time reversal operation acting
on the original Neumann data is employed.

However, to extend Theorem 4.16 to other inverse obstacle problems whose
governing equations are not like the wave equation, it seems that the time reversal invariance together with the property (4.94) 
become a strong obstruction.  However, reconsidering the role of the support of the final time data, we can take an
extremely simple way.
In \cite{EVI} we have found the idea of directly taking the support of the initial data of (4.83) given
by a spherical shell which yields also a similar formula to those of Theorem 4.16.
One can expect that the method therein 
shall cover a broad class of inverse obstacle problems governed by various partial differential equations in time domain.

Note that, in \cite{IHW, EHR} the idea of Theorems 4.15 and 4.16  has also been applied to inverse obstacle problems governed by the heat equation.
Therein to generate a necessary heat flux a solution of a wave equation with propagation speed depending on
a large parameter or its time reverse has been used.



$$\quad$$

\centerline{{\bf Acknowledgment}}

The author was partially supported by Grant-in-Aid for
Scientific Research (C)(No. 17K05331) and (B)(No. 18H01126) of Japan  Society for
the Promotion of Science.

$$\quad$$

\vskip1cm
\noindent
e-mail address

ikehataprobe@gmail.com

\end{document}